\documentclass[a4paper,11pt]{amsart}
\usepackage[margin=28mm]{geometry}
\usepackage{enumitem,microtype}
\usepackage{amssymb,mathrsfs,amsfonts,amsthm,amsmath,amsxtra,dsfont,mathabx}
\usepackage[bookmarksnumbered, colorlinks, linkcolor=blue, anchorcolor=blue, citecolor=blue, plainpages]{hyperref}

\newtheorem{theorem}{Theorem}[section]
\newtheorem{corollary}[theorem]{Corollary}
\newtheorem{lemma}[theorem]{Lemma}
\newtheorem{proposition}[theorem]{Proposition}
\theoremstyle{definition}
\newtheorem{definition}[theorem]{Definition}
\theoremstyle{remark}

\numberwithin{equation}{section}
\newcommand{\CC}{\mathbb C}
\newcommand{\NN}{\mathbb N}
\newcommand{\RR}{\mathbb R}
\newcommand{\ZZ}{\mathbb Z}
\newcommand{\BS}{\mathbb S}

\newcommand{\calB}{\mathcal B}
\newcommand{\calD}{\mathcal D}
\newcommand{\calK}{\mathcal K}
\newcommand{\calL}{\mathcal L}

\newcommand{\calS}{\mathcal S}
\newcommand{\dW}{\dot W}
\newcommand{\Tr}{\operatorname{Tr}}
\newcommand{\rank}{\operatorname{rank}}

\newcommand{\ii}{\mathrm i}
\newcommand{\re}{\operatorname{Re}}
\newcommand{\bcap}{\bigcap}

\DeclareRobustCommand{\dbar}{{\mkern3mu\mathchar'26\mkern-12mu d}}

\title[Endpoint Schatten class properties]{Endpoint Schatten class properties of commutators on quantum Euclidean spaces}

\author[C. Pang]{Changbao Pang}
\address{%
School of Mathematical Sciences\\
Shanxi Normal University\\
Taiyuan, 030031\\
China}
\email{cbpangmath@sxnu.edu.cn}

\author[M. Wang]{Maofa Wang}
\address{%
School of Mathematics and Statistics\\
Wuhan University\\
Wuhan, 430072\\
China}
\email{mfwang.math@whu.edu.cn}

\author[H. Zhang]{Hao Zhang}
\address{%
Instituto de Ciencias Matemáticas\\
Consejo Superior de Investigaciones Científicas\\
C/ Nicolás Cabrera 13-15, 28049, Madrid, Spain}
\email{hao.zhang@icmat.es}

\subjclass[2020]{Primary 46L87, 47B10; Secondary 46L52, 47G30}

\keywords{Quantum Euclidean spaces, Sobolev spaces, fractional Laplacian, commutators, weak Schatten class, Cwikel estimates}

\date{November 3, 2026}

\thanks{C. Pang was supported by National Natural Science Foundation of China (Grant No. 12301153). M. Wang was supported by National Natural Science Foundation of China (Grant No. 12571145).}

\begin{document}

\begin{abstract}
We study commutators $[(-\Delta_\theta)^{\frac{\epsilon}{2}},M_x]$ involving fractional Laplacians on quantum Euclidean spaces associated with arbitrary real antisymmetric matrices $\theta$. For $d\geq2$ and $-\frac d2<\epsilon\neq 0<1$, we characterize critical weak Schatten class membership of $[(-\Delta_\theta)^{\frac{\epsilon}{2}},M_x]$ for any $x$ as a tempered distribution by Sobolev spaces, with equivalent corresponding norms. We also establish precise asymptotics of singular values for $[(-\Delta_\theta)^{\frac{\epsilon}{2}},M_x]$. For this purpose, we extend the Cwikel estimates and the required asymptotics of singular values to arbitrary real antisymmetric matrices, and we also establish a systematic approximation theory of tempered distributions on quantum Euclidean spaces.

\end{abstract}

\maketitle
\tableofcontents

\section{Introduction}
Quantized calculus was introduced by Connes as a noncommutative analogue of the classical differential calculus \cite{C1985} and subsequently related to action functionals and trace formulae in noncommutative geometry \cite{C1988}; see also \cite{C1994}. The language of spectral triples provides a natural framework for formulating quantized calculus. Recall that a spectral triple $(\mathscr A,H,D)$ consists of a Hilbert space $H$, a pre-$C^*$-algebra $\mathscr A$ faithfully represented on $H$ and a densely defined self-adjoint operator $D$ on $H$ with compact resolvent such that every $a\in\mathscr A$ preserves the domain of $D$ and the commutator $[D,a]$ extends to a bounded operator for every $a\in\mathscr A$. Then the quantized differential of $a\in\mathscr A$ is defined as the following commutator:
$$\dbar a=\ii[\operatorname{sgn}(D),a].$$
In the framework of Connes' quantized differential calculus on Euclidean spaces $\RR^d$, the operator $D$ is precisely the Dirac operator, and so $\operatorname{sgn}(D)$ exactly becomes Riesz transforms. 

From the perspective of harmonic analysis, the study of commutators involving Riesz transforms on Euclidean spaces (or more generally, singular integral operators) dates back to the pioneering work of Coifman, Rochberg and Weiss \cite{CRW1976}, and has been studied extensively since then. On Euclidean spaces, among interesting problems on commutators are boundedness, compactness and Schatten class properties. We refer to the interested reader to \cite{CPSZ2019,HLW2016,HLW2017,HW2018,Hyt2021,JP1988,LOR1,LOR2,P2000,PS2004,WZ2026} for more information on boundedness and compactness. In addition, Schatten class properties of commutators have been studied for scalar-valued, matrix-valued, operator-valued and weighted symbols, as well as for spaces of homogeneous type \cite{FWXZZ,Hyt2026,HW2025,JW1982,LLW2024,LLWW2025,LLSW2026,WZ2024}. Beyond Euclidean spaces, related questions have also been considered for Neumann and Bessel operators, Riesz transforms on Heisenberg and stratified Lie groups, and fractional Laplacian commutators \cite{FLL2023,FLLVW2025,FLMSZ2024,FLSZ2024,Hyt2026a,JTZZ2026,LXY2024}. 

From the viewpoint of quantized differential calculus, Schatten class properties of commutators have also been studied extensively on Euclidean spaces $\RR^d$. For instance, Connes, Sullivan and Teleman \cite{CST1994} gave necessary and sufficient conditions for the Schatten--Lorentz classes of quantum derivatives (also called as ``commutators"). Lord, McDonald, Sukochev and Zanin \cite{LMSZ2017} obtained equivalent characterizations for the weak Schatten $\calL_{d,\infty}$-quasi-norm of commutators, and established a Dixmier trace formula. Later, Frank, Sukochev and Zanin \cite{FSZ2023} further established asymptotics of singular values for commutators on $\RR^d$. In addition, Gong, Li and Wick \cite{GLW2023} established the strong Schatten classes and endpoint weak Schatten classes for commutators involving Riesz transforms on weighted $L_2$-spaces. 

In the noncommutative setting, McDonald, Sukochev and Xiong developed quantum differentiability on quantum tori \cite{MSX2019} and quantum Euclidean spaces \cite{MSX2020}, which are indeed the first truly noncommutative examples. Subsequently, they gave a unified analytical approach and strengthened the results on quantum Euclidean spaces in \cite{MSX2023}. Recently, Sukochev, Xiong and Zanin \cite{SXZ2023} obtained asymptotics of singular values for quantized derivatives on quantum tori and refined the corresponding quantum integration formula. More recently, Tian \cite{T2025} established asymptotics of singular values for quantized derivatives on quantum Euclidean spaces $L_\infty(\RR^d_\theta)$ under the nondegeneracy assumption $\det\theta\ne0$. For other related developments in quantized calculus, we refer to \cite{KSZ2024,LZ2026,MSZ2019,MSZ2022,SZcmp2018,SZjot2018} and references therein.

We would like to remark that quantum Euclidean spaces, which originate in the phase-space formulation of quantum mechanics \cite{G1946,M1949}, are now fundamental objects in mathematical physics. These spaces also have played an important role in the development of noncommutative analysis; for example, Fourier multipliers, singular integrals, Sobolev inequalities, nonlinear equations and quantum differentiability on quantum Euclidean spaces have been extensively studied in \cite{GJM2022,GJP2021,JMPX2021,LSZ2020,M2024,MSX2020,MSX2023,T2025}.

In the remarkable paper \cite{FSZ2024}, Frank, Sukochev and Zanin proved that for $-\frac{d}{2}<\epsilon\ne0<1$ and $f\in L_{1}^{\mathrm{loc}}(\RR^d)$, $[(-\Delta)^{\varepsilon/2},M_f]\in \mathcal{L}_{\frac{d}{1-\varepsilon},\infty}(L_2(\RR^d))$ if and only if $f\in \dW^{1}_{\frac{d}{1-\epsilon}}(\RR^d)$, and they also determined the corresponding asymptotics of singular values. It is worthwhile to mention that such commutators involve fractional powers of the Laplacian, and have a similar form to the quantized differential. More recently, Jiao, Tian, Zanin and Zhou \cite{JTZZ2026} studied the strong and weak Schatten class properties of such commutators involving the fractional Laplacian on weighted $L_2$-spaces, and also obtained the corresponding asymptotics of singular values. McDonald, Xiong and Zhang \cite{MXZ2026} proved Weyl's laws for operator-valued pseudo-differential operators and obtained corresponding asymptotic formulae for operator-valued fractional-Laplacian commutators. In addition, Zhang \cite{Z2026QT} characterized weak Schatten class membership of fractional-Laplacian commutators on quantum tori and established the corresponding asymptotics of singular values.

Inspired by the remarkable work of Frank, Sukochev and Zanin \cite{FSZ2024}, the purpose of this paper is to establish the analogous results on quantum Euclidean spaces for arbitrary real antisymmetric matrices $\theta$ within a general setting, which includes the degenerate case $\det\theta=0$ and symbols that may be tempered distributions.

\

Now we set up some notation for the statement of our main results. For detailed definitions, see Section \ref{s:pre} and Section \ref{sec:qtd}. Throughout this paper, unless stated otherwise, $d\ge2$ and $\theta=\{\theta_{j,k}\}_{1\le j,k\le d}$ is a fixed antisymmetric real $d\times d$ matrix, which may be degenerate, i.e. $\det\theta=0$. Denote by $L_\infty(\RR^d_\theta)$ the quantum Euclidean space associated with $\theta$. Let $\calS(\RR^d_\theta)$ denote the Schwartz space and $L_p(\RR^d_\theta)$ the noncommutative $L_p$-spaces on $L_\infty(\RR^d_\theta)$. Denote by $\calS'(\RR^d_\theta)$ the space of tempered distributions. Let $\dW^1_{p}(\RR^d_\theta)$ be the homogeneous Sobolev space associated with $L_\infty(\RR^d_\theta)$. Let $-\Delta_\theta=\sum_{k=1}^d\partial_k^2$. For $\epsilon\in(-\frac d2,1)$ and $x\in\calS'(\RR^d_\theta)$, we define 
$$[(-\Delta_\theta)^{\frac{\epsilon}{2}}, M_x]=(-\Delta_\theta)^{\frac{\epsilon}{2}} M_x-M_x(-\Delta_\theta)^{\frac{\epsilon}{2}}:\calS(\RR^d_\theta)\to\calS'(\RR^d_\theta).$$
For $x\in L_1(\RR^d_\theta)+L_\infty(\RR^d_\theta)$, $M_x$ is just the left multiplication on $\calS(\RR^d_\theta)$, i.e., $M_xy=xy$ for $y\in \calS(\RR^d_\theta)$. For general $x\in\calS'(\RR^d_\theta)$, Proposition \ref{prop:dw1} implies that the two operators $(-\Delta_\theta)^{\frac{\epsilon}{2}}M_x$ and $M_x(-\Delta_\theta)^{\frac{\epsilon}{2}}$ are indeed $\calS'(\RR^d_\theta)$-valued, and so well-defined. More precisely, for $a,b\in\calS(\RR^d_\theta)$ the two distributions are determined respectively by the pairings
$$\big\langle(-\Delta_\theta)^{\frac{\epsilon}{2}}M_xa,b^*\big\rangle=\big\langle x,a\big((-\Delta_\theta)^{\frac{\epsilon}{2}}b\big)^*\big\rangle$$
and
$$\big\langle M_x(-\Delta_\theta)^{\frac{\epsilon}{2}}a,b^*\big\rangle
=\big\langle x,\big((-\Delta_\theta)^{\frac{\epsilon}{2}}a\big)b^*\big\rangle.$$
These formulas are well-defined because the two expressions on the right-hand-side belong to $\calS(\RR^d_\theta)$; see Proposition \ref{prop:dw1}.

Our main result characterizes the endpoint weak Schatten class membership of $[(-\Delta_\theta)^{\frac{\epsilon}{2}}, M_x]$, and determines the corresponding asymptotics of singular values.

\begin{theorem}\label{thm:m1}
Let $-\frac d2<\epsilon\neq 0<1$ and $x\in\calS'(\RR^d_\theta)$. The following conditions are equivalent:
\begin{description}
\item[(i)] $[(-\Delta_\theta)^{\frac{\epsilon}{2}},M_x]\in\calL_{\frac{d}{1-\epsilon},\infty}(L_2(\RR^d_\theta))$;
\item[(ii)] $x\in\dW^1_{\frac{d}{1-\epsilon}}(\RR^d_\theta)$.
\end{description}
Moreover, in this case,
$$
\|[(-\Delta_\theta)^{\frac{\epsilon}{2}},M_x]\|_{\calL_{\frac{d}{1-\epsilon},\infty}(L_2(\RR^d_\theta))}
\asymp_{d,\epsilon}\|x\|_{\dW^1_{\frac{d}{1-\epsilon}}(\RR^d_\theta)}
$$
and
\begin{align*}
&\lim_{t\to\infty}t^{\frac{1-\epsilon}{d}}\mu\big(t,[(-\Delta_\theta)^{\frac{\epsilon}{2}},M_x]\big)\\
=&|\epsilon|(d(2\pi)^d)^{-\frac{1-\epsilon}{d}}
\Big(\int_{\BS^{d-1}}\tau_\theta\Big(
\Big|\sum_{j=1}^dt_j\partial_jx\Big|^{\frac{d}{1-\epsilon}}
\Big)dS(t)\Big)^{\frac{1-\epsilon}{d}}.
\end{align*}
\end{theorem}

Theorem \ref{thm:m1} extends the classical results of \cite{FSZ2024} to quantum Euclidean spaces associated with arbitrary real antisymmetric matrices. We would like to stress that the symbol $x$ can also be allowed to be a tempered distribution, whereas \cite{FSZ2024} assumed $f\in L_{1,\mathrm{loc}}(\RR^d)$ and \cite{Z2026QT} assumed $x\in L_1(\mathbb T_\theta^d)$. In other words, no a priori integrability condition is imposed on symbols in our Theorem \ref{thm:m1}. More explicitly, starting only from the weak Schatten class assumption on the commutators, we prove that the distributional derivatives $\partial_jx$ belong to $L_{\frac{d}{1-\epsilon}}(\RR^d_\theta)$, and hence this yields that $x\in\dW^1_{\frac{d}{1-\epsilon}}(\RR^d_\theta)$.

The ideas underlying our proofs, including tempered distributions on quantum Euclidean spaces, double operator integrals, Cwikel estimates and spectral asymptotics, are mainly drawn from \cite{FSZ2024}. However, extending these arguments to arbitrary real antisymmetric matrices $\theta$ and symbols in $\calS'(\RR^d_\theta)$ presents additional difficulties. First, the Cwikel estimates and spectral asymptotics needed in the proof must be established without the nondegeneracy assumption on $\theta$. Moreover, since we only assume $x\in\calS'(\RR^d_\theta)$ here, while \cite{FSZ2024} worked with locally integrable functions, our proofs rely more heavily on the study of $\calS'(\RR^d_\theta)$. This also forces us to develop a systematic and self-contained theory of tempered distributions on quantum Euclidean spaces, which is of independent interest.

On quantum Euclidean spaces, we need the Cwikel estimates for arbitrary real antisymmetric matrices $\theta$. While the upper bound in Theorem \ref{thm:m1} only relies on the Cwikel estimates for $2\leq p<\infty$, the lower bound requires the Cwikel estimates for $1\leq p<2$. The Cwikel estimates, conjectured by Simon \cite{S1976}, were solved by Cwikel \cite{Cwikel1977} (see also \cite{BS1977}). Levitina, Sukochev and Zanin \cite{LSZ2020} revisited these estimates in an abstract framework based on generalized singular values and symmetric interpolation spaces, and established the corresponding strong and weak Schatten estimates on quantum Euclidean spaces in nondegenerate cases. We would like to remark that the Cwikel estimates of \cite{LSZ2020} for the strong Schatten class ($2\leq p<\infty$) and the weak Schatten class ($2<p<\infty$) remain valid without assuming nondegeneracy. For the range $1\leq p<2$, building on \cite{LSZ2020} with appropriate modifications, we obtain the corresponding version of the Cwikel estimates for arbitrary real antisymmetric matrices $\theta$.

\

Now we illustrate the detailed methodology. At first, we use the transference method to pass from $L_2(\RR^d_\theta)$ to $L_2(\RR^d)$, and so we need to make essential use of the isometric isomorphism between $L_2(\RR^d)$ and $L_2(\RR^d_\theta)$. Recall from \cite{MSX2020} that this isometric isomorphism is given as
$$V=(2\pi)^{-d/2}U_\theta:L_2(\RR^d)\longrightarrow L_2(\RR^d_\theta).$$
Under this $V$, $-\Delta_\theta$ is transferred to an operator on $L_2(\RR^d)$, of which the detailed definition is presented as follows.

For $1\leq k\leq d$, let $\calD_k$ be the coordinate multiplication operator on $L_2(\RR^d)$, defined on
$$\operatorname{Dom}(\calD_k)=\{f\in L_2(\RR^d):\xi_kf(\xi)\in L_2(\RR^d)\}$$
by
$$(\calD_kf)(\xi)=\xi_kf(\xi),\quad \xi\in\RR^d.$$
Let $\nabla=(-\ii\calD_1,\ldots,-\ii\calD_d)$, and denote by $E_\nabla$ the joint spectral measure of the commuting self-adjoint tuple $\ii\nabla=(\calD_1,\ldots,\calD_d)$. For a Borel function $g:\RR^d\to\CC$, define $$M_g=g(\ii\nabla)=\int_{\RR^d}g(t) dE_\nabla(t).$$
Thus $M_g$ is the multiplication operator by $g$, with the domain
$$\operatorname{Dom}(M_g)=\{f\in L_2(\RR^d):gf\in L_2(\RR^d)\}.$$
Define $-\Delta=\sum_{k=1}^d\calD_k^2$. Hence, by virtue of $V$, for $x\in\calS(\RR^d_\theta)$,
$$V^*[(-\Delta_\theta)^{\frac{\epsilon}{2}},M_x]V=[(-\Delta)^{\frac{\epsilon}{2}},x].$$
With the help of this transference, the commutator $[(-\Delta)^{\frac{\epsilon} {2}},x]$ for $x\in\calS'(\RR^d_\theta)$, is defined as an map from $\calS(\RR^d)$ to $\calS'(\RR^d)$ as follows:
$$\big\langle[(-\Delta)^{\frac{\epsilon}{2}},x]f,\overline g\big\rangle=\langle x,\Phi_{\epsilon;Vf,Vg}\rangle,\quad f,g\in\calS(\RR^d),$$
where
$$\Phi_{\epsilon;Vf,Vg}=(2\pi)^{-d}\Big(U_\theta(f)U_\theta(|\cdot|^\epsilon g)^*-U_\theta(|\cdot|^\epsilon f)U_\theta(g)^*\Big).$$
Therefore, it suffices to consider Schatten class properties of $[(-\Delta)^{\frac{\epsilon}{2}},x]$ for $x\in\calS'(\RR^d_\theta)$. The proof will be divided into two parts, respectively; one is the upper bound, and the other one is the lower bound. In addition, we would like to mention that in this paper, we adopt the notation $(-\Delta)^{\epsilon/2}$ merely as a symbol for the transferred fractional Laplacian; it is not the classical operator used in \cite{FSZ2024}. In \cite{FSZ2024}, $(-\Delta)^{\epsilon/2}$ is the Fourier multiplier with symbol $|\xi|^\epsilon$.

Regarding the upper bound, we first consider a Schwartz symbol $x\in\calS(\RR^d_\theta)$, where the commutator can be expressed in terms of double operator integrals and then bounded by the Cwikel estimates. This gives
\begin{align}\label{eq:dxu1}
\|[(-\Delta)^{\frac{\epsilon}{2}},x]\|_{\calL_{\frac{d}{1-\epsilon},\infty}(L_2(\RR^d))}\lesssim_{d,\epsilon}\|x\|_{\dW^1_{\frac{d}{1-\epsilon}}(\RR^d_\theta)}.
\end{align}
As established in \cite{MSX2023}, $\calS(\RR^d_\theta)$ is dense in
$\dW^1_{\frac{d}{1-\epsilon}}(\RR^d_\theta)$, so the estimate for general $x$ follows by approximation. 

For the lower bound, asymptotics of singular values of $[(-\Delta)^{\epsilon/2},x]$ will be needed. To establish asymptotics of singular values of $[(-\Delta)^{\epsilon/2},x]$, we first prove that
\begin{align}\label{eq:dx1}
[(-\Delta)^{\frac{\epsilon}{2}},x]-\epsilon\sum_{j=1}^d(\partial_jx)\pi_2(\varLambda_j)J^{-(1-\epsilon)}\in\big(\calL_{\frac{d}{1-\epsilon},\infty}(L_2(\RR^d))\big)_0,
\end{align}
where $J=(1-\Delta)^{1/2}$, $\varLambda_j(s)=s_j$ and $\pi_2(g)$ denotes multiplication by $g(\xi/|\xi|)$; see Section \ref{s:spec} for detailed definitions. Thus, this difference does not affect asymptotics of singular values, so it remains to find asymptotics of finite sums $\sum_jx_j\pi_2(g_j)J^{-\alpha}$. In Theorem \ref{thm:allp}, we establish
\begin{align}\label{eq:tx1}
\lim_{t\to\infty}t^{\frac{\alpha}{d}}\mu\Big(t,\sum_jx_j\pi_2(g_j)J^{-\alpha}\Big)=(d(2\pi)^d)^{-\frac{\alpha}{d}}\Big(\int_{\BS^{d-1}}\tau_\theta\Big(\Big|\sum_jg_j(s)x_j\Big|^{\frac{d}{\alpha}}\Big) dS(s)\Big)^{\frac{\alpha}{d}}
\end{align}
for $0<\alpha\leq d$, $x_j\in\calS(\RR^d_\theta)$, and $g_j\in C(\BS^{d-1})$.
To establish \eqref{eq:tx1} , we derive the trace formula directly from the Hilbert--Schmidt identity, so that the abstract result of \cite[Theorem 6.4]{T2025} applies without the nondegeneracy assumption. We also develop smooth functional calculus on $\calS(\RR^d_\theta)$; see Lemma \ref{lem:root}.

By virtue of \eqref{eq:dx1} and \eqref{eq:tx1}, we show that for $x\in \dW^1_{\frac{d}{1-\epsilon}}(\RR^d_\theta)$,
\begin{align}\label{eq:xd1}
\|x\|_{\dW^1_{\frac{d}{1-\epsilon}}(\RR^d_\theta)}
\lesssim_{d,\epsilon}\|[(-\Delta)^{\frac{\epsilon}{2}},x]\|_{\calL_{\frac{d}{1-\epsilon},\infty}(L_2(\RR^d))},
\end{align}
%However, for general $x\in\calS'(\RR^d_\theta)$, we utilize the approximation theory for tempered distributions on quantum Euclidean spaces, which is established in Section \ref{s:appr}. 
For general $x\in\calS'(\RR^d_\theta)$, since $\partial_jx\in L_{\frac{d}{1-\epsilon}}(\RR^d_\theta)$ does not hold in general, so approximation in the homogeneous Sobolev spaces cannot be used directly. To this end, we utilize approximation properties for tempered distributions on quantum Euclidean spaces established in Section \ref{s:appr}. More precisely, by using the heat kernel $h_t=(4\pi t)^{-d/2}e^{-|\cdot|^2/(4t)}$ for $t>0$, we define $y_t=h_t*x$, which satisfy
$$\widehat{y_t}=e^{-t|\cdot|^2}\widehat x,\quad y_t\to x\ \text{in }\calS'(\RR^d_\theta) \ \text{as} \ t\to 0^+ ,$$
and
$$ [(-\Delta)^{\epsilon/2},y_t]=h_t*[(-\Delta)^{\epsilon/2},x],$$
so that the weak Schatten quasi-norm of $[(-\Delta)^{\epsilon/2},y_t]$ can be controlled uniformly by that of $[(-\Delta)^{\epsilon/2},x]$. Then we construct a net of self-adjoint Schwartz elements $z_\delta=z_\delta^*\in\calS(\RR^d_\theta)$ with $\|z_\delta\|_{L_\infty(\RR^d_\theta)}\leq1$ such that $z_\delta y z_\delta\to y$ in $\calS(\RR^d_\theta)$ for every Schwartz element $y$. Indeed, we will show a stronger result that $z_\delta (\partial_jy_t) z_\delta$ also belongs to $\calS(\RR^d_\theta)$. Moreover, similar asymptotics of singular values for $z_\delta[(-\Delta)^{\epsilon/2},y_t]z_\delta$ shall be established as well. Therefore, by applying asymptotics of singular values for Schwartz symbols to $z_\delta[(-\Delta)^{\epsilon/2},y_t]z_\delta$, we obtain the following uniform bounds
$$\sum_{j=1}^d\|z_\delta(\partial_jy_t)z_\delta\|_{L_{\frac{d}{1-\epsilon}}(\RR^d_\theta)}\lesssim_{d,\epsilon}\|[(-\Delta)^{\epsilon/2},x]\|_{\calL_{\frac{d}{1-\epsilon},\infty}(L_2(\RR^d))}.$$
Finally, by means of the approximation theory for distributions established in Section \ref{s:appr}, we obtain \eqref{eq:xd1} for general $x\in\calS'(\RR^d_\theta)$.

The paper is organized as follows. In Section \ref{s:pre}, we present some notation and recall several basic properties of Schatten classes and quantum Euclidean spaces. In Section \ref{sec:qtd}, we introduce Sobolev and Schwartz spaces together with the commutator $[(-\Delta_\theta)^{\epsilon/2},M_x]$ defined for tempered distributions, and in Section \ref{s:appr} we develop the approximation theory of tempered distributions on quantum Euclidean spaces. In Section \ref{s:cw}, we establish the Cwikel estimates for arbitrary real antisymmetric matrices $\theta$, without the nondegeneracy assumption. In Section \ref{s:doi}, we use double operator integrals to prove the weak Schatten class bound \eqref{eq:dxu1}. In Section \ref{s:prin}, we show \eqref{eq:dx1}. In Section \ref{s:spec}, we establish \eqref{eq:tx1}, namely asymptotics of singular values of $\sum_j\pi_1(x_j)\pi_2(g_j)J^{-\alpha}$. In Section \ref{s:comm}, we complete the proof of Theorem \ref{thm:m1}.

\bigskip

\section{Preliminaries}\label{s:pre}
In this section, we fix the notation used throughout the paper and recall some basic facts about the Schatten classes $\calL_p$ and the weak Schatten classes $\calL_{p,\infty}$. We also recall the quantum Euclidean spaces $L_\infty(\RR^d_\theta)$.

\

\noindent\textbf{Notation:} We write $A\lesssim B$, or equivalently $B\gtrsim A$, if $A\le CB$ where $C$ is a constant independent of the quantities being compared, and we write $A\asymp B$ if both $A\lesssim B$ and $B\lesssim A$. Subscripts list the dependence on the permitted parameters of the implicit constants; for example, the constant in $A\lesssim_{\alpha,\beta,\gamma}B$ may depend only on $\alpha$, $\beta$, and $\gamma$. We write $\NN_0=\{0,1,2,\ldots\}$ and $\NN=\{1,2,\ldots\}$. Denote by $\langle\xi,\eta\rangle$ the Euclidean inner product of $\xi,\eta\in\RR^d$, and by $\ii$ the imaginary unit such that $\ii^2=-1$. For $1\le p\le\infty$, let $p'$ denote the H\"older conjugate exponent such that $\frac1p+\frac1{p'}=1$. For two operators $A,B$, denote by $[A,B]=AB-BA$ their commutator. For every integer $k\geq1$ and $\psi\in L_1(\RR^k)$, the Fourier transform of $\psi$ is given by
$$\widehat\psi(\xi)=\int_{\RR^k}e^{-\ii\langle s,\xi\rangle}\psi(s) ds,\quad \xi\in\RR^k.$$
For $\psi\in\calS(\RR^k)$, the Fourier inversion formula is
$$\psi(s)=\frac{1}{(2\pi)^k}\int_{\RR^k}e^{\ii\langle s,\xi\rangle}\widehat\psi(\xi) d\xi,\quad s\in\RR^k.$$
Given multi-indices $\alpha=(\alpha_1,\dots,\alpha_d)$ and $\beta=(\beta_1,\dots,\beta_d)$, the notation $\alpha\le\beta$ means that $\alpha_i\le\beta_i$ for every $i=1,\ldots,d$. For a multi-index $\alpha=(\alpha_1,\dots,\alpha_d)$ and $\xi=(\xi_1,\ldots,\xi_d)\in\RR^d$, we write $\xi^\alpha$ for $\xi_1^{\alpha_1}\cdots\xi_d^{\alpha_d}$. Let
$$Q=[0,1)^d,\quad Q_n=n+Q\quad \text{and}\quad P_n=\mathds{1}_{Q_n}(\ii\nabla)=M_{\mathds{1}_{Q_n}}, \quad n\in\ZZ^d.$$

\subsection{Schatten class}\label{ssec:Sc1}
We present definitions and some properties of Schatten classes, weak Schatten classes and singular values; see \cite{LSZ2012} and \cite{S2005} for more information on related subjects. Let $H$ be a complex separable Hilbert space. Denote by $\calB(H)$ and $\calK(H)$ the sets of bounded and compact linear operators on $H$, respectively. For a normed space $X$, the operator norm is denoted by $\|\cdot\|_{\calB(X)}$. For $T\in\calK(H)$, define its singular value function as
$$\mu(t,T) = \inf\{\|T-R\|_{\mathcal B(H)} : \operatorname{rank}(R) \le t\}. $$
Denote by $\mu(T)$ the sequence $\{\mu(n,T)\}_{n=0}^{\infty}$. For $0<p<\infty$, the weak Schatten class, denoted by $\calL_{p,\infty}(H)$, is the set of all operators $T\in\calK(H)$ such that $\mu(T)$ is in the weak-$L_p$ space $\ell_{p,\infty}$, endowed with the following norm:
$$\|T\|_{\calL_{p,\infty}(H)}=\sup_{n\ge 0}(n+1)^{1/p}\mu(n,T)<\infty.$$
If $1<p<\infty$, $\calL_{p,\infty}(H)$ admits one equivalent Banach norm given by
$$\|T\|_{\calL_{p,\infty}(H)}^{\#}
:=\sup_{N\geq1}N^{1/p-1}\sum_{n=0}^{N-1}\mu(n,T).$$
Moreover, for $1<p<\infty$ and any $T\in \calL_{p,\infty}(H)$,
\begin{equation}\label{eq:bn1}
\|T\|_{\calL_{p,\infty}(H)}\leq\|T\|_{\calL_{p,\infty}(H)}^{\#}
\leq\frac p{p-1}\|T\|_{\calL_{p,\infty}(H)}.
\end{equation}
The following H\"older-type inequality will be used:
$$\|TS\|_{\calL_{r,\infty}(H)}\le2^{\frac{1}{r}}\|T\|_{\calL_{p,\infty}(H)}\|S\|_{\calL_{q,\infty}(H)}, \quad \frac{1}{r}=\frac{1}{p}+\frac{1}{q}.$$
For $0<p<\infty$, the Schatten class $\calL_p(H)$ consists of all $T\in\calK(H)$ such that $\mu(T)$ belongs to $\ell_p$, endowed with the follwoing norm:
$$\|T\|_{\calL_p(H)}:=\|\mu(T)\|_{\ell_p(\NN_0)}
=\big(\sum_{n=0}^\infty \mu(n,T)^p\big)^{1/p}.$$
Denote by $\Tr$ the canonical semifinite trace on the positive operators in $\calB(H)$ and its linear extension to $\calL_1(H)$. For $0<p<\infty$, denote by $(\mathcal{L}_{p,\infty}(H))_{0}$ the set of $T\in \mathcal{L}_{p,\infty}(H)$ satisfying
$$\lim_{n\to\infty}(n+1)^{1/p}\mu(n,T)=0.$$
For $0<p<q<\infty$, the following embedding properties clearly follow from the definition:
\begin{description}
\item[(i)] $\calL_{p,\infty}(H)\hookrightarrow\calL_{q,\infty}(H)$;
\item[(ii)] $\calL_{p,\infty}(H)\subseteq(\calL_{q,\infty}(H))_0$;
\item[(iii)] $\calL_{p,\infty}(H)\hookrightarrow\calL_q(H)$; \item[(iv)] $\calL_p(H)\hookrightarrow\calL_{p,\infty}(H)$;
\item[(v)] $\calL_p(H)\subseteq(\calL_{p,\infty}(H))_0$ when $1\le p<\infty$.
\end{description}
All embeddings above are continuous. We refer the interested reader to \cite{LSZ2012,S2005} for more information on Schatten classes. 

We also need to use the following properties of singular values:
$$\mu(t+s,A+B)\leq \mu(t,A)+\mu(s,B),\quad t,s\geq 0,$$
and, for every $0<p<\infty$,
$$\mu(t,A)=\mu(t,A^{*})=\mu(t,|A|^{p})^{\frac{1}{p}},\quad t\geq 0.$$
The interested reader is referred to \cite[Section 2.3]{LSZ2012} for detailed proofs. The following Fatou property for $\calL_p(H)$ and $\calL_{p,\infty}(H)$ will be used later as well, and e.g., see \cite[Theorem 2.7(d)]{S2005}.

\begin{lemma}\label{lem:fatou}
Let $0<p<\infty$. Suppose that $T_n,T\in\calB(H)$, $T_n\to T$ in the weak operator topology and $\sup_n\|T_n\|_{\calL_{p,\infty}(H)}<\infty$. Then $T\in\calL_{p,\infty}(H)$ and
\begin{equation}\label{eq:fatou}
\|T\|_{\calL_{p,\infty}(H)}
\leq\liminf_{n\to\infty}\|T_n\|_{\calL_{p,\infty}(H)}.
\end{equation}
In addition, the analogous lower-semicontinuity assertion holds for $\calL_p(H)$.
\end{lemma}

The following perturbation result is from \cite[Lemma 3.1]{FSZ2023}.
\begin{lemma}\label{lem:asst1}
Let $0<p<\infty$, $A\in\mathcal{L}_{p,\infty}(H)$ and $B\in(\mathcal{L}_{p,\infty}(H))_0$. Then
$$\limsup_{t\to\infty} t^{\frac{1}{p}}\mu(t,A+B) = \limsup_{t\to\infty} t^{\frac{1}{p}}\mu(t,A),$$
and
$$\liminf_{t\to\infty} t^{\frac{1}{p}}\mu(t,A+B) = \liminf_{t\to\infty} t^{\frac{1}{p}}\mu(t,A).$$
In particular,
$$\lim_{t\to\infty} t^{\frac{1}{p}}\mu(t,A+B) = \lim_{t\to\infty} t^{\frac{1}{p}}\mu(t,A),$$
provided that the limit on the right-hand-side exists.
\end{lemma}

In order to implement approximation arguments, we need to use the following asymptotic stability theorem in its singular-value form; see \cite[Theorem 4.1 and Remark 4.2]{BS1970}.

\begin{lemma}\label{lem:asapp}
Let $0<p<\infty$ and $T,T_n\in\calL_{p,\infty}(H)$. Suppose that
$$\lim_{n\to\infty}\limsup_{t\to\infty}
t^{1/p}\mu(t,T_n-T)=0.$$
If the limits $\lim_{t\to\infty}t^{1/p}\mu(t,T_n)$ exist and there exists $a\geq 0$ such that
$$\lim_{n\to\infty}\lim_{t\to\infty}t^{1/p}\mu(t,T_n)=a,$$
then
$$\lim_{t\to\infty}t^{1/p}\mu(t,T)=a.$$
In particular, the conclusion applies when $\|T_n-T\|_{\calL_{p,\infty}(H)}\to0$.
\end{lemma}

The following lemma is a trace formulation of Ky Fan's variational principle; see \cite[Theorem 3.4.1]{HJ1991} for the matrix version. For completeness, we include a proof for compact operators on any complex separable Hilbert space $H$. Recall that if $A\in B(H)$ and $A\geq0$, then for any orthonormal basis $(e_j)$ of $H$, 
$$\Tr(A)=\sum_j\langle Ae_j,e_j\rangle_H\in[0,\infty].$$

\begin{lemma}\label{lem:mTr1}
Let $H$ be a complex separable Hilbert space and $A$ be a compact operator on $H$. Then for every integer $N\geq1$,
$$\sum_{j=0}^{N-1}\mu(j,A)
=\sup_{\substack{F\in\calB(H),\ \rank F\leq N\\
\|F\|_{\calB(H)}\leq1}}|\Tr(FA)|.$$
The sequence of singular values is extended by zeros when necessary.
\end{lemma}
\begin{proof}
For every $F\in\calB(H)$ with $\rank F\leq N$ and $\|F\|_{\calB(H)}\leq1$, it is clear that
$$|\Tr(FA)|
\leq\|FA\|_{\calL_1(H)}
\leq\|F\|_{\calB(H)}\sum_{j=0}^{N-1}\mu(j,A)
\leq\sum_{j=0}^{N-1}\mu(j,A).$$
For the converse inequality, let $A=U|A|$ be the polar decomposition of $A$. By the spectral theorem, we can choose mutually orthonormal eigenvectors corresponding to the positive singular values:
\begin{equation}\label{eq:Ax1}
|A|\xi_j=\mu(j,A)\xi_j,\quad \mu(j,A)>0.
\end{equation}
Let $P_N$ be the projection onto the span of $\xi_j$ for all $0\leq j<N$. Then $\rank P_N\leq N$. Since each selected eigenvector $\xi_j$ belongs to $\operatorname{Ran}|A|$, we deduce $P_N\leq U^*U$. Then
$$\Tr(P_N|A|)=\sum_{\substack{0\leq j<N\\\mu(j,A)>0}}\mu(j,A)
=\sum_{j=0}^{N-1}\mu(j,A).$$
The operator $F_N=P_NU^*$ satisfies $\rank F_N\leq N$, $\|F_N\|_{\calB(H)}\leq1$ and $F_NA=P_N|A|$. Thus one can choose $F_N$ to attain the required supremum.
\end{proof}

We shall use the following characterization of duality for Schatten classes by finite-rank operators. Its proof is standard, and so we omit it.

\begin{lemma}\label{lem:dual}
Let $H$ be a complex separable Hilbert space and $T\in\calB(H)$. 
\begin{description}
\item[(I)] Let $1< p<\infty$. If
$$C_p:=\sup\{|\Tr(BT)|:B\in\calB(H),\ \rank B<\infty, \|B\|_{\calL_{p'}(H)}\leq1\}<\infty,$$
then $T\in\calL_p(H)$ and
\begin{align}\label{eq:Tp}
\|T\|_{\calL_p(H)}=C_p.
\end{align}
\item[(ii)] if
$$C_1:=\sup\{|\Tr(BT)|:B\in\calB(H),\ \rank B<\infty, \|B\|_{\calB(H)}\leq1\}<\infty,$$
then $T\in\calL_1(H)$ and
\begin{align}\label{eq:T1}
\|T\|_{\calL_1(H)}=C_1.
\end{align}
\end{description}
\end{lemma}

\subsection{Quantum Euclidean spaces}\label{ssec:QEs1}
We follow the same approach as in \cite{LSZ2020} to present quantum Euclidean spaces; equivalent formulations can be found in \cite[Sect. 5.2.2.2]{CGRS2014}, \cite{GGSV2004}, \cite{GJP2021} and \cite[Chapter 14]{Ha2013}. Associated with $\theta$ is the strongly continuous family of unitary operators
\begin{equation}\label{eq:wun1}
(U_\theta(s)\xi)(t)
=e^{-\frac{\ii}{2}\langle s,\theta t\rangle}\xi(t-s),
\quad s,t\in\RR^d,\xi\in L_2(\RR^d).
\end{equation}
The family $\{U_\theta(t)\}_{t\in\RR^d}$ satisfies the Weyl form of the canonical commutation relations:
$$U_\theta(s)U_\theta(t)
=e^{-\frac{\ii}{2}\langle s,\theta t\rangle}U_\theta(s+t).$$
The von Neumann subalgebra of $\calB(L_2(\RR^d))$ generated by $\{U_\theta(t)\}_{t\in\RR^d}$ is called \emph{quantum Euclidean space} denoted by $L_\infty(\RR^d_\theta)$. Notably, Gonz\'alez-P\'erez, Junge and Parcet \cite{GJP2021} adopted the term \emph{quantum Euclidean spaces} for the associated von Neumann algebras.

For $f\in L_1(\RR^d)$, the operator $U_\theta(f)$ is defined by
\begin{equation}\label{eq:defufh1}
U_\theta(f)= \int_{\RR^d} f(t) U_\theta(t) dt.
\end{equation}
The above integral is absolutely convergent in the Bochner sense and defines a bounded linear operator $U_\theta(f):L_2(\RR^d)\to L_2(\RR^d)$ belonging to $L_\infty(\RR^d_\theta)$; see \cite[Lemma 2.3]{MSX2020}. Denote by $\calS(\RR^d)$ the classical Schwartz space on $\RR^d$. The Schwartz space associated with $L_\infty(\RR^d_\theta)$ is defined as
$$\calS(\RR^d_\theta) = \{ x \in L_\infty(\RR^d_\theta) : x = U_\theta(f) \text{ for some } f \in \calS(\RR^d) \}.$$
For $x,y\in\calS(\RR^d_\theta)$, it is clear that $xy\in\calS(\RR^d_\theta)$. We write $x=U_\theta(f)\in\calS(\RR^d_\theta)$ to mean $f\in\calS(\RR^d)$. For $N\in\NN_0$ and $f\in\calS(\RR^d)$, define
$$\mathfrak p_N(f)=\max_{|\alpha|\leq N}\sup_{s\in\RR^d}(1+|s|)^N|\partial_s^\alpha f(s)|.
$$
We equip $\calS(\RR^d_\theta)$ with the topology induced by the canonical Fr\'echet topology on $\calS(\RR^d)$ via the topological linear isomorphism $U_\theta$, that is 
$$\mathfrak p_N(x)=\max_{|\alpha|\leq N}\sup_{s\in\RR^d}(1+|s|)^N|\partial_s^\alpha f(s)|,\quad x=U_\theta(f)\in\calS(\RR^d_\theta),
$$
Here we still use $\mathfrak p_N$ to denote the seminorms on $\calS(\RR^d_\theta)$ as this will not cause any confusion.

For $f\in\calS(\RR^d)$, a functional $\tau_\theta:\calS(\RR^d_\theta)\to\CC$ is defined by
\begin{equation}\label{eq:fun1}
\tau_\theta( U_\theta(f)) = \tau_\theta\Big( \int_{\RR^d} f(\eta) U_\theta(\eta) d\eta \Big) = (2\pi)^d f(0).
\end{equation}
This functional extends to a normal faithful semifinite trace on $L_\infty(\RR^d_\theta)$; see \cite{GJP2021}. The inverse map $U_\theta^{-1}:\calS(\RR^d_\theta)\to\calS(\RR^d)$ is given by
\begin{equation}\label{eq:inveru1}U_\theta^{-1}(x)=\hat{x},\quad
\hat{x}(s)=(2\pi)^{-d}\tau_\theta\big(xU_\theta(s)^*\big),
\quad s\in\mathbb{R}^d,
\end{equation}
where $\widehat{x}$ denotes the Fourier transform of $x$, following
\cite{GJP2021}.

For $1\le p<\infty$, the space $L_p(\RR^d_\theta)$ is the completion of $\{x\in L_\infty(\RR^d_\theta):\tau_\theta(|x|^p)<\infty\}$ with respect to the norm $\|x\|_{L_p(\RR^d_\theta)}=\tau_\theta(|x|^p)^{\frac{1}{p}}$. By \cite{GJP2021}, $\calS(\RR^d_\theta)$ is $w^*$ dense in $L_\infty(\RR^d_\theta)$, and also dense in $L_p(\RR^d_\theta)$ with respect to the $L_p$-norm for $1\le p<\infty$. As explained in \cite{GJP2021}, if $\rank\theta=2r$ and $k=d-2r$, then
$$L_\infty(\RR^d_\theta)\simeq L_\infty(\RR^k) \overline\otimes\calB(L_2(\RR^r)),$$
When $\theta=0$, $L_p(\RR^d_\theta)$ is just the usual Lebesgue space $L_p(\RR^d)$, whereas for $\det\theta\ne0$ it can be identified with the Schatten class $\calL_p(L_2(\RR^{\frac{d}{2}}))$. By \cite{MSX2020} or \cite[Chapter 2, Lemma 3.1]{T2008}, one has the following proposition.

\begin{proposition}\label{prop:iso1}
The map
$$V:=(2\pi)^{-\frac{d}{2}}U_\theta: L_2(\RR^d)\to L_2(\RR^d_\theta)$$
is an isometric isomorphism.
\end{proposition}

We will use the noncommutative Hausdorff--Young inequality on quantum Euclidean spaces; see \cite[Proposition 2.10]{MSX2020}.
\begin{proposition}\label{prop:iso2}
Let $1\le p\le 2$. Then for every $f\in L_p(\RR^d)\bcap L_1(\RR^d)$, $U_\theta(f)\in L_{p'}(\RR^d_\theta)$ and
$$\|U_\theta(f)\|_{L_{p'}(\RR^d_\theta)}\le (2\pi)^{\frac{d}{2}}\|f\|_{L_p(\RR^d)},$$
and hence $U_\theta$ has continuous extension from $L_p(\RR^d)$ to $L_{p'}(\RR^d_\theta)$.
\end{proposition}

By Proposition \ref{prop:iso2} and \cite[Proposition 1.12]{GJP2021}, we have the following continuous inclusion:
$$\calS(\RR^d_\theta)\hookrightarrow L_p(\RR^d_\theta),
\quad 1\leq p\leq\infty.$$

\bigskip

\section{Sobolev spaces, commutators and Schwartz spaces}\label{sec:qtd}
In this section, following \cite{MSX2020}, we first recall the definition of the homogeneous Sobolev spaces $\dW^1_p(\RR^d_\theta)$ in terms of tempered distributions. Then we define the commutator $[(-\Delta_\theta)^{\epsilon/2},M_x]$ for tempered distributions by duality, which can be transferred to $[(-\Delta)^{\epsilon/2},x]$ on $L_2(\RR^d)$ under the unitary identification $V$. We also develop smooth functional calculus on $\calS(\RR^d_\theta)$, showing that $F(x)\in\calS(\RR^d_\theta)$ whenever $F(0)=0$ and $F$ is smooth near the spectrum of a self-adjoint Schwartz element $x$. These results will be used in the proof of asymptotics of singular values.

\subsection{Sobolev spaces}\label{ssec:SS1}
To define Sobolev spaces on quantum Euclidean spaces, we need to use the notion of tempered distributions. Following \cite{GJP2021}, we present tempered distributions on quantum Euclidean spaces; see \cite{GJP2021} and \cite{MSX2020}.
\begin{definition}
The space $\calS'(\RR^d_\theta)$ consists of all continuous linear functionals on $\calS(\RR^d_\theta)$, which is called the space of tempered distributions on quantum Euclidean spaces. A symbol $x\in\calS'(\RR^d_\theta)$ is called a tempered distribution or simply called a distribution.
\end{definition}

In the classical setting, $\langle T,\phi\rangle$ usually denotes the pairing of $T \in \calS'(\RR^d)$ and $\phi \in \calS(\RR^d)$. As in the classical setting, the pairing of $x \in \calS'(\RR^d_\theta)$ and $y \in \calS(\RR^d_\theta)$ is denoted by $\langle x,y\rangle$. As in \cite{MSX2020}, every $x\in L_1(\RR^d_\theta)+L_\infty(\RR^d_\theta)$ defines a tempered distribution through
$$
\langle x,a\rangle:=\tau_\theta(xa),
\quad a\in\calS(\RR^d_\theta).
$$
Indeed, if $x=x_1+x_\infty$ with $x_1\in L_1(\RR_\theta^d)$ and $x_\infty\in L_\infty(\RR_\theta^d)$, then H\"older's inequality implies
\begin{align*}
|\tau_\theta(xa)|&\leq \|x_1\|_{L_1(\RR_\theta^d)}\|a\|_{L_\infty(\RR_\theta^d)}+\|x_\infty\|_{L_\infty(\RR_\theta^d)}\|a\|_{L_1(\RR_\theta^d)}.
\end{align*}
The continuous inclusions $\mathcal{S}(\RR_\theta^d)\hookrightarrow L_1(\RR_\theta^d)\bcap L_\infty(\RR_\theta^d)$ show that this functional is continuous.

Given $x\in L_\infty(\RR^d_\theta)$, whenever $[\calD_k,x]$ admits a bounded extension to $L_2(\RR^d)$, denote  by $\partial_kx$ this extension, which is called the partial derivatives of $x$. Then for $x=U_\theta(f)\in\calS(\RR^d_\theta)$, one has $$\partial_kx=U_\theta(\xi_kf(\xi))\in\calS(\RR^d_\theta).$$
For $\alpha\in\NN_0^d$, define mixed derivatives by successive applications of these derivations, whenever the intermediate commutators admit bounded extensions:
$$\partial^\alpha x
=\operatorname{ad}(\calD_1)^{\alpha_1}\cdots
\operatorname{ad}(\calD_d)^{\alpha_d}(x),
\quad \operatorname{ad}(\calD_j)(y)=[\calD_j,y],
\quad \operatorname{ad}(\calD_j)^0(y)=y.$$
For $x\in L_\infty(\RR^d_\theta)$, it is called smooth if $\partial^\alpha(x)$ is bounded for all $\alpha\in\NN_0^d$. In particular,
$$\partial^\alpha U_\theta(f)=U_\theta(\xi^\alpha f(\xi)),
\quad f\in\calS(\RR^d).$$

Now we introduce partial derivatives of distributions. Given any multi-index $\alpha\in\mathbb{N}_0^d$ and $T\in\mathcal{S}'(\mathbb{R}_\theta^d)$, the distributional derivative $\partial^\alpha T\in \mathcal{S}'(\mathbb{R}_\theta^d)$ is defined by
$$\langle\partial^\alpha T,\phi\rangle= (-1)^{|\alpha|}\langle T,\partial^\alpha \phi\rangle, \quad \forall \phi\in\mathcal{S}(\mathbb{R}_\theta^d).$$
In particular, the partial derivatives of bounded elements arise as a special case. Let $x\in L_\infty(\RR^d_\theta)$ viewed as the tempered distribution via $\langle x,a\rangle=\tau_\theta(xa)$ for $a\in\calS(\RR^d_\theta)$. Whenever $\partial^\alpha x$ admits a bounded extension to $L_2(\RR^d)$, this bounded element coincides with the distributional derivative. Indeed, for $1\le j\le d$,
\begin{align*}
\tau_\theta((\partial_j x)\phi)
&=\tau_\theta([\calD_j,x]\phi)\\
&=\tau_\theta(x\phi\calD_j)-\tau_\theta(x\calD_j\phi)\\
&=-\tau_\theta(x[\calD_j,\phi])
=-\tau_\theta(x \partial_j\phi).
\end{align*}
Applying this identity once for each factor of $\partial^\alpha=\partial_1^{\alpha_1}\cdots\partial_d^{\alpha_d}$ yields directly
$$\tau_\theta((\partial^\alpha x)\phi)=(-1)^{|\alpha|}\tau_\theta(x \partial^\alpha\phi)=\langle\partial^\alpha x,\phi\rangle.$$
Hence, the definitions of partial derivatives via commutators and distributions agree for every multi-index $\alpha$, if the former partial derivatives exist. For a detailed account of the above, we refer the reader to \cite{LSZ2020,MSX2020}.

Given $T\in\calS'(\RR^d_\theta)$ and $1< p<\infty$, by using the trace-pairing embedding in \cite{MSX2020}, we say $T\in L_p(\RR^d_\theta)$ if there exists $y\in L_p(\RR^d_\theta)$ such that
$$\langle T,a\rangle=\tau_\theta(ya),\quad a\in\calS(\RR^d_\theta).$$
If such $y$ exists, then it is unique as $\calS(\RR^d_\theta)$ is dense in $L_{p'}(\RR^d_\theta)$. 
\iffalse
Indeed, let $y_1,y_2\in L_p(\RR^d_\theta)$ both satisfy
$$\tau_\theta(y_1a)=\tau_\theta(y_2a)=\langle T,a\rangle,\quad a\in\calS(\RR^d_\theta).$$
Since $\calS(\RR^d_\theta)$ is dense in $L_{p'}(\RR^d_\theta)$ and the map $a\mapsto\tau_\theta((y_1-y_2)a)$ is continuous on $L_{p'}(\RR^d_\theta)$ by H\"older's inequality, it follows that
$$\tau_\theta((y_1-y_2)a)=0,\quad a\in L_{p'}(\RR^d_\theta).$$
Let $y_1-y_2=u|y_1-y_2|$ be the polar decomposition and take $a=|y_1-y_2|^{p-1}u^*\in L_{p'}(\RR^d_\theta)$ (for $p=1$, take $a=u^*$). Then
\begin{align*}
0=\tau_\theta((y_1-y_2)a)
&=\tau_\theta\big(u|y_1-y_2|\cdot|y_1-y_2|^{p-1}u^*\big)	=\tau_\theta(u|y_1-y_2|^p u^*)\\
&=\tau_\theta(u^*u|y_1-y_2|^p)
=\tau_\theta(|y_1-y_2|^p)
=\|y_1-y_2\|_{L_p(\RR^d_\theta)}^p,
\end{align*}
which implies $y_1=y_2$. 
\fi
In this case, we define
$$\|T\|_{L_p(\RR^d_\theta)}=\tau_\theta(|y|^p)^{1/p},\quad 1< p<\infty.$$
If such $y$ does not exist, we set $\|T\|_{L_p(\RR^d_\theta)}=+\infty$. All statements that a distributional derivative belongs to $L_p(\RR^d_\theta)$ and all $L_p(\RR^d_\theta)$ norms of such derivatives in the sequel are understood in this sense. 

%Now we define homogeneous and inhomogeneous Sobolev spaces on quantum Euclidean spaces.
McDonald, Sukochev and Xiong first gave the following distributional definitions of both homogeneous and inhomogeneous Sobolev spaces in \cite{MSX2020}.

\begin{definition}
Let $1\le p<\infty$ and $m\in\mathbb N$. The homogeneous Sobolev space $\dot{W}_p^m(\mathbb{R}_\theta^d)$ consists of those $x\in\mathcal{S}'(\mathbb{R}_\theta^d)$ such that every partial derivative of $x$ of order $m$ belongs to $L_p(\mathbb{R}_\theta^d)$, endowed with the following seminorm:
$$\|x\|_{\dot{W}_p^m(\RR^d_\theta)}
=\sum_{|\alpha|=m}\|\partial^\alpha x\|_{L_p(\RR^d_\theta)}.$$
The inhomogeneous quantum Sobolev space is defined as
$$W^m_p(\RR^d_\theta)
=\big\{x\in L_p(\RR^d_\theta):
\partial^\beta x\in L_p(\RR^d_\theta),\ |\beta|\leq m\big\},$$
endowed with the following norm
$$\|x\|_{W^m_p(\RR^d_\theta)}=\sum_{|\beta|\leq m}\|\partial^\beta x\|_{L_p(\RR^d_\theta)}.$$
\end{definition}

Recall the following approximation proposition from \cite[Proposition 4.9]{MSX2023}.
\begin{proposition}
\label{prop:hd1}
Let $m\in\NN$ and $1<p<\infty$. For every $x\in\dW^m_p(\RR^d_\theta)$, there exists a sequence $(x_n)_{n\geq1}\subseteq\calS(\RR^d_\theta)$ such that
\begin{equation}\label{eq:hd1}
\lim_{n\to\infty}\|x_n-x\|_{\dW^m_p(\RR^d_\theta)}=0.
\end{equation}
\end{proposition}

We also need the following inhomogeneous approximation proposition, which follows from \cite[Proposition 3.14]{MSX2020}.

\begin{proposition}\label{prop:id1}
Let $m\in\NN_0$ and $1\leq p<\infty$. For every $x\in W^m_p(\RR^d_\theta)$, there exists a sequence $(x_n)_{n\geq1}\subseteq\calS(\RR^d_\theta)$ such that
\begin{equation}\label{eq:id1}
\lim_{n\to\infty}\|x_n-x\|_{W^m_p(\RR^d_\theta)}=0.
\end{equation}
\end{proposition}

The following lemma concerns the Fatou property for the homogeneous Sobolev spaces $\dW^1_p(\RR^d_\theta)$, which is quite helpful in the proof of Theorem \ref{thm:m1}.

\begin{lemma}\label{lem:wgrad}
Let $1<p<\infty$ and $(x_n)_{n\geq 1}\subseteq\dW^1_p(\RR^d_\theta)$. Suppose that $x_n\to x$ in $\calS'(\RR^d_\theta)$. If $\liminf_{n\to\infty}\|x_n\|_{\dW^1_p(\RR^d_\theta)}<\infty$, then $x\in\dW^1_p(\RR^d_\theta)$ and
$$
\|x\|_{\dW^1_p(\RR^d_\theta)}
\leq\liminf_{n\to\infty}\|x_n\|_{\dW^1_p(\RR^d_\theta)}.
$$
\end{lemma}
\begin{proof}
We pass to a subsequence along which the Sobolev seminorms converge to their lower limit. By the reflexivity of $L_p(\RR^d_\theta)$, without loss of generality, we assume that $\partial_jx_n\rightharpoonup v_j$ with respect to the $w^*$ topology in $L_p(\RR^d_\theta)$ for $1\leq j\leq d$. For $a\in\calS(\RR^d_\theta)\subseteq L_{p'}(\RR^d_\theta)$,
\begin{align*}
\tau_\theta(v_ja)
=&\lim_{n\to\infty}\tau_\theta((\partial_jx_n)a)
=-\lim_{n\to\infty}\langle x_n,\partial_ja\rangle\\
=&-\langle x,\partial_ja\rangle
=\langle\partial_jx,a\rangle.
\end{align*}
Thus $v_j$ represents $\partial_jx$. For every $j=1,\ldots,d$ and $a\in L_{p'}(\RR^d_\theta)$ with $\|a\|_{L_{p'}(\RR^d_\theta)}\leq1$, the $w^*$ convergence and H\"older's inequality yield
\begin{align*}
|\tau_\theta((\partial_jx)a)|
=&\lim_{n\to\infty}|\tau_\theta((\partial_jx_n)a)|\\
\leq&\liminf_{n\to\infty}\|\partial_jx_n\|_{L_p(\RR^d_\theta)}\|a\|_{L_{p'}(\RR^d_\theta)}\\
\leq&\liminf_{n\to\infty}\|\partial_jx_n\|_{L_p(\RR^d_\theta)}.
\end{align*}
Taking the supremum over all such $a$, we obtain
\begin{align*}
\|\partial_jx\|_{L_p(\RR^d_\theta)}
\leq\liminf_{n\to\infty}\|\partial_jx_n\|_{L_p(\RR^d_\theta)}.
\end{align*}
Summing over $j$, we conclude that $x\in\dW^1_p(\RR^d_\theta)$ and
\begin{align*}
\|x\|_{\dW^1_p(\RR^d_\theta)}
\leq\liminf_{n\to\infty}\|x_n\|_{\dW^1_p(\RR^d_\theta)},
\end{align*}
as desired.
\end{proof}

%Next we prove that the spherical gradient norm and the homogeneous Sobolev seminorm are equivalent. 
Next, we compare the $L_p$-norm over the sphere of $\sum_{j=1}^d s_j\partial_jx$ with the homogeneous Sobolev seminorms. Let $dS$ denote the surface measure on $\BS^{d-1}$.

\begin{lemma}\label{lem:ge1}
For $1<p<\infty$ and $x\in\dW^1_p(\RR^d_\theta)$,
$$\Big\|\sum_{j=1}^d(\partial_jx)s_j\Big\|_{L_p(L_\infty(\RR^d_\theta)\overline\otimes L_\infty(\BS^{d-1}))}\asymp_{d,p}\|x\|_{\dW^1_p(\RR^d_\theta)}.$$
\end{lemma}

\begin{proof}
By the triangle inequality, we have
\begin{align*}
\|\sum_{j=1}^d(\partial_jx)s_j\Big\|_{L_p(L_\infty(\RR^d_\theta)\overline\otimes L_\infty(\BS^{d-1}))}&\leq\sum_{j=1}^d
\|\partial_jx\|_{L_p(\RR^d_\theta)}\|s_j\|_{L_p(\BS^{d-1})}\\
&\leq|\BS^{d-1}|^{1/p}\|x\|_{\dW^1_p(\RR^d_\theta)}.
\end{align*}
For the converse direction, note that
$$\int_{\BS^{d-1}}s_ks_j dS
 =\frac{|\BS^{d-1}|}{d}\delta_{jk},\quad k,j\in\{1,\ldots,d\}.$$
Therefore, one has
$$ \partial_kx=\frac d{|\BS^{d-1}|}\sum_{j=1}^d
 \int_{\BS^{d-1}}(\partial_jx)s_ks_j dS(s).$$
By Minkowski's and H\"older's inequalities,
\begin{align*}
 \|\partial_kx\|_{L_p(\RR^d_\theta)}
 \leq&\frac d{|\BS^{d-1}|}
 \int_{\BS^{d-1}}|s_k|\Big\|\sum_{j=1}^d(\partial_jx)s_j\Big\|_{L_p(\RR^d_\theta)} dS\\
 \leq&\frac d{|\BS^{d-1}|}
 \|s_k\|_{L_{p'}(\BS^{d-1})}\Big\|\sum_{j=1}^d(\partial_jx)s_j\Big\|_{L_p(L_\infty(\RR^d_\theta)\overline\otimes L_\infty(\BS^{d-1}))}(x).
\end{align*}
This completes the proof.
\end{proof}

\subsection{The definition of $[(-\Delta_\theta)^{\frac{\epsilon}{2}},M_x]$}\label{ssec:dcom1}
We first introduce the definition of $(-\Delta_\theta)^{\frac{\epsilon}{2}}$ via functional calculus. Define $\nabla_\theta=(-\ii\partial_1,\ldots,-\ii\partial_d)$. Let $E_{\nabla_\theta}$ be the joint spectral measure of $(\partial_1,\ldots,\partial_d)$. For a Borel function $g:\RR^d\to\CC$, define
$$g(\ii\nabla_\theta)=g(\partial_1,\ldots,\partial_d)
=\int_{\RR^d}g(t) dE_{\nabla_\theta}(t),$$
with the domain
$$\operatorname{Dom}(g(\ii\nabla_\theta))
=\Big\{x\in L_2(\RR^d_\theta):
\int_{\RR^d}|g(t)|^2 d\langle E_{\nabla_\theta}(t)x,x\rangle<\infty\Big\}.$$
We have the following proposition.

\begin{proposition}\label{prop:LLt1}
Let $g:\RR^d\to\CC$ be a Borel function. Then
$$\operatorname{Dom}(g(\ii\nabla_\theta))
=V\operatorname{Dom}(g(\ii\nabla))$$
and
$$V^*g(\ii\nabla_\theta)Vf=g(\ii\nabla)f,\quad
f\in\operatorname{Dom}(g(\ii\nabla)).$$
\end{proposition}
\begin{proof}
Since $\partial_j=V\calD_jV^*$, the joint spectral measures satisfy $E_{\nabla_\theta}(B)=VE_\nabla(B)V^*$ for every Borel set $B$. Thus for $f\in L_2(\RR^d)$,
$$\int_{\RR^d}|g(t)|^2 d\langle E_{\nabla_\theta}(t)Vf,Vf\rangle
=\int_{\RR^d}|g(t)f(t)|^2 dt.$$
This proves the first identity involving domain. On that domain, spectral integration yields
$$V^*g(\ii\nabla_\theta)Vf
=\int_{\RR^d}g(t) dE_\nabla(t)f=g(\ii\nabla)f.$$
This completes the proof.
\end{proof}

By Proposition \ref{prop:LLt1}, for $\epsilon\in(-\frac{d}{2},1)$ and for any $f,|\cdot|^\epsilon f\in L_2(\RR^d)$, one has
$$(-\Delta_\theta)^{\frac{\epsilon}{2}}U_\theta(f)
=U_\theta(|\cdot|^\epsilon f).$$
In what follows, we give the definition of the commutator $[(-\Delta_\theta)^{\frac{\epsilon}{2}}, M_x ]$ for $x\in\calS'(\RR^d_\theta)$ via tempered distributions on quantum Euclidean spaces. We would like to point out that fractional powers of the Laplacians need not preserve the Schwartz spaces. Thus, for tempered distributions, the two terms in the commutator require a definition by duality. The following proposition shows that the products in the corresponding pairings belong to $\calS(\RR^d_\theta)$.
%and also proves compatibility with approximation in the homogeneous Sobolev seminorm.

\begin{proposition}\label{prop:dw1}
Let $-\frac d2<\epsilon\neq 0<1$.
\begin{enumerate}[label=\textup{(\roman*)}]
\item For every $a,b\in\calS(\RR^d_\theta)$, one has
$$a\big((-\Delta_\theta)^{\frac{\epsilon}{2}}b\big)^*,\big((-\Delta_\theta)^{\frac{\epsilon}{2}}a\big)b^*\in\calS(\RR^d_\theta),$$
and the two maps
$$(a,b)\mapsto a\big((-\Delta_\theta)^{\frac{\epsilon}{2}}b\big)^* \quad \text{or} \quad \big((-\Delta_\theta)^{\frac{\epsilon}{2}}a\big)b^*$$
are both continuous from $\calS(\RR^d_\theta)\times\calS(\RR^d_\theta)$ to $\calS(\RR^d_\theta)$.

\item Let $x\in\calS'(\RR^d_\theta)$. For any $a\in\calS(\RR^d_\theta)$, there exist two unique distributions, denoted by $(-\Delta_\theta)^{\frac{\epsilon}{2}}M_xa$ and $M_x(-\Delta_\theta)^{\frac{\epsilon}{2}}a$, such that for all $b\in\calS(\RR^d_\theta)$ one has
\begin{equation}\label{eq:drb1}
\big\langle (-\Delta_\theta)^{\frac{\epsilon}{2}} M_xa,b^*\big\rangle
=\big\langle x,a\big((-\Delta_\theta)^{\frac{\epsilon}{2}}b\big)^*\big\rangle
\end{equation}
and
\begin{equation}\label{eq:dlb1}
\big\langle M_x(-\Delta_\theta)^{\frac{\epsilon}{2}} a,b^*\big\rangle
=\big\langle x,\big((-\Delta_\theta)^{\frac{\epsilon}{2}}a\big)b^*\big\rangle.
\end{equation}

\item If $x_n,x\in L_\infty(\RR^d_\theta)$ and $x_n\to x$ with respect to the $w^*$-topology in $L_\infty(\RR^d_\theta)$, then for every fixed $a\in\calS(\RR^d_\theta)$,
$$(-\Delta_\theta)^{\frac{\epsilon}{2}} M_{x_n}a\to (-\Delta_\theta)^{\frac{\epsilon}{2}} M_xa$$
and
$$M_{x_n}(-\Delta_\theta)^{\frac{\epsilon}{2}} a\to M_x(-\Delta_\theta)^{\frac{\epsilon}{2}} a$$
in $\calS'(\RR^d_\theta)$.

\item Let $1<p<\infty$ and $x\in\dW^1_p(\RR^d_\theta)$. Suppose that a sequence $(x_n)_{n\geq 1}\subseteq\calS(\RR^d_\theta)$ satisfies $\|x_n-x\|_{\dW^1_p(\RR^d_\theta)}\to0$. Then for every $a\in\calS(\RR^d_\theta)$,
\begin{equation}\label{eq:dgc1}
[(-\Delta_\theta)^{\frac{\epsilon}{2}}, M_{x_n}]a\to
[(-\Delta_\theta)^{\frac{\epsilon}{2}}, M_x]a
\end{equation}
in $\calS'(\RR^d_\theta)$.
\end{enumerate}
\end{proposition}

\begin{proof}
We first verify the mapping property needed to define the products above. For $a=U_\theta(f)$ with $f\in\calS(\RR^d)$, recall that $(-\Delta_\theta)^{\frac{\epsilon}{2}} a=U_\theta(| \cdot |^\epsilon f)$. Choose $N>d+\max\{\epsilon,0\}$. Since $\epsilon>-\frac{d}{2}>-d$,
\begin{align*}
\|| \cdot |^\epsilon f\|_{L_1(\RR^d)}
\leq&
\|f\|_{L_\infty(\RR^d)}\int_{|t|\leq1}|t|^\epsilon dt
+
\sup_t(1+|t|)^N|f(t)|
\int_{|t|>1}|t|^\epsilon(1+|t|)^{-N} dt\\
\lesssim&_{d,\epsilon,N}
\Big(\|f\|_{L_\infty(\RR^d)}+
\sup_t(1+|t|)^N|f(t)|\Big).
\end{align*}
By Proposition \ref{prop:iso2}, 
$$(-\Delta_\theta)^{\frac{\epsilon}{2}}:\calS(\RR^d_\theta)\to
L_\infty(\RR^d_\theta)$$
is continuous. Then
$$a\big((-\Delta_\theta)^{\frac{\epsilon}{2}}b\big)^*,\big((-\Delta_\theta)^{\frac{\epsilon}{2}}a\big)b^*\in L_1(\RR^d_\theta)\bcap L_\infty(\RR^d_\theta),\quad a,b\in\calS(\RR^d_\theta).$$
We show (i), (ii), (iii) and (iv) respectively.

\

\noindent (i): Write $a=U_\theta(f)$ and $b=U_\theta(g)$, where $f,g\in \calS(\RR^d)$. Let
$$I_1(s)=\int_{\RR^d}|t|^\epsilon
e^{-\frac{\ii}{2}\langle s,\theta t\rangle}
f(s-t)\overline{g(-t)} dt$$
and
$$I_2(s)=\int_{\RR^d}|t|^\epsilon
e^{\frac{\ii}{2}\langle s,\theta t\rangle}
f(t)\overline{g(t-s)} dt.$$
For every multi-index $\beta$, taking differentiation yields
$$\partial_s^\beta I_1(s)
=\sum_{\gamma\leq\beta}
\binom{\beta}{\gamma}
\Big(-\frac{\ii}{2}\Big)^{|\gamma|}
\int_{\RR^d}|t|^\epsilon(\theta t)^\gamma
e^{-\frac{\ii}{2}\langle s,\theta t\rangle}
(\partial^{\beta-\gamma}f)(s-t)
\overline{g(-t)} dt$$
and
$$\partial_s^\beta I_2(s)
=\sum_{\gamma\leq\beta}
\binom{\beta}{\gamma}
\Big(\frac{\ii}{2}\Big)^{|\gamma|}
(-1)^{|\beta-\gamma|}
\int_{\RR^d}|t|^\epsilon(\theta t)^\gamma
e^{\frac{\ii}{2}\langle s,\theta t\rangle}
f(t)(\partial^{\beta-\gamma}\overline g)(t-s) dt.$$
For $N\in\NN_0$, choose an integer $L>N+|\beta|+d+|\epsilon|$. Since
$$1+|s|\leq(1+|s-t|)(1+|t|),
\quad\int_{\RR^d}|t|^\epsilon(1+|t|)^{N+|\beta|-L} dt<\infty,$$
one has
\begin{align*}
&\sum_{\ell=1}^2\sup_{s\in\RR^d}
(1+|s|)^N|\partial_s^\beta I_\ell(s)|\\
\lesssim&_{N,\beta,d,\theta}
\Big(\max_{|\kappa|\leq|\beta|}\sup_{u\in\RR^d}
(1+|u|)^L|\partial^\kappa f(u)|\Big)
\Big(\max_{|\kappa|\leq|\beta|}\sup_{u\in\RR^d}
(1+|u|)^L|\partial^\kappa g(u)|\Big)\\
&\times\int_{\RR^d}|t|^\epsilon(1+|t|)^{N+|\beta|-L} dt.
\end{align*}
This implies that $I_1,I_2\in\calS(\RR^d)$ and that both maps $(f,g)\mapsto I_1,I_2$ are continuous from $\calS(\RR^d)\times\calS(\RR^d)$ to $\calS(\RR^d)$. By \eqref{eq:inveru1}, $I_1$ and $I_2$ are the Fourier transform of $a\big((-\Delta_\theta)^{\frac{\epsilon}{2}}b\big)^*$ and $\big((-\Delta_\theta)^{\frac{\epsilon}{2}}a\big)b^*$, respectively. Since $U_\theta$ is a topological isomorphism, the maps $\varphi_1(a,b)=a\big((-\Delta_\theta)^{\frac{\epsilon}{2}}b\big)^*$ and $\varphi_2(a,b)=\big((-\Delta_\theta)^{\frac{\epsilon}{2}}a\big)b^*$ are continuous from $\calS(\RR^d_\theta)\times\calS(\RR^d_\theta)$ to $\calS(\RR^d_\theta)$. This proves (i).

\

\noindent (ii): Fix $a\in\calS(\RR^d_\theta)$. By (i), the maps $b\mapsto\langle x,\varphi_j(a,b)\rangle$ for $j\in \{1,2\}$ are continuous and conjugate-linear. Therefore, they define unique distributions $(-\Delta_\theta)^{\frac{\epsilon}{2}}M_xa$ and $M_x(-\Delta_\theta)^{\frac{\epsilon}{2}}a$ through their pairings with $b^*$, respectively. This proves (ii). 

\

\noindent (iii): This just follows immediately from (i) and (ii). 

\

\noindent (iv): If $\Phi=U_\theta(F)\in\calS(\RR^d_\theta)$ and $F(0)=0$, then there exist $F_j\in \calS(\RR^d)$ such that
$$F(s)=\sum_{j=1}^ds_jF_j(s).$$
To this end, choose $\rho\in C_c^\infty(\RR^d)$ with $\rho=1$ near the origin and put
$$F_j(s)=\rho(s)\int_0^1(\partial_j F)(rs) dr +(1-\rho(s))\frac{s_j}{|s|^2}F(s).$$
Then $F_j\in\calS(\RR^d)$ and $\sum_{j=1}^ds_jF_j(s)=F(s)$ because $F(0)=0$.
Hence, with $\eta_j=U_\theta(F_j)$,
\begin{equation}\label{eq:div0}
\Phi=\sum_{j=1}^d\partial_j\eta_j,
\quad \eta_j\in\calS(\RR^d_\theta).
\end{equation}
Let $a,b\in\calS(\RR^d_\theta)$ with $a=U_\theta(f)$ and $b=U_\theta(g)$. Then
$$I_1(0)=\int_{\RR^d}|t|^\epsilon f(-t)\overline{g(-t)} dt
=\int_{\RR^d}|u|^\epsilon f(u)\overline{g(u)} du
=I_2(0).$$
Therefore, \eqref{eq:div0} can be applied to
$$\Phi_{\epsilon;a,b}
=a\big((-\Delta_\theta)^{\frac{\epsilon}{2}}b\big)^*-\big((-\Delta_\theta)^{\frac{\epsilon}{2}}a\big)b^*.$$
By (ii),
\begin{align*}
&\Big|\big\langle
[(-\Delta_\theta)^{\frac{\epsilon}{2}}, M_{x_n}]a,b^*
\big\rangle-\big\langle
[(-\Delta_\theta)^{\frac{\epsilon}{2}}, M_x]a,b^*
\big\rangle\Big|\\
=&\Big|\langle x_n-x,\Phi_{\epsilon;a,b}\rangle\Big|
=\Big|\sum_{j=1}^d\langle\partial_jx_n-\partial_jx,-\eta_j\rangle\Big|\\
\leq&\sum_{j=1}^d
\|\partial_jx_n-\partial_jx\|_{L_p(\RR^d_\theta)}
\|\eta_j\|_{L_{p'}(\RR^d_\theta)}
\to0.
\end{align*}
This proves (iv).
\end{proof}

Indeed, Proposition \ref{prop:dw1} implies that for $x\in\calS(\RR^d_\theta)$, the operator $[(-\Delta_\theta)^{\frac{\epsilon}{2}}, M_x]$ defined on $\calS(\RR^d_\theta)$ is precisely the standard commutator, namely
\begin{align}\label{eq:dmx1}
[(-\Delta_\theta)^{\frac{\epsilon}{2}}, M_x]=(-\Delta_\theta)^{\frac{\epsilon}{2}} M_x-M_x(-\Delta_\theta)^{\frac{\epsilon}{2}},
\end{align}
with $M_x$ the left multiplication operator given by $M_xy=xy$. For $f,g\in \calS(\RR^d)$, we will use the multiplication identity (see \cite{MSX2020})
\begin{equation}\label{eq:mul1}
U_\theta(f) U_\theta(g) = U_\theta(f *_\theta g),
\end{equation}
with the $\theta$-convolution given by
\begin{equation}\label{eq:con1}
f *_\theta g(t) = \int_{\RR^d} e^{-\frac{\ii}{2} \langle t, \theta s\rangle} f(t- s) g(s) ds.
\end{equation}
When $x\in\calS(\RR^d_\theta)$, the following proposition implies that the operator $M_x$ is unitarily equivalent to $x$, the latter being viewed as a bounded linear operator on $L_2(\RR^d)$.

\begin{proposition}\label{prop:umu1}
Let $x=U_\theta(h)\in\calS(\RR^d_\theta)$. Then for any $f\in\calS(\RR^d)$,
\begin{equation}\label{eq:vmv1}
(V^*M_x V)f = h *_\theta f=xf.
\end{equation}
\end{proposition}
\begin{proof}
Let $f\in\calS(\RR^d)$. By \eqref{eq:mul1} and the definition of $V$, we have $V^*M_xVf=h*_\theta f$. Using \eqref{eq:wun1} and \eqref{eq:defufh1}, we obtain
\begin{align*}
(xf)(s)
=&\int_{\RR^d}h(t)e^{-\frac{\ii}{2}\langle t,\theta s\rangle}f(s-t) dt\\
=&\int_{\RR^d}e^{-\frac{\ii}{2}\langle s,\theta u\rangle}h(s-u)f(u) du\\
=&(h*_\theta f)(s).
\end{align*}
This proves \eqref{eq:vmv1}.
\end{proof}

%We shall use the following invariance of Schatten norms under the unitary map $V$.

%\begin{lemma}\label{lem:UTU1}
%For $0<p<\infty$ and $T\in\calK(L_2(\RR^d_\theta))$, we have
%$$\|V^*TV\|_{\calL_{p}(L_2(\RR^d))}=\|T\|_{\calL_{p}(L_2(\RR^d_\theta))}$$
%and
%$$\|V^*TV\|_{\calL_{p,\infty}(L_2(\RR^d))}=\|T\|_{\calL_{p,\infty}%(L_2(\RR^d_\theta))}.$$
%\end{lemma}
%\begin{proof}
%Since unitary conjugation preserves rank and operator norm, we have
%$$\mu(t,V^*TV)
%=\inf_{\rank A\leq t}\|V^*(T-A)V\|_{\calB(L_2(\RR^d))}
%=\mu(t,T),\quad t\geq0.$$
%Then both assertions follow from the definitions of the Schatten quasi-norms.
%\end{proof}

By \eqref{eq:dmx1} and Proposition \ref{prop:umu1}, we obtain the following proposition.

%By Proposition \ref{prop:LLt1}, Proposition \ref{prop:umu1}, and Lemma \ref{lem:UTU1}, we obtain the following proposition.

\begin{proposition}\label{prop:tr1}
Let $-\frac d2<\epsilon\neq 0<1$ and $x\in\calS(\RR^d_\theta)$. Then
$$[(-\Delta)^{\frac{\epsilon}{2}}, x ]=V^*[(-\Delta_\theta)^{\frac{\epsilon}{2}}, M_x ]V.$$
%\begin{description}
%\item[(i)] $[(-\Delta)^{\frac{\epsilon}{2}}, x ]=V^*[(-\Delta_\theta)^{\frac{\epsilon}{2}}, M_x ]V$.
%\item[(ii)] Either commutator has a bounded extension in its indicated weak Schatten ideal if and only if the other does; for these extensions,
%$$\|[(-\Delta_\theta)^{\frac{\epsilon}{2}},M_x]\|_{\calL_{\frac{d}{1-\epsilon},\infty}(L_2(\RR^d_\theta))}
%=\|[(-\Delta)^{\frac{\epsilon}{2}},x]\|_{\calL_{\frac{d}{1-\epsilon},\infty}(L_2(\RR^d))}.$$
%\end{description}
\end{proposition}

For any arbitrary $x\in\calS'(\RR^d_\theta)$, we will define the transferred commutators by the distributional pairing in \eqref{eq:cform0}. Proposition \ref{prop:dw1} guarantees that this definition is well-defined and
agrees with the preceding commutator when $x\in\calS(\RR^d_\theta)$.

Let $-\frac d2<\epsilon\neq 0<1$. For $a,b\in\calS(\RR^d_\theta)$, define
\begin{equation}\label{eq:Phi}
\Phi_{\epsilon;a,b}
=a\big((-\Delta_\theta)^{\frac{\epsilon}{2}}b\big)^*-\big((-\Delta_\theta)^{\frac{\epsilon}{2}}a\big)b^*.
\end{equation}
By Proposition \ref{prop:dw1} (i), $\Phi_{\epsilon;a,b}\in\calS(\RR^d_\theta)$. Using the bilinear pairing between $\calS'(\RR^d)$ and $\calS(\RR^d)$, define the corresponding distributional commutator for every $x\in\calS'(\RR^d_\theta)$ by
\begin{equation}\label{eq:cform0}
\big\langle[(-\Delta)^{\frac{\epsilon}{2}},x]f,\overline g\big\rangle
=\langle x,\Phi_{\epsilon;Vf,Vg}\rangle,\quad f,g\in\calS(\RR^d),
\end{equation}
where
\begin{equation}\label{eq:cform1}
\Phi_{\epsilon;Vf,Vg}
=(2\pi)^{-d}\Big(
U_\theta(f)U_\theta(|\cdot|^\epsilon g)^*
-U_\theta(|\cdot|^\epsilon f)U_\theta(g)^*\Big).
\end{equation}
For $x\in\calS(\RR^d_\theta)$, the definition given by \eqref{eq:cform0} agrees with the definitions established in Proposition \ref{prop:tr1} and Proposition \ref{prop:dw1}. In addition, by Proposition \ref{prop:dw1} (ii) and \eqref{eq:cform0}, for every $x\in\calS'(\RR^d_\theta)$ and $f,g\in\calS(\RR^d)$,
\begin{align*}
\big\langle[(-\Delta_\theta)^{\frac{\epsilon}{2}},M_x]Vf,(Vg)^*\big\rangle
=&\langle x,\Phi_{\epsilon;Vf,Vg}\rangle\\
=&\big\langle[(-\Delta)^{\frac{\epsilon}{2}},x]f,\overline g\big\rangle.
\end{align*}
Since $V\calS(\RR^d)=\calS(\RR^d_\theta)$ is dense in $L_2(\RR^d_\theta)$ and $V$ is unitary, either commutator has a bounded extension if and only if the other one does. For these extensions, the preceding equality in Proposition \ref{prop:tr1} implies
\begin{equation}\label{eq:transfer7}
[(-\Delta_\theta)^{\frac{\epsilon}{2}},M_x]
=V[(-\Delta)^{\frac{\epsilon}{2}},x]V^*.
\end{equation}
In paticular, if $x\in \dW^1_{\frac{d}{1-\epsilon}}(\RR^d_\theta)$, $(x_n)_{n\geq 1}\subseteq\calS(\RR^d_\theta)$ and $\|x_n-x\|_{\dW^1_{\frac{d}{1-\epsilon}}(\RR^d_\theta)}\to0$, then by Proposition \ref{prop:dw1} and Proposition \ref{prop:tr1}, we have
\begin{align*}
\big\langle[(-\Delta)^{\frac{\epsilon}{2}},x]f,\overline g\big\rangle
=&\big\langle[(-\Delta_\theta)^{\frac{\epsilon}{2}},M_x]Vf,(Vg)^*\big\rangle\\
=&\lim_{n\to\infty}\langle[(-\Delta_\theta)^{\frac{\epsilon}{2}},M_{x_n}]Vf, (Vg)^*\rangle\\
=&\lim_{n\to\infty}\big\langle[(-\Delta)^{\frac{\epsilon}{2}},x_n]f,\overline g\big\rangle.
\end{align*}

\subsection{Functional calculus on $\calS(\RR^d_\theta)$}\label{ssec:Fcs1}
We will use the functional calculus of some special functions that sends Schwartz elements to Schwartz elements. In \cite{T2025}, square roots of positive Schwartz elements were used to verify the commutator estimates required for the spectral asymptotic formula. Under the assumption $\det\theta\ne0$, the matrix description of $\calS(\RR^d_\theta)$ ensured that the square root of a positive Schwartz element was again a Schwartz element; see \cite[Proposition 2.13 and Lemma 2.14]{T2025}. In the present setting, applying an arbitrary positive power or taking square root, need not generate a new Schwartz element, even when the symbol is strictly positive. To see this, consider the following strictly positive Schwartz function
$$f(t)=e^{-t^2}\sin^2t+e^{-2t^2},\quad t\in\RR,$$
for which
$$(f^{\frac12})''(n\pi)=1+(4n^2\pi^2-2)e^{-n^2\pi^2}\to1, \quad \text{as} \ n\to\infty;$$
this implies that $f^{1/2}$ is smooth but not Schwartz. Therefore, we use smooth functional calculus to show that $F(x)\in\calS(\RR^d_\theta)$ whenever $F(0)=0$ and $F$ is smooth near the spectrum of $x$, and then construct positive Schwartz elements with Schwartz square roots, which approximate $|x|$ through their integer powers. These two facts are shown in Lemma \ref{lem:fc} and Lemma \ref{lem:root}, respectively.

For $N\in\NN_0$ and $f\in\calS(\RR^d)$, define
$$
q_N(f)=\sum_{|\alpha|+|\beta|\leq N}\|s^\alpha\partial_s^\beta f(s)\|_{L_2(\RR^d)}.
$$
The two families of seminorms $(q_N)$ and $(\mathfrak p_N)$ are equivalent and define the same topology on $\mathcal S(\mathbb R^d)$. This is because the estimate $\mathfrak p_N\lesssim q_{N'}$ for some $N'$ depending on $N$ is a consequence of the Sobolev embedding \cite[Theorem 12.46]{L2017}.

The following lemma follows directly by elementary calculation, and so we include its detailed proof in Appendix \ref{sec:app1}.
\begin{lemma}\label{lem:exp}
Let $x=x^*=U_\theta(f)\in\calS(\RR^d_\theta)$. For every $t\in\RR$, there exists $u_t\in\calS(\RR^d)$ such that
$$
U_\theta(u_t)=e^{\ii tx}-1.
$$
Moreover, for every $N\in\NN_0$,
$$
q_N(u_t)\lesssim_{N,f,\theta,d}|t|(1+|t|)^N.
$$
\end{lemma}

Now we establish the stability of functional calculus on $\calS(\RR^d_\theta)$.

\begin{lemma}\label{lem:fc}
Let $x=x^*=U_\theta(f)\in\calS(\RR^d_\theta)$. Denote by $\sigma(x)$ its spectrum. Let $F$ be smooth on a neighborhood of $\sigma(x)\cup\{0\}$ with $F(0)=0$. Then $F(x)\in\calS(\RR^d_\theta)$.
\end{lemma}

\begin{proof}
Let $O\subset\RR$ be an open neighborhood of $\sigma(x)\cup\{0\}$ on which $F$ is smooth. Choose $\chi\in C_c^\infty(O)$ that is equal to one near $\sigma(x)\cup\{0\}$, and extend $\chi \cdot F$ to the whole $\mathbb{R}$ by zero outside $O$. Denote this extension by $\Phi$. Then
$$\Phi\in C_c^\infty(\RR),\quad
\Phi(0)=0\quad\text{and}\quad
\Phi(x)=F(x), \ \forall x\in O.$$
Recall that the one-dimensional Fourier transform of $\Phi(x)$ is given by
$$\widehat\Phi(t)=\int_\RR e^{-\ii t\lambda}\Phi(\lambda) d\lambda.$$
Let $u_t$ be as in Lemma \ref{lem:exp}. For $R>0$, let
$$g_R=\frac1{2\pi}\int_{-R}^R\widehat\Phi(t)u_t dt.$$
For $R'>R$ and $N\in\NN_0$,
\begin{align*}
q_N(g_{R'}-g_R)
\leq&\frac1{2\pi}\int_{R<|t|\leq R'}|\widehat\Phi(t)|q_N(u_t) dt\\
\leq&\frac{C_{N,f,\theta,d}}{2\pi}\int_{|t|>R}|\widehat\Phi(t)| |t|(1+|t|)^N dt\\
\leq&\frac{C_{N,f,\theta,d}}{2\pi}\sup_{t\in\RR}\big((1+|t|)^{N+3}|\widehat\Phi(t)|\big)\int_{|t|>R}\frac{|t|}{(1+|t|)^3} dt\\
\leq&\frac{C_{N,f,\theta,d}}{2\pi}\sup_{t\in\RR}\big((1+|t|)^{N+3}|\widehat\Phi(t)|\big)\int_{|t|>R}(1+|t|)^{-2} dt\\
=&\frac{C_{N,f,\theta,d}}{\pi(1+R)}\sup_{t\in\RR}\big((1+|t|)^{N+3}|\widehat\Phi(t)|\big)\to0.
\end{align*}
Since the family of the seminorms $(q_N)_{N\in\NN_0}$ defines the topology on $\calS(\RR^d)$, there exists $g\in\calS(\RR^d)$ such that $q_N(g-g_R)\to0$ for every $N\in\NN_0$. Denote this limit by $g=\frac1{2\pi}\int_\RR\widehat\Phi(t)u_t dt$. Moreover, by Lemma \ref{lem:exp},
\begin{align*}
q_N(g)=&\lim_{R\to\infty}q_N(g_R)\\
\leq&\frac1{2\pi}\int_\RR|\widehat\Phi(t)|q_N(u_t) dt\\
\leq&\frac{C_{N,f,\theta,d}}{2\pi}\int_\RR|\widehat\Phi(t)| |t|(1+|t|)^N dt
<\infty.
\end{align*}
By the continuity of $U_\theta$, we obtain
\begin{align*}
U_\theta(g)=&\frac1{2\pi}\int_\RR\widehat\Phi(t)U_\theta(u_t) dt\\
=&\frac1{2\pi}\int_\RR\widehat\Phi(t)(e^{\ii tx}-1) dt\\
=&\Phi(x)-\Phi(0)\\
=&F(x).
\end{align*}
This implies $F(x)\in\calS(\RR^d_\theta)$. This completes the proof.
\end{proof}

By Lemma \ref{lem:fc}, we obtain the following lemma.

\begin{lemma}\label{lem:root}
Let $x\in\calS(\RR^d_\theta)$ and $p\in\NN$. For every $\delta>0$, the positive operators
$$
|x|^2(|x|^2+\delta)^{\frac{1}{4p}-1},
\quad
\big(|x|^2(|x|^2+\delta)^{\frac{1}{4p}-1}\big)^2
\quad\text{and}\quad
\big(|x|^2(|x|^2+\delta)^{\frac{1}{4p}-1}\big)^{2p}
$$
belong to $\calS(\RR^d_\theta)$. Moreover,
\begin{equation}\label{eq:sroot}
\big(|x|^2(|x|^2+\delta)^{\frac{1}{4p}-1}\big)^{2p}
=
|x|
\Big(\frac{|x|^2}{|x|^2+\delta}\Big)^{2p-\frac12}.
\end{equation}
When $\delta\to0$,
$$
\big(|x|^2(|x|^2+\delta)^{\frac{1}{4p}-1}\big)^{2p}
\to |x|
$$
in every $L_r(\RR^d_\theta)$ $(0<r<\infty)$ containing $|x|$. Moreover, if $|x|P_0\in\calL_r(L_2(\RR^d))$ for some $0<r<\infty$, then
\begin{equation}\label{eq:slconv}
\Big\|
\Big(
\big(|x|^2(|x|^2+\delta)^{\frac{1}{4p}-1}\big)^{2p}
-|x|
\Big)P_0
\Big\|_{\calL_r(L_2(\RR^d))}
\to0.
\end{equation}
\end{lemma}

\begin{proof}
Since $x^*x=|x|^2\in\calS(\RR^d_\theta)$, we apply Lemma \ref{lem:fc} to $t\mapsto t(t+\delta)^{\frac{1}{4p}-1}$ on $(-\delta,\infty)$ and obtain
$$
0\leq |x|^2(|x|^2+\delta)^{\frac{1}{4p}-1}
\in\calS(\RR^d_\theta).
$$
Then $$\big(|x|^2(|x|^2+\delta)^{\frac{1}{4p}-1}\big)^2
\quad\text{and}\quad
\big(|x|^2(|x|^2+\delta)^{\frac{1}{4p}-1}\big)^{2p}\in \calS(\RR^d_\theta).$$
For $s\geq0$, note that
$$
\big(s^2(s^2+\delta)^{\frac{1}{4p}-1}\big)^{2p}
=s\Big(\frac{s^2}{s^2+\delta}\Big)^{2p-\frac12},
$$
and this immediately implies \eqref{eq:sroot} by functional calculus.

For $0<r<\infty$ and $x\in L_r(\RR^d_\theta)$, define the positive Borel measure associated with $x$ as
$$
\nu_r(B)=
\tau_\theta\big(|x|^r1_B(|x|)\big),
$$
where $B$ is a Borel subset of $[0,\infty)$. Since $\nu_r([0,\infty))=\tau_\theta(|x|^r)<\infty$, one has from  \eqref{eq:sroot} that as $\delta\to0$,
\begin{align*}
\Big\|\big(|x|^2(|x|^2+\delta)^{\frac{1}{4p}-1}\big)^{2p}-|x|
\Big\|_{L_r(\RR^d_\theta)}^r
=&\int_{[0,\infty)}\Big[1-\Big(\frac{s^{2/r}}{s^{2/r}+\delta}\Big)^{2p-\frac12}\Big]^r d\nu_r(s)
\to0.
\end{align*}
It remains to show \eqref{eq:slconv}. By the spectral theorem, one has
$$\Big[1-\Big(\frac{|x|^2}{|x|^2+\delta}\Big)^{2p-\frac12}
\Big]\mathds{1}_{(0,\infty)}(|x|)\to0\quad\text{strongly},$$
and
$$\Big\|1-\Big(\frac{|x|^2}{|x|^2+\delta}\Big)^{2p-\frac12}
\Big\|_{\calB(L_2(\RR^d))}\leq1.$$
Since $|x|P_0\in\calL_r(L_2(\RR^d))$, one has
$$\sum_{k=0}^\infty \mu(k,|x|P_0)^r=\big\||x|P_0\big\|_{\calL_r(L_2(\RR^d))}^r<\infty.$$
Hence, for every $n\geq1$,
$$\rank\Big(\mathds{1}_{(1/n,\infty)}\big(|x|P_0\big)\Big)
=\#\Big\{k:\mu(k,|x|P_0)>\frac1n\Big\}<\infty.$$
By \eqref{eq:sroot},
$$\Big(\big(|x|^2(|x|^2+\delta)^{\frac{1}{4p}-1}\big)^{2p}-|x|\Big)P_0
=-\Big[1-\Big(\frac{|x|^2}{|x|^2+\delta}\Big)^{2p-\frac12}\Big]|x|P_0.$$
Then by the triangle inequality or the $r$-norm inequality for $0<r<1$, one has
\begin{align}\label{eq:tterm}
\begin{split}
&\Big\|\Big(\big(|x|^2(|x|^2+\delta)^{\frac{1}{4p}-1}\big)^{2p}
-|x|\Big)P_0\Big\|_{\calL_r(L_2(\RR^d))}^{\min\{1,r\}}\\
\leq&\Big\|\Big[1-\Big(\frac{|x|^2}{|x|^2+\delta}\Big)^{2p-\frac12}\Big]|x|P_0\mathds{1}_{[0,1/n]}\big(|x|P_0\big)\Big\|_{\calL_r(L_2(\RR^d))}^{\min\{1,r\}}\\
&+\Big\|\Big[1-\Big(\frac{|x|^2}{|x|^2+\delta}\Big)^{2p-\frac12}\Big]|x|P_0\mathds{1}_{(1/n,\infty)}\big(|x|P_0\big)\Big\|_{\calL_r(L_2(\RR^d))}^{\min\{1,r\}}.
\end{split}
\end{align}
The first term in the summand of \eqref{eq:tterm} is bounded by
$$\Big\||x|P_0\mathds{1}_{[0,1/n]}\big(|x|P_0\big)\Big\|_{\calL_r(L_2(\RR^d))}^{\min\{1,r\}}.$$
Note that
\begin{align*}
\Big\||x|P_0\mathds{1}_{[0,1/n]}\big(|x|P_0\big)\Big\|_{\calL_r(L_2(\RR^d))}^r
=&\Tr\Big(\big(|x|P_0\big)^r\mathds{1}_{(0,1/n]}\big(|x|P_0\big)\Big)\\
=&\sum_{\substack{k\geq0\\0<\mu(k,|x|P_0)\leq1/n}}\mu(k,|x|P_0)^r.
\end{align*}
Since $|x|P_0\in\calL_r(L_2(\RR^d))$, one has that as $n\to\infty$,
$$\sum_{\substack{k\geq0\\0<\mu(k,|x|P_0)\leq1/n}}
\mu(k,|x|P_0)^r=\Tr\Big(\big(|x|P_0\big)^r
\mathds{1}_{(0,1/n]}\big(|x|P_0\big)\Big)
\to0. $$
Now we estimate the second term in the summand of \eqref{eq:tterm}. Note that $\mathds{1}_{(1/n,\infty)}(|x|P_0)$ has finite rank. By using the following elementary inequality
$$\|TR\|_{\calL_r(H)}\leq\rank(R)^{1/r}\|TR\|_{\calB(H)},$$
we obtain
\begin{align*}
&\Big\|\Big[1-\Big(\frac{|x|^2}{|x|^2+\delta}\Big)^{2p-\frac12}\Big]|x|P_0\mathds{1}_{(1/n,\infty)}\big(|x|P_0\big)\Big\|_{\calL_r(L_2(\RR^d))}^{\min\{1,r\}}\\
\leq&\Big[\rank\Big(\mathds{1}_{(1/n,\infty)}\big(|x|P_0\big)\Big)\Big]^{\frac{\min\{1,r\}}{r}}\cdot\Big\|\Big[1-\Big(\frac{|x|^2}{|x|^2+\delta}\Big)^{2p-\frac12}\Big]|x|P_0\mathds{1}_{(1/n,\infty)}\big(|x|P_0\big)\Big\|_{\calB(L_2(\RR^d))}^{\min\{1,r\}}.
\end{align*}
The operator $|x|P_0 \mathds{1}_{(1/n,\infty)}\big(|x|P_0\big)$ has finite rank, and its range is contained in $\mathds{1}_{(0,\infty)}(|x|)L_2(\RR^d)$. Let $e_1,\ldots,e_m$ be an orthonormal basis of this range. For a unit vector $\xi\in L_2(\RR^d)$, write
$$|x|P_0 \mathds{1}_{(1/n,\infty)}\big(|x|P_0\big)\xi
=\sum_{j=1}^m c_j e_j. $$
Then by Parseval's identity,
\begin{align*}
\sum_j|c_j|^2=&\Big\||x|P_0\mathds{1}_{(1/n,\infty)}\big(|x|P_0\big)\xi\Big\|_{L_2(\RR^d)}^2\\
\leq&\Big\||x|P_0\mathds{1}_{(1/n,\infty)}\big(|x|P_0\big)\Big\|_{\calB(L_2(\RR^d))}^2.
\end{align*}
Note that for each $j$, as $\delta \to 0$,
$$\Big\|\Big[1-\Big(\frac{|x|^2}{|x|^2+\delta}\Big)^{2p-\frac12}\Big]e_j\Big\|_{L_2(\RR^d)}
\to0.$$
Hence, the Cauchy--Schwarz inequality yields
\begin{align*}
&\Big\|\Big[1-\Big(\frac{|x|^2}{|x|^2+\delta}\Big)^{2p-\frac12}\Big]|x|P_0 \mathds{1}_{(1/n,\infty)}\big(|x|P_0\big)\xi\Big\|_{L_2(\RR^d)}\\
\leq&\Big(\sum_j|c_j|^2\Big)^{1/2}\Big(\sum_j\Big\|\Big[1-\Big(\frac{|x|^2}{|x|^2+\delta}\Big)^{2p-\frac12}\Big]e_j\Big\|_{L_2(\RR^d)}^2\Big)^{1/2}
\to0,
\end{align*}
uniformly over any unit vectors $\xi$. Hence, as $\delta \to 0$,
$$\Big\|\Big[1-\Big(\frac{|x|^2}{|x|^2+\delta}\Big)^{2p-\frac12}\Big]|x|P_0 \mathds{1}_{(1/n,\infty)}\big(|x|P_0\big)\Big\|_{\calB(L_2(\RR^d))}\to0.$$
After taking the limit first as $\delta\to0$ and then as $n\to\infty$ in \eqref{eq:tterm}, we obtain \eqref{eq:slconv}.
\end{proof}

\bigskip

\section{Approximation of tempered distributions}\label{s:appr}
%In this section, we establish several approximation properties needed in the proof of the necessity part of Theorem \ref{thm:m1} for arbitrary tempered distributions. More precisely, this section is devoted to the introduction of convolution by heat kernels and approximation by Schwartz elements. We proceed with the definition of convolution of tempered distributions.

In this section, we establish the approximation properties needed in the proof of the necessity part of Theorem \ref{thm:m1} for arbitrary tempered distributions. The required approximation properties must preserve the control of the weak Schatten quasi-norms of the commutators and allow us to apply the asymptotic formula to Schwartz elements. We first use convolution by heat kernels to obtain the required bound for the commutators. Then we show that multiplication by Schwartz elements turns the distributional derivatives of $h_t*x$ into Schwartz elements. This will allow us to prove that the distributional derivatives of $x$ belong to $L_{\frac{d}{1-\epsilon}}(\RR^d_\theta)$ in Section \ref{s:comm}. We proceed with the definition of convolution of tempered distributions.

\subsection{Convolution of tempered distributions and commutators}
For $\eta,t\in\RR^d$, $f\in \calS(\RR^d)$ and $x\in\calS(\RR^d_\theta)$, define
$$(m_\eta f)(t)=e^{-\ii\langle\eta,t\rangle}f(t)\quad\text{and}\quad
\sigma_\eta(x)=m_\eta x m_\eta^*.$$
Then for any $f\in L_2(\RR^d)$,
\begin{align}\label{eq:mum1}
\sigma_\eta(U_\theta(f))=U_\theta(e^{-\ii\langle\eta,\cdot\rangle}f)=U_\theta(m_\eta f).
\end{align}
For $\eta\in\RR^d$ and $x\in\calS'(\RR^d_\theta)$, define $\sigma_\eta(x)$ by
$$\langle\sigma_\eta(x),a\rangle
=\langle x,\sigma_{-\eta}(a)\rangle,
\quad a\in\calS(\RR^d_\theta).$$
%This definition is well posed by Proposition \ref{prop:mod}. Indeed, for each %$x\in\calS'(\RR^d_\theta)$, we choose $N\in\NN_0$ such that
%$$|\langle x,U_\theta(f)\rangle|\lesssim_x\mathfrak p_N(f),\quad f\in\calS(\RR^d).$$
%Then
%\begin{align*}
%|\langle\sigma_\eta(x),U_\theta(f)\rangle|
%=&|\langle x,U_\theta(m_\eta^*f)\rangle|\\
%\lesssim&_x\mathfrak p_N(m_\eta^*f)\\
%\lesssim&_{x,N}(1+|\eta|)^N\mathfrak p_N(f).
%\end{align*}
%Thus $\sigma_\eta(x)\in\calS'(\RR^d_\theta)$.
This definition is consistent with the special case $x\in\calS(\RR^d_\theta)$. Then $\sigma_\eta$ is a trace-preserving normal $*$-automorphism on $L_\infty(\RR^d_\theta)$, and $\sigma_\eta$ is an isometry on every $L_p(\RR^d_\theta)$ for all $0<p\leq\infty$. 
%For more details, see \cite{GJP2021,MSX2020}. 
Note that Gonz\'alez-P\'erez, Junge and Parcet described the corresponding action
through a corepresentation in \cite{GJP2021}. Here, we use the argument in \cite{MSX2020}.

Next, we show that, for $\psi\in\calS(\RR^d)$ and $x\in\calS(\RR^d_\theta)$, the integral
$$\int_{\RR^d}\psi(\eta)\sigma_{\eta}(x)  d\eta$$
still belongs to $\calS(\RR^d_\theta)$. The explicit estimates of Schwartz seminorms are also shown in the following proposition.

\begin{proposition}\label{prop:mod}
The group of unitary operators $(m_\eta)_{\eta\in\RR^d}$ is strongly continuous on $L_2(\RR^d)$. For every $f\in\calS(\RR^d)$ and $N\in\NN_0$,
$$\mathfrak p_N(m_\eta f)\lesssim_N
(1+|\eta|)^N\mathfrak p_N(f), \quad \eta\in\RR^d,$$
and the map $\eta\mapsto m_\eta f$ is continuous from $\RR^d$ to $\calS(\RR^d)$. Moreover, for every $\psi\in\calS(\RR^d)$ and $x\in\calS(\RR^d_\theta)$, the integral $\int_{\RR^d}\psi(\eta)\sigma_{\eta}(x) d\eta$ belongs to $\calS(\RR^d_\theta)$.
\end{proposition}

\begin{proof}
Note that as $\eta \to \eta_0$,
$$
\|(m_\eta-m_{\eta_0})f\|_{L_2(\RR^d)}^2
=\int_{\RR^d}|e^{-\ii\langle\eta-\eta_0,t\rangle}-1|^2|f(t)|^2 dt\to0.
$$
The implies that $(m_\eta)_{\eta\in\RR^d}$ is strongly continuous on $L_2(\RR^d)$. 

Let $N\in\NN_0$. For $f\in\calS(\RR^d)$, we use the Leibniz rule to get
$$\partial_s^\alpha(m_\eta f)(s)=e^{-\ii\langle\eta,s\rangle}\sum_{\beta\leq\alpha}\binom{\alpha}{\beta}(-\ii\eta)^\beta\partial_s^{\alpha-\beta}f(s).$$
Hence, we obtain
\begin{align}
\begin{split}\label{eq:pnf1}
\mathfrak p_N(m_\eta f)
\leq&\max_{|\alpha|\leq N}\sum_{\beta\leq\alpha}
\binom{\alpha}{\beta}|\eta|^{|\beta|}\sup_{s\in\RR^d}(1+|s|)^N|\partial_s^{\alpha-\beta}f(s)|\\
\lesssim&_N(1+|\eta|)^N\mathfrak p_N(f).
\end{split}
\end{align}

Next, we show the continuity of the map $\eta\mapsto m_\eta f$. It is clear that
$$|e^{-\ii\langle\eta,s\rangle}-e^{-\ii\langle\eta_0,s\rangle}|\leq|\eta-\eta_0||s|$$
and for $|\beta|\geq1$,
$$|\eta^\beta-\eta_0^\beta|\lesssim_\beta|\eta-\eta_0|(1+|\eta|+|\eta_0|)^{|\beta|-1}.$$
For $|\alpha|\leq N$, we write
\begin{align*}
\partial_s^\alpha\big((m_\eta-m_{\eta_0})f\big)(s)
=&\sum_{\beta\leq\alpha}\binom{\alpha}{\beta}(-\ii\eta)^\beta\big(e^{-\ii\langle\eta,s\rangle}-e^{-\ii\langle\eta_0,s\rangle}
\big)\partial_s^{\alpha-\beta}f(s)\\
&+e^{-\ii\langle\eta_0,s\rangle}\sum_{\substack{\beta\leq\alpha\\|\beta|\geq1}}\binom{\alpha}{\beta}\big((-\ii\eta)^\beta-(-\ii\eta_0)^\beta\big)\partial_s^{\alpha-\beta}f(s).
\end{align*}
Multiplying by $(1+|s|)^N$ and taking the supremum over $s$, we obtain
\begin{align}\label{eq:pmf1}
\mathfrak p_N\big((m_\eta-m_{\eta_0})f\big)\lesssim_N|\eta-\eta_0|(1+|\eta|+|\eta_0|)^N\mathfrak p_{N+1}(f)\to0,
\end{align}
as $\eta \to \eta_0$. This implies that the map $\eta\mapsto m_\eta f$ is continuous from $\RR^d$ to $\calS(\RR^d)$. 

Let $x=U_\theta(f)$ with $f\in\calS(\RR^d)$.
%We have $\sigma_{-\eta}(x)=U_\theta(m_\eta^*f)$. For every $N\in\NN_0$,
%\begin{align*}
%\int_{\RR^d}\mathfrak p_N\big(\psi(\eta)m_\eta^*f\big) d\eta
%=&\int_{\RR^d}|\psi(\eta)|\mathfrak p_N(m_\eta^*f) d\eta\\
%\lesssim&_N\mathfrak p_N(f)\int_{\RR^d}|\psi(\eta)|(1+|\eta|)^N d\eta\\
%\leq&\mathfrak p_N(f)\mathfrak p_{N+d+1}(\psi)\int_{\RR^d}(1+|\eta|)^{-d-1} d\eta
%<\infty.
%\end{align*}
By \eqref{eq:pnf1}, for $R'>R>0$,
\begin{align*}
&\mathfrak p_N\Big(\int_{[-R',R']^d}\psi(\eta)m_\eta f d\eta
-\int_{[-R,R]^d}\psi(\eta)m_\eta f d\eta
\Big)\\
\leq&
\int_{[-R',R']^d\setminus[-R,R]^d}
|\psi(\eta)|\mathfrak p_N(m_\eta f) d\eta\\
\lesssim&_N
\mathfrak p_N(f)
\int_{|\eta|>R}
|\psi(\eta)|(1+|\eta|)^N d\eta
\to0
\quad(R\to\infty).
\end{align*}
Thus the truncated integrals are a Cauchy sequence for every seminorm $\mathfrak p_N$. By completeness,
$$
\int_{\RR^d}\psi(\eta)m_\eta f d\eta
\in\calS(\RR^d).
$$
Hence, we have
$$\int_{\RR^d}\psi(\eta)\sigma_{\eta}(x) d\eta=U_\theta\Big(\int_{\RR^d}\psi(\eta)m_\eta f d\eta\Big)\in\calS(\RR^d_\theta).$$
This completes the proof.
\end{proof}

The proposition below establishes the relation between the commutator and $\sigma_\eta$. Recall that $\Phi_{\epsilon;a, b}$ is defined in \eqref{eq:Phi}.

\begin{proposition}
Let $-\frac d2<\epsilon\neq 0<1$. For every $\eta\in\RR^d$ and $f,g\in\calS(\RR^d)$,
\begin{equation}\label{eq:cov}
\Phi_{\epsilon;Vm_\eta^*f,Vm_\eta^*g}
=
\sigma_{-\eta}(\Phi_{\epsilon;Vf,Vg}).
\end{equation}
Moreover, for every $x\in\calS'(\RR^d_\theta)$,
$$
\big\langle[(-\Delta)^{\frac{\epsilon}{2}},\sigma_\eta(x)]f,\overline g\big\rangle
=\big\langle[(-\Delta)^{\frac{\epsilon}{2}},x]m_\eta^*f,\overline{m_\eta^*g}
\big\rangle.
$$
\end{proposition}

\begin{proof}
Let $f,g\in\calS(\RR^d)$. By the definitions of $V$, $(-\Delta_\theta)^{\frac{\epsilon}{2}}$ and $\sigma_{-\eta}$,
\begin{align}\label{eq:v1}
Vm_\eta^*f
=(2\pi)^{-\frac d2}U_\theta(m_\eta^*f)
=\sigma_{-\eta}\big((2\pi)^{-\frac d2}U_\theta(f)\big)
=\sigma_{-\eta}(Vf)
\end{align}
and
\begin{align*}
V\big(|\cdot|^\epsilon m_\eta^*f\big)=V\big(m_\eta^*|\cdot|^\epsilon f\big)=\sigma_{-\eta}(-\Delta_\theta)^{\frac{\epsilon}{2}}V(f),
\end{align*}
where we use the fact that $|\cdot|^\epsilon f\in L_2(\RR^d)$. Then, by \eqref{eq:cform1} and $\sigma_\eta$ is $*$-automorphism,
\begin{align*}
\Phi_{\epsilon;Vm_\eta^*f,Vm_\eta^*g}
=&\Big(V(m_\eta^*f)V(|\cdot|^\epsilon m_\eta^*g)^*
-V(|\cdot|^\epsilon m_\eta^*f)V(m_\eta^*g)^*\Big)\\
=&\sigma_{-\eta}(Vf)\sigma_{-\eta}\Big(\big((-\Delta_\theta)^{\frac{\epsilon}{2}}Vg\big)^*\Big)-\sigma_{-\eta}\Big(
(-\Delta_\theta)^{\frac{\epsilon}{2}}Vf\Big)
\sigma_{-\eta}\big((Vg)^*\big)\\
=&\sigma_{-\eta}\Big(Vf\big((-\Delta_\theta)^{\frac{\epsilon}{2}}Vg\big)^*-\big((-\Delta_\theta)^{\frac{\epsilon}{2}}Vf\big)(Vg)^*
\Big)\\
=&\sigma_{-\eta}(\Phi_{\epsilon;Vf,Vg}).
\end{align*}
By Proposition \ref{prop:dw1} (i) and Proposition \ref{prop:mod}, both sides belong to $\calS(\RR^d_\theta)$, so this proves \eqref{eq:cov}.

For $x\in\calS'(\RR^d_\theta)$, by using \eqref{eq:cov} and the definitions of $[(-\Delta)^{\frac{\epsilon}{2}},\sigma_\eta(x)]$ and $[(-\Delta)^{\frac{\epsilon}{2}},x]$, we obtain
\begin{align*}
\big\langle
[(-\Delta)^{\frac{\epsilon}{2}},\sigma_\eta(x)]f,
\overline g\big\rangle
=&\langle\sigma_\eta(x),\Phi_{\epsilon;Vf,Vg}\rangle\\
=&\langle x,\sigma_{-\eta}(\Phi_{\epsilon;Vf,Vg})\rangle\\
=&\langle x,\Phi_{\epsilon;Vm_\eta^*f,Vm_\eta^*g}\rangle\\
=&\big\langle[(-\Delta)^{\frac{\epsilon}{2}},x]m_\eta^*f,
\overline{m_\eta^*g}
\big\rangle.
\end{align*}
This completes the proof.
\end{proof}

For $x\in L_p(\RR^d_\theta)$, $1\leq p<\infty$, and $\psi\in L_1(\RR^d)$, convolution has already been defined in \cite{MSX2020}. We use the same formula for Schwartz elements and extend this definition to arbitrary $x\in\calS'(\RR^d_\theta)$ by duality.

\begin{definition}
For $\psi\in\calS(\RR^d)$ and $x\in\calS(\RR^d_\theta)$, we use the convolution formula in \cite[Definition 3.9]{MSX2020} to define
%For $\psi\in\calS(\RR^d)$ and $x\in\calS(\RR^d_\theta)$, define
$$\psi*x=\int_{\RR^d}\psi(\eta)\sigma_\eta(x) d\eta,$$
where the integral converges in the topology of $\calS(\RR^d_\theta)$ as guaranteed by Proposition \ref{prop:mod}. For $x\in\calS'(\RR^d_\theta)$, define $\psi*x\in\calS'(\RR^d_\theta)$ by
$$\langle\psi*x,a\rangle
=\int_{\RR^d}\psi(\eta)\langle x,\sigma_{-\eta}(a)\rangle d\eta, \quad \forall a\in\calS(\RR^d).$$ 
\end{definition}

We need to verify that $\psi*x$ is well-defined for $x\in\calS'(\RR^d_\theta)$. Indeed, by the continuity of $x$ on $\calS(\RR^d_\theta)$, there exist $N\in\NN_0$ and $C_x>0$ such that
$$|\langle x,\sigma_{-\eta}(a)\rangle|\le C_x\mathfrak p_N(\sigma_{-\eta}(a)),\quad a\in\calS(\RR^d_\theta).$$
By Proposition \ref{prop:mod} together with $\sigma_{-\eta}(U_\theta(f))=U_\theta(m_\eta^*f)$, gives
$$\mathfrak p_N(\sigma_{-\eta}(a))\lesssim_N(1+|\eta|)^N\mathfrak p_N(a).$$
Since $\psi\in\calS(\RR^d)$, one has
$$|\psi(\eta)|(1+|\eta|)^N\le \mathfrak p_{N+d+1}(\psi)(1+|\eta|)^{-d-1}.$$
Note that $\eta\mapsto\psi(\eta)\langle x,\sigma_{-\eta}(a)\rangle$ is measurable, since $\eta\mapsto\sigma_{-\eta}(a)$ is continuous in $\calS(\RR^d_\theta)$ by Proposition \ref{prop:mod}. Therefore, we obtain
\begin{align*}
\int_{\RR^d}|\psi(\eta)| |\langle x,\sigma_{-\eta}(a)\rangle| d\eta
&\lesssim_{x,N}\mathfrak p_N(a)\int_{\RR^d}|\psi(\eta)|(1+|\eta|)^N d\eta\\
&\lesssim_{x,N,d}\mathfrak p_{N+d+1}(\psi)\mathfrak p_N(a)<\infty.
\end{align*}
This implies that $\eta\mapsto\psi(\eta)\langle x,\sigma_{-\eta}(a)\rangle$ is integrable. The same estimate shows that $a\mapsto\langle\psi*x,a\rangle$ is a continuous linear functional on $\calS(\RR^d_\theta)$, hence $\psi*x\in\calS'(\RR^d_\theta)$. Next, we verify that these two definitions agree when $x\in\calS(\RR^d_\theta)$. For $x,a\in\calS(\RR^d_\theta)$, we have
\begin{align*}
\Big\langle\int_{\RR^d}\psi(\eta)\sigma_\eta(x) d\eta,a\Big\rangle
=&\int_{\RR^d}\psi(\eta)\tau_\theta(\sigma_\eta(x)a) d\eta\\
=&\int_{\RR^d}\psi(\eta)\tau_\theta\big(\sigma_\eta(x\sigma_{-\eta}(a))\big) d\eta\\
=&\int_{\RR^d}\psi(\eta)\tau_\theta(x\sigma_{-\eta}(a)) d\eta\\
=&\int_{\RR^d}\psi(\eta)\langle x,\sigma_{-\eta}(a)\rangle d\eta.
\end{align*}

In the following lemma, we show that $\eta\mapsto m_\eta Tm_\eta^*$ is continuous with respect to the operator norm for every compact operator $T$.

\begin{lemma}\label{lem:orbit}
Let $T$ be a compact operator on $L_2(\RR^d)$. Then the map $\eta\mapsto m_\eta Tm_\eta^*$ is continuous with respect to the operator norm. More precisely, for every $\eta_0\in\RR^d$,
$$
\|m_\eta Tm_\eta^*-m_{\eta_0}Tm_{\eta_0}^*\|_{\calB(L_2(\RR^d))}\to0\quad \text{as} \quad \eta\to\eta_0.
$$
\end{lemma}

\begin{proof}
Af first, since $T$ is compact, we can choose a sequence $(F_n)_{n\geq 1}$ of finite-rank operators such that as $n\to \infty$,
$$
\|T-F_n\|_{\calB(L_2(\RR^d))}\to0.
$$
For each $n$, we write
$$
F_n\xi=\sum_{j=1}^{k_n}\langle\xi,v_{n,j}\rangle_{L_2(\RR^d)}u_{n,j},\quad \forall \xi\in L_2(\RR^d),
$$
where $u_{n,j},v_{n,j}\in L_2(\RR^d)$. Here, $\langle \cdot,\cdot\rangle_{L_2(\RR^d)}$ is the inner product of $L_2(\RR^d)$. It is immediate that one has
$$
m_\eta F_nm_\eta^*\xi=\sum_{j=1}^{k_n}\langle\xi,m_\eta v_{n,j}\rangle_{L_2(\RR^d)}m_\eta u_{n,j}.
$$
Therefore, we obtain
\begin{align*}
(m_\eta F_nm_\eta^*-m_{\eta_0}F_nm_{\eta_0}^*)\xi
=&\sum_{j=1}^{k_n}\langle\xi,m_\eta v_{n,j}\rangle_{L_2(\RR^d)}(m_\eta-m_{\eta_0})u_{n,j}\\
&+\sum_{j=1}^{k_n}\langle\xi,(m_\eta-m_{\eta_0})v_{n,j}\rangle_{L_2(\RR^d)}m_{\eta_0}u_{n,j}.
\end{align*}
By the Cauchy--Schwarz inequality, one has
\begin{align*}
&\|(m_\eta F_nm_\eta^*-m_{\eta_0}F_nm_{\eta_0}^*)\xi\|_{L_2(\RR^d)}\\
\leq&\|\xi\|_{L_2(\RR^d)}\sum_{j=1}^{k_n}\|(m_\eta-m_{\eta_0})u_{n,j}\|_{L_2(\RR^d)}\|v_{n,j}\|_{L_2(\RR^d)}\\
&+\|\xi\|_{L_2(\RR^d)}\sum_{j=1}^{k_n}\|u_{n,j}\|_{L_2(\RR^d)}\|(m_\eta-m_{\eta_0})v_{n,j}\|_{L_2(\RR^d)}.
\end{align*}
After taking the supremum over all $\|\xi\|_{L_2(\RR^d)}\leq1$ and using the strong continuity of $(m_\eta)_{\eta\in\RR^d}$ established in Proposition \ref{prop:mod}, we obtain that for each fixed $n$, as $\eta\to\eta_0$,
\begin{align*}
&\|m_\eta F_nm_\eta^*-m_{\eta_0}F_nm_{\eta_0}^*\|_{\calB(L_2(\RR^d))}\\
\leq&\sum_{j=1}^{k_n}\|(m_\eta-m_{\eta_0})u_{n,j}\|_{L_2(\RR^d)}\|v_{n,j}\|_{L_2(\RR^d)}\\
&+\sum_{j=1}^{k_n}\|u_{n,j}\|_{L_2(\RR^d)}\|(m_\eta-m_{\eta_0})v_{n,j}\|_{L_2(\RR^d)}\to 0.
\end{align*}
By the triangle inequality, it follows that
\begin{align*}
&\|m_\eta Tm_\eta^*-m_{\eta_0}Tm_{\eta_0}^*\|_{\calB(L_2(\RR^d))}\\
\leq&\|m_\eta(T-F_n)m_\eta^*\|_{\calB(L_2(\RR^d))}+\|m_{\eta_0}(T-F_n)m_{\eta_0}^*\|_{\calB(L_2(\RR^d))}\\
&+\|m_\eta F_nm_\eta^*-m_{\eta_0}F_nm_{\eta_0}^*\|_{\calB(L_2(\RR^d))}\\
=&2\|T-F_n\|_{\calB(L_2(\RR^d))}+\|m_\eta F_nm_\eta^*-m_{\eta_0}F_nm_{\eta_0}^*\|_{\calB(L_2(\RR^d))}.
\end{align*}
Thus, for every $n\geq1$,
$$0\leq\limsup_{\eta\to\eta_0}\|m_\eta Tm_\eta^*-m_{\eta_0}Tm_{\eta_0}^*\|_{\calB(L_2(\RR^d))}
\leq2\|T-F_n\|_{\calB(L_2(\RR^d))}.$$
Letting $n\to\infty$, we conclude that
$$\lim_{\eta\to\eta_0}\|m_\eta Tm_\eta^*-m_{\eta_0}Tm_{\eta_0}^*\|_{\calB(L_2(\RR^d))}=0.$$
This completes the proof.
\end{proof}

\begin{definition}\label{def:45}
Let $-\frac d2<\epsilon\neq 0<1$, $x\in\calS'(\RR^d_\theta)$ and $\psi\in\calS(\RR^d)$. Define $\psi*[(-\Delta)^{\frac{\epsilon}{2}},x]$ as
$$\big\langle\big(\psi*[(-\Delta)^{\frac{\epsilon}{2}},x]\big)f,\overline g\big\rangle=\int_{\RR^d}\psi(\eta)\big\langle[(-\Delta)^{\frac{\epsilon}{2}},x]m_\eta^*f,\overline{m_\eta^*g}\big\rangle d\eta,
\quad \forall f,g\in\calS(\RR^d).$$
\end{definition}

We would like to remark that the integral in Definition \ref{def:45} is integrable because of the following proposition. In addition, the following proposition is also quite helpful when approximating commutators by convolution.

\begin{proposition}\label{prop:psix1}
Let $-\frac d2<\epsilon\neq 0<1$, $x\in\calS'(\RR^d_\theta)$ and $\psi\in\calS(\RR^d)$. 
\begin{description}
\item[(i)] for fixed $f,g\in\calS(\RR^d)$ the map
$$\eta\longmapsto\psi(\eta)\big\langle[(-\Delta)^{\frac{\epsilon}{2}},x]m_\eta^*f,\overline{m_\eta^*g}\big\rangle$$
is continuous on $\RR^d$ and absolutely integrable;
\item[(ii)] $\psi*[(-\Delta)^{\frac{\epsilon}{2}},x]$ defines a continuous linear map from $\calS(\RR^d)$ to $\calS'(\RR^d)$, and
\begin{equation}\label{eq:avform}
[(-\Delta)^{\frac{\epsilon}{2}},\psi*x]=\psi*[(-\Delta)^{\frac{\epsilon}{2}},x];
\end{equation}
\item[(iii)] If $[(-\Delta)^{\frac{\epsilon}{2}},x]$ is compact on $L_2(\RR^d)$, then
$$\psi*[(-\Delta)^{\frac{\epsilon}{2}},x]=\int_{\RR^d}\psi(\eta)m_\eta[(-\Delta)^{\frac{\epsilon}{2}},x]m_\eta^* d\eta.$$
\end{description}
\end{proposition}

\begin{proof}
(i): We choose $N\in\NN_0$ such that
$$
|\langle x,U_\theta(F)\rangle|\lesssim_x\mathfrak p_N(F),\quad \forall F\in\calS(\RR^d).
$$
By Proposition \ref{prop:dw1} (i) and the continuity of $V$ and $U_\theta^{-1}$ on the corresponding Schwartz spaces, we can choose $M\in\NN_0$ such that
$$
\mathfrak p_N\big(
U_\theta^{-1}(\Phi_{\epsilon;Vf,Vg})
\big)
\lesssim_{N,\epsilon,d,\theta}
\mathfrak p_M(f)\mathfrak p_M(g),
\quad f,g\in\calS(\RR^d).
$$
Then for $u,v\in\calS(\RR^d)$, one has
\begin{align*}
\big|\langle[(-\Delta)^{\frac{\epsilon}{2}},x]u,\overline v\rangle\big|
=&|\langle x,\Phi_{\epsilon;Vu,Vv}\rangle|\\
\lesssim&_x\mathfrak p_N\big(U_\theta^{-1}(\Phi_{\epsilon;Vu,Vv})\big)\\
\lesssim&_{x,\epsilon,d,\theta}\mathfrak p_M(u)\mathfrak p_M(v).
\end{align*}
Fix $\eta_0\in\RR^d$. Hence, we obtain
\begin{align*}
&\Big|\psi(\eta)\big\langle[(-\Delta)^{\frac{\epsilon}{2}},x]m_\eta^*f,\overline{m_\eta^*g}\big\rangle-\psi(\eta_0)\big\langle[(-\Delta)^{\frac{\epsilon}{2}},x]m_{\eta_0}^*f,\overline{m_{\eta_0}^*g}\big\rangle\Big|\\
\leq&|\psi(\eta)-\psi(\eta_0)|\cdot\big|\big\langle[(-\Delta)^{\frac{\epsilon}{2}},x]m_\eta^*f,\overline{m_\eta^*g}\big\rangle\big|\\
&+|\psi(\eta_0)|\cdot\big|\big\langle[(-\Delta)^{\frac{\epsilon}{2}},x](m_\eta^*-m_{\eta_0}^*)f,\overline{m_\eta^*g}\big\rangle\big|\\
&+|\psi(\eta_0)|\cdot\big|\big\langle[(-\Delta)^{\frac{\epsilon}{2}},x]m_{\eta_0}^*f,\overline{(m_\eta^*-m_{\eta_0}^*)g}\big\rangle\big|\\
\lesssim&_{x,\epsilon,d,\theta}|\psi(\eta)-\psi(\eta_0)|\mathfrak p_M(m_\eta^*f)\mathfrak p_M(m_\eta^*g)\\
&+|\psi(\eta_0)|\mathfrak p_M\big((m_\eta^*-m_{\eta_0}^*)f\big)\mathfrak p_M(m_\eta^*g)\\
&+|\psi(\eta_0)|\mathfrak p_M(m_{\eta_0}^*f)\mathfrak p_M\big((m_\eta^*-m_{\eta_0}^*)g\big).
\end{align*}
Suppose that $|\eta-\eta_0|\leq1$. Then by the triangle inequality,
$$
1+|\eta|\leq2+|\eta_0|\quad\text{and}\quad1+|\eta|+|\eta_0|\leq2+2|\eta_0|.
$$
Thus by Proposition \ref{prop:mod} and \eqref{eq:pmf1}, for every $h\in\calS(\RR^d)$,
\begin{align*}
\mathfrak p_M(m_\eta^*h)+\mathfrak p_M(m_{\eta_0}^*h)
\lesssim&_M\big((1+|\eta|)^M+(1+|\eta_0|)^M\big)\mathfrak p_M(h)\\
\lesssim&_{M,\eta_0}\mathfrak p_M(h)
\end{align*}
and
\begin{align*}
\mathfrak p_M\big((m_\eta^*-m_{\eta_0}^*)h\big)
\lesssim&_M|\eta-\eta_0|(1+|\eta|+|\eta_0|)^M\mathfrak p_{M+1}(h)\\
\lesssim&_{M,\eta_0}|\eta-\eta_0|\mathfrak p_{M+1}(h).
\end{align*}
Plugging into these estimates and using $\mathfrak p_M\leq\mathfrak p_{M+1}$, we obtain that as $\eta\to\eta_0$,
\begin{align*}
&\Big|\psi(\eta)\big\langle[(-\Delta)^{\frac{\epsilon}{2}},x]m_\eta^*f,\overline{m_\eta^*g}\big\rangle-\psi(\eta_0)\big\langle[(-\Delta)^{\frac{\epsilon}{2}},x]m_{\eta_0}^*f,\overline{m_{\eta_0}^*g}\big\rangle\Big|\\
\lesssim&_{x,\epsilon,d,\theta,M,\eta_0}|\psi(\eta)-\psi(\eta_0)|\mathfrak p_M(f)\mathfrak p_M(g)\\
&+|\psi(\eta_0)||\eta-\eta_0|\Big(\mathfrak p_{M+1}(f)\mathfrak p_M(g)+\mathfrak p_M(f)\mathfrak p_{M+1}(g)\Big)\\
\lesssim&_{x,\epsilon,d,\theta,M,\eta_0}\Big(|\psi(\eta)-\psi(\eta_0)|+|\psi(\eta_0)||\eta-\eta_0|\Big)\mathfrak p_{M+1}(f)\mathfrak p_{M+1}(g)
\to0.
\end{align*}
Thus the map $\eta\mapsto\psi(\eta)\big\langle[(-\Delta)^{\frac{\epsilon}{2}},x]m_\eta^*f,\overline{m_\eta^*g}\big\rangle$ is continuous on $\RR^d$.

By using \eqref{eq:cov}, \eqref{eq:v1} and Proposition \ref{prop:mod}, we obtain
\begin{align}\label{eq:mfg1}
\begin{split}
\big|\big\langle[(-\Delta)^{\frac{\epsilon}{2}},x]m_\eta^*f,\overline{m_\eta^*g}\big\rangle\big|
=&\big|\langle x,\Phi_{\epsilon;Vm_\eta^*f,Vm_\eta^*g}\rangle\big|\\
=&\big|\langle x,\sigma_{-\eta}(\Phi_{\epsilon;Vf,Vg})\rangle\big|\\
=&\big|\big\langle x,U_\theta\big(m_\eta^*U_\theta^{-1}(\Phi_{\epsilon;Vf,Vg})\big)
\big\rangle\big|\\
\lesssim&_{x}\mathfrak p_N\big(m_\eta^*U_\theta^{-1}(\Phi_{\epsilon;Vf,Vg})\big)\\
\lesssim&_{x,N}(1+|\eta|)^N\mathfrak p_N\big(U_\theta^{-1}(\Phi_{\epsilon;Vf,Vg})\big)\\
\lesssim&_{x,\epsilon,d,\theta}(1+|\eta|)^N\mathfrak p_M(f)\mathfrak p_M(g).
\end{split}
\end{align}
Since $\psi\in\calS(\RR^d)$,
\begin{align*}
\int_{\RR^d}|\psi(\eta)|(1+|\eta|)^N d\eta
\leq&
\mathfrak p_{N+d+1}(\psi)
\int_{\RR^d}(1+|\eta|)^{-d-1}d\eta\\
\lesssim&_d\mathfrak p_{N+d+1}(\psi)
<\infty.
\end{align*}
Then by \eqref{eq:mfg1}, one has
\begin{align*}
&\int_{\RR^d}|\psi(\eta)|
\big|
\big\langle
[(-\Delta)^{\frac{\epsilon}{2}},x]m_\eta^*f,
\overline{m_\eta^*g}
\big\rangle
\big|d\eta\\
\lesssim&_{x,\epsilon,d,\theta}
\mathfrak p_M(f)\mathfrak p_M(g)
\int_{\RR^d}|\psi(\eta)|(1+|\eta|)^N d\eta\\
\lesssim&_{x,\epsilon,d,\theta}
\mathfrak p_{N+d+1}(\psi)
\mathfrak p_M(f)\mathfrak p_M(g).
\end{align*}
Thus the defining integral is absolutely convergent and gives a continuous sesquilinear form in $f,g$. In particular, it defines a continuous linear map from $\calS(\RR^d)$ to $\calS'(\RR^d)$.

\

\noindent (ii): By the definition of $\psi*x$ and \eqref{eq:cov}, we have
\begin{align*}
\big\langle
[(-\Delta)^{\frac{\epsilon}{2}},\psi*x]f,\overline g
\big\rangle
=&\langle\psi*x,\Phi_{\epsilon;Vf,Vg}\rangle\\
=&\int_{\RR^d}\psi(\eta)
\langle x,\sigma_{-\eta}(\Phi_{\epsilon;Vf,Vg})\rangle d\eta\\
=&\int_{\RR^d}\psi(\eta)
\langle x,\Phi_{\epsilon;Vm_\eta^*f,Vm_\eta^*g}\rangle d\eta\\
=&\int_{\RR^d}\psi(\eta)
\big\langle
[(-\Delta)^{\frac{\epsilon}{2}},x]m_\eta^*f,
\overline{m_\eta^*g}
\big\rangle d\eta\\
=&\big\langle
\big(\psi*[(-\Delta)^{\frac{\epsilon}{2}},x]\big)f,
\overline g
\big\rangle.
\end{align*}
Since $f,g\in\calS(\RR^d)$ are arbitrary, we obtain \eqref{eq:avform}.

\

\noindent (iii): Now suppose that $[(-\Delta)^{\frac{\epsilon}{2}},x]$ is compact on $L_2(\RR^d)$. By Lemma \ref{lem:orbit}, the map $\eta\mapsto m_\eta[(-\Delta)^{\frac{\epsilon}{2}},x]m_\eta^*$ is continuous with respect to the operator norm. 
\iffalse
Note that
\begin{align*}
&\int_{\RR^d}|\psi(\eta)|
\|m_\eta[(-\Delta)^{\frac{\epsilon}{2}},x]m_\eta^*\|_{\calB(L_2(\RR^d))}
d\eta\\
=&\int_{\RR^d}|\psi(\eta)|
\|[(-\Delta)^{\frac{\epsilon}{2}},x]\|_{\calB(L_2(\RR^d))}
d\eta\\
=&\|\psi\|_{L_1(\RR^d)}
\|[(-\Delta)^{\frac{\epsilon}{2}},x]\|_{\calB(L_2(\RR^d))}
<\infty,
\end{align*}
this gives convergence of the integral in operator norm. 
\fi
Then for $f,g\in\calS(\RR^d)$, we obtain
\begin{align*}
&\Big\langle
\Big(
\int_{\RR^d}\psi(\eta)
m_\eta[(-\Delta)^{\frac{\epsilon}{2}},x]m_\eta^*d\eta
\Big)f,\overline g
\Big\rangle\\
=&\int_{\RR^d}\psi(\eta)
\big\langle
m_\eta[(-\Delta)^{\frac{\epsilon}{2}},x]m_\eta^*f,
\overline g
\big\rangle d\eta\\
=&\int_{\RR^d}\psi(\eta)
\big\langle
[(-\Delta)^{\frac{\epsilon}{2}},x]m_\eta^*f,
\overline{m_\eta^*g}
\big\rangle d\eta\\
=&\big\langle
\big(\psi*[(-\Delta)^{\frac{\epsilon}{2}},x]\big)f,
\overline g
\big\rangle.
\end{align*}
This completes the proof.
\end{proof}

In general, for $\psi\in L_1(\RR^d)$ and a compact operator $T$ on $L_2(\RR^d)$, we can also define
$$\psi*T=\int_{\RR^d}\psi(\eta)m_\eta Tm_\eta^* d\eta \in \calB(L_2(\RR^d)).$$
It is clear that the convolution above is well-defined as the integral is convergent. The next lemma concerns the fact that the $\calL_{p,\infty}(L_2(\RR^d))$-norm of $\psi*T$ is bounded above by the corresponding norm of $T$ for $1<p<\infty$. We refer to \cite[Lemma 18]{LMSZ2017} for a similar result on $\RR$, and give a different proof by duality of the result we need below.

\begin{lemma}\label{lem:avg}
Let $1<p<\infty$, $\psi\in L_1(\RR^d)$ and $T\in\calL_{p,\infty}(L_2(\RR^d))$. Then
\begin{equation}\label{eq:avg}
\|\psi*T\|_{\calL_{p,\infty}(L_2(\RR^d))}^{\#}
\leq\|\psi\|_{L_1(\RR^d)}\|T\|_{\calL_{p,\infty}(L_2(\RR^d))}^{\#}.
\end{equation}
\end{lemma}

\begin{proof}
\iffalse
By Lemma \ref{lem:orbit}, the map $\eta\mapsto m_\eta[(-\Delta)^{\frac{\epsilon}{2}},x]m_\eta^*$ is continuous in operator norm. Note that 
$$
\int_{\RR^d}|\psi(\eta)|
\|m_\eta Tm_\eta^*\|_{\calB(L_2(\RR^d))} d\eta
=\|\psi\|_{L_1(\RR^d)}\|T\|_{\calB(L_2(\RR^d))}<\infty,
$$
we may take the Bochner integral in $\calK(L_2(\RR^d))$. In particular, $\psi*T$ is compact.
\fi
Fix $N\geq1$ and a finite-rank operator $F$ with $\rank F\leq N$ and $\|F\|_{\calB(L_2(\RR^d))}\leq1$. For $A,B\in\calB(L_2(\RR^d))$, we have
$$
|\Tr(FA)-\Tr(FB)|=|\Tr(F(A-B))|\leq N\|A-B\|_{\calB(L_2(\RR^d))}.
$$
Thus $A\mapsto\Tr(FA)$ is a continuous linear functional on $\calB(L_2(\RR^d))$. Therefore, we obtain
$$
\Tr\Big(F\int_{\RR^d}\psi(\eta)m_\eta Tm_\eta^* d\eta\Big)
=\int_{\RR^d}\psi(\eta)\Tr(Fm_\eta Tm_\eta^*) d\eta.
$$
Note that the scalar integral is absolutely convergent, since
\begin{align*}
\int_{\RR^d}|\psi(\eta)||\Tr(Fm_\eta Tm_\eta^*)| d\eta
\leq&N\int_{\RR^d}|\psi(\eta)|\|m_\eta Tm_\eta^*\|_{\calB(L_2(\RR^d))} d\eta\\
=&N\|\psi\|_{L_1(\RR^d)}\|T\|_{\calB(L_2(\RR^d))}<\infty.
\end{align*}
By Lemma \ref{lem:mTr1} and the unitary invariance of singular values,
\begin{align*}
|\Tr(F(\psi*T))|
=&\Big|\int_{\RR^d}\psi(\eta)
\Tr(Fm_\eta Tm_\eta^*) d\eta\Big|\\
\leq&\int_{\RR^d}|\psi(\eta)|
\sum_{j=0}^{N-1}\mu(j,m_\eta Tm_\eta^*) d\eta\\
=&\|\psi\|_{L_1(\RR^d)}
\sum_{j=0}^{N-1}\mu(j,T).
\end{align*}
By taking the supremum over $F$ in Lemma \ref{lem:mTr1}, we obtain
\begin{equation}\label{eq:NT1}
\begin{split}
N^{\frac{1}{p}-1}\sum_{j=0}^{N-1}\mu(j,\psi*T)
\leq&\|\psi\|_{L_1(\RR^d)}
N^{\frac{1}{p}-1}\sum_{j=0}^{N-1}\mu(j,T)\\
\leq&\|\psi\|_{L_1(\RR^d)}
\|T\|_{\calL_{p,\infty}(L_2(\RR^d))}^{\#}.
\end{split}
\end{equation}
After taking the supremum over $N\geq1$, the proof of \eqref{eq:avg} is finished.
\end{proof}

\subsection{Approximation by heat kernels}
For $x\in\calS(\RR^d_\theta)$, recall that the Fourier transform of $x$ is given by \eqref{eq:inveru1}. Now we come to the introduction of the Fourier transform on $\calS'(\RR^d_\theta)$ by duality.

\begin{definition}
For $x\in\calS'(\RR^d_\theta)$, define its Fourier transform $\widehat x\in\calS'(\RR^d)$ by
$$
\langle\widehat x,a\rangle
=(2\pi)^{-d}\langle x,U_\theta(a(-\cdot))\rangle,
\quad a\in\calS(\RR^d).
$$ 
\end{definition}
We verify that this definition agrees with the usual Fourier transform of $x\in \calS(\RR^d_\theta)$. Let $x=U_\theta(f)$ with $f\in\calS(\RR^d)$, whose Fourier transform in \eqref{eq:inveru1} is $U_\theta^{-1}(x)=f$. For every $a\in\calS(\RR^d)$, the trace pairing, \eqref{eq:mul1}, \eqref{eq:fun1} and \eqref{eq:con1} give
\begin{align*}
\langle\widehat x,a\rangle
=&(2\pi)^{-d}\langle U_\theta(f),U_\theta(a(-\cdot))\rangle\\
=&(2\pi)^{-d}\tau_\theta\big(U_\theta(f)U_\theta(a(-\cdot))\big)\\
=&(2\pi)^{-d}\tau_\theta\big(U_\theta(f*_\theta(a(-\cdot)))\big)\\
=&\big(f*_\theta(a(-\cdot))\big)(0)\\
=&\int_{\RR^d}f(u)a(u) du,
\end{align*}
which implies $\widehat{x}=f$. 

In the classical setting, the multiplication of a tempered distribution $T\in \calS'(\RR^d)$ by $\phi\in\calS(\RR^d)$ is defined through
$$
\langle \phi T,f\rangle=\langle T\phi,f\rangle=\langle T,f\phi\rangle,
\quad f\in\calS(\RR^d).
$$
Convolution on $L_p(\RR_\theta^d)$ and convergence under an approximate identity were given in \cite[Definition 3.9 and Theorem 3.10]{MSX2020}. However, we need to require convergence in $\calS'(\RR^d_\theta)$ and estimates for the corresponding commutators. The lemma below, employing the heat kernels as the convolution kernel, establishes a helpful approximation property for tempered distributions on quantum Euclidean spaces.

\begin{lemma}\label{lem:heat}
For $t>0$, define
$$
h_t(\eta)=(4\pi t)^{-\frac{d}{2}}e^{-\frac{|\eta|^2}{4t}}, \quad \forall \eta\in \mathbb{R}^d.
$$
Let $-\frac{d}{2}<\epsilon\neq 0<1$ and $x\in\calS'(\RR^d_\theta)$. Let $y_t=h_t*x$. Then one has the following.
\begin{description}
\item[(i)] The tempered distribution $y_t$ has the Fourier transform $\widehat{y_t}=e^{-t|\cdot|^2}\widehat x$. Moreover, $y_t\to x$ in $\calS'(\RR^d_\theta)$ as $t\to0^+$;
\item[(ii)] For $t>0$,
\begin{equation}\label{eq:heatid}
[(-\Delta)^{\frac{\epsilon}{2}},y_t]
=h_t*[(-\Delta)^{\frac{\epsilon}{2}},x];
\end{equation}
\item[(iii)] If $[(-\Delta)^{\frac{\epsilon}{2}},x]\in\calL_{\frac{d}{1-\epsilon},\infty}(L_2(\RR^d))$, then
\begin{equation}\label{eq:heatc}
\|[(-\Delta)^{\frac{\epsilon}{2}},y_t]\|_{\calL_{\frac{d}{1-\epsilon},\infty}(L_2(\RR^d))}
\lesssim_{d,\epsilon}
\|[(-\Delta)^{\frac{\epsilon}{2}},x]\|_{\calL_{\frac{d}{1-\epsilon},\infty}(L_2(\RR^d))}.
\end{equation}
\end{description}
\end{lemma}

\begin{proof}
(i): By Proposition \ref{prop:mod}, for every $a\in\calS(\RR^d)$ and $N\in\NN_0$,
\begin{align*}
\int_{\RR^d}h_t(\eta)\mathfrak p_N(m_\eta a) d\eta
\lesssim&_N\mathfrak p_N(a)\int_{\RR^d}h_t(\eta)(1+|\eta|)^N d\eta
<\infty.
\end{align*}
Thus the integral converges in $\calS(\RR^d)$. For every $s=(s_1,\ldots,s_d)\in\RR^d$, by the definition of $m_\eta$ and Fubini's theorem, we have
\begin{align*}
\Big(\int_{\RR^d}h_t(\eta)m_\eta a d\eta\Big)(s)
=&\int_{\RR^d}h_t(\eta)(m_\eta a)(s) d\eta\\
=&a(s)\int_{\RR^d}h_t(\eta)e^{-\ii\langle\eta,s\rangle} d\eta\\
=&a(s)\widehat h_t(s)\\
=&e^{-t|s|^2}a(s).
\end{align*}
This implies
$$
\int_{\RR^d}h_t(\eta)m_\eta a d\eta=\widehat h_t a=e^{-t|\cdot|^2}a
\quad\text{in }\calS(\RR^d).
$$
Since $\widehat x$ is a continuous linear functional on $\calS(\RR^d)$, we obtain
\begin{align}\label{eq:htm1}
\int_{\RR^d}h_t(\eta)\langle\widehat x,m_\eta a\rangle d\eta
=\Big\langle\widehat x,\int_{\RR^d}h_t(\eta)m_\eta a d\eta\Big\rangle=\langle\widehat x,e^{-t|\cdot|^2}a\rangle.
\end{align}
For every $a\in\calS(\RR^d)$, note that
$$
\sigma_{-\eta}(U_\theta(a(-\cdot)))
=U_\theta(e^{\ii\langle\eta,\cdot\rangle}a(-\cdot))
=U_\theta((m_\eta a)(-\cdot)).
$$
Then by \eqref{eq:htm1},
\begin{align*}
\langle\widehat{y_t},a\rangle
=&(2\pi)^{-d}\langle h_t*x,U_\theta(a(-\cdot))\rangle\\
=&(2\pi)^{-d}\int_{\RR^d}h_t(\eta)\langle x,\sigma_{-\eta}(U_\theta(a(-\cdot)))\rangle d\eta\\
=&(2\pi)^{-d}\int_{\RR^d}h_t(\eta)\langle x,U_\theta((m_\eta a)(-\cdot))\rangle d\eta\\
=&\int_{\RR^d}h_t(\eta)\langle\widehat x,m_\eta a\rangle d\eta\\
=&\langle\widehat x,e^{-t|\cdot|^2}a\rangle\\
=&\langle e^{-t|\cdot|^2}\widehat x,a\rangle.
\end{align*}
Since $a\in\calS(\RR^d)$ is arbitrary, we conclude that
$$
\widehat{y_t}=e^{-t|\cdot|^2}\widehat x
\quad\text{in }\calS'(\RR^d).
$$

By using the estimate \eqref{eq:pmf1} in the proof of Proposition \ref{prop:mod} with $\eta_0=0$, one has
$$
\mathfrak p_N((m_\eta-1)a)\lesssim_N|\eta|(1+|\eta|)^N\mathfrak p_{N+1}(a).
$$
Since $\int_{\RR^d}h_t(\eta) d\eta=1$, we have
$$
(e^{-t|\cdot|^2}-1)a=\int_{\RR^d}h_t(\eta)(m_\eta a-a) d\eta
\quad\text{in }\calS(\RR^d).
$$
For $0<t\leq1$, the change of variables $\eta=\sqrt t u$ gives
\begin{align*}
\mathfrak p_N\big((e^{-t|\cdot|^2}-1)a\big)
\leq&\int_{\RR^d}h_t(\eta)\mathfrak p_N((m_\eta-1)a) d\eta\\
\lesssim&_N\mathfrak p_{N+1}(a)\int_{\RR^d}h_t(\eta)|\eta|(1+|\eta|)^N d\eta\\
=&(4\pi)^{-\frac d2} t^{1/2} \mathfrak p_{N+1}(a)\int_{\RR^d}e^{-\frac{|u|^2}{4}}|u|(1+\sqrt t |u|)^N du\\
\leq&(4\pi)^{-\frac d2} t^{1/2} \mathfrak p_{N+1}(a)\int_{\RR^d}e^{-\frac{|u|^2}{4}}|u|(1+|u|)^N du\\
\lesssim&_N t^{1/2} \mathfrak p_{N+1}(a)\to0\quad \text{as} \ t\to0^+.
\end{align*}
This yields
$$
e^{-t|\cdot|^2}a\to a\quad\text{in }\calS(\RR^d).
$$
Let $\mathcal A$ be a bounded subset of $\calS(\RR^d)$. By definition,
$$
\sup_{a\in\mathcal A}\mathfrak p_{N+1}(a)<\infty,\quad N\in\NN_0.
$$
For every $N\in\NN_0$ and $0<t\leq1$, the preceding estimate gives
\begin{align*}
\sup_{a\in\mathcal A}\mathfrak p_N\big(e^{-t|\cdot|^2}a-a\big)
\lesssim&_N t^{1/2} \sup_{a\in\mathcal A}\mathfrak p_{N+1}(a)
\to0\quad \text{as} \ t\to0^+.
\end{align*}

Since $\widehat x\in\calS'(\RR^d)$, we can choose some $N\in\NN_0$ such that
$$
|\langle\widehat x,a\rangle|\lesssim_x\mathfrak p_N(a),\quad a\in\calS(\RR^d).
$$
By the definition of the Fourier transform for tempered distributions on quantum Euclidean spaces, we obtain
\begin{align*}
\langle y_t-x,U_\theta(a(-\cdot))\rangle
=&(2\pi)^d\langle\widehat{y_t}-\widehat x,a\rangle\\
=&(2\pi)^d\langle\widehat x,(e^{-t|\cdot|^2}-1)a\rangle.
\end{align*}
Therefore, the previous estimate gives, for $0<t\leq1$,
\begin{align*}
\sup_{a\in\mathcal A}|\langle y_t-x,U_\theta(a(-\cdot))\rangle|
\lesssim&_x\sup_{a\in\mathcal A}\mathfrak p_N\big((e^{-t|\cdot|^2}-1)a\big)\\
\lesssim&_{x,N} t^{1/2} \sup_{a\in\mathcal A}\mathfrak p_{N+1}(a)
\to0\quad(t\to0^+).
\end{align*}
The map $a\mapsto U_\theta(a(-\cdot))$ is a topological isomorphism from $\calS(\RR^d)$ onto $\calS(\RR^d_\theta)$. Hence every bounded subset of $\calS(\RR^d_\theta)$ is the image of a bounded subset of $\calS(\RR^d)$ under this map. We conclude that $y_t\to x$ uniformly on bounded subsets of $\calS(\RR^d_\theta)$, that is, with respect to the strong topology of $\calS'(\RR^d_\theta)$. This proves (i).

\

\noindent (ii): It is immediate that (ii) is obtained from \eqref{eq:avform} with $\psi=h_t$. 

\

\noindent (iii): Since $1<\frac{d}{1-\epsilon}<\infty$, by \eqref{eq:bn1} and Lemma \ref{lem:avg}, one has
\begin{align*}
\|[(-\Delta)^{\frac{\epsilon}{2}},y_t]\|_{\calL_{\frac{d}{1-\epsilon},\infty}(L_2(\RR^d))}
\leq&\|h_t*[(-\Delta)^{\frac{\epsilon}{2}},x]\|_{\calL_{\frac{d}{1-\epsilon},\infty}(L_2(\RR^d))}^{\#}\\
\leq&\|[(-\Delta)^{\frac{\epsilon}{2}},x]\|_{\calL_{\frac{d}{1-\epsilon},\infty}(L_2(\RR^d))}^{\#}\\
\lesssim&_{d,\epsilon}
\|[(-\Delta)^{\frac{\epsilon}{2}},x]\|_{\calL_{\frac{d}{1-\epsilon},\infty}(L_2(\RR^d))}.
\end{align*}
This proves (iii).
\end{proof}

\subsection{Approximation by Schwartz elements}

We present the definition of multiplication for distributions on quantum Euclidean spaces. The multiplication of a tempered distribution $u\in \calS'(\RR^d_\theta)$ by $z\in\calS(\RR^d_\theta)$ is defined through
$$
\langle zu,a\rangle=\langle u,az\rangle,
\quad
\langle uz,a\rangle=\langle u,za\rangle,
\quad a\in\calS(\RR^d_\theta).
$$
With this notation in hand, we need the following lemma to handle approximation of distributions, where the net $(z_\delta)_{\delta>0}\subseteq\calS(\RR^d_\theta)$ is needed in the proof of Theorem \ref{thm:m1}.

%This lemma plays an important role in the proof of the main result.

\begin{lemma}\label{lem:ztest}
There exists a net $(z_\delta)_{\delta>0}\subseteq\calS(\RR^d_\theta)$ with $z_\delta^*=z_\delta$ and $\|z_\delta\|_{L_\infty(\RR^d_\theta)}\le1$ such that $z_\delta xz_\delta\to x$ in $\calS(\RR^d_\theta)$ as $\delta\to0$ for any $x\in \calS(\RR^d_\theta)$.
\end{lemma}

\begin{proof}
We choose
$$\varphi(s)=\pi^{-\frac d2}e^{-|s|^2},\quad\varphi_\delta(s)=\delta^{-d}\varphi\big(\frac{s}{\delta}\big),\quad z_\delta=U_\theta(\varphi_\delta),\quad\delta>0.$$
By \eqref{eq:defufh1}, one has
$$z_\delta^*=U_\theta(\overline{\varphi_\delta(-\cdot)})=U_\theta(\varphi_\delta)=z_\delta$$
and
$$\|z_\delta\|_{L_\infty(\RR^d_\theta)}\leq\|\varphi_\delta\|_{L_1(\RR^d)}=\|\varphi\|_{L_1(\RR^d)}=1.$$
It is clear that $z_\delta\in\calS(\RR^d_\theta)$. Let $x=U_\theta(f)$ with $f\in\calS(\RR^d)$. By \eqref{eq:mul1}, for every $y\in\calS(\RR^d_\theta)$, $\widehat{z_\delta y}=\varphi_\delta *_\theta \widehat y$. In particular,
$$\widehat{z_\delta x}=\varphi_\delta *_\theta f,\quad
\widehat{xz_\delta}=f *_\theta \varphi_\delta.$$
It follows from \eqref{eq:con1} that for $t\in\RR^d$,
\begin{align*}
(\varphi_\delta *_\theta f)(t)
=\int_{\RR^d}e^{\frac{\ii\delta}{2}\langle t,\theta u\rangle}
\varphi(u)f(t-\delta u) du
\end{align*}
and
\begin{align*}
(f *_\theta \varphi_\delta)(t)
=\int_{\RR^d}e^{-\frac{\ii\delta}{2}\langle t,\theta u\rangle}
\varphi(u)f(t-\delta u) du.
\end{align*}

Let $N\in\NN_0$. We show that for $0<\delta\le1$ and $h\in\calS(\RR^d)$,
\begin{equation}\label{eq:pleft}
\mathfrak p_N(\varphi_\delta *_\theta h)\lesssim_{N,\theta}\mathfrak p_N(h).
\end{equation}
Indeed, by the Leibniz rule, for $|\alpha|\le N$,
\begin{align*}
\partial_t^\alpha(\varphi_\delta *_\theta h)(t)
=&\int_{\RR^d}\varphi(u)\sum_{\gamma\le\alpha}\binom{\alpha}{\gamma}\Big(\frac{\ii\delta}{2}\Big)^{|\gamma|}(\theta u)^\gamma
e^{\frac{\ii\delta}{2}\langle t,\theta u\rangle}
\partial^{\alpha-\gamma}h(t-\delta u) du.
\end{align*}
For $0<\delta\le1$, since
$$(1+|t|)^N\le(1+|t-\delta u|)^N(1+|\delta u|)^N
\le(1+|u|)^N(1+|t-\delta u|)^N,$$
one has \eqref{eq:pleft}.

For $|\alpha|\le N$, note that by the Leibniz rule again,
\begin{align*}
&\partial_t^\alpha\big(e^{\pm\frac{\ii\delta}{2}\langle t,\theta u\rangle}f(t-\delta u)-f(t)\big)\\
=&\delta\int_0^1\sum_{\gamma\le\alpha}\binom{\alpha}{\gamma}
\Big(\pm\frac{\ii v\delta}{2}\Big)^{|\gamma|}(\theta u)^\gamma
e^{\pm\frac{\ii v\delta}{2}\langle t,\theta u\rangle}\times\Big[\pm\frac{\ii}{2}\langle t,\theta u\rangle \partial^{\alpha-\gamma}f(t-v\delta u)\\
&\pm\frac{\ii}{2}\sum_{j=1}^d(\alpha-\gamma)_j(\theta u)_j\partial^{\alpha-\gamma-e_j}f(t-v\delta u)-\sum_{j=1}^d u_j\partial^{\alpha-\gamma+e_j}f(t-v\delta u)\Big] dv,
\end{align*}
where $e_j=(0,\ldots,0,1,0,\ldots,0)$ is the $j$-th coordinate unit vector of $\RR^d$. Then for every $|\alpha|\le N$,
\begin{align*}
&(1+|t|)^N\big|\partial_t^\alpha(\varphi_\delta *_\theta f-f)(t)\big|\\
\lesssim&_{N,\theta}\delta\mathfrak p_{N+1}(f)
\int_{\RR^d}\varphi(u)(1+|u|)^{N+1}
\sum_{\gamma\le\alpha}|\theta u|^{|\gamma|}
\big(1+|\theta u|+|u|\big) du\\
\lesssim&_{N,\theta}\delta\mathfrak p_{N+1}(f).
\end{align*}
Hence, we have
\begin{align*}
\mathfrak p_N(\varphi_\delta *_\theta f-f)
\lesssim_{N,\theta}\delta\mathfrak p_{N+1}(f).
\end{align*}
In the same way, we obtain
\begin{align*}
\mathfrak p_N(f *_\theta \varphi_\delta-f)
\lesssim_{N,\theta}\delta\mathfrak p_{N+1}(f).
\end{align*}
By using \eqref{eq:pleft} with $h=\widehat{xz_\delta-x}$, 
\begin{align*}
\mathfrak p_N(\widehat{z_\delta xz_\delta}-f)
\le&\mathfrak p_N\big(\widehat{z_\delta(xz_\delta-x)}\big)
+\mathfrak p_N(\widehat{z_\delta x}-f)\\
=&\mathfrak p_N\big(\varphi_\delta *_\theta(\widehat{xz_\delta-x})\big)
+\mathfrak p_N(\varphi_\delta *_\theta f-f)\\
\lesssim&_{N,\theta}\mathfrak p_N(\widehat{xz_\delta}-f)
+\mathfrak p_N(\widehat{z_\delta x}-f)\\
\lesssim&_{N,\theta}\delta\mathfrak p_{N+1}(f)\to0
\quad\text{as }\delta\to0.
\end{align*}
This completes the proof.
\end{proof}

By virtue of the heat kernels, we can construct the mollification of symbols for tempered distributions on quantum Euclidean spaces.

\begin{lemma}\label{lem:hprod}
Let $x\in\calS'(\RR^d_\theta)$ and $y_t=h_t*x$ for $t>0$. For every $z\in\calS(\RR^d_\theta)$ and every multi-index $\alpha\in\NN_0^d$, one has
$$
z\partial^\alpha y_t,\ (\partial^\alpha y_t)z\in\calS(\RR^d_\theta).
$$
For fixed $t,z,\alpha$, both maps from $x$ to these products are continuous from the strong topology of $\calS'(\RR^d_\theta)$ to $\calS(\RR^d_\theta)$.
\end{lemma}

\begin{proof}
For every $a\in\calS(\RR^d)$,
\begin{align*}
\partial^\alpha U_\theta(a(-\cdot))
=&U_\theta((\cdot)^\alpha a(-\cdot))\\
=&(-1)^{|\alpha|}U_\theta\big(((\cdot)^\alpha a)(-\cdot)\big).
\end{align*}
Then one has
\begin{align*}
\langle\widehat{\partial^\alpha y_t},a\rangle
=&(2\pi)^{-d}\langle\partial^\alpha y_t,U_\theta(a(-\cdot))\rangle\\
=&(-1)^{|\alpha|}(2\pi)^{-d}\langle y_t,\partial^\alpha U_\theta(a(-\cdot))\rangle\\
=&(2\pi)^{-d}\big\langle y_t,U_\theta\big(((\cdot)^\alpha a)(-\cdot)\big)\big\rangle\\
=&\langle\widehat{y_t},(\cdot)^\alpha a\rangle\\
=&\langle\widehat x,e^{-t|\cdot|^2}(\cdot)^\alpha a\rangle\\
=&\langle(\cdot)^\alpha \big(e^{-t|\cdot|^2}\widehat x\big),a\rangle.
\end{align*}
Since $a\in\calS(\RR^d)$ is arbitrary, we obtain
$$
\widehat{\partial^\alpha y_t}=(\cdot)^\alpha e^{-t|\cdot|^2}\widehat x
\quad\text{in }\calS'(\RR^d).
$$
Thus the Fourier transform of $\partial^\alpha y_t$ is $s^\alpha e^{-t|s|^2}\widehat x(s)$. 

Write $z=U_\theta(f)$ with $f\in\calS(\RR^d)$ and define
$$
K_\pm(r,s)=e^{-t|s|^2}f(r-s)e^{\pm\frac{\ii}{2}\langle r,\theta s\rangle},
\quad r,s\in\RR^d.
$$
It is clear that $K_\pm\in\calS(\RR^{2d})$. For each $r\in\RR^d$, the function $(\cdot)^\alpha K_\pm(r,\cdot)$ belongs to $\calS(\RR^d)$. Therefore, we can define
$$
k_{\pm,\alpha}(r)=\langle\widehat x,(\cdot)^\alpha K_\pm(r,\cdot)\rangle,
$$
where the distributional pairing is taken in the second variable of $K_\pm$ with $r$ fixed. Since $K_\pm\in\calS(\RR^{2d})$, for every $M\in\NN_0$ and multi-index $\beta$,
$$
\sup_{r\in\RR^d}\mathfrak p_M\big((\cdot)^\alpha\partial_r^\beta K_\pm(r,\cdot)\big)<\infty,
$$
where $\mathfrak p_M$ is taken in the variable $s$. Let $e_j$ denote the $j$-th coordinate vector. For $h\in\RR\setminus\{0\}$, note that
\begin{align*}
&s^\alpha\Big(\frac{\partial_r^\beta K_\pm(r+he_j,s)-\partial_r^\beta K_\pm(r,s)}{h}-\partial_r^{\beta+e_j}K_\pm(r,s)\Big)\\
=&h\int_0^1(1-u)s^\alpha\partial_r^{\beta+2e_j}K_\pm(r+uhe_j,s) du.
\end{align*}
Consequently, as $h\to0$,
\begin{align*}
&\mathfrak p_M\Big((\cdot)^\alpha\Big(\frac{\partial_r^\beta K_\pm(r+he_j,\cdot)-\partial_r^\beta K_\pm(r,\cdot)}{h}-\partial_r^{\beta+e_j}K_\pm(r,\cdot)\Big)\Big)\\
\leq&\frac{|h|}{2}\sup_{v\in\RR^d}\mathfrak p_M\big((\cdot)^\alpha\partial_r^{\beta+2e_j}K_\pm(v,\cdot)\big)
\to0.
\end{align*}
Thus the difference quotient converges in $\calS(\RR^d)$ in the variable $s$. By the continuity of $\widehat x$ on $\calS(\RR^d)$,
\begin{align*}
&\frac{\partial}{\partial r_j}\big\langle\widehat x,s^\alpha\partial_r^\beta K_\pm(r,s)\big\rangle\\
=&\lim_{h\to0}\Big\langle\widehat x,s^\alpha\frac{\partial_r^\beta K_\pm(r+he_j,s)-\partial_r^\beta K_\pm(r,s)}{h}\Big\rangle\\
=&\big\langle\widehat x,s^\alpha\partial_r^{\beta+e_j}K_\pm(r,s)\big\rangle.
\end{align*}
These pairings are continuous in $r$. Indeed, as $r'\to r$,
\begin{align*}
&\mathfrak p_M\big((\cdot)^\alpha(\partial_r^\beta K_\pm(r',\cdot)-\partial_r^\beta K_\pm(r,\cdot))\big)\\
\leq&\sum_{j=1}^d|r'_j-r_j|\sup_{v\in\RR^d}\mathfrak p_M\big((\cdot)^\alpha\partial_r^{\beta+e_j}K_\pm(v,\cdot)\big)
\to0.
\end{align*}
Applying these identities successively, we obtain $k_{\pm,\alpha}\in C^\infty(\RR^d)$ and
$$
\partial_r^\beta k_{\pm,\alpha}(r)=\big\langle\widehat x,s^\alpha\partial_r^\beta K_\pm(r,s)\big\rangle
$$
for every multi-index $\beta$. However, note that there exist $M\in\NN_0$ and $C_x>0$ such that, for every $N\in\NN_0$ and multi-index $\beta$,
\begin{align*}
&\sup_{r\in\RR^d}(1+|r|)^N
|\partial_r^\beta k_{\pm,\alpha}(r)|\\
\leq&C_x\sum_{|\gamma|\leq M}
\sup_{r,s\in\RR^d}(1+|r|)^N(1+|s|)^M
\big|\partial_s^\gamma\partial_r^\beta(s^\alpha K_\pm(r,s))\big|
<\infty.
\end{align*}
This implies $k_{\pm,\alpha}\in\calS(\RR^d)$.

Let $a\in\calS(\RR^d)$. By the definition of multiplication on distributions and \eqref{eq:mul1},
\begin{align*}
\langle z\partial^\alpha y_t,U_\theta(a)\rangle
=&\langle\partial^\alpha y_t,U_\theta(a)z\rangle\\
=&(2\pi)^d\Big\langle\widehat x,s^\alpha e^{-t|s|^2}
\int_{\RR^d}f(r-s)e^{-\frac{\ii}{2}\langle r,\theta s\rangle}a(-r) dr\Big\rangle\\
=&(2\pi)^d\int_{\RR^d}k_{-,\alpha}(r)a(-r) dr.
\end{align*}
Similarly, in the same way, one has
$$
\langle(\partial^\alpha y_t)z,U_\theta(a)\rangle
=(2\pi)^d\int_{\RR^d}k_{+,\alpha}(r)a(-r) dr.
$$
We will identify these products by their pairings with elements in $\calS(\RR^d_\theta)$. Since $k_{\pm,\alpha}\in\calS(\RR^d)$, the multiplication formula \eqref{eq:mul1} and the trace normalization \eqref{eq:fun1} imply that for every $a\in\calS(\RR^d)$,
\begin{align*}
\langle U_\theta(k_{\pm,\alpha}),U_\theta(a)\rangle
=&\tau_\theta\big(U_\theta(k_{\pm,\alpha})U_\theta(a)\big)\\
=&\tau_\theta\big(U_\theta(k_{\pm,\alpha}*_\theta a)\big)\\
=&(2\pi)^d(k_{\pm,\alpha}*_\theta a)(0)\\
=&(2\pi)^d\int_{\RR^d}k_{\pm,\alpha}(-s)a(s) ds\\
=&(2\pi)^d\int_{\RR^d}k_{\pm,\alpha}(r)a(-r) dr.
\end{align*}
Combining this identity with the preceding pairing formulae, we obtain
\begin{align*}
\langle z\partial^\alpha y_t-U_\theta(k_{-,\alpha}),U_\theta(a)\rangle=0
\end{align*}
and
\begin{align*}
\langle(\partial^\alpha y_t)z-U_\theta(k_{+,\alpha}),U_\theta(a)\rangle=0.
\end{align*}
Hence, we conclude that
$$
z\partial^\alpha y_t=U_\theta(k_{-,\alpha})\in\calS(\RR^d_\theta),
\quad
(\partial^\alpha y_t)z=U_\theta(k_{+,\alpha})\in\calS(\RR^d_\theta).
$$

To prove continuity, suppose that $x_n\to x$ strongly in $\calS'(\RR^d_\theta)$. For every bounded subset $\mathcal A$ of $\calS(\RR^d)$, the set $\{U_\theta(a(-\cdot)):a\in\mathcal A\}$ is also bounded in $\calS(\RR^d_\theta)$. Then
\begin{align*}
\sup_{a\in\mathcal A}|\langle\widehat{x_n}-\widehat x,a\rangle|
=&(2\pi)^{-d}\sup_{a\in\mathcal A}|\langle x_n-x,U_\theta(a(-\cdot))\rangle|
\to0.
\end{align*}
Thus $\widehat{x_n}\to\widehat x$ uniformly on every bounded subset of $\calS(\RR^d)$.

We fix $N\in\NN_0$ and a multi-index $\beta$. Since $K_\pm\in\calS(\RR^{2d})$, for every $M\in\NN_0$,
\begin{align*}
&\sup_{r\in\RR^d}\mathfrak p_M\big((1+|r|)^N\partial_r^\beta((\cdot)^\alpha K_\pm(r,\cdot))\big)\\
=&\max_{|\gamma|\leq M}\sup_{r,s\in\RR^d}(1+|r|)^N(1+|s|)^M\big|\partial_s^\gamma\partial_r^\beta(s^\alpha K_\pm(r,s))\big|
<\infty.
\end{align*}
Hence the family
$$
\big\{s\mapsto(1+|r|)^N\partial_r^\beta(s^\alpha K_\pm(r,s)):r\in\RR^d\big\}
$$
is bounded in $\calS(\RR^d)$. Applying the preceding uniform convergence to this family, we obtain that as $n\to \infty$,
\begin{align*}
&\sup_{r\in\RR^d}(1+|r|)^N\big|\partial_r^\beta(k_{\pm,\alpha,n}-k_{\pm,\alpha})(r)\big|\\
=&\sup_{r\in\RR^d}\big|\big\langle\widehat{x_n}-\widehat x,(1+|r|)^N\partial_r^\beta(s^\alpha K_\pm(r,s))\big\rangle\big|
\to0.
\end{align*}
Since there are only finitely many multi-indices $\beta$ with $|\beta|\leq N$, we conclude that
$$
\mathfrak p_N(k_{\pm,\alpha,n}-k_{\pm,\alpha})
=\max_{|\beta|\leq N}\sup_{r\in\RR^d}(1+|r|)^N\big|\partial_r^\beta(k_{\pm,\alpha,n}-k_{\pm,\alpha})(r)\big|
\to0.
$$
As $N$ is arbitrary, $k_{\pm,\alpha,n}\to k_{\pm,\alpha}$ in $\calS(\RR^d)$ as $n\to\infty$. This proves the desired continuity.
\end{proof}

Using the preceding lemma, we obtain the following result in the two-sided case.
\begin{proposition}\label{prop:hprod}
Let $x\in\calS'(\RR^d_\theta)$ and $y_t=h_t*x$ for $t>0$. For every $w,z\in\calS(\RR^d_\theta)$ and every multi-index $\alpha\in\NN_0^d$, one has
$$
w(\partial^\alpha y_t) z\in\calS(\RR^d_\theta).
$$
For fixed $t,w,z,\alpha$, the map $x\mapsto w\partial^\alpha y_t z$ is continuous from the strong topology of $\calS'(\RR^d_\theta)$ to $\calS(\RR^d_\theta)$.
\end{proposition}

\begin{proof}
By Lemma \ref{lem:hprod}, we have $w(\partial^\alpha y_t) z\in\calS(\RR^d_\theta)$. To prove continuity, suppose $x_n\to x$ strongly in $\calS'(\RR^d_\theta)$. Let $y_{t,n}=h_t*x_n$. Write $z=U_\theta(f)$ with $f\in\calS(\RR^d)$. By \eqref{eq:mul1},
\begin{align*}
w\partial^\alpha y_{t,n} z
=\big(w\partial^\alpha y_{t,n}\big)z
=U_\theta\Big(U_\theta^{-1}\big(w\partial^\alpha y_{t,n}\big)*_\theta f\Big),
\end{align*}
and
\begin{align*}
w\partial^\alpha y_t z
=\big(w\partial^\alpha y_t\big)z
=U_\theta\Big(U_\theta^{-1}\big(w\partial^\alpha y_t\big)*_\theta f\Big).
\end{align*}
Note that for every $h\in\calS(\RR^d)$ and every multi-index $\beta$,
\begin{align*}
\partial_r^\beta(h*_\theta f)(r)
=&\sum_{\gamma\leq\beta}\binom{\beta}{\gamma}
\Big(-\frac{\ii}{2}\Big)^{|\gamma|}
\int_{\RR^d}(\theta s)^\gamma
e^{-\frac{\ii}{2}\langle r,\theta s\rangle}
(\partial^{\beta-\gamma}h)(r-s)f(s) ds.
\end{align*}
By using $1+|r|\leq(1+|r-s|)(1+|s|)$ and choosing $L>N+|\beta|+d$, we obtain that for every $N\in\NN_0$, 
\begin{align*}
\sup_{r\in\RR^d}(1+|r|)^N
\big|\partial_r^\beta(h*_\theta f)(r)\big|
\lesssim_{N,\beta,f,\theta}\mathfrak p_L(h),
\quad h\in\calS(\RR^d).
\end{align*}
Thus the map $h\mapsto h*_\theta f$ is continuous on $\calS(\RR^d)$. By Lemma \ref{lem:hprod},
$$
w\partial^\alpha y_{t,n}\to w\partial^\alpha y_t
\quad\text{in }\calS(\RR^d_\theta).
$$
Then the map $x\mapsto w\partial^\alpha y_t z$ is continuous from the strong topology of $\calS'(\RR^d_\theta)$ to $\calS(\RR^d_\theta)$. This completes the proof.
\end{proof}

Finally, we establish the following quite important approximation property for tempered distributions by terms of $\calS(\RR^d_\theta)$. We would like to remark that the following lemma is much stronger than Lemma \ref{lem:heat} as one can approximate $\calS'(\RR^d_\theta)$ by virtue of $\calS(\RR^d_\theta)$ with the help of the following lemma.

\begin{lemma}\label{lem:sapp}
For every $x\in\calS'(\RR^d_\theta)$, there exists a sequence $(x_n)_{n\geq 1}\subseteq\calS(\RR^d_\theta)$ such that $x_n\to x$ uniformly on bounded subsets of $\calS(\RR^d_\theta)$ with respect to the strong topology in $\calS'(\RR^d_\theta)$.
\end{lemma}

\begin{proof}
Choose real even functions $\rho,\chi\in C_c^\infty(\RR^d)$ such that $\int_{\RR^d}\rho=1$ and $\chi=1$ on a neighborhood of zero. Let
$$
\rho_n(s)=n^d\rho(ns),
\quad
F_n(s)=\chi\big(\frac{s}{n}\big)(\rho_n*\widehat x)(s),
\quad\text{and}\quad
x_n=U_\theta(F_n).
$$
Here the convolution is the usual convolution of a test function with a tempered distribution in $\calS'(\RR^d)$, defined by
$$
(\rho_n*\widehat x)(s)=\langle\widehat x,\rho_n(s-\cdot)\rangle,\quad s\in\RR^d.
$$
Since $\rho_n$ is even, for every $a\in\calS(\RR^d)$, we have
\begin{align*}
\langle\rho_n*\widehat x,a\rangle
=&\int_{\RR^d}a(s)\langle\widehat x,\rho_n(s-\cdot)\rangle ds\\
=&\Big\langle\widehat x,\int_{\RR^d}a(s)\rho_n(s-\cdot) ds\Big\rangle\\
=&\Big\langle\widehat x,\int_{\RR^d}\rho_n(\cdot-s)a(s) ds\Big\rangle\\
=&\langle\widehat x,\rho_n*a\rangle.
\end{align*}
Since $\widehat x\in\calS'(\RR^d)$, we can choose $M\in\NN_0$ such that
$$
|\langle\widehat x,a\rangle|\lesssim_x\mathfrak p_M(a),\quad a\in\calS(\RR^d).
$$
The map $s\mapsto\rho_n(s-\cdot)$ is smooth with values in $\calS(\RR^d)$. Differentiating under the distributional pairing, we obtain
$$
\partial_s^\beta(\rho_n*\widehat x)(s)=\langle\widehat x,(\partial^\beta\rho_n)(s-\cdot)\rangle,\quad \beta\in\NN_0^d.
$$
By using $1+|u|\leq(1+|s|)(1+|s-u|)$, one has
\begin{align*}
|\partial_s^\beta(\rho_n*\widehat x)(s)|
\lesssim&_x\mathfrak p_M\big((\partial^\beta\rho_n)(s-\cdot)\big)\\
=&\max_{|\gamma|\leq M}\sup_{u\in\RR^d}(1+|u|)^M|\partial^{\beta+\gamma}\rho_n(s-u)|\\
\leq&(1+|s|)^M\max_{|\gamma|\leq M}\sup_{v\in\RR^d}(1+|v|)^M|\partial^{\beta+\gamma}\rho_n(v)|\\
\leq&(1+|s|)^M\mathfrak p_{M+|\beta|}(\rho_n).
\end{align*}
Thus $\rho_n*\widehat x$ is smooth, and every derivative has at most polynomial growth. By Leibniz's rule,
$$
\partial^\beta F_n(s)=\sum_{\gamma\leq\beta}\binom{\beta}{\gamma}n^{-|\gamma|}(\partial^\gamma\chi)\big(\tfrac{s}{n}\big)\langle\widehat x,(\partial^{\beta-\gamma}\rho_n)(s-\cdot)\rangle.
$$
Together with
$$
\operatorname{supp}F_n\subseteq n\operatorname{supp}\chi,
$$
we obtain
$$
F_n\in C_c^\infty(\RR^d)\subseteq\calS(\RR^d).
$$
For every $a\in\calS(\RR^d)$, direct calculation yields
\begin{align*}
\langle F_n,a\rangle
=&\int_{\RR^d}(\rho_n*\widehat x)(s)\chi\big(\frac{s}{n}\big)a(s) ds\\
=&\Big\langle\rho_n*\widehat x,\chi\big(\frac{\cdot}{n}\big)a\Big\rangle\\
=&\Big\langle\widehat x,\rho_n*\Big(\chi\big(\frac{\cdot}{n}\big)a\Big)\Big\rangle.
\end{align*}

All implicit constants below may depend on the fixed functions $\rho,\chi$ and on $d$. Since $\int_{\RR^d}\rho(u) du=1$, for every $|\alpha|\leq N$,
\begin{align*}
\partial_s^\alpha(\rho_n*a-a)(s)
=&\int_{\RR^d}\rho(u)\Big(\partial^\alpha a\big(s-\frac{u}{n}\big)-\partial^\alpha a(s)\Big) du\\
=&-\frac1n\sum_{j=1}^d\int_{\RR^d}\rho(u)u_j\int_0^1\partial^{\alpha+e_j}a\big(s-\frac{vu}{n}\big) dv du.
\end{align*}
For $n\geq1$ and $0\leq v\leq1$, by the triangle inequality,
$$
1+|s|\leq\Big(1+\big|s-\frac{vu}{n}\big|\Big)(1+|u|).
$$
Consequently,
\begin{align*}
\mathfrak p_N(\rho_n*a-a)
\leq&\frac{\sqrt d}{n}\mathfrak p_{N+1}(a)\int_{\RR^d}|\rho(u)||u|(1+|u|)^N du\\
\lesssim&_N\frac1n\mathfrak p_{N+1}(a).
\end{align*}
By Leibniz's rule,
\begin{align*}
\mathfrak p_{N+1}\Big(\chi\big(\frac{\cdot}{n}\big)a\Big)
\leq&\max_{|\alpha|\leq N+1}\sum_{\beta\leq\alpha}\binom{\alpha}{\beta}n^{-|\beta|}\|\partial^\beta\chi\|_{L_\infty(\RR^d)}\mathfrak p_{N+1}(a)\\
\lesssim&_N\mathfrak p_{N+1}(a).
\end{align*}
Applying the convolution estimate to $\chi(\frac{\cdot}{n})a$, we obtain
\begin{align*}
\mathfrak p_N\Big(\rho_n*\Big(\chi\big(\frac{\cdot}{n}\big)a\Big)-\chi\big(\frac{\cdot}{n}\big)a\Big)
\lesssim&_N\frac1n\mathfrak p_{N+1}\Big(\chi\big(\frac{\cdot}{n}\big)a\Big)\\
\lesssim&_N\frac1n\mathfrak p_{N+1}(a).
\end{align*}
Choose $c>0$ such that $\chi(s)=1$ whenever $|s|\leq c$. For $|\alpha|\leq N$,
\begin{align*}
\partial_s^\alpha\Big(\Big(\chi\big(\frac{s}{n}\big)-1\Big)a(s)\Big)
=&\Big(\chi\big(\frac{s}{n}\big)-1\Big)\partial^\alpha a(s)\\
&+\sum_{\substack{\beta\leq\alpha\\|\beta|\geq1}}\binom{\alpha}{\beta}n^{-|\beta|}(\partial^\beta\chi)\big(\frac{s}{n}\big)\partial^{\alpha-\beta}a(s).
\end{align*}
Every term vanishes when $|s|<cn$. Hence
\begin{align*}
\mathfrak p_N\Big(\Big(\chi\big(\frac{\cdot}{n}\big)-1\Big)a\Big)
\lesssim&_N\max_{|\alpha|\leq N}\sup_{|s|\geq cn}(1+|s|)^N|\partial^\alpha a(s)|\\
\leq&\mathfrak p_{N+1}(a)\sup_{|s|\geq cn}\frac1{1+|s|}\\
=&\frac1{1+cn}\mathfrak p_{N+1}(a)\\
\lesssim&_c\frac1n\mathfrak p_{N+1}(a).
\end{align*}
Combining the two estimates, we conclude that
\begin{align*}
\mathfrak p_N\Big(\rho_n*\Big(\chi\big(\frac{\cdot}{n}\big)a\Big)-a\Big)
\leq&\mathfrak p_N\Big(\rho_n*\Big(\chi\big(\frac{\cdot}{n}\big)a\Big)-\chi\big(\frac{\cdot}{n}\big)a\Big)\\
&+\mathfrak p_N\Big(\Big(\chi\big(\frac{\cdot}{n}\big)-1\Big)a\Big)\\
\lesssim&_{N,c}\frac1n\mathfrak p_{N+1}(a).
\end{align*}

Let $\mathcal A$ be a bounded subset of $\calS(\RR^d)$. Then
$$
\sup_{a\in\mathcal A}\mathfrak p_{N+1}(a)<\infty,\quad N\in\NN_0,
$$
and therefore as $n\to \infty$,
$$
\sup_{a\in\mathcal A}\mathfrak p_N\Big(\rho_n*\Big(\chi\big(\frac{\cdot}{n}\big)a\Big)-a\Big)
\lesssim_{N,c}\frac1n\sup_{a\in\mathcal A}\mathfrak p_{N+1}(a)\to0.
$$
By using the preceding formula for $\langle F_n,a\rangle$, we obtain
\begin{align*}
\sup_{a\in\mathcal A}|\langle F_n-\widehat x,a\rangle|
=&\sup_{a\in\mathcal A}\Big|\Big\langle\widehat x,\rho_n*\Big(\chi\big(\frac{\cdot}{n}\big)a\Big)-a\Big\rangle\Big|\\
\lesssim&_x\sup_{a\in\mathcal A}\mathfrak p_M\Big(\rho_n*\Big(\chi\big(\frac{\cdot}{n}\big)a\Big)-a\Big)\\
\lesssim&_{x,M,c}\frac1n\sup_{a\in\mathcal A}\mathfrak p_{M+1}(a)\to0.
\end{align*}
Thus $F_n\to\widehat x$ strongly in $\calS'(\RR^d)$. Since $x_n=U_\theta(F_n)$, one has
\begin{align*}
\sup_{a\in\mathcal A}|\langle x_n-x,U_\theta(a(-\cdot))\rangle|
=&(2\pi)^d\sup_{a\in\mathcal A}|\langle F_n-\widehat x,a\rangle|
\to0.
\end{align*}
The map $a\mapsto U_\theta(a(-\cdot))$ is a topological isomorphism between the two Schwartz spaces. We conclude that $x_n\to x$ strongly in $\calS'(\RR^d_\theta)$.
\end{proof}

\bigskip

\section{Cwikel estimates}\label{s:cw}
In this section, we establish the Cwikel estimates on quantum Euclidean spaces associated with arbitrary antisymmetric real matrices. To prove the upper bound in Theorem \ref{thm:m1}, we need the Cwikel estimates of \cite{LSZ2020} for the strong Schatten class ($2\leq p<\infty$) and the weak Schatten class ($2<p<\infty$). We state below the version needed here, following \cite{LSZ2020,MSX2020}. The nondegeneracy condition required in \cite{LSZ2020} is removed here, and for the reader's convenience the proof is given in Appendix \ref{sec:app2}.

\begin{theorem}\label{thm:cw-hi}
Let $E$ be a symmetric interpolation space for $(L_2(0,\infty),L_\infty(0,\infty))$. Let $x\otimes g\in E(L_\infty(\RR^d_\theta)\overline{\otimes}L_\infty(\RR^d))$. Then $xM_g\in E(L_2(\RR^d))$ and
$$\|xM_g\|_{E(L_2(\RR^d))} \lesssim_E\|x\otimes g\|_{E(L_\infty(\RR^d_\theta)\overline{\otimes}L_\infty(\RR^d))}.$$
In particular, for $2\leq p<\infty$,
\begin{equation}\label{eq:cw-q}
\|xM_g\|_{\calL_p(L_2(\RR^d))}
\lesssim_p\|x\|_{L_p(\RR^d_\theta)}
\|g\|_{L_p(\RR^d)},
\end{equation}
and, for $2<p<\infty$,
\begin{equation}\label{eq:cw-wq}
\|xM_g\|_{\calL_{p,\infty}(L_2(\RR^d))}
\lesssim_p\|x\|_{L_p(\RR^d_\theta)}
\|g\|_{L_{p,\infty}(\RR^d)}.
\end{equation}
\end{theorem}

Recall that
$$Q=[0,1)^d,\quad Q_n=n+Q,\quad \text{and}\quad P_n=\mathds{1}_{Q_n}(\ii\nabla)=M_{\mathds{1}_{Q_n}}, \quad  \ \forall n\in\ZZ^d.$$
Let
$$\|g\|_{\ell_p(L_\infty(\RR^d);Q)} =\Big(\sum_{n\in\ZZ^d} \|g\mathds{1}_{n+Q}\|_{L_\infty(\RR^d)}^p\Big)^{1/p}.$$
Denote by $a_n(g)=\|g\mathds{1}_{n+Q}\|_{L_\infty(\RR^d)}$ and let $(a_j^*(g))_{j\geq0}$ be its decreasing rearrangement. Then define
$$\|g\|_{\ell_{p,\infty}(L_\infty(\RR^d);Q)}=\sup_{j\geq0}(j+1)^{1/p}a_j^*(g) =\sup_{\lambda>0}\lambda \#\{n\in\ZZ^d:a_n(g)>\lambda\}^{1/p}.$$

However, the proof of the lower bound in Theorem \ref{thm:m1} requires the Cwikel estimates for $1\leq p<2$. For $1\leq p<2$, we obtain the following Cwikel estimates. 
%Our work builds on this framework by removing the nondegeneracy assumption from the required localization arguments.
The main difficulty is to obtain the local Schatten estimate and the unitary equivalence of cube localizations without assuming that $\theta$ is nondegenerate. The proof below establishes these two facts for arbitrary $\theta$, so that the argument of \cite{LSZ2020} can also apply in the degenerate case.

\begin{theorem}\label{thm:cwl1}
Let $1\leq p<2$. 
\begin{description}
\item[(i)] Suppose that $x\in L_2(\RR^d_\theta)$ and $xP_0\in\calL_p(L_2(\RR^d))$. Then
\begin{equation}\label{eq:cw-lo}
\|xM_g\|_{\calL_p(L_2(\RR^d))}
\lesssim_p
\|xP_0\|_{\calL_p(L_2(\RR^d))}
\|g\|_{\ell_p(L_\infty(\RR^d);Q)},
\quad g\in\ell_p(L_\infty(\RR^d);Q),
\end{equation}
and
\begin{equation}\label{eq:cw-wlo}
\|xM_g\|_{\calL_{p,\infty}(L_2(\RR^d))}
\lesssim_p\|xP_0\|_{\calL_p(L_2(\RR^d))}
\|g\|_{\ell_{p,\infty}(L_\infty(\RR^d);Q)},
\quad g\in\ell_{p,\infty}(L_\infty(\RR^d);Q).
\end{equation}
\item[(ii)] If $x\in W^d_p(\RR^d_\theta)$, then
$$
\|xM_g\|_{\calL_p(L_2(\RR^d))}
\lesssim_{d,p}\|x\|_{W^d_p(\RR^d_\theta)}
\|g\|_{\ell_p(L_\infty(\RR^d);Q)}$$
and
$$\|xM_g\|_{\calL_{p,\infty}(L_2(\RR^d))}
\lesssim_{d,p}\|x\|_{W^d_p(\RR^d_\theta)}
\|g\|_{\ell_{p,\infty}(L_\infty(\RR^d);Q)}.$$
\end{description}

\end{theorem}

For $1\leq p<2$, we extend the framework of Levitina, Sukochev and Zanin \cite[Theorems 7.6--7.7]{LSZ2020} to arbitrary antisymmetric real matrices. First, we combine Dirac-power factorization, coordinate commutators and the high-exponent Cwikel estimates with Schwartz approximation to obtain
$$\|xP_0\|_{\calL_p(L_2(\RR^d))}
\lesssim_{d,p}\|x\|_{W^d_p(\RR^d_\theta)},
\quad x\in W^d_p(\RR^d_\theta),$$
which extends \cite[Lemma 7.4]{LSZ2020} to arbitrary antisymmetric real matrices. Then we replace the normal-form construction of \cite[Lemma 7.5]{LSZ2020} with the weighted right translations
$$(V_\theta(b)\xi)(t)=e^{\frac{\ii}{2}\langle b,\theta t\rangle}\xi(t-b).$$
We verify that these unitaries commute with $L_\infty(\RR^d_\theta)$ and obtain
$$xP_n=V_\theta(n)(xP_0)V_\theta(n)^*,$$
recovering the unitary equivalence of cube localizations for arbitrary $\theta$; see Lemma \ref{lem:shift}.

\subsection{Proof of Theorem \ref{thm:cwl1}}\label{ssec:pfcwl1}
At first, we record the Hilbert--Schmidt identity for arbitrary $\theta$ from \cite[Lemma 7.1]{LSZ2020}. Although $\theta$ is assumed to be nondegenerate in \cite{LSZ2020}, the identity itself does not require such assumption. For the reader's convenience, we include the detailed proof in Appendix \ref{sec:app2}.

\begin{lemma}\label{lem:hs}
For every $x\in L_2(\RR^d_\theta)$ and $g\in L_2(\RR^d)$, 
\begin{equation}\label{exMg1}
\|xM_g\|_{\calL_2(L_2(\RR^d))}=(2\pi)^{-\frac{d}{2}}\|x\|_{L_2(\RR^d_\theta)} \|g\|_{L_2(\RR^d)}.
\end{equation}
\end{lemma}

Then the following lemma follows from the polarization identity.

\begin{lemma}\label{lem:hs2}
Let $x, y \in L_2(\mathbb{R}_\theta^d)$ and $g, h \in L_2(\mathbb{R}^d)$. Then
$$\operatorname{Tr}\Big((x M_g)^* y M_h\Big)
= (2\pi)^{-d} \tau_\theta(x^* y) \int_{\mathbb{R}^d} \overline{g(t)} h(t)  dt.
$$
\end{lemma}
\begin{proof}
Note that
\begin{align*}
\operatorname{Tr}\big((xM_g)^*yM_g\big)
=\frac14\sum_{k=0}^{3} i^{-k}\operatorname{Tr}
\big(((x+i^k y)M_g)^*(x+i^k y)M_g\big).
\end{align*}
By Lemma \ref{lem:hs},
\begin{align*}
\operatorname{Tr}\big((xM_g)^*yM_g\big)
=&(2\pi)^{-d}\int_{\mathbb{R}^d}|g(t)|^2 dt\cdot\Big[\frac14
\sum_{k=0}^{3}i^{-k}\tau_\theta\big((x+i^k y)^*(x+i^k y)\big)\Big]\\
=&(2\pi)^{-d}\tau_\theta(x^*y)\int_{\mathbb{R}^d}|g(t)|^2 dt.
\end{align*}
On the other hand,
\begin{align*}
\operatorname{Tr}\big((xM_g)^*yM_h\big)
=\frac14\sum_{k=0}^{3} i^{-k}\operatorname{Tr}
\big((xM_{g+i^k h})^*yM_{g+i^k h}\big).
\end{align*}
Hence, this implies
$$\operatorname{Tr}\Big((x M_g)^* y M_h\Big)
= (2\pi)^{-d} \tau_\theta(x^* y) \int_{\mathbb{R}^d} \overline{g(t)} h(t)  dt.
$$
This completes the proof.
\end{proof}

By using Lemma \ref{lem:hs} and the general framework of Cwikel estimates developed in \cite{LSZ2020}, one can prove Theorem \ref{thm:cw-hi}; for the reader's convenience, the proof is provided in Appendix \ref{sec:app2}.

Now we extend the low-exponent estimates of \cite[Theorems 7.6--7.7]{LSZ2020} to arbitrary antisymmetric real matrices $\theta$ without the nondegeneracy assumption in the range $1\leq p<2$. We need to adapt \cite[Lemma 7.4]{LSZ2020} and \cite[Lemma 7.5]{LSZ2020} to our setting so that the analogous conclusions hold for degenerate $\theta$, i.e. $\det\theta=0$.

Let $N_d=2^{\lfloor \frac{d}{2}\rfloor}$. Choose self-adjoint matrices $\gamma_1,\ldots,\gamma_d\in M_{N_d}$ such that
$$\gamma_j\gamma_k+\gamma_k\gamma_j=2\delta_{jk}1_{N_d},$$
where $1_{N_d}$ is the identity of $M_{N_d}$. On $\CC^{N_d}\otimes L_2(\RR^d)$, define $\mathscr D=\sum_{j=1}^d\gamma_j\otimes \calD_j$, $J=(1-\Delta)^{\frac{1}{2}}$, $R=\mathscr D-\ii$ and
$$\widetilde J=(1+\mathscr D^2)^{\frac{1}{2}}=1_{N_d}\otimes J;$$
see \cite[Definition 6.17]{LSZ2020}.

\begin{lemma}\label{lem:rp}
For $s\in\RR$, there exists a unitary $\mathscr U$ commuting with $\widetilde J$ such that
\begin{equation}\label{eq:rp}
R^s=\mathscr U^s\widetilde J^s,
\quad
\operatorname{Dom}(R^s)=\operatorname{Dom}(\widetilde J^s).
\end{equation}
Moreover, for $a,b\in\RR$ and $y\in\calS(\RR^d_\theta)$,
\begin{equation}\label{eq:rj}
R^a(1_{N_d}\otimes y)R^b
=\mathscr U^a
\big[\widetilde J^a(1_{N_d}\otimes y)\widetilde J^b\big]
\mathscr U^b.
\end{equation}
\end{lemma}

\begin{proof}
Since $R^*R=RR^*=\widetilde J^2$, we have the following unique polar decomposition of $R$:
$$R=\mathscr U\widetilde J,\quad \mathscr U=R\widetilde J^{-1}.$$
Let $E_{\mathscr D}$ denote the spectral measure of $\mathscr D$ and $\beta(\lambda)=\arctan\lambda-\frac{\pi}{2}$. Then $\lambda-\ii=(1+\lambda^2)^{\frac12}e^{\ii\beta(\lambda)}$, and we deduce that
$$\mathscr U=R\widetilde J^{-1}
=\int_{\RR}\frac{\lambda-\ii}{(1+\lambda^2)^{\frac12}} dE_{\mathscr D}(\lambda)
=e^{\ii\beta(\mathscr D)}$$
is unitary and commutes with $\widetilde J$. It is also clear that by spectral calculus, one has

\begin{align*}
R^s=&\int_{\RR}\exp\big(s\operatorname{Log}(\lambda-\ii)\big)
dE_{\mathscr D}(\lambda)\\
=&\int_{\RR}e^{\ii s\beta(\lambda)}(1+\lambda^2)^{s/2}
dE_{\mathscr D}(\lambda)
=\mathscr U^s\widetilde J^s.
\end{align*}
Moreover,
$$\operatorname{Dom}(R^s)
=\Big\{\xi:\int_{\RR}(1+\lambda^2)^s
d\langle E_{\mathscr D}(\lambda)\xi,\xi\rangle<\infty\Big\}
=\operatorname{Dom}(\widetilde J^s).$$
This proves \eqref{eq:rp}. We obtain \eqref{eq:rj} by substituting $R^a=\mathscr U^a\widetilde J^a$ and $R^b=\widetilde J^b\mathscr U^b$.
\end{proof}

We use Lemma \ref{lem:rp} to prove the following estimate of Schatten classes by induction.

\begin{lemma}\label{lem:rtk1}
Under the same notation as in Lemma \ref{lem:rp}, for $k\in\NN_0$, define
$$T_k(y)=R^{k-1-\frac{d}{2}}(1_{N_d}\otimes y)R^{-k-1-\frac{d}{2}}, \quad y\in\calS(\RR^d_\theta).$$
If $1\leq p<\infty$, then for any $k\in\NN_0$,
\begin{equation}\label{eq:tk}
\|T_k(y)\|_{\calL_p(\CC^{N_d}\otimes L_2(\RR^d))}
\lesssim_{d,p,k}\|y\|_{W^k_p(\RR^d_\theta)}.
\end{equation}

\end{lemma}
\begin{proof}
Let $y\in\calS(\RR^d_\theta)$. For $k=0$, we put $\alpha=1+\frac{d}{2}$. By \eqref{eq:rj},
$$\|T_0(y)\|_{\calL_p(\CC^{N_d}\otimes L_2(\RR^d))}
=N_d^{\frac1p}\|J^{-\alpha}yJ^{-\alpha}\|_{\calL_p(L_2(\RR^d))}.$$
Let $y=u|y|$ be the polar decomposition of $y$, and we factor
$$J^{-\alpha}yJ^{-\alpha}
=\big(J^{-\alpha}u|y|^{\frac12}\big)
\big(|y|^{\frac12}J^{-\alpha}\big).$$
The function $g_\alpha(t)=(1+|t|^2)^{-\frac{\alpha}{2}}$ belongs to $L_{2p}(\RR^d)$ because $2\alpha p>d$. By Theorem \ref{thm:cw-hi},
\begin{align*}
\|J^{-\alpha}u|y|^{\frac12}\|_{\calL_{2p}(L_2(\RR^d))}
\lesssim_{d,p}\||y|^{\frac12}\|_{L_{2p}(\RR^d_\theta)}
\end{align*}
and
\begin{align*}
\||y|^{\frac12}J^{-\alpha}\|_{\calL_{2p}(L_2(\RR^d))}
\lesssim_{d,p}\||y|^{\frac12}\|_{L_{2p}(\RR^d_\theta)}.
\end{align*}
By H\"older's inequality, it yields \eqref{eq:tk} for $k=0$.

Now assume that \eqref{eq:tk} holds for $k$. Let $a_k=k-1-\frac d2$, $b_k=-k-1-\frac d2$. Since $a_{k+1}=a_k+1$ and $b_{k+1}=b_k-1$,
\begin{align*}
T_{k+1}(y)
=&R^{a_k+1}(1_{N_d}\otimes y)R^{b_k-1}\\
=&R^{a_k}(1_{N_d}\otimes y)R^{b_k}+R^{a_k}
\big(R(1_{N_d}\otimes y)-(1_{N_d}\otimes y)R\big)R^{b_k-1}
\\
=&T_k(y)+R^{a_k}[R,1_{N_d}\otimes y]R^{b_k-1}.
\end{align*}
By using $R=\mathscr D-\ii$ and $\partial_jy=[\calD_j,y]$, we obtain
$$[R,1_{N_d}\otimes y]=\sum_{j=1}^d
\big[(\gamma_j\otimes \calD_j),(1_{N_d}\otimes y)\big]
=\sum_{j=1}^d\gamma_j\otimes\partial_jy.$$
By Lemma \ref{lem:rp},
\begin{align*}
W_{j,k}:=&R^{a_k}(\gamma_j\otimes1)R^{-a_k}\\
=&\mathscr U^{a_k}\widetilde J^{a_k}
(\gamma_j\otimes1)\widetilde J^{-a_k}\mathscr U^{-a_k}\\
=&\mathscr U^{a_k}(\gamma_j\otimes1)\mathscr U^{-a_k}.
\end{align*}
Thus $W_{j,k}$ is unitary. Note also that
\begin{align*}
R^{a_k}(\gamma_j\otimes\partial_jy)R^{b_k-1}
=&
R^{a_k}(\gamma_j\otimes1)R^{-a_k}
R^{a_k}(1_{N_d}\otimes\partial_jy)R^{b_k}R^{-1}\\
=&W_{j,k}T_k(\partial_jy)R^{-1}.
\end{align*}
Since $\|R^{-1}\|_{\calB(\CC^{N_d}\otimes L_2(\RR^d))}\leq1$, the induction hypothesis implies
\begin{align*}
\|T_{k+1}(y)\|_{\calL_p(\CC^{N_d}\otimes L_2(\RR^d))}
\leq&
\|T_k(y)\|_{\calL_p(\CC^{N_d}\otimes L_2(\RR^d))}
+\sum_{j=1}^d
\|T_k(\partial_jy)\|_{\calL_p(\CC^{N_d}\otimes L_2(\RR^d))}\\
\lesssim&_{d,p,k}
\|y\|_{W^k_p(\RR^d_\theta)}
+\sum_{j=1}^d\|\partial_jy\|_{W^k_p(\RR^d_\theta)}\\
\lesssim&_{d,p,k}\|y\|_{W^{k+1}_p(\RR^d_\theta)}.
\end{align*}
This completes the proof.
\end{proof}

%Recall that for $n\in\ZZ^d$,
%$$Q=[0,1)^d,\quad Q_n=n+Q\quad\text{and } %P_n=\mathds{1}_{Q_n}.%%(\ii\nabla)=M_{\mathds{1}_{Q_n}}. $$
Next we use Lemma \ref{lem:rp} and Lemma \ref{lem:rtk1} to prove the following Schatten bound of $xP_0$.

\begin{lemma}\label{lem:loc}
Let $1\leq p\leq2$. If $x\in W^d_p(\RR^d_\theta)$, then $x\in L_2(\RR^d_\theta)$ and
\begin{equation}\label{eq:loc}
\|xP_0\|_{\calL_p(L_2(\RR^d))}
\lesssim_{d,p}\|x\|_{W^d_p(\RR^d_\theta)}.
\end{equation}
\end{lemma}
\begin{proof}
First, suppose $x\in\calS(\RR^d_\theta)$ and choose
$m=\lfloor\frac{d+3}{2}\rfloor$, $a=m-1-\frac d2$.
Then $a=0$ for even $d$, $a=\frac12$ for odd $d$ and $m\leq d$. By Lemma \ref{lem:rp} and Lemma \ref{lem:rtk1},
$$\|\widetilde J^a(1_{N_d}\otimes x)
\widetilde J^{-m-1-\frac{d}{2}}\|_{\calL_p(\CC^{N_d}\otimes L_2(\RR^d))}\lesssim_{d,p}\|x\|_{W^m_p(\RR^d_\theta)}
\leq\|x\|_{W^d_p(\RR^d_\theta)}.$$
We use the factorization
$$1_{N_d}\otimes xP_0
=\widetilde J^{-a}
\big[\widetilde J^a(1_{N_d}\otimes x)
\widetilde J^{-m-1-\frac{d}{2}}\big]
\widetilde J^{m+1+\frac{d}{2}}(1_{N_d}\otimes P_0).$$
Note that
\begin{align*}
\|\widetilde J^{-a}\|_{\calB(\CC^{N_d}\otimes L_2(\RR^d))}
\leq1
\end{align*}
and
\begin{align*}
\|\widetilde J^{m+1+\frac{d}{2}}(1_{N_d}\otimes P_0)\|_{\calB(\CC^{N_d}\otimes L_2(\RR^d))}
\leq(1+d)^{\frac{m+1+\frac{d}{2}}{2}}.
\end{align*}
Together with
$\|1_{N_d}\otimes xP_0\|_{\calL_p(\CC^{N_d}\otimes L_2(\RR^d))}
=N_d^{\frac1p}\|xP_0\|_{\calL_p(L_2(\RR^d))}$,
this proves \eqref{eq:loc} for $x\in\calS(\RR^d_\theta)$.

For $x\in W^d_p(\RR^d_\theta)$, by Proposition \ref{prop:id1}, there exists a sequence $(x_\nu)_{\nu\geq 1}\subseteq\calS(\RR^d_\theta)$ converging to $x$ in $W^d_p(\RR^d_\theta)$. By the previous estimate, $\{x_\nu P_0\}$ is a Cauchy sequence in $\calL_p(L_2(\RR^d))$. By Lemma \ref{lem:hs} and $|Q|=1$,
\begin{align*}
\|x_\nu-x_\mu\|_{L_2(\RR^d_\theta)}
=&(2\pi)^{d/2}
\|(x_\nu-x_\mu)P_0\|_{\calL_2(L_2(\RR^d))}\\
\leq&(2\pi)^{d/2}
\|(x_\nu-x_\mu)P_0\|_{\calL_p(L_2(\RR^d))}
\to0.
\end{align*}
Thus $x_\nu\to z$ in $L_2(\RR^d_\theta)$ for some $z$. Note that $x_\nu\to x$ in $L_p(\RR^d_\theta)$ as well. Both convergences imply convergence in measure, so $z=x$ from uniqueness of the measure limit. In particular, $x\in L_2(\RR^d_\theta)$. Then by Lemma \ref{lem:hs} again, one has
$$\|(x_\nu-x)P_0\|_{\calL_2(L_2(\RR^d))}
=(2\pi)^{-\frac d2}\|x_\nu-x\|_{L_2(\RR^d_\theta)}\to0.$$
This implies that the limit of $\{x_\nu P_0\}$ in $\calL_p(L_2(\RR^d))$ is still $xP_0$, since both convergences imply the convergence with respect to the operator norm. Passing to the limit, we proves \eqref{eq:loc}.
\end{proof}

To pass the Schatten bound of $xP_0$ in Lemma \ref{lem:loc} to that of $xP_n$, the following weighted translation operator shall be used.

\begin{lemma}\label{lem:shift}
For $b\in\RR^d$, define
$$(V_\theta(b)\xi)(t) = e^{\frac{\ii}{2}\langle b,\theta t\rangle}\xi(t-b), \quad \xi\in L_2(\RR^d).$$
Then $V_\theta(b)$ is unitary and commutes with every element of $L_\infty(\RR^d_\theta)$, and
\begin{equation}\label{eq:shift}
V_\theta(b)\calD_jV_\theta(b)^*=\calD_j-M_{b_j},\quad 1\leq j\leq d.
\end{equation}
Consequently, for every $n\in\ZZ^d$ and every operator $x$ affiliated with $L_\infty(\RR^d_\theta)$, we have
\begin{equation}\label{eq:peq}
V_\theta(n)P_0V_\theta(n)^*=P_n\quad\text{and}\quad
xP_n=V_\theta(n)(xP_0)V_\theta(n)^*.
\end{equation}
\end{lemma}

\begin{proof}
Translation and multiplication by a unimodular function are unitary on $L_2(\RR^d)$, so $V_\theta(b)$ is unitary. For $r\in\RR^d$ and $\xi\in L_2(\RR^d)$, by the antisymmetry of $\theta$, one has
\begin{align*}
(V_\theta(b)U_\theta(r)\xi)(t)
=&e^{\frac{\ii}{2}\langle b,\theta t\rangle
-\frac{\ii}{2}\langle r,\theta(t-b)\rangle}\xi(t-b-r)\\
=&(U_\theta(r)V_\theta(b)\xi)(t).
\end{align*}
Hence, $V_\theta(b)$ commutes with every element of $L_\infty(\RR^d_\theta)$. Its adjoint is given by
$$(V_\theta(b)^*\xi)(t)=e^{-\frac{\ii}{2}\langle b,\theta t\rangle}\xi(t+b).$$
Therefore, we obtain $(V_\theta(b)\calD_jV_\theta(b)^*\xi)(t)=(t_j-b_j)\xi(t)$ on $\calS(\RR^d)$, and so this proves \eqref{eq:shift}. Direct calculation yields $V_\theta(n)P_0V_\theta(n)^*=P_n$. Since $x$ is affiliated with $L_\infty(\RR^d_\theta)$, one has
$V_\theta(n)\operatorname{Dom}(x)=\operatorname{Dom}(x)$ and $xV_\theta(n)\xi=V_\theta(n)x\xi$ for $\xi\in\operatorname{Dom}(x)$. Then
\begin{align*}
\operatorname{Dom}(xP_n)
=&\{\xi\in L_2(\RR^d):P_0V_\theta(n)^*\xi\in\operatorname{Dom}(x)\}\\
=&V_\theta(n)\operatorname{Dom}(xP_0).
\end{align*}
On this domain, $xP_n\xi=V_\theta(n)(xP_0)V_\theta(n)^*\xi$. This completes the proof.
\end{proof}

For nonnegative decreasing functions $f,g$ on $(0,\infty)$, the notation $f\prec\prec g$ means
$$\int_0^t f(s) ds\leq\int_0^t g(s) ds,\quad t>0.$$
In addition, we write $f\prec g$, if $f\prec\prec g$ and
$$\int_0^\infty f(s) ds=\int_0^\infty g(s) ds.$$

\begin{lemma}\label{lem:bmaj}
Let $1\leq p<2$ and $T\in\calL_2(L_2(\RR^d))$. Suppose that $T=\sum_{n\in\ZZ^d}B_n\in\calL_2(L_2(\RR^d))$,
where the summands satisfy
$$B_nB_m^*=0,\quad n,m\in\ZZ^d,\quad n\neq m.$$
On $\mathscr K=\ell_2(\ZZ^d)\otimes L_2(\RR^d)$, define $B=\bigoplus_{n\in\ZZ^d}B_n\in B(\mathscr K)$. Then
\begin{equation}\label{eq:bmaj}
\mu( \cdot ,B)^2\prec\mu( \cdot ,T)^2.
\end{equation}
Consequently,
\begin{equation}\label{eq:bstr}
\|T\|_{\calL_p(L_2(\RR^d))}
\leq\|B\|_{\calL_p(\mathscr K)},
\quad\text{if }B\in\calL_p(\mathscr K),
\end{equation}
and
\begin{equation}\label{eq:bweak}
\|T\|_{\calL_{p,\infty}(L_2(\RR^d))}
\le 2^{1/p}\Big(\frac p{2-p}\Big)^{\frac{1}{2}}\|B\|_{\calL_{p,\infty}(\mathscr K)},
\quad\text{if }B\in\calL_{p,\infty}(\mathscr K).
\end{equation}
\end{lemma}
\begin{proof}
\eqref{eq:bmaj} follows directly from \cite[Lemma 2.9(i)]{LSZ2020}. If $B\in\calL_p(\mathscr K)$, then by \cite[Proposition 2.7(i)]{LSZ2020} and \eqref{eq:bmaj}, we get \eqref{eq:bstr}. It remains to show \eqref{eq:bweak}. By \eqref{eq:bmaj}, we obtain
$$\int_t^\infty\mu(s,T)^2 ds \leq \int_t^\infty\mu(s,B)^2 ds.$$
Since both singular value functions are decreasing and $\mu(s,B)\leq s^{-1/p}\|B\|_{\calL_{p,\infty}(\mathscr K)}$, we obtain for $t>0$,
\begin{align*}
t \mu(2t,T)^2
\leq&\int_t^{2t}\mu(s,T)^2 ds\leq\int_t^\infty\mu(s,T)^2 ds\\
\leq& \|B\|_{\calL_{p,\infty}(\mathscr K)}^2\int_t^\infty s^{-2/p} ds=\frac p{2-p}\|B\|_{\calL_{p,\infty}(\mathscr K)}^2t^{1-2/p}.
\end{align*}
This implies
$$t^{1/p}\mu(t,T) \leq 2^{1/p}\Big(\frac p{2-p}\Big)^{\frac{1}{2}}\|B\|_{\calL_{p,\infty}(\mathscr K)}.$$
Taking the supremum over $t>0$, we obtain \eqref{eq:bweak}. This completes the proof.
\end{proof}

Now we are in a position to prove Theorem \ref{thm:cwl1}.

\begin{proof}[Proof of Theorem \ref{thm:cwl1}]
First we show (i). Suppose that $x\in L_2(\RR^d_\theta)$ and $xP_0\in\calL_p(L_2(\RR^d))$. For $n\in\ZZ^d$, denote by $a_n=\|g\mathds{1}_{Q_n}\|_{L_\infty(\RR^d)}$. If $g\in\ell_p(L_\infty(\RR^d);Q)$, then
$$\sum_{n\in\ZZ^d}a_n^2
\leq\Big(\sum_{n\in\ZZ^d}a_n^p\Big)^{\frac2p}
=\|g\|_{\ell_p(L_\infty(\RR^d);Q)}^2<\infty.$$
If $g\in\ell_{p,\infty}(L_\infty(\RR^d);Q)$, suppose that $(a_j^*)_{j\geq1}$ is the decreasing rearrangement of $(a_n)$, and then
$$\sum_{n\in\ZZ^d}a_n^2
=\sum_{j=1}^\infty(a_j^*)^2
\leq\|(a_n)\|_{\ell_{p,\infty}(\ZZ^d)}^2
\sum_{j=1}^\infty j^{-2/p}<\infty.$$
Consequently, in both cases, we obtain
$$\|g\|_{L_2(\RR^d)}^2\leq\sum_{n\in\ZZ^d}a_n^2<\infty.$$
By Lemma \ref{lem:hs}, $xM_g\in\calL_2(L_2(\RR^d))$. Note that $xM_{g\mathds{1}_{Q_n}}=(xM_{g1_{Q_n}})P_n$. For $n\in\ZZ^d$ and $t>0$, by Lemma \ref{lem:shift},
\begin{equation}\label{eq:bmu}
\mu(t,xM_{g\mathds{1}_{Q_n}})
\leq a_n\mu(t,xP_n)=a_n\mu(t,xP_0).
\end{equation}
In addition, 
$$(xM_{g\mathds{1}_{Q_n}})(xM_{g\mathds{1}_{Q_m}})^*
=(xM_{g\mathds{1}_{Q_n}})P_nP_m(xM_{g1_{Q_m}})^*
=0,\quad n\neq m.$$
By Lemma \ref{lem:hs}, for every finite $F\subseteq\ZZ^d$,
$$\Big\|xM_g-\sum_{n\in F}xM_{g\mathds{1}_{Q_n}}\Big\|_{\calL_2(L_2(\RR^d))}
=(2\pi)^{-\frac{d}{2}}\|x\|_{L_2(\RR^d_\theta)}
\Big\|g\mathds{1}_{\RR^d\setminus\bigcup_{n\in F}Q_n}
\Big\|_{L_2(\RR^d)}
\to0$$
as $F\to\ZZ^d$. Thus the hypotheses in Lemma \ref{lem:bmaj} are satisfied. Put $\mathscr K=\ell_2(\ZZ^d)\otimes L_2(\RR^d)$. If $(a_n)\in\ell_p(\ZZ^d)$, then
\begin{align*}
\Big\|\bigoplus_{n\in\ZZ^d}xM_{g\mathds{1}_{Q_n}}
\Big\|_{\calL_p(\mathscr K)}^p
=&\sum_{n\in\ZZ^d}
\|xM_{g\mathds{1}_{Q_n}}\|_{\calL_p(L_2(\RR^d))}^p\\
\leq&\|xP_0\|_{\calL_p(L_2(\RR^d))}^p
\sum_{n\in\ZZ^d}a_n^p,
\end{align*}
where in the last inequality we use \eqref{eq:bmu}. By \eqref{eq:bstr} and the preceding estimate,
\begin{align*}
\|xM_g\|_{\calL_p(L_2(\RR^d))}
\leq&\Big\|\bigoplus_{n\in\ZZ^d}xM_{g\mathds{1}_{Q_n}}
\Big\|_{\calL_p(\mathscr K)}\\
\leq&\|xP_0\|_{\calL_p(L_2(\RR^d))}
\Big(\sum_{n\in\ZZ^d}a_n^p\Big)^{1/p}\\
=&\|xP_0\|_{\calL_p(L_2(\RR^d))}
\|g\|_{\ell_p(L_\infty(\RR^d);Q)}.
\end{align*}
This proves \eqref{eq:cw-lo}. Next suppose that $(a_n)\in\ell_{p,\infty}(\ZZ^d)$. By \eqref{eq:bmu},
$$\mu\Big(t,\bigoplus_{n\in\ZZ^d}xM_{g\mathds{1}_{Q_n}}\Big)
\leq\mu\Big(t,\bigoplus_{n\in\ZZ^d}a_nxP_0\Big),\quad t>0.$$
Therefore, one can apply \cite[Lemma 2.3]{LSZ2020} to obtain
\begin{align*}
\Big\|\bigoplus_{n\in\ZZ^d}xM_{g\mathds{1}_{Q_n}}
\Big\|_{\calL_{p,\infty}(\mathscr K)}
\leq&\|\operatorname{diag}(a_n)\otimes(xP_0)
\|_{\calL_{p,\infty}(\mathscr K)}\\
\leq&\|(a_n)\|_{\ell_{p,\infty}(\ZZ^d)}
\|xP_0\|_{\calL_p(L_2(\RR^d))}.
\end{align*}
Combining this estimate with \eqref{eq:bweak}, we obtain \eqref{eq:cw-wlo}.

For (ii), if $x\in W^d_p(\RR^d_\theta)$, Lemma \ref{lem:loc} gives $x\in L_2(\RR^d_\theta)$ and
$$\|xP_0\|_{\calL_p(L_2(\RR^d))}
\lesssim_{d,p}\|x\|_{W^d_p(\RR^d_\theta)}.$$
Substituting this estimate into \eqref{eq:cw-lo} and \eqref{eq:cw-wlo}, the proof is finished.
\end{proof}

\subsection{Further discussion on Schatten classes}\label{ssec:fdsc1}
In this subsection, we will show the following three propositions.

\begin{proposition}\label{prop:col2}
Let $x=U_\theta(f)\in\calS(\RR^d_\theta)$. Suppose that $g\in L_2(\RR^d)$ has compact support. Then
$$xM_g\in\bigcap_{0<p<\infty}\calL_p(L_2(\RR^d)).$$
\end{proposition}

\begin{proposition}\label{prop:sdw1}
Let $x=U_\theta(f)\in\calS(\RR^d_\theta)$. Suppose that $a,b\geq0$, $a+b\ne0$ and $1\leq p<\infty$. If
\begin{equation}\label{eq:sanc1}
p(a+b)>d,
\end{equation}
then
\begin{equation}\label{eq:sandwich}
J^{-a}xJ^{-b}\in\calL_p(L_2(\RR^d)).
\end{equation}
\end{proposition}

\begin{proposition}\label{prop:cube}
Let $x=U_\theta(f)\in\calS(\RR^d_\theta)$. Suppose that $h\in C^\infty(\RR^d)$ satisfies that for some fixed $\gamma\geq0$,
$$\sup_{t\in\RR^d}(1+|t|^2)^{\frac{\gamma+|\nu|}{2}}
|\partial^\nu h(t)|<\infty\quad\text{for every }\nu\in\NN_0^d.$$
If $a,b\geq0$ and $\rho=a+b+\gamma+1$, then
$$J^{-a}[x,h(\ii\nabla)]J^{-b}\in\calL_{\frac d\rho,\infty}(L_2(\RR^d)).$$
\end{proposition}

The following proposition is a direct consequence of Proposition \ref{prop:cube}.

\begin{proposition}\label{prop:lower}
Let $x=U_\theta(f)\in\calS(\RR^d_\theta)$. Suppose that $a,b\geq0$ and $p>1$. If $p(a+b+1)>d$, then
$$[J^{-a},x]J^{-b}\in\calL_p(L_2(\RR^d)).$$
\end{proposition}
\begin{proof}
If $a=0$, the assertion follows at once. For $a>0$, by Proposition \ref{prop:cube} with $h(t)=(1+|t|^2)^{-a/2}$, $$[J^{-a},x]J^{-b}\in\calL_{\frac{d}{a+b+1},\infty}(L_2(\RR^d))\subseteq\calL_p(L_2(\RR^d)).$$
This completes the proof.
\end{proof}

The remaining part of this subsection is devoted the proofs of Propsition \ref{prop:col2}, Proposition \ref{prop:sdw1} and Proposition \ref{prop:cube}, respectively. The proofs require some properties of rank-one operators. Let $H$ be a Hilbert space. For $a,b\in H$, denote by $a\otimes b$ the rank-one operator on $H$ defined by $a\otimes b(\xi)=\langle\xi,b\rangle_H a$. Then 
$$\mu(0,a\otimes b)=\|a\|_H\|b\|_H,
\quad\mu(j,a\otimes b)=0
\quad(j\geq1)$$
and for every $0<p<\infty$,
$$\|a\otimes b\|_{\calL_p(H)}=\|a\otimes b\|_{\calL_{p,\infty}(H)}=\|a\|_H\|b\|_H.$$

Fix $\eta\in C_c^\infty((-2,2)^d)$ such that $\eta=1$ on a neighbourhood of $\overline Q$. For $m\in\ZZ^d$ and $u\in(-2,2)^d$, define $e_m(u)=4^{-\frac{d}{2}}e^{\frac{\pi\ii}{2}\langle m,u\rangle}$. We shall need the following properties of Fourier series, and we refer to Appendix \ref{sec:app2} for their detailed proofs.

\begin{lemma}\label{lem:coef}
Let $f\in\calS(\RR^d)$. For $r,m\in\ZZ^d$, $v\in Q$ and $u\in(-2,2)^d$, define
\begin{description}
\item[(i)] $h_{r,v}(u)
=\eta(u)f(r+u-v)
e^{\frac{\ii}{2}\langle v,\theta u\rangle}$;
\item[(ii)] $c_{m,r}(v)
=\int_{(-2,2)^d}h_{r,v}(u)\overline{e_m(u)} du$.
\end{description}
Then for any $M,N\in\NN_0$,
$$\|c_{m,r}\|_{L_\infty(Q)} \lesssim_{M,N,f,\theta} (1+|m|)^{-M}(1+|r|)^{-N}, \quad m,r\in\ZZ^d.$$
If $M>d$, then the Fourier series
$$f(r+u-v)e^{\frac{\ii}{2}\langle v,\theta u\rangle} =\sum_{m\in\ZZ^d}c_{m,r}(v)e_m(u)$$
converges absolutely and uniformly for $(u,v)\in Q\times Q$.
\end{lemma}

The lemma below is elementary but quite helpful, and will be frequently used throughout this subsection.

\begin{lemma}\label{lem:absum}
Let $H$ be a separable Hilbert space and $0<p<\infty$. Let $(T_k)_{k\in I}$ be a countable family of operators on $H$ satisfying the following conditions:
\begin{description}
\item[(i)] $\sum_{k\in I}\|T_k\|_{\calL_{\min\{1,p\}}(H)}^{\min\{1,p\}}<\infty$;
\item[(ii)] there exists $S\in\calL_2(H)$ such that the net of finite partial sums converges to $S$ in $\calL_2(H)$, that is,
$$\lim_{F\to I}\Big\|S-\sum_{k\in F}T_k\Big\|_{\calL_2(H)}=0,$$
where $F$ runs over the finite subsets of $I$.
\end{description}
Then $S\in\calL_p(H)$ and
$$\|S\|_{\calL_p(H)}\le\Big(\sum_{k\in I}\|T_k\|_{\calL_{\min\{1,p\}}(H)}^{\min\{1,p\}}\Big)^{\frac{1}{\min\{1,p\}}}.$$
\end{lemma}

\begin{proof}
By (i), we can write $T=\sum_{k\in I}T_k \in\calL_{\min\{1,p\}}(H)$. It follows from the the continuous inclusion
$$\calL_{\min\{1,p\}}(H)\hookrightarrow\calL_2(H)$$
that $T=\sum_{k\in I}T_k\in\calL_2(H)$. By (ii) and uniqueness of the Hilbert--Schmidt limit, we have $T=S$ and
\begin{align*}
\|S\|_{\calL_p(H)}^{\min\{1,p\}}\le&\|S\|_{\calL_{\min\{1,p\}}(H)}^{\min\{1,p\}}\leq \lim_{F\to I}\Big\|\sum_{k\in F}T_k\Big\|_{\calL_{\min\{1,p\}}(H)}^{\min\{1,p\}}\\
\le&\sum_{k\in I}\|T_k\|_{\calL_{\min\{1,p\}}(H)}^{\min\{1,p\}}.
\end{align*}
This completes the proof.
\end{proof}

The following lemma provides a localized $\calL_p(L_2(\RR^d))$-norm estimate for $xM_g$ in the case $x\in \calS(\RR^d_\theta)$ and $g\in L_2(\RR^d)$.

\begin{lemma}\label{lem:bl1}
Let $x=U_\theta(f)\in\calS(\RR^d_\theta)$ and $g\in L_2(\RR^d)$. Then for every $l,n\in\ZZ^d$ and $0<p<\infty$,
$$P_lxM_gP_n\in\calL_p(L_2(\RR^d)).$$
Moreover, for each $N\in\NN_0$,
$$\|P_lxM_gP_n\|_{\calL_p(L_2(\RR^d))}
\lesssim_{p,N,f,\theta}(1+|l-n|)^{-N}\|g\mathds{1}_{Q_n}\|_{L_2(\RR^d)}.$$
\end{lemma}

\begin{proof}
Recall that the kernel of $xM_g$ is
$$(t,s)\mapsto f(t-s)e^{\frac{\ii}{2}\langle s,\theta t\rangle}g(s).$$
Fix $l,n\in\ZZ^d$ and write $t=l+u$, $s=n+v$, where $u,v\in Q$. We apply Lemma \ref{lem:coef} to obtain that for almost every $(t,s)\in Q_l\times Q_n$,
\begin{align*}
f(t-s)e^{\frac{\ii}{2}\langle s,\theta t\rangle}g(s)
=&e^{\frac{\ii}{2}\langle n,\theta t\rangle}
e^{\frac{\ii}{2}\langle v,\theta l\rangle}g(s)
\sum_{m\in\ZZ^d}c_{m,l-n}(v)e_m(u).
\end{align*}
For $m\in\ZZ^d$, define
$$a_{m,l,n}(t)=\mathds{1}_{Q_l}(t)e^{\frac{\ii}{2}\langle n,\theta t\rangle}e_m(t-l)$$
and
$$b_{m,l,n}(s)=\mathds{1}_{Q_n}(s)e^{-\frac{\ii}{2}\langle s-n,\theta l\rangle}\overline{c_{m,l-n}(s-n)g(s)}.$$
Therefore, one has
\begin{align*}
\|a_{m,l,n}\otimes b_{m,l,n}\|_{\calL_{\min\{1,p\}}(L_2(\RR^d))}
=&\|a_{m,l,n}\|_{L_2(\RR^d)}
\|b_{m,l,n}\|_{L_2(\RR^d)}\\
=&4^{-\frac d2}
\Big(\int_Q|c_{m,l-n}(v)|^2|g(n+v)|^2 dv\Big)^{\frac{1}{2}}\\
\leq&4^{-\frac d2}\|c_{m,l-n}\|_{L_\infty(Q)}
\|g\mathds{1}_{Q_n}\|_{L_2(\RR^d)}.
\end{align*}
Now choose $M\in\NN$ depending only on $d,p$ such that $M\min\{1,p\}>d$. Fix $0<p<\infty$ and $N\in\NN_0$. By Lemma \ref{lem:coef}, we obtain
\begin{align*}
&\sum_{m\in\ZZ^d}
\|a_{m,l,n}\otimes b_{m,l,n}\|_{\calL_{\min\{1,p\}}(L_2(\RR^d))}^{\min\{1,p\}}\\
\leq&4^{-\frac{d}{2}\min\{1,p\}}\|g\mathds{1}_{Q_n}\|_{L_2(\RR^d)}^{\min\{1,p\}}\sum_{m\in\ZZ^d}\|c_{m,l-n}\|_{L_\infty(Q)}^{\min\{1,p\}}\\
\lesssim&_{s,N,f,\theta}(1+|l-n|)^{-N\min\{1,p\}}\|g\mathds{1}_{Q_n}\|_{L_2(\RR^d)}^{\min\{1,p\}}\sum_{m\in\ZZ^d}(1+|m|)^{-M\min\{1,p\}}\\
\lesssim&_{p,N,f,\theta}(1+|l-n|)^{-N\min\{1,p\}}\|g\mathds{1}_{Q_n}\|_{L_2(\RR^d)}^{\min\{1,p\}}
<\infty.
\end{align*}

For every finite set $F\subseteq\ZZ^d$, by Lemma \ref{lem:coef},
\begin{align*}
&\Big\|P_lxM_gP_n
-\sum_{m\in F}a_{m,l,n}\otimes b_{m,l,n}
\Big\|_{\calL_2(L_2(\RR^d))}\\
=&\Big(
\int_Q\int_Q
|g(n+v)|^2
\Big|\sum_{m\notin F}c_{m,l-n}(v)e_m(u)\Big|^2
 du dv
\Big)^{\frac{1}{2}}\\
\leq&4^{-\frac d2}\|g\mathds{1}_{Q_n}\|_{L_2(\RR^d)}
\sum_{m\notin F}\|c_{m,l-n}\|_{L_\infty(Q)}\\
\lesssim&_{s,N,f,\theta}
(1+|l-n|)^{-N}\|g\mathds{1}_{Q_n}\|_{L_2(\RR^d)}
\sum_{m\notin F}(1+|m|)^{-M}
\to0
\end{align*}
as $F\to\ZZ^d$, since $M>d$. By Lemma \ref{lem:absum},
\begin{align*}
\|P_lxM_gP_n\|_{\calL_p(L_2(\RR^d))}
\leq&\Big(
\sum_{m\in\ZZ^d}
\|a_{m,l,n}\otimes b_{m,l,n}\|_{\calL_{\min\{1,p\}}(L_2(\RR^d))}^{\min\{1,p\}}
\Big)^{\frac1{\min\{1,p\}}}\\
\lesssim&_{p,N,f,\theta}
(1+|l-n|)^{-N}\|g\mathds{1}_{Q_n}\|_{L_2(\RR^d)}.
\end{align*}
This completes the proof.
\end{proof}

Then we come to the proof of Proposition \ref{prop:col2}.

\begin{proof}[Proof of Proposition \ref{prop:col2}]
By Lemma \ref{lem:hs}, we have $xM_g\in\calL_2(L_2(\RR^d))$. Since $g$ has compact support, one has
$$
\#\{n\in\ZZ^d:\|g\mathds{1}_{Q_n}\|_{L_2(\RR^d)}>0\}<\infty.
$$
Fix $0<p<\infty$ and choose $N\in\NN$ such that $N\min\{1,p\}>d$. By Lemma \ref{lem:bl1}, we obtain
\begin{align*}
&\sum_{n\in\ZZ^d}\sum_{l\in\ZZ^d}
\|P_lxM_gP_n\|_{\calL_{\min\{1,p\}}(L_2(\RR^d))}^{\min\{1,p\}}\\
\lesssim&_{p,N,f,\theta}
\sum_{n\in\ZZ^d}
\|g\mathds{1}_{Q_n}\|_{L_2(\RR^d)}^{\min\{1,p\}}
\sum_{l\in\ZZ^d}(1+|l-n|)^{-N\min\{1,p\}}\\
=&\Big(\sum_{l\in\ZZ^d}(1+|l|)^{-N\min\{1,p\}}\Big)
\sum_{n\in\ZZ^d}
\|g\mathds{1}_{Q_n}\|_{L_2(\RR^d)}^{\min\{1,p\}}
<\infty.
\end{align*}
Since only finitely many functions $g\mathds{1}_{Q_n}$ are nonzero, for every $R>0$,
\begin{align*}
\sum_{n\in\ZZ^d}\sum_{\substack{l\in\ZZ^d\\|l|\leq R}}
P_lxM_gP_n
=\sum_{\substack{l\in\ZZ^d\\|l|\leq R}}
P_lxM_{\sum_{n\in\ZZ^d}g\mathds{1}_{Q_n}}
=\Big(\sum_{\substack{l\in\ZZ^d\\|l|\leq R}}P_l\Big)xM_g.
\end{align*}
This implies
\begin{align*}
&\Big\|xM_g-
\sum_{n\in\ZZ^d}\sum_{\substack{l\in\ZZ^d\\|l|\leq R}}
P_lxM_gP_n
\Big\|_{\calL_2(L_2(\RR^d))}^2\\
=&\int_{\RR^d\setminus\bigcup_{|l|\leq R}Q_l}
\int_{\RR^d}|f(t-r)|^2|g(r)|^2 dr dt
\to0
\end{align*}
as $R\to\infty$. By Lemma \ref{lem:absum}, $xM_g\in\calL_p(L_2(\RR^d))$. Since $p>0$ is arbitrary,
$$xM_g\in\bigcap_{0<p<\infty}\calL_p(L_2(\RR^d)).$$
This completes the proof.
\end{proof}

Recall that $J=(1-\Delta)^{\frac{1}{2}}$. The following condition ensures that $xJ^{-(a+b)}\in\calL_p(L_2(\RR^d))$ for $x\in\calS(\RR^d_\theta)$.

\begin{lemma}\label{lem:jend}
Let $x\in\calS(\RR^d_\theta)$, $1\leq p<\infty$ and $a,b\geq0$ with $a+b\ne0$. If $p(a+b)>d$, then
$$xJ^{-(a+b)}\in\calL_p(L_2(\RR^d)).$$
\end{lemma}

\begin{proof}
For $2\leq p<\infty$, one has
\begin{align*}
\|(1+|\cdot|^2)^{-\frac{a+b}{2}}\|_{L_p(\RR^d)}^p
=\int_{\RR^d}(1+|t|^2)^{-\frac{p(a+b)}2} dt<\infty.
\end{align*}
For $1\leq p<2$ and $t\in n+Q$,
$$1+|n|\leq1+|t|+\sqrt d\leq\sqrt2(1+\sqrt d)(1+|t|^2)^{\frac{1}{2}}.$$
This yields
\begin{align*}
\|(1+|\cdot|^2)^{-\frac{a+b}{2}}
\|_{\ell_p(L_\infty(\RR^d);Q)}^p
=&\sum_{n\in\ZZ^d}
\|(1+|\cdot|^2)^{-\frac{a+b}{2}}\mathds{1}_{Q_n}
\|_{L_\infty(\RR^d)}^p\\
\lesssim&_{a,b,d,p}
\sum_{n\in\ZZ^d}(1+|n|)^{-p(a+b)}
<\infty.
\end{align*}
Since $x\in\calS(\RR^d_\theta)\subseteq L_p(\RR^d_\theta)\bcap W^d_p(\RR^d_\theta)$, one can apply \eqref{eq:cw-q} when $p\geq2$ and Theorem \ref{thm:cwl1} when $p<2$ to obtain
$$\|xJ^{-(a+b)}\|_{\calL_p(L_2(\RR^d))}<\infty.$$
This completes the proof.
\end{proof}

The proof of the following lemma is elementary, and is put into Appendix \ref{sec:app2}.

\begin{lemma}\label{lem:jtrace}
Let $x\in\calS(\RR^d_\theta)$ and $a,b\geq0$ with $a+b\ne0$. Let $B\in\calB(L_2(\RR^d))$ be of finite rank. Then
$$z\mapsto
\Tr\big(BJ^{-(a+b)z}xJ^{-(a+b)(1-z)}\big)$$
is holomorphic on $0<\Re z<1$ and continuous on $0\leq\Re z\leq1$. On the closed strip, we have
$$\big|\Tr\big(BJ^{-(a+b)z}xJ^{-(a+b)(1-z)}\big)\big|
\leq\|B\|_{\calL_1(L_2(\RR^d))}
\|x\|_{\calB(L_2(\RR^d))}.$$
\end{lemma}

We apply the three-line theorem to the following trace estimate.

\begin{lemma}\label{lem:jtest}
Let $x\in\calS(\RR^d_\theta)$, $a,b>0$ and $1\leq p<\infty$ with $p(a+b)>d$. Let $B\in\calB(L_2(\RR^d))$ be of finite rank and satisfy $\|B\|_{\calL_{p'}(L_2(\RR^d))}\leq1$. Then
$$\big|\Tr(BJ^{-a}xJ^{-b})\big|
\leq\|xJ^{-(a+b)}\|_{\calL_p(L_2(\RR^d))}^{\frac{b}{a+b}}
\|x^*J^{-(a+b)}\|_{\calL_p(L_2(\RR^d))}^{\frac{a}{a+b}}.$$
\end{lemma}

\begin{proof}
It follows from H\"older's inequality that for $r\in \mathbb R$,
\begin{align*}
\big|\Tr\big(BJ^{-\ii(a+b)r}xJ^{-(a+b)(1-\ii r)}\big)\big|
\leq\|xJ^{-(a+b)}\|_{\calL_p(L_2(\RR^d))}
\end{align*}
and
\begin{align*}
\big|\Tr\big(BJ^{-(a+b)(1+\ii r)}xJ^{\ii(a+b)r}\big)\big|
\leq&\|J^{-(a+b)}x\|_{\calL_p(L_2(\RR^d))}\\
=&\|x^*J^{-(a+b)}\|_{\calL_p(L_2(\RR^d))}.
\end{align*}
By Lemma \ref{lem:jend} and Lemma \ref{lem:jtrace}, one can apply the three-line theorem at $z=\frac{a}{a+b}$ to obtain the desired estimate.
\end{proof}

Now we can show Proposition \ref{prop:sdw1}.

\begin{proof}[Proof of Proposition \ref{prop:sdw1}]
Lemma \ref{lem:jend} implies directly \eqref{eq:sandwich} when $a=0$. When $b=0$, we apply Lemma \ref{lem:jend} again to $x^*$ and take adjoint to get \eqref{eq:sandwich}. For $a,b>0$, we combine Lemma \ref{lem:jtest} with Lemma \ref{lem:dual} to obtain
$$\|J^{-a}xJ^{-b}\|_{\calL_p(L_2(\RR^d))}
\leq\|xJ^{-(a+b)}\|_{\calL_p(L_2(\RR^d))}^{\frac{b}{a+b}}
\|x^*J^{-(a+b)}\|_{\calL_p(L_2(\RR^d))}^{\frac{a}{a+b}}
<\infty.$$
This proves \eqref{eq:sandwich}.
\end{proof}

We derive the following lemma from the Fatou property stated in Lemma \ref{lem:fatou}.

\begin{lemma}\label{lem:ws}
Let $p>0$, $\kappa=p/(p+1)$ and $T_j\in\calL_{p,\infty}(L_2(\RR^d))$. If
$$\sum_{j=1}^\infty
\|T_j\|_{\calL_{p,\infty}(L_2(\RR^d))}^{\kappa}<\infty,$$
then $\sum_jT_j$ converges in $\calL_{p,\infty}(L_2(\RR^d))$. Moreover, in this case,
$$\Big\|\sum_{j=1}^\infty T_j\Big\|_{\calL_{p,\infty}(L_2(\RR^d))}
\leq\Big(\sum_{j=1}^\infty\|T_j\|_{\calL_{p,\infty}(L_2(\RR^d))}^{\kappa}\Big)^{1/\kappa}.$$
\end{lemma}

\begin{proof}
Since $0<\kappa<1$,
$$\sum_j\|T_j\|_{\calB(L_2(\RR^d))}
\leq\sum_j\|T_j\|_{\calL_{p,\infty}(L_2(\RR^d))}
\leq\Big(\sum_j\|T_j\|_{\calL_{p,\infty}(L_2(\RR^d))}^\kappa\Big)^{1/\kappa}<\infty.$$
Thus $\sum_jT_j$ converges in $\calB(L_2(\RR^d))$. If all $T_j$ vanish, the assertion is immediate. Otherwise, we ignore the zero summands and take a finite nonempty set $F\subseteq\{j\in\NN:T_j\neq0\}$. Let
$$A_F=\sum_{j\in F}\|T_j\|_{\calL_{p,\infty}(L_2(\RR^d))}^\kappa
\quad\text{and}\quad
\delta_j=\frac{t\|T_j\|_{\calL_{p,\infty}(L_2(\RR^d))}^\kappa}{A_F},
\quad t>0.$$
Then $A_F>0$, $\delta_j>0$ for every $j\in F$ and $\sum_{j\in F}\delta_j=t$. Hence, one has
\begin{align*}
\mu\Big(t,\sum_{j\in F}T_j\Big)
\leq&\sum_{j\in F}\mu(\delta_j,T_j)
\leq\sum_{j\in F}\|T_j\|_{\calL_{p,\infty}(L_2(\RR^d))}\delta_j^{-1/p}\\
=&t^{-1/p}A_F^{1/p}
\sum_{j\in F}\|T_j\|_{\calL_{p,\infty}(L_2(\RR^d))}^{1-\kappa/p}\\
=&t^{-1/p}A_F^{1/\kappa}.
\end{align*}
Consequently,
$$\Big\|\sum_{j\in F}T_j\Big\|_{\calL_{p,\infty}(L_2(\RR^d))}
\leq\Big(\sum_j\|T_j\|_{\calL_{p,\infty}(L_2(\RR^d))}^\kappa\Big)^{1/\kappa}.$$
By passing to the limit and Lemma \ref{lem:fatou}, this yields the desired result.
\end{proof}

Define $\omega(t)=(1+|t|^2)^{\frac{1}{2}}$ for $t\in \mathbb{R}^d$. Recall that the functions $e_m$ were defined before Lemma \ref{lem:coef}. The reader is referred to Appendix \ref{sec:app2} for the detailed proof of the following lemma, as the proof is just quite standard but tedious.

\begin{lemma}\label{lem:w11}
Let $f\in\calS(\RR^d)$ and $b,c,\gamma$ be arbitrary nonnegative real numbers. Let $\sigma\in C^\infty(\RR^d)$ satisfy
$$\sup_{t\in\RR^d}\omega(t)^{\gamma+|\nu|}|\partial^\nu\sigma(t)|
<\infty\quad\text{for every }\nu\in\NN_0^d.$$
Let $\rho=b+c+\gamma+1$. For $r,n\in\ZZ^d$, $u\in(-2,2)^d$ and $v\in Q$, define
\begin{align*}
h_{r,n,v}(u)=&\eta(u)\omega(n+r+u)^{-b}
f(r+u-v)e^{\frac{\ii}{2}\langle v,\theta u\rangle}\\
&\cdot\big(\sigma(n+v)-\sigma(n+r+u)\big)\omega(n+v)^{-c},
\end{align*}
and
$$c_{m,r,n}(v)=\int_{(-2,2)^d} h_{r,n,v}(u)\overline{e_m(u)} du,
\quad m\in\ZZ^d.$$
For every $M,N\geq0$, one has
$$\|c_{m,r,n}\|_{L_\infty(Q)}\lesssim_{M,N,b,c,\gamma,d,f,\sigma,\theta,\eta}(1+|m|)^{-M}(1+|r|)^{-N}(1+|n|)^{-\rho}.$$
Moreover,
$$h_{r,n,v}(u)=\sum_{m\in\ZZ^d}c_{m,r,n}(v)e_m(u),
\quad (u,v)\in Q\times Q,$$
converges absolutely and uniformly for $(u,v)\in Q\times Q$. The implicit constant in the coefficient estimate is independent of $m,r,n$; its dependence on $f$ and $\sigma$ involves only finitely many Schwartz and symbol seminorms, respectively.
\end{lemma}

Now we can prove Proposition \ref{prop:cube}.

\begin{proof}[Proof of Proposition \ref{prop:cube}]
Let $T=J^{-a}[x,M_h]J^{-b}$. Since $\gamma\geq0$,
$$\|h\|_{L_\infty(\RR^d)}\leq\sup_{t\in\RR^d}\omega(t)^\gamma|h(t)|<\infty.$$
Since $J^{-a}$ and $J^{-b}$ are contractions, we have
$$\|T\|_{\calB(L_2(\RR^d))}\leq2\|x\|_{\calB(L_2(\RR^d))}
\|h\|_{L_\infty(\RR^d)}<\infty.$$
Let $\kappa=d/(d+\rho)$ and choose $M,N>0$ such that $M\kappa>d$ and $N\kappa>d$. In particular, $M,N>d$. The kernel of $T$ is given by
$$K(t,s)=\omega(t)^{-a}f(t-s)e^{\frac{\ii}{2}\langle s,\theta t\rangle}\big(h(s)-h(t)\big)\omega(s)^{-b}.$$
For $t=n+r+u$ and $s=n+v$, with $u,v\in Q$, we have
$$\langle s,\theta t\rangle=\langle n,\theta t\rangle+\langle v,\theta(n+r)\rangle+\langle v,\theta u\rangle.$$
By Lemma \ref{lem:w11},
$$K(t,s)=e^{\frac{\ii}{2}\langle n,\theta t\rangle}e^{\frac{\ii}{2}\langle s-n,\theta(n+r)\rangle}\sum_{m\in\ZZ^d}c_{m,r,n}(s-n)e_m(t-n-r),\quad (t,s)\in Q_{n+r}\times Q_n.$$
Define
\begin{description}
\item[(i)] $\phi_{m,r,n}(t)=\mathds{1}_{Q_{n+r}}(t)e^{\frac{\ii}{2}\langle n,\theta t\rangle}e_m(t-n-r)$;
\item[(ii)] $\psi_{m,r,n}(s)=\mathds{1}_{Q_n}(s)e^{-\frac{\ii}{2}\langle s-n,\theta(n+r)\rangle}\overline{c_{m,r,n}(s-n)}$.
%\item[(iii)] $B_{m,r,n}\xi=\phi_{m,r,n}\int_{\RR^d}\overline{\psi_{m,r,n}(s)}\xi(s) ds$ for $\xi\in L_2(\RR^d)$.
\end{description}
Let $B_{m,r,n}=\mathcal \phi_{m,r,n}\otimes\psi_{m,r,n}$. By Lemma \ref{lem:w11},
\begin{align*}
\alpha_{m,r,n}:=&\mu(0,B_{m,r,n})\\
=&\|\phi_{m,r,n}\|_{L_2(\RR^d)}\|\psi_{m,r,n}\|_{L_2(\RR^d)}\\
=&4^{-\frac{d}{2}}\|c_{m,r,n}\|_{L_2(Q)}\\
\leq&A_{m,r}(1+|n|)^{-\rho},
\end{align*}
where
$$A_{m,r}=C_{M,N}(1+|m|)^{-M}(1+|r|)^{-N}.$$
As in Lemma \ref{lem:w11}, $C_{M,N}>0$ may also depend on $a,b,\gamma,d,f,h,\theta,\eta$. For fixed $r,n$, if a sequence of finite sets $F\subseteq\ZZ^d$ increase to $\ZZ^d$, then
\begin{align}\label{eq:PTP1}
\Big\|P_{n+r}TP_n-\sum_{m\in F}B_{m,r,n}\Big\|_{\calL_2(L_2(\RR^d))}
\leq4^{-\frac{d}{2}}\sum_{m\notin F}\|c_{m,r,n}\|_{L_\infty(Q)}
\to0.
\end{align}

Fix $m,r$. Note that
$$\operatorname{supp}\phi_{m,r,n}\subseteq\overline{Q_{n+r}},
\quad\operatorname{supp}\psi_{m,r,n}\subseteq\overline{Q_n},$$
and this implies that for $n\neq n'$,
$$\langle\phi_{m,r,n},\phi_{m,r,n'}\rangle_{L_2(\RR^d)}=
\langle\psi_{m,r,n},\psi_{m,r,n'}\rangle_{L_2(\RR^d)}=0$$
and
$$B_{m,r,n}^*B_{m,r,n'}=B_{m,r,n}B_{m,r,n'}^*=0.$$
For every finite nonempty set $F\subseteq\ZZ^d$, it yields
$$\Big\|\sum_{n\in F}B_{m,r,n}\Big\|_{\calB(L_2(\RR^d))}
=\max_{n\in F}\alpha_{m,r,n}.$$
Hence, one has
\begin{align*}
\sup_{\substack{G\subseteq\ZZ^d\setminus F\\G\text{ finite}}}
\Big\|\sum_{n\in G}B_{m,r,n}\Big\|_{\calB(L_2(\RR^d))}
\leq&\sup_{n\notin F}\alpha_{m,r,n}\\
\leq&A_{m,r}\sup_{n\notin F}(1+|n|)^{-\rho}
\to0,
\end{align*}
as $F$ increases to $\ZZ^d$. We conclude that $S_{m,r}=\sum_nB_{m,r,n}$ converges in $\calB(L_2(\RR^d))$. Since its finite partial sums have finite rank, $S_{m,r}$ is compact.

The singular values of $S_{m,r}$ are the decreasing rearrangement of $(\alpha_{m,r,n})_n$. If $\lambda\geq A_{m,r}$, then
$$\Tr\big(\mathds{1}_{(\lambda,\infty)}(|S_{m,r}|)\big)
=\#\{n:\alpha_{m,r,n}>\lambda\}=0.$$
For $0<\lambda<A_{m,r}$, one has
\begin{align*}
\Tr\big(\mathds{1}_{(\lambda,\infty)}(|S_{m,r}|)\big)
=&\#\{n:\alpha_{m,r,n}>\lambda\}\\
\leq&\#\{n:A_{m,r}(1+|n|)^{-\rho}>\lambda\}\\
=&\#\Big\{n:1+|n|<(A_{m,r}/\lambda)^{1/\rho}\Big\}\\
\leq&\Big(2(A_{m,r}/\lambda)^{1/\rho}+1\Big)^d\\
\leq&3^d(A_{m,r}/\lambda)^{d/\rho}.
\end{align*}
Therefore, we obtain
\begin{align*}
\|S_{m,r}\|_{\calL_{\frac d\rho,\infty}(L_2(\RR^d))}
=\sup_{\lambda>0}
\lambda\Big[\Tr\big(\mathds{1}_{(\lambda,\infty)}(|S_{m,r}|)\big)\Big]^{\frac\rho d}\leq3^\rho A_{m,r}.
\end{align*}
By our choice of $M,N$,
\begin{align*}
&\sum_{m,r\in\ZZ^d}
\|S_{m,r}\|_{\calL_{\frac d\rho,\infty}(L_2(\RR^d))}^{\kappa}\\
\leq&3^{\rho\kappa}C_{M,N}^{\kappa}
\Big(\sum_{m\in\ZZ^d}(1+|m|)^{-M\kappa}\Big)
\Big(\sum_{r\in\ZZ^d}(1+|r|)^{-N\kappa}\Big)
<\infty.
\end{align*}
Since $\kappa=(d/\rho)/(d/\rho+1)$, Lemma \ref{lem:ws} yields
\begin{align}\label{eq:Smr2}
S=\sum_{m,r\in\ZZ^d}S_{m,r}
\in\calL_{\frac d\rho,\infty}(L_2(\RR^d)).
\end{align}
Moreover,
\begin{align}\label{eq:Smr1}
\begin{split}
\sum_{m,r\in\ZZ^d}\|S_{m,r}\|_{\calB(L_2(\RR^d))}
\leq&\sum_{m,r\in\ZZ^d}
\|S_{m,r}\|_{\calL_{\frac d\rho,\infty}(L_2(\RR^d))}\\
\leq&\Big(
\sum_{m,r\in\ZZ^d}
\|S_{m,r}\|_{\calL_{\frac d\rho,\infty}(L_2(\RR^d))}^{\kappa}
\Big)^{1/\kappa}
<\infty.
\end{split}
\end{align}

For $l,n\in\ZZ^d$, note that
$$P_lS_{m,r}P_n=
\begin{cases}
B_{m,r,n},&r=l-n,\\
0,&r\neq l-n.
\end{cases}
$$
Using \eqref{eq:Smr1}, we obtain
$$P_lSP_n=\sum_{m,r\in\ZZ^d}P_lS_{m,r}P_n
=\sum_{m\in\ZZ^d}B_{m,l-n,n}
\quad\text{in }\calB(L_2(\RR^d)).$$
On the other hand, by \eqref{eq:PTP1}, as $F$ increases to $\ZZ^d$,
\begin{align*}
\Big\|P_lTP_n-\sum_{m\in F}B_{m,l-n,n}
\Big\|_{\calB(L_2(\RR^d))}\leq 
\Big\|P_lTP_n-\sum_{m\in F}B_{m,l-n,n}
\Big\|_{\calL_2(L_2(\RR^d))}
\to0.
\end{align*}
By uniqueness of the operator-norm limit, this implies that $P_lSP_n=P_lTP_n$.

Let $R\ge0$, and define $E_R=\sum_{|n|\leq R}P_n$. Then  $E_R\to1$ strongly. Note that
$$E_R(S-T)E_R=\sum_{|l|\leq R, |n|\leq R}P_l(S-T)P_n
=0,$$
and this implies that $T=S$.
\iffalse
\begin{align*}
\|(S-T)\xi\|_{L_2(\RR^d)}
=&\|(S-T-E_R(S-T)E_R)\xi\|_{L_2(\RR^d)}\\
\leq&\|(1-E_R)(S-T)\xi\|_{L_2(\RR^d)}\\
&+\|S-T\|_{\calB(L_2(\RR^d))}\|(1-E_R)\xi\|_{L_2(\RR^d)}
\to0.
\end{align*}
\fi
Hence, we conclude from \eqref{eq:Smr2} that $T=S\in\calL_{\frac d\rho,\infty}(L_2(\RR^d))$. This completes the proof.
\end{proof}

We end this section with the following corollary.

\begin{corollary}\label{cor:ang1}
Let $x\in\calS(\RR^d_\theta)$, $g\in C^\infty(\BS^{d-1})$ and $a,b,\gamma\geq0$. Let $\rho=a+b+\gamma+1$. This
$$J^{-a}[x,\pi_2(g)J^{-\gamma}]J^{-b}
\in\calL_{d/\rho,\infty}(L_2(\RR^d)).$$
In particular, if $\rho>d$, the above operator is of trace class.
\end{corollary}

\begin{proof}
Choose $\chi\in C_c^\infty(\RR^d)$ such that it is equal to one near zero, and split $\sigma=\widetilde g \omega^{-\gamma}$ into the sum $\chi\sigma+(1-\chi)\sigma$. Since $\chi\sigma$ is bounded and compactly supported, Proposition \ref{prop:col2} yields that $xM_{\chi\sigma},M_{\chi\sigma}x\in\calL_s(L_2(\RR^d))$ for every $s>0$. Since $J^{-a}$ and $J^{-b}$ are bounded, we obtain
$$J^{-a}[x, M_{\chi\sigma}]J^{-b}
\in\calL_{d/\rho,\infty}(L_2(\RR^d)).$$
The remaining symbol is smooth and satisfies
$$|\partial^\nu((1-\chi)\sigma)(t)|
\lesssim_{\nu,g,\gamma,\chi}\omega(t)^{-\gamma-|\nu|}.$$
By Proposition \ref{prop:cube}, we obtain 
$$J^{-a}[x, M_{(1-\chi)\sigma}]J^{-b}
\in\calL_{d/\rho,\infty}(L_2(\RR^d)).$$
This implies that 
$$J^{-a}[x,\pi_2(g)J^{-\gamma}]J^{-b}
\in\calL_{d/\rho,\infty}(L_2(\RR^d)).$$
When $\rho>d$, $\calL_{d/\rho,\infty}(L_2(\RR^d))$ is contained in $\calL_1(L_2(\RR^d))$. This completes the proof.
\end{proof}

\bigskip

\section{The estimate for the upper bound of $[(-\Delta)^{\frac{\epsilon}{2}},x]$}\label{s:doi}
 
In this section, we prove the upper bound of the weak Schatten class for $[(-\Delta)^{\frac{\epsilon}{2}},x]$ of any Schwartz symbol by expressing it through double operator integrals and applying the Cwikel estimates. Then we extend the estimate to every element in the homogeneous Sobolev space $\dW^1_{\frac{d}{1-\epsilon}}(\RR^d_\theta)$ by standard approximation. The argument of double operator integral follows from \cite{FSZ2024}. 
%The estimates for quantum Euclidean symbols and the passage from Schwartz elements to $\dW^1_{\frac{d}{1-\epsilon}}(\RR_\theta^d)$ are proved here using the Cwikel estimates and the approximation results established above.

At first, we establish the desired weak Schatten bound for $x \in \calS(\RR^d_\theta)$. If $\mathbf{A} = (A_1,\dots,A_d)$ is a $d$-tuple of mutually commuting self-adjoint operators on a Hilbert space and $\Phi$ is a sufficiently regular function on $\mathbb{R}^d\times\mathbb{R}^d$, denote by $T_\Phi^{\mathbf{A}}$ the corresponding double operator integral (DOI). We refer to \cite{CSZ2018, KPSS2012} for more information.

For $\epsilon\in\RR$, define the function $\phi_\epsilon$ on $\RR^d \times \RR^d$ by
\begin{align}\label{eq:phi1}
\phi_\epsilon(\lambda,\mu) =
\begin{cases}
\frac{|\lambda|^\epsilon - |\mu|^\epsilon}{|\lambda| - |\mu|} \cdot |\lambda|^{\frac{1-\epsilon}{2}} |\mu|^{\frac{1-\epsilon}{2}}, & \text{ if } |\lambda| \neq |\mu|,\\
\epsilon, & \text{ if } |\lambda| = |\mu|.
\end{cases}
\end{align}
Next, we record the Schatten boundedness of the double operator integral with the symbol $\phi_\epsilon$ from \eqref{eq:phi1}. Since the proof of Lemma \ref{lem:TLL1} is completely analogous to that of \cite[Lemma 4.2]{FSZ2024}, the details are omitted.

\begin{lemma}\label{lem:TLL1}
Let $\epsilon \in [-1,1]$ and $\phi_\epsilon$ be as in \eqref{eq:phi1}. Then for any $1 < p < \infty$, the double operator integral $T_{\phi_\epsilon}^{\ii\nabla}$ is bounded on $\calL_p(L_2(\RR^d))$ and on $\calL_{p,\infty}(L_2(\RR^d))$.
\end{lemma}

Define the function $\psi_k$ $(1 \le k \le d)$ on $\RR^d \times \RR^d$ by
\begin{align}\label{eq:psi1}
\psi_k(\lambda,\mu) :=
\begin{cases}
\displaystyle \frac{|\lambda| - |\mu|}{|\lambda - \mu|^2} \cdot (\lambda_k - \mu_k), & \text{ if } \lambda \neq \mu,\\
0, & \text{ if } \lambda = \mu.
\end{cases}
\end{align}
The following lemma concerns the Schatten bound of $T_{\psi_k}^{\ii\nabla}$, which follows from \cite[Corollary 5.2]{CSZ2018} (see also \cite[Theorem 5.1]{CMPS2014}) by spectral truncation.
%, and interpolation gives the weak Schatten bound.
\begin{lemma}\label{lem:TLL2}
Let $1\le k\le d$. For any $1 < p < \infty$, the double operator integral $T_{\psi_k}^{\ii\nabla}$ is bounded on $\calL_p(L_2(\RR^d))$ and on $\calL_{p,\infty}(L_2(\RR^d))$.
\end{lemma}

The following lemma will be needed later, and the proof relies on Theorem \ref{thm:cw-hi}.

\begin{lemma}\label{lem:DDM1}
Let $\epsilon \in (-\frac{d}{2},1)$ and $x \in \calS(\RR^d_\theta)$.
\begin{description}
\item[(i)] We have
$$\| (-\Delta)^{\frac{\epsilon-1}{4}} x (-\Delta)^{\frac{\epsilon-1}{4}}\|_{\calL_{\frac{d}{1-\epsilon},\infty}(L_2(\RR^d))}\lesssim_{d,\epsilon}\|x\|_{L_{\frac{d}{1-\epsilon}}(\RR^d_\theta)}.$$
\item[(ii)] Let $n \in \mathbb{N}$ satisfy $\epsilon \in [-2n - 1, -2n + 1)$. Then
$$\| (-\Delta)^{\frac{\epsilon+k-1}{2}} x (-\Delta)^{-\frac{k}{2}} \|_{\calL_{\frac{d}{1-\epsilon},\infty}(L_2(\RR^d))}
\lesssim_{d,\epsilon,k}\|x\|_{L_{\frac{d}{1-\epsilon}}(\RR^d_\theta)},\quad 1 \le k \le n.$$
\end{description}
\end{lemma}
\begin{proof}
(i) Let $x=u|x|$ be the polar decomposition of $x$. By H\"older's inequality, one has
\begin{align*}
&\| (-\Delta)^{\frac{\epsilon-1}{4}} x (-\Delta)^{\frac{\epsilon-1}{4}}\|_{\calL_{\frac{d}{1-\epsilon},\infty}(L_2(\RR^d))}\\
\lesssim_{d,\epsilon}&\| (-\Delta)^{\frac{\epsilon-1}{4}} u|x|^{\frac12}\|_{\calL_{\frac{2d}{1-\epsilon},\infty}(L_2(\RR^d))}
\||x|^{\frac12} (-\Delta)^{\frac{\epsilon-1}{4}}\|_{\calL_{\frac{2d}{1-\epsilon},\infty}(L_2(\RR^d))}.
\end{align*} 
Since $t\mapsto |t|^{\frac{\epsilon-1}{2}}$ belongs to $L_{\frac{2d}{1-\epsilon},\infty}(\RR^d)$ and $\frac{2d}{1-\epsilon}>2$, by Theorem \ref{thm:cw-hi}, each factor is dominated by $\|x\|_{L_{\frac{d}{1-\epsilon}}(\RR^d_\theta)}^{\frac12}$. This proves (i).

\noindent (ii) Note that $d>-2\epsilon>1-\epsilon>2n\ge 2k\ge2$. By H\"older's inequality, one has
\begin{align*}
&\| (-\Delta)^{\frac{\epsilon+k-1}{2}} x (-\Delta)^{-\frac{k}{2}} \|_{\calL_{\frac{d}{1-\epsilon},\infty}(L_2(\RR^d))}\\
\lesssim&_{k,d,\epsilon}\|(-\Delta)^{\frac{\epsilon+k-1}{2}} u|x|^{1-\frac{k}{1-\epsilon}} \|_{\calL_{\frac{d}{1-\epsilon-k},\infty}(L_2(\RR^d))}
\||x|^{\frac{k}{1-\epsilon}} (-\Delta)^{-\frac{k}{2}}\|_{\calL_{\frac{d}{k},\infty}(L_2(\RR^d))}.
\end{align*} 
Note that $\frac{d}{1-\epsilon-k}>2$ and $\frac{d}{k}>2$. By Theorem \ref{thm:cw-hi} again, we complete the proof of (ii).
\end{proof}

Lemma \ref{lem:DDM1} can be applied directly to obtain the following lemma if $x$ is replaced by $\partial_k x$.
\begin{lemma}\label{lem:XY1}
Let $\epsilon \in (-\frac{d}{2},1)$ and $n \in \NN\cup\{0\}$ with $\epsilon \in [-2n - 1, -2n + 1)$. For $x \in \calS(\RR^d_\theta)$, define
$$X_k = (-\Delta)^{\frac{\epsilon-1}{4}} (\partial_k x) (-\Delta)^{\frac{\epsilon-1}{4}},\quad 1 \le k \le d,$$
and
$$Y_k = \sum_{j=1}^n ( (-\Delta)^{\frac{\epsilon+j-1}{2}} (\partial_k x) (-\Delta)^{-\frac{j}{2}}
+ (-\Delta)^{-\frac{j}{2}} (\partial_k x) (-\Delta)^{\frac{\epsilon+j-1}{2}} ),\quad 1 \le k \le d.$$
Here, for $n = 0$, $Y_k$ is assumed to be zero. Then for any $1 \le k \le d$,
$$\|X_k\|_{\calL_{\frac{d}{1-\epsilon},\infty}(L_2(\RR^d))} + \|Y_k\|_{\calL_{\frac{d}{1-\epsilon},\infty}(L_2(\RR^d))}
\lesssim_{n,d,\epsilon}\|x\|_{\dW^1_{\frac{d}{1-\epsilon}}(\RR^d_\theta)}.$$
\end{lemma}

We shall need the following approximation properties.

\begin{lemma}\label{lem:epmg0}
Let $-\frac{d}{2}<\epsilon\neq 0<1$ and $x\in L_\infty(\RR^d_\theta)$.
For $m\ge1$, let $P_m^\theta=\mathds{1}_{(1/m,m)}(-\Delta_\theta)$. Then for any $\phi\in\calS(\RR^d_\theta)$,
\begin{equation}\label{eq:PDM1}
(P_m^\theta(-\Delta_\theta)^{\frac{\epsilon}{2}})M_xP_m^\theta\phi
\to(-\Delta_\theta)^{\frac{\epsilon}{2}}x\phi
\end{equation}
and
\begin{equation}\label{eq:PDM2}
P_m^\theta x(P_m^\theta(-\Delta_\theta)^{\frac{\epsilon}{2}})\phi
\to x(-\Delta_\theta)^{\frac{\epsilon}{2}}\phi
\end{equation}
in $\calS'(\RR^d_\theta)$ as $m\to\infty$.
\end{lemma}

\begin{proof}
By Proposition \ref{prop:LLt1}, for $f\in\calS(\RR^d)$, one has
$$P_m^\theta U_\theta(f)=U_\theta\big(\mathds{1}_{(\frac1m,m)}(|\xi|^2)f\big)$$
and
$$P_m^\theta(-{\Delta_\theta})^{\frac{\epsilon}{2}}U_\theta(f)=U_\theta(\mathds{1}_{(\frac1m,m)}(|\xi|^2)|\xi|^\epsilon f).$$
It is clear that the operators $P_m^\theta$ and $P_m^\theta(-{\Delta_\theta})^{\frac{\epsilon}{2}}$ are bounded on $L_2(\RR^d_\theta)$.
For every $y=U_\theta(h)\in\calS(\RR^d_\theta)$ with $h\in\calS(\RR^d)$, Proposition \ref{prop:iso2} and the dominated convergence theorem yield
\begin{equation}\label{eq:epmg0-cutoff}
\|(I-P_m^\theta)y\|_{L_2(\RR^d_\theta)}
\lesssim_d\big\|\mathds{1}_{\{1/m<|\xi|^2<m\}^{\mathrm c}} h\big\|_{L_2(\RR^d)}
\to 0, \quad m\to \infty.
\end{equation}
By Proposition \ref{prop:iso1}, one has that as $m\to \infty$,
\begin{equation}\label{eq:epmg0-l2}
\|(I-P_m^\theta)(-{\Delta_\theta})^{\frac{\epsilon}{2}}y\|_{L_2(\RR^d_\theta)}^2
\asymp_d \int_{\{1/m<|\xi|^2<m\}^{\mathrm c}}|\xi|^{2\epsilon}|h(\xi)|^2 d\xi
\to 0.
\end{equation}

For all $\phi,\psi\in\calS(\RR^d_\theta)$, since both
$$(P_m^\theta(-{\Delta_\theta})^{\frac{\epsilon}{2}})M_xP_m^\theta\phi
\quad\text{and}\quad
P_m^\theta M_x(P_m^\theta(-{\Delta_\theta})^{\frac{\epsilon}{2}})\phi$$
belong to $L_2(\RR^d_\theta)$, we have, in the sense of distributions,
$$\big\langle (P_m^\theta(-{\Delta_\theta})^{\frac{\epsilon}{2}})M_xP_m^\theta\phi,\psi^*\big\rangle=\tau_\theta\Big((xP_m^\theta\phi) (P_m^\theta(-{\Delta_\theta})^{\frac{\epsilon}{2}}\psi)^*\Big)$$
and
$$\big\langle P_m^\theta M_x(P_m^\theta(-{\Delta_\theta})^{\frac{\epsilon}{2}})\phi,\psi^*\big\rangle=\tau_\theta\Big((P_m^\theta(-{\Delta_\theta})^{\frac{\epsilon}{2}}\phi)( x^*P_m^\theta\psi)^*\Big).$$
It remains to show the convergence of the pairings for every fixed $\psi\in\calS(\RR^d_\theta)$. Using \eqref{eq:drb1}, we decompose
\begin{align*}
&\big\langle (P_m^\theta(-{\Delta_\theta})^{\frac{\epsilon}{2}})M_xP_m^\theta\phi-(-{\Delta_\theta})^{\frac{\epsilon}{2}}M_x\phi,\psi^*\big\rangle\\
=&\tau_\theta\Big(x(P_m^\theta-I)\phi \big(P_m^\theta(-{\Delta_\theta})^{\frac{\epsilon}{2}}\psi\big)^*\Big)
+\tau_\theta\Big(x\phi \big((P_m^\theta-I)(-{\Delta_\theta})^{\frac{\epsilon}{2}}\psi\big)^*\Big).
\end{align*}
By H\"older's inequality, \eqref{eq:epmg0-cutoff} and \eqref{eq:epmg0-l2}, as $m\to \infty$,
\begin{align*}
&\Big|\tau_\theta\Big(x(P_m^\theta-I)\phi \big(P_m^\theta(-{\Delta_\theta})^{\frac{\epsilon}{2}}\psi\big)^*\Big)\Big|\\
\le& \|x\|_{L_\infty(\RR^d_\theta)}\|(P_m^\theta-I)\phi\|_{L_2(\RR^d_\theta)}
\|P_m^\theta(-{\Delta_\theta})^{\frac{\epsilon}{2}}\psi\|_{L_2(\RR^d_\theta)}\\
\le& \|x\|_{L_\infty(\RR^d_\theta)}\|(P_m^\theta-I)\phi\|_{L_2(\RR^d_\theta)}
\|(-{\Delta_\theta})^{\frac{\epsilon}{2}}\psi\|_{L_2(\RR^d_\theta)}\to 0.
\end{align*}
By H\"older's inequality and \eqref{eq:epmg0-l2}, as $m\to \infty$,
\begin{align*}
&\Big|\tau_\theta\Big(x\phi \big((P_m^\theta-I)(-{\Delta_\theta})^{\frac{\epsilon}{2}}\psi\big)^*\Big)\Big|\\
\le& \|x\phi\|_{L_2(\RR^d_\theta)}
\|(P_m^\theta-I)(-{\Delta_\theta})^{\frac{\epsilon}{2}}\psi\|_{L_2(\RR^d_\theta)}\\
\le& \|x\|_{L_\infty(\RR^d_\theta)}\|\phi\|_{L_2(\RR^d_\theta)}
\|(P_m^\theta-I)(-{\Delta_\theta})^{\frac{\epsilon}{2}}\psi\|_{L_2(\RR^d_\theta)}
\to 0.
\end{align*}
This proves \eqref{eq:PDM1}. In addition, \eqref{eq:PDM2} can be shown in the same way.
\end{proof}

Now we can transfer the convergence in Lemma \ref{lem:epmg0} to the convergence in $\calS'(\RR^d)$.

\begin{lemma}\label{lem:epmg1}
Let $-\frac{d}{2}<\epsilon\neq 0<1$ and $x\in \calS(\RR^d_\theta)$. For $m\ge1$, let $P_m=\mathds{1}_{(1/m,m)}(-\Delta)$. Then for any $\phi\in\calS(\RR^d)$,
$$(P_m(-\Delta)^{\frac{\epsilon}{2}})xP_m\phi-P_mx(P_m(-\Delta)^{\frac{\epsilon}{2}})\phi\to[(-\Delta)^{\frac{\epsilon}{2}},x]\phi$$
in $\calS'(\RR^d)$ as $m\to\infty$.
\end{lemma}
\begin{proof}
Let $P_m^\theta=\mathds{1}_{(1/m,m)}(-\Delta_\theta)$. By Proposition \ref{prop:LLt1} and Proposition \ref{prop:umu1},
$$(P_m(-\Delta)^{\frac{\epsilon}{2}})xP_m=V^*(P_m^\theta(-\Delta_\theta)^{\frac{\epsilon}{2}})M_xP_m^\theta V$$
and
$$P_mx(P_m(-\Delta)^{\frac{\epsilon}{2}})=V^*P_m^\theta M_x(P_m^\theta(-\Delta_\theta)^{\frac{\epsilon}{2}})V.$$
We conclude the desired result by Lemma \ref{lem:epmg0}.
\end{proof}

We will represent the commutator by double operator integrals in terms of $X_k$ and $Y_k$.

\begin{lemma}\label{lem:XY2}
Let $-\frac{d}{2}<\epsilon\neq0<1$. Choose $n \in \NN\cup\{0\}$ such that $\epsilon \in [-2n-1,-2n+1)$. For every $x \in \calS(\RR^d_\theta)$, we have
$$[(-{\Delta})^{\frac{\epsilon}{2}}, x] = \sum_{k=1}^d T_{\psi_k}^{\ii\nabla} (T_{\phi_{\epsilon+2n}}^{\ii\nabla}(X_k)) - \sum_{k=1}^d T_{\psi_k}^{\ii\nabla}(Y_k).$$
Here, $X_k$ and $Y_k$ are given in Lemma \ref{lem:XY1}, $\phi_\epsilon$ is given in \eqref{eq:phi1} and $\psi_k$ is given in \eqref{eq:psi1}.
\end{lemma}
\begin{proof}
For a tuple $\mathbf{A}=(A_1,\dots,A_d)$, we write
$$|\mathbf{A}|=(A_1^2+\dots+A_d^2)^{\frac12}.$$
Fix $m \in \mathbb{N}$. Let $P_m = \mathds{1}_{(\frac{1}{m},m)}(-{\Delta})$. On the Hilbert space $P_m(L_2(\RR^d))$, consider the tuple $\mathbf{A} = {\ii\nabla}\cdot P_m$ and $B = P_m x P_m$. Clearly, $\mathbf{A}$ is bounded and $|\mathbf{A}|$ is bounded below. Let $x=u|x|$ be the polar decomposition of $x$. We write $$B=P_mu|x|^{\frac{1}{2}}|x|^{\frac{1}{2}}P_m.$$
By Lemma \ref{lem:hs},
$|x|^{\frac{1}{2}}P_m,P_mu|x|^{\frac{1}{2}}\in \calL_2(L_2(\RR^d))$, and hence $B \in \calL_1(L_2(\RR^d))$. Write $\delta:= \epsilon + 2n$. By applying \cite[Lemma 4.6]{FSZ2024} (for $n \in \mathbb{N}$) or \cite[Lemma 4.5]{FSZ2024} (for $n=0$), we obtain
\begin{align*}
T_m :=& (P_m(-{\Delta})^{\frac{\epsilon}{2}})  x P_m - P_m  x (P_m(-{\Delta})^{\frac{\epsilon}{2}}) = [|\mathbf{A}|^\epsilon, B] \\
=& \sum_{k=1}^d T_{\psi_k}^{\mathbf{A}} ( T_{\phi_{\delta}}^{\mathbf{A}} ( |\mathbf{A}|^{\frac{\epsilon-1}{2}} [A_k,B] |\mathbf{A}|^{\frac{\epsilon-1}{2}} ) ) \\
&- \sum_{k=1}^d T_{\psi_k}^{\mathbf{A}} ( \sum_{j=1}^n |\mathbf{A}|^{\epsilon+j-1}[A_k,B]|\mathbf{A}|^{-j} + \sum_{j=1}^n |\mathbf{A}|^{-j}[A_k,B]|\mathbf{A}|^{\epsilon+j-1} ).
\end{align*}
Since
$$
|\mathbf{A}|^{\frac{\epsilon-1}{2}} [A_k,B] |\mathbf{A}|^{\frac{\epsilon-1}{2}} = P_m X_k P_m,
$$
$$
\sum_{j=1}^n |\mathbf{A}|^{\epsilon+j-1}[A_k,B]|\mathbf{A}|^{-j} + \sum_{j=1}^n |\mathbf{A}|^{-j}[A_k,B]|\mathbf{A}|^{\epsilon+j-1} = P_m Y_k P_m,
$$
we obtain
$$
T_m = \sum_{k=1}^d T_{\psi_k}^{\ii\nabla} (T_{\phi_{\delta}}^{\ii\nabla} (P_m  X_k  P_m)) - \sum_{k=1}^d T_{\psi_k}^{\ii\nabla} (P_m  Y_k  P_m).
$$
By Lemma \ref{lem:XY1}, one has
$$
X_k, Y_k \in \calL_{\frac{d}{1-\epsilon},\infty}(L_2(\RR^d)) \subseteq \calL_p(L_2(\RR^d)),\quad \frac{d}{1-\epsilon} < p < \infty.
$$
We pass to the limit $m \to \infty$. Since $P_m \to 1$ strongly, we have
$$
P_m  X_k  P_m \to X_k,\quad P_m  Y_k  P_m \to Y_k \quad \text{ in } \calL_p(L_2(\RR^d)).
$$
The operators $T_{\phi_{\delta}}^{\ii\nabla}$ and $T_{\psi_k}^{\ii\nabla}$ are bounded on $\calL_p(L_2(\RR^d))$ by Lemma \ref{lem:TLL1} and Lemma \ref{lem:TLL2}, respectively.
It follows that
$$
T_{\psi_k}^{\ii\nabla}(T_{\phi_{\delta}}^{\ii\nabla}(P_m  X_k  P_m)) \to T_{\psi_k}^{\ii\nabla}(T_{\phi_{\delta}}^{\ii\nabla}(X_k)),
\quad T_{\psi_k}^{\ii\nabla}(P_m  Y_k  P_m) \to T_{\psi_k}^{\ii\nabla}(Y_k)
$$
in $\calL_p(L_2(\RR^d))$. Thus
$$
\calL_p(L_2(\RR^d)) \ni T_\infty := \lim_{m\to\infty} T_m
= \sum_{k=1}^d T_{\psi_k}^{\ii\nabla}(T_{\phi_{\delta}}^{\ii\nabla}(X_k)) - \sum_{k=1}^d T_{\psi_k}^{\ii\nabla}(Y_k).
$$
It remains to show that $T_\infty$ is the unique bounded extension of $[(-{\Delta})^{\frac{\epsilon}{2}}, x]$ to $L_2(\RR^d)$. By Lemma \ref{lem:epmg1}, for every $\phi \in \calS(\RR^d)$, one has
$$
T_m \phi \to (-{\Delta})^{\frac{\epsilon}{2}} x \phi - x (-{\Delta})^{\frac{\epsilon}{2}} \phi = [(-{\Delta})^{\frac{\epsilon}{2}}, x]\phi
$$
in $\calS'(\RR^d)$. This implies
$$
T_\infty \phi = [(-{\Delta})^{\frac{\epsilon}{2}}, x]\phi,\quad \phi \in \calS(\RR^d),
$$
which proves the claim.
\end{proof}

We can combine these estimates to obtain the upper bound for Schwartz symbols.

\begin{proposition}\label{prop:DW1}
Let $\epsilon \in (-\frac{d}{2},1)$. For every $x \in \calS(\RR^d_\theta)$, we have
$$\| [(-\Delta)^{\frac{\epsilon}{2}}, x ] \|_{\calL_{\frac{d}{1-\epsilon},\infty}(L_2(\RR^d))}
\lesssim_{d,\epsilon} \|x\|_{\dW^1_{\frac{d}{1-\epsilon}}(\RR^d_\theta)}.$$
\end{proposition}

\begin{proof}
By applying the triangle inequality to Lemma \ref{lem:XY2} and using Lemma \ref{lem:TLL1} and Lemma \ref{lem:TLL2}, we obtain
\begin{align*}
&\| [(-\Delta)^{\frac{\epsilon}{2}}, x ] \|_{\calL_{\frac{d}{1-\epsilon},\infty}(L_2(\RR^d))}\\
\lesssim_{d,\epsilon}&\sum_{k=1}^d\Big( \|X_k\|_{\calL_{\frac{d}{1-\epsilon},\infty}(L_2(\RR^d))}
+ \|Y_k\|_{\calL_{\frac{d}{1-\epsilon},\infty}(L_2(\RR^d))}
\Big).
\end{align*}
We conclude the desired result by Lemma \ref{lem:XY1}.
\end{proof}

By Proposition \ref{prop:hd1}, we obtain the following:

\begin{proposition}\label{prop:clapp}
Let $-\frac d2<\epsilon\neq0<1$ and $x\in\dW^1_{\frac{d}{1-\epsilon}}(\RR^d_\theta)$. Then
\begin{equation}\label{eq:clapp}
\begin{split}
\|[(-\Delta)^{\frac{\epsilon}{2}},x]\|_{\calL_{\frac{d}{1-\epsilon},\infty}(L_2(\RR^d))}
\lesssim_{d,\epsilon}\|x\|_{\dW^1_{\frac{d}{1-\epsilon}}(\RR^d_\theta)}.
\end{split}
\end{equation}
\end{proposition}

\begin{proof}
Let $x\in \dW^1_{\frac{d}{1-\epsilon}}(\RR^d_\theta)$. By Proposition \ref{prop:hd1}, there exists a sequence $\{x_n\}\subseteq\calS(\RR^d_\theta)$ satisfying
$$\|x_n-x\|_{\dW^1_{\frac{d}{1-\epsilon}}(\RR^d_\theta)}\to0.$$
By Proposition \ref{prop:tr1} and Proposition \ref{prop:DW1},
\begin{align*}
&\|[(-\Delta)^{\frac{\epsilon}{2}},x_n]-[(-\Delta)^{\frac{\epsilon}{2}},x_m]\|_{\calL_{\frac{d}{1-\epsilon},\infty}(L_2(\RR^d))}\\
=&\|[(-\Delta)^{\frac{\epsilon}{2}},x_n-x_m]\|_{\calL_{\frac{d}{1-\epsilon},\infty}(L_2(\RR^d))}\\
\quad\lesssim&_{d,\epsilon} \|x_n-x_m\|_{\dW^1_{\frac{d}{1-\epsilon}}(\RR^d_\theta)} \leq
\|x_n-x\|_{\dW^1_{\frac{d}{1-\epsilon}}(\RR^d_\theta)}
+\|x_m-x\|_{\dW^1_{\frac{d}{1-\epsilon}}(\RR^d_\theta)}\to0.
\end{align*}
In particular, $[(-\Delta)^{\frac{\epsilon}{2}},x_n]$ converges in $\mathcal B (L_2(\RR^d))$. Denote by $T$ the limit operator. Instead, we consider $V[(-\Delta)^{\frac{\epsilon}{2}},x_n]V^*=[(-\Delta_\theta)^{\frac{\epsilon}{2}},M_{x_n}]$. By Proposition \ref{prop:dw1} (iv), one has
$$\langle [(-\Delta_\theta)^{\frac{\epsilon}{2}},M_{x_n}-M_x]Vf, Vg\rangle\to0,\quad f,g\in\calS(\RR^d).$$
By the density of $\calS(\RR^d_\theta)$ in $L_2(\RR^d_\theta)$, the sequence $[(-\Delta_\theta)^{\frac{\epsilon}{2}},M_{x_n}]$ converges to $[(-\Delta_\theta)^{\frac{\epsilon}{2}},M_x]$ in $\mathcal{B}(L_2(\RR^d_\theta))$, which gives
$$ VTV^*=[(-\Delta_\theta)^{\frac{\epsilon}{2}},M_x].$$
This implies that $T=[(-\Delta)^{\frac{\epsilon}{2}},x]$. By \eqref{eq:bn1} and Proposition \ref{prop:DW1},
\begin{align*}
\|[(-\Delta)^{\frac{\epsilon}{2}},x]\|_{\calL_{\frac{d}{1-\epsilon},\infty}(L_2(\RR^d))}
=&\lim_{n\to\infty}
\|[(-\Delta)^{\frac{\epsilon}{2}},x_n]\|_{\calL_{\frac{d}{1-\epsilon},\infty}(L_2(\RR^d))}\\
\lesssim&_{d,\epsilon}\lim_{n\to\infty}
\|x_n\|_{\dW^1_{\frac{d}{1-\epsilon}}(\RR^d_\theta)}
=\|x\|_{\dW^1_{\frac{d}{1-\epsilon}}(\RR^d_\theta)}.
\end{align*}
This completes the proof of \eqref{eq:clapp}.
\end{proof}

\bigskip

\section{The principal part of $[(-\Delta)^{\frac{\epsilon}{2}},x]$}\label{s:prin}
In this section, we aim to show the following proposition, which is vital in the proofs of Theorem \ref{thm:m1} and Theorem \ref{thm:allp}. Let $\varLambda_j(t)=t_j$ denote the coordinate functions on $\BS^{d-1}$.

\begin{proposition}\label{prop:hom}
Let $x\in\calS(\RR^d_\theta)$ and $-\frac{d}{2}<\epsilon\ne0<1$. Then
\begin{equation}\label{eq:hom}
[(-\Delta)^{\frac{\epsilon}{2}},x]-\epsilon\sum_{j=1}^d(\partial_jx)\pi_2(\varLambda_j)J^{-(1-\epsilon)}\in\big(\calL_{\frac{d}{1-\epsilon},\infty}(L_2(\RR^d))\big)_0.
\end{equation}
\end{proposition}

The remainder of this section is devoted to the proof of Proposition \ref{prop:hom}. To this end, we introduce some notation. Let
$$A=J^2=1-\Delta,\quad A_{1,x}=[A,x],\quad A_{2,x}=[A,A_{1,x}].$$
The following identities will be used later, which follow directly from definition.
\begin{lemma}\label{lem:iA1}
Let $x\in\calS(\RR^d_\theta)$. Then
\begin{description}
\item[(i)] $[\calD_k^2,x]=2(\partial_kx)\calD_k+\partial_k^2x$;
\item[(ii)] $A_{1,x}=2\sum_{j=1}^d(\partial_jx)\calD_j+\sum_{j=1}^d\partial_j^2x$;
\item[(iii)] $A_{2,x}=4\sum_{j,k=1}^d(\partial_j\partial_kx)\calD_j\calD_k
+4\sum_{j,k=1}^d(\partial_j\partial_k^2x)\calD_j
+\sum_{j,k=1}^d\partial_j^2\partial_k^2x$.
\end{description}
\end{lemma}

\begin{proof}
(i) Note that \begin{align*}
\big[\mathcal{D}_k^2,  x\big]
&= \mathcal{D}_k \big[\mathcal{D}_k,  x\big] + \big[\mathcal{D}_k,  x\big] \mathcal{D}_k \\
&= \mathcal{D}_k (\partial_k x) + (\partial_k x) \mathcal{D}_k \\
&= (\partial_k x) \mathcal{D}_k + \big[\mathcal{D}_k,  \partial_k x\big] + (\partial_k x) \mathcal{D}_k \\
&= 2 (\partial_k x) \mathcal{D}_k + \partial_k^2 x.
\end{align*}

\noindent (ii) By (i), one has
\begin{align*}
A_1(x)
=&\sum_{k=1}^d[\mathcal D_k^2,x]=2\sum_{k=1}^d(\partial_kx)\mathcal D_k
+\sum_{k=1}^d\partial_k^2x.
\end{align*}

\noindent (iii) We use (i) and the Leibniz rule to obtain
\begin{align*}
A_{2,x}=[-\Delta,A_{1,x}]
=&2\sum_{j=1}^d[-\Delta,(\partial_jx)\mathcal D_j]
+\big[-\Delta,\sum_{j=1}^d\partial_j^2x\big]\\
=&2\sum_{j=1}^d([-\Delta,\partial_jx]\mathcal D_j
+\partial_jx[-\Delta,\mathcal D_j])
+\big[-\Delta,\sum_{j=1}^d\partial_j^2x\big]\\
=&2\sum_{j=1}^d[-\Delta,\partial_jx]\mathcal D_j
+\big[-\Delta,\sum_{j=1}^d\partial_j^2x\big].
\end{align*}
By using (i) to $\partial_jx$ and
$\sum_{j=1}^d\partial_j^2x$ respectively, we obtain
$$[-\Delta,\partial_jx]=2\sum_{k=1}^d(\partial_k\partial_jx)\mathcal D_k+\sum_{k=1}^d\partial_k^2\partial_jx$$
and
$$\big[-\Delta,\sum_{j=1}^d\partial_j^2x\big]
=2\sum_{k=1}^d\partial_k\Big(\sum_{j=1}^d\partial_j^2x\Big)\mathcal D_k
+\sum_{k=1}^d\partial_k^2\sum_{j=1}^d\partial_j^2x.$$
Since $\partial_k^2\partial_j=\partial_j\partial_k^2$,
$$A_{2,x}=4\sum_{j,k=1}^d(\partial_j\partial_kx)\mathcal D_j\mathcal D_k
+4\sum_{j,k=1}^d(\partial_j\partial_k^2x)\mathcal D_j
+\sum_{j,k=1}^d\partial_j^2\partial_k^2x.$$
This completes the proof.
\end{proof}

Instead of considering $[(-\Delta)^{\frac{\epsilon}{2}},x]$ directly, we deal with $[J^\epsilon,x]$, since in this case $J^{-1}$ is bounded and much easier to handle. For the principal part of $[J^\epsilon,x]$, we follow the proof of \cite[Theorem 5.1]{FSZ2024}. To pass from $J^\epsilon$ to $(-\Delta)^{\frac{\epsilon}{2}}$, we use the argument in the proof of \cite[Lemma 7.1]{FSZ2024}. The required remainder estimates for $x\in\calS(\RR^d_\theta)$ follow from Propositions \ref{prop:col2}, \ref{prop:sdw1} and \ref{prop:lower}. At first, we give the principal term of $[J^\epsilon,x]$ for $0<\epsilon<1$.

\begin{lemma}\label{lem:pres}
Let $x\in\calS(\RR^d_\theta)$ and $0<\epsilon<1$. Then
\begin{equation}\label{eq:posres}
[J^\epsilon,x]=\frac{\epsilon}{2} A_{1,x}J^{\epsilon-2}-\frac{\sin(\pi\frac{\epsilon}{2})}{\pi}\int_0^\infty\lambda^{\frac{\epsilon}{2}}(A+\lambda)^{-1}A_{2,x}(A+\lambda)^{-2} d\lambda.
\end{equation}
For every $p\in(d,\infty)$, the remainder integral converges in $\calL_p(L_2(\RR^d))$.
\end{lemma}

\begin{proof}
By functional calculus, one can write
$$J^\epsilon
=\frac{\sin(\pi\frac{\epsilon}{2})}{\pi}\int_0^\infty
A(A+\lambda)^{-1}
\lambda^{\frac{\epsilon}{2}-1} d\lambda.$$
Note that
\begin{align}\label{eq:Ax2}
\begin{split}
[(A+\lambda)^{-1},x]
=&-(A+\lambda)^{-1}[A+\lambda,x](A+\lambda)^{-1}\\
=&-(A+\lambda)^{-1}[A,x](A+\lambda)^{-1}\\
=&-(A+\lambda)^{-1}A_{1,x}(A+\lambda)^{-1},
\end{split}
\end{align}
and
$$
A(A+\lambda)^{-1}
=1-\lambda(A+\lambda)^{-1}.
$$
Hence, we have
\begin{align*}
[J^\epsilon,x]
=&\frac{\sin(\pi\frac{\epsilon}{2})}{\pi}\int_0^\infty
\big[A(A+\lambda)^{-1},x\big]
\lambda^{\frac{\epsilon}{2}-1} d\lambda\\
=&-\frac{\sin(\pi\frac{\epsilon}{2})}{\pi}\int_0^\infty
\big[(A+\lambda)^{-1},x\big]
\lambda^{\frac{\epsilon}{2}} d\lambda\\
=&\frac{\sin(\pi\frac{\epsilon}{2})}{\pi}\int_0^\infty
(A+\lambda)^{-1}A_{1,x}(A+\lambda)^{-1}
\lambda^{\frac{\epsilon}{2}} d\lambda.
\end{align*}
Since
\begin{align}\label{eq:Ax3}
\begin{split}
(A+\lambda)^{-1}A_{1,x}
=&A_{1,x}(A+\lambda)^{-1}
+\big[(A+\lambda)^{-1},A_{1,x}\big]\\
=&A_{1,x}(A+\lambda)^{-1}
-(A+\lambda)^{-1}[A,A_{1,x}](A+\lambda)^{-1}\\
=&A_{1,x}(A+\lambda)^{-1}
-(A+\lambda)^{-1}A_{2,x}(A+\lambda)^{-1},
\end{split}
\end{align}
we get
\begin{align*}
[J^\epsilon,x]
=&\frac{\sin(\pi\frac{\epsilon}{2})}{\pi}\int_0^\infty \lambda^{\frac{\epsilon}{2}}\big[A_{1,x}(A+\lambda)^{-2}-(A+\lambda)^{-1}A_{2,x}(A+\lambda)^{-2}\big]d\lambda\\
=&\frac{\epsilon}{2}A_{1,x}J^{\epsilon-2}-\frac{\sin(\pi\frac{\epsilon}{2})}{\pi}\int_0^\infty\lambda^{\frac{\epsilon}{2}}(A+\lambda)^{-1}A_{2,x}(A+\lambda)^{-2}d\lambda.
\end{align*}
By Lemma \ref{lem:iA1},
\begin{align*}
A_{2,x}J^{-3}
=&4\sum_{j,k}\Big((\partial_j\partial_kx)J^{-1}\Big)(J^{-2}\calD_j\calD_k)\\
&+4\sum_{j,k}\Big((\partial_j\partial_k^2x)J^{-1}\Big)J^{-1}(\calD_jJ^{-1})
+\sum_{j,k}\Big((\partial_j^2\partial_k^2x)J^{-1}\Big)J^{-2}.
\end{align*}
By Proposition \ref{prop:sdw1}, $A_{2,x}J^{-3}\in\calL_p(L_2(\RR^d))$ for $p>d$. Thus
\begin{align*}
&\|(A+\lambda)^{-1}A_{2,x}(A+\lambda)^{-2}\|_{\calL_p(L_2(\RR^d))}\\
=&\|(A+\lambda)^{-1}(A_{2,x}J^{-3})J^3(A+\lambda)^{-2}\|_{\calL_p(L_2(\RR^d))}\\
\lesssim&(1+\lambda)^{-3/2}\|A_{2,x}J^{-3}\|_{\calL_p(L_2(\RR^d))}.
\end{align*}
Since $\int_0^\infty\lambda^{\frac{\epsilon}{2}}(1+\lambda)^{-3/2}d\lambda<\infty$, we conclude that
$$\frac{\sin(\pi\frac{\epsilon}{2})}{\pi}\int_0^\infty\lambda^{\frac{\epsilon}{2}}(A+\lambda)^{-1}A_{2,x}(A+\lambda)^{-2} d\lambda$$
converges in $\calL_p(L_2(\RR^d))$. This completes the proof.
\end{proof}

Next, we turn to the principal part of $[J^{-\epsilon},x]$ for $0<\epsilon<1$. Let
$$\mathcal R_{\epsilon,x}=\frac{\sin(\pi\frac{\epsilon}{2})}{\pi}\int_0^\infty\lambda^{-\frac{\epsilon}{2}}(A+\lambda)^{-1}A_{2,x}(A+\lambda)^{-2} d\lambda.$$
In the same way as in Lemma \ref{lem:pres}, the above integral converges in $\calL_p(L_2(\RR^d))$ for $d<p<\infty$.

\begin{lemma}\label{lem:nres}
Let $x\in\calS(\RR^d_\theta)$ and $0<\epsilon<1$. Then the following statements hold.
\begin{description}
\item[(i)] One has
\begin{equation}\label{eq:negres}
[J^{-\epsilon},x]=-\frac\epsilon2A_{1,x}J^{-\epsilon-2}+\mathcal R_{\epsilon,x}.
\end{equation}
\item[(ii)] If $a,b\geq0$, $p>1$, and $p(a+b+2)>d$, then $J^{-a}\mathcal R_{\epsilon,x}J^{-b}\in\calL_p(L_2(\RR^d))$.
\end{description}
\end{lemma}

\begin{proof}
For $0<\epsilon<1$, by functional calculus, we get
$$
J^{-\epsilon}
=\frac{\sin(\pi\frac{\epsilon}{2})}{\pi}\int_0^\infty
\lambda^{-\frac{\epsilon}{2}}(A+\lambda)^{-1}d\lambda.
$$
By \eqref{eq:Ax2} and \eqref{eq:Ax3},
\begin{align*}
[J^{-\epsilon},x]=&\frac{\sin(\pi\frac{\epsilon}{2})}{\pi}\int_0^\infty
\lambda^{-\frac{\epsilon}{2}}[(A+\lambda)^{-1},x]d\lambda\\
=&-\frac{\sin(\pi\frac{\epsilon}{2})}{\pi}\int_0^\infty
\lambda^{-\frac{\epsilon}{2}}(A+\lambda)^{-1}A_{1,x}(A+\lambda)^{-1} d\lambda\\
=&-\frac{\sin(\pi\frac{\epsilon}{2})}{\pi}A_{1,x}
\int_0^\infty\lambda^{-\frac{\epsilon}{2}}(A+\lambda)^{-2} d\lambda+\mathcal R_{\epsilon,x}\\
=&-\frac\epsilon2A_{1,x}J^{-\epsilon-2}+\mathcal R_{\epsilon,x}.
\end{align*}
By Lemma \ref{lem:iA1}, we obtain
\begin{align*}
J^{-a}A_{2,x}J^{-b-4}
=&4\sum_{j,k}
\Big(J^{-a}(\partial_j\partial_kx)J^{-b-2}\Big)
(\mathcal D_j\mathcal D_kJ^{-2})\\
&+4\sum_{j,k}
\Big(J^{-a}(\partial_j\partial_k^2x)J^{-b-2}\Big)J^{-1}
(\mathcal D_jJ^{-1})
+\sum_{j,k}
\Big(J^{-a}(\partial_j^2\partial_k^2x)J^{-b-2}\Big)J^{-2}.
\end{align*}
Then by Proposition \ref{prop:sdw1},
$$
J^{-a}A_{2,x}J^{-b-4}\in\calL_p(L_2(\RR^d))
\quad\text{whenever}\quad p(a+b+2)>d.
$$
This completes the proof.
\end{proof}

Now we can give the principal part of $[J^{\epsilon},x]$ for the full range of $\epsilon$.

\begin{proposition}\label{prop:bessel}
Let $-\frac{d}{2}<\epsilon\neq 0<1$ and $x\in\calS(\RR^d_\theta)$. Then
\begin{equation}\label{eq:bessel}
[J^\epsilon,x]-\frac\epsilon2A_{1,x}J^{\epsilon-2}\in\big(\calL_{\frac{d}{1-\epsilon},\infty}(L_2(\RR^d))\big)_0.
\end{equation}
\end{proposition}

\begin{proof}
For $0<\epsilon<1$, choose $d<p<\frac{d}{1-\epsilon}$. We apply Lemma \ref{lem:pres} and the inclusion $\calL_p(L_2(\RR^d))\subseteq(\calL_{\frac{d}{1-\epsilon},\infty}(L_2(\RR^d)))_0$ to get \eqref{eq:bessel}. It remains to consider the case $-\frac{d}{2}<\epsilon<0$. Choose $m\in\NN$ so that $\delta=-\frac\epsilon m\in(0,1)$. Since $\frac{d}{1-\epsilon}>1$ and $2-\epsilon-\delta>1-\epsilon$, we can choose
\begin{equation}\label{eq:qneg}
\max\Big\{1,\frac d{2-\epsilon-\delta}\Big\}<p<\frac{d}{1-\epsilon}.
\end{equation}
By the Leibniz identity,
\begin{align}\label{eq:J1}
[J^\epsilon,x]=\sum_{l=0}^{m-1}J^{-l\delta}[J^{-\delta},x]J^{-(m-1-l)\delta}.
\end{align}
For $a=l\delta$ and $b=(m-1-l)\delta$, we have $a+b+2=2-\epsilon-\delta$. By Lemma \ref{lem:nres} and \eqref{eq:qneg},
$$J^{-l\delta}\mathcal R_{\delta,x}J^{-(m-1-l)\delta}\in\calL_p(L_2(\RR^d)),\quad 0\leq l<m.$$
From Lemma \ref{lem:iA1}, one can write
\begin{align}\label{eq:JAAJE1}
J^{-l\delta}A_{1,x}J^{-(m-l)\delta-2}
=A_{1,x}J^{-m\delta-2}+E_l,
\end{align}
where
\begin{align*}
E_l:=&2\sum_j\Big([J^{-l\delta},\partial_jx]J^{-(m-l)\delta-1}\Big)(\calD_jJ^{-1})\\
&+\sum_j\Big([J^{-l\delta},\partial_j^2x]J^{-(m-l)\delta-1}\Big)J^{-1}.
\end{align*}
By Proposition \ref{prop:lower},
$$[J^{-l\delta},\partial_jx]J^{-(m-l)\delta-1},\quad
[J^{-l\delta},\partial_j^2x]J^{-(m-l)\delta-1}\in\calL_p(L_2(\RR^d)),$$
so $E_l\in\calL_p(L_2(\RR^d))$. By \eqref{eq:negres}, \eqref{eq:J1} and \eqref{eq:JAAJE1},
$$[J^\epsilon,x]=\frac{\epsilon}{2}A_{1,x}J^{\epsilon-2}-\frac{\delta}{2}\sum_{l=0}^{m-1}E_l+\sum_{l=0}^{m-1}J^{-l\delta}\mathcal R_{\delta,x}J^{-(m-1-l)\delta}.$$
This implies 
\begin{align*}
[J^\epsilon,x]-\frac\epsilon2A_{1,x}J^{\epsilon-2}\in\calL_p(L_2(\RR^d)).
\end{align*}
Since $p<\frac{d}{1-\epsilon}$, we get \eqref{eq:bessel}. This completes the proof.
\end{proof}

Finally, we come to the proof of Proposition \ref{prop:hom}.

\begin{proof}[Proof of Proposition \ref{prop:hom}]
Recall that $\omega(t)=(1+|t|^2)^{\frac{1}{2}}$, and choose a real $\chi\in C_c^\infty(\RR^d)$ such that it is equal to one near zero. Fix
$$
\max\Big\{1,\frac d{2-\epsilon}\Big\}<p<\frac{d}{1-\epsilon}.
$$
By Proposition \ref{prop:sdw1}, one can see that if $y\in\calS(\RR^d_\theta)$, $p\beta>d$ and $|g(t)|\lesssim \omega(t)^{-\beta}$, then
\begin{equation}\label{eq:decay}
yM_g, M_gy\in\calL_p(L_2(\RR^d)).
\end{equation}
Let $h_\epsilon(t)=|t|^\epsilon-\omega(t)^\epsilon$, $h_0=\chi h_\epsilon$ and $h_\infty=(1-\chi)h_\epsilon$. Note that $|h_\infty(t)|\lesssim \omega(t)^{-(2-\epsilon)}$. Moreover, $h_0$ is compactly supported and belongs to $L_2(\RR^d)$, and it is because for $\epsilon>0$ it is bounded, and for $\epsilon<0$ this follows from $\int_0^1r^{2\epsilon+d-1} dr<\infty$. By Proposition \ref{prop:col2} and \eqref{eq:decay},
$$xM_{h_0},\ M_{h_0}x,\ xM_{h_\infty},\ M_{h_\infty}x\in\calL_p(L_2(\RR^d)).$$
Then
\begin{equation}\label{eq:hom-bessel}
[(-\Delta)^{\frac{\epsilon}{2}},x]-[J^\epsilon,x]=[M_{h_\epsilon},x]\in\calL_p(L_2(\RR^d)).
\end{equation}

For $1\leq j\leq d$, let $b_j(0)=0$ and for $t\ne0$,
$$b_j(t)=t_j\omega(t)^{\epsilon-2}-\frac{t_j}{|t|}\omega(t)^{\epsilon-1}=-\frac{t_j\omega(t)^{\epsilon-2}}{|t|(\omega(t)+|t|)}.$$
This function is bounded near zero and satisfies $|(1-\chi(t))b_j(t)|\lesssim \omega(t)^{-(3-\epsilon)}$. By using Proposition \ref{prop:col2} to $\chi b_j$ and \eqref{eq:decay} to $(1-\chi)b_j$, one has $(\partial_jx)M_{b_j}\in\calL_p(L_2(\RR^d))$. Note that $\calD_jJ^{\epsilon-2}-\pi_2(\varLambda_j)J^{\epsilon-1}=M_{b_j}$, and Proposition \ref{prop:sdw1} implies $(\partial_j^2x)J^{\epsilon-2}\in\calL_p(L_2(\RR^d))$. By using the formula for $A_{1,x}$ in Lemma \ref{lem:iA1}, we obtain
\begin{align*}
&\frac\epsilon2 A_{1,x}J^{\epsilon-2}-\epsilon\sum_j(\partial_jx)\pi_2(\varLambda_j)J^{\epsilon-1}\\
=&\epsilon\sum_{j=1}^d(\partial_jx)M_{b_j}+\frac\epsilon2\sum_{j=1}^d(\partial_j^2x)J^{\epsilon-2}
\in\calL_p(L_2(\RR^d)).
\end{align*}
Combining \eqref{eq:hom-bessel}, Proposition \ref{prop:bessel} and the inclusion $\calL_p(L_2(\RR^d))\subset\Big(\calL_{\frac{d}{1-\epsilon},\infty}(L_2(\RR^d))\Big)_0$, we obtain
\begin{align*}
&[(-\Delta)^{\frac{\epsilon}{2}},x]-\epsilon\sum_{j=1}^d(\partial_jx)\pi_2(\varLambda_j)J^{-(1-\epsilon)}\\
=&[(-\Delta)^{\frac{\epsilon}{2}},x]-[J^\epsilon,x]+[J^\epsilon,x]-\frac\epsilon2A_{1,x}J^{\epsilon-2}+\frac\epsilon2 A_{1,x}J^{\epsilon-2}-\epsilon\sum_j(\partial_jx)\pi_2(\varLambda_j)J^{\epsilon-1}\\
\in&\Big(\calL_{\frac{d}{1-\epsilon},\infty}(L_2(\RR^d))\Big)_0.
\end{align*}
This completes the proof.
\end{proof}

\bigskip

\section{Asymptotics of singular values}\label{s:spec}
In this section, we are devoted to establishing asymptotics of singular values for arbitrary antisymmetric real matrices $\theta$. For quantum Euclidean spaces with $\det\theta\neq0$, Tian \cite{T2025} established asymptotics of singular values for quantized derivatives. For arbitrary antisymmetric real $d\times d$ matrices $\theta$, we establish asymptotics of singular values for finite sums, namely Theorem \ref{thm:allp}. For $g\in L_\infty(\BS^{d-1})$, define
\begin{equation}\label{eq:ghat}
\widetilde g(t)=g\Big(\frac{t}{|t|}\Big)\quad(t\neq0), \quad\text{and}\quad \widetilde g(0)=0.
\end{equation}
Define
$$\pi_1(x)\xi=x\xi \quad\text{and}\quad
\pi_2(g)\xi=M_{\widetilde g}\xi=\widetilde g(\ii\nabla)\xi,
\quad \xi\in L_2(\RR^d),\ x\in L_\infty(\RR^d_\theta).$$
Thus $g$ is a function on the sphere, while $\widetilde g$ denotes its extension to Euclidean spaces. Let $J=(1-\Delta)^{\frac{1}{2}}$. Recall that $dS$ denotes the standard non-normalized surface measure on the unit sphere $\BS^{d-1}$ in $\RR^d$. This section is devoted to the proof of the following theorem.

\begin{theorem}\label{thm:allp}
Let $\alpha\in(0,d]$. For every finite family $x_j\in\calS(\RR^d_\theta)$ and $g_j\in C(\BS^{d-1})$, one has
$$\lim_{t\to\infty}t^{\frac{\alpha}{d}}
\mu\Big(t,\sum_j\pi_1(x_j)\pi_2(g_j)J^{-\alpha}\Big)
=(d(2\pi)^d)^{-\frac{\alpha}{d}}
\|\sum_jx_j\otimes g_j\|_{L_{\frac{d}{\alpha}}(L_\infty(\RR^d_\theta)\overline\otimes L_\infty(\BS^{d-1}))}.$$
\end{theorem}

Tian \cite[Theorem 7.1]{T2025} established asymptotics of singular values for finite sums with $J^{-1}$ under the assumption $\det\theta\ne0$. Here, we consider $J^{-\alpha}$, $0<\alpha\leq d$, for arbitrary real antisymmetric matrices $\theta$. The proof follows \cite[Sections 6--7]{T2025} by passing from positive elements to finite sums. 
%We derive the trace identity in Lemma \ref{lem:zeta} from the Hilbert--Schmidt identity and verify the hypotheses of Theorem \ref{thm:tr1} using the Cwikel and commutator estimates established before.
To carry out this argument when $\theta$ is degenerate, we derive the trace identity in Lemma \ref{lem:zeta} directly from the Hilbert--Schmidt identity and verify the hypotheses of Theorem \ref{thm:tr1} by using the Cwikel and commutator estimates established before.

%We proceed with the following useful lemma.

\begin{lemma}\label{lem:zeta}
Let $0<\beta<1$, $0\leq x\in\calS(\RR^d_\theta)$ and $0\leq h\in C^\infty(\BS^{d-1})$. Let $G=\pi_2(h)J^{-\beta}$. For $\Re z>\frac{d}{\beta}$, we have
\begin{equation}\label{eq:szeta}
\Tr(G^zx^z)
=(2\pi)^{-d}\frac{\Gamma(\frac{d}{2})\Gamma((\beta z-d)/2)}{2\Gamma(\beta z/2)}\tau_\theta(x^z)
\int_{\BS^{d-1}}h(s)^z dS.
\end{equation}
Moreover, the function $\Tr(G^zx^z)-\frac{c}{z-\frac{d}{\beta}}$ extends continuously to $\Re z\geq \frac{d}{\beta}$, where
\begin{equation}\label{eq:sres}
c=\frac{(2\pi)^{-d}}{\beta}\tau_\theta(x^{\frac{d}{\beta}})
\int_{\BS^{d-1}}h(s)^{\frac{d}{\beta}} dS.
\end{equation}

\end{lemma}
\begin{proof}
Let $z_0=\frac{d}{\beta}>2$. For $z\in\CC$, we define $u^z=e^{z\log u}$ for $u>0$ and $0^z=0$ for $\Re z>0$. Recall that $\widetilde h(t)$ is defined in \eqref{eq:ghat}. For $s:=\Re z>z_0$, let
$$g_z(t)=\widetilde h(t)^z(1+|t|^2)^{-\frac{\beta z}{2}}.$$
By Proposition \ref{prop:col2}, $xP_0\in\calL_1(L_2(\RR^d))$.
For $t\in n+Q$, we have $1+|n|\lesssim_d(1+|t|^2)^{\frac12}$, so
$$\|g_{\bar z}\|_{\ell_1(L_\infty(\RR^d);Q)}
\lesssim_{d,\beta,s}\|h\|_{L_\infty(\BS^{d-1})}^{s}
\sum_{n\in\ZZ^d}(1+|n|)^{-\beta s}<\infty.$$
It follows from \eqref{eq:cw-lo} and $\|x^{\bar z}P_0\|_{\calL_1(L_2(\RR^d))}\leq\|x\|_{L_\infty(\RR^d_\theta)}^{s-1}\|xP_0\|_{\calL_1(L_2(\RR^d))}$ that
$$\|G^zx^z\|_{\calL_1(L_2(\RR^d))}
\lesssim\|x^{\bar z}P_0\|_{\calL_1(L_2(\RR^d))}
\|g_{\bar z}\|_{\ell_1(L_\infty(\RR^d);Q)}<\infty.$$
Let $B_R=\{t\in\RR^d:|t|<R\}$ and $P_R=M_{\mathds{1}_{B_R}}$. By Lemma \ref{lem:hs2},
\begin{align*}
\Tr(P_RG^zx^zP_R)=&\Tr\Big((x^{\frac{\bar z}{2}}G^{\bar z}P_R)^*(x^{\frac{z}{2}}P_R)\Big)\\
=&(2\pi)^{-d}\tau_\theta\big((x^{\bar z/2})^*x^{z/2}\big)
\int_{B_R}\overline{g_{\bar z}(t)} dt\\
=&(2\pi)^{-d}\tau_\theta(x^z)\int_{B_R}g_z(t) dt.
\end{align*}
\iffalse
We write $T=G^zx^z=v|T|$ for the polar decomposition. Since $P_R\to1$ strongly, H\"older's inequality gives
\begin{align*}
&\|T-P_RTP_R\|_{\calL_1(L_2(\RR^d))}\\
\leq&\|(1-P_R)T\|_{\calL_1(L_2(\RR^d))}
+\|T(1-P_R)\|_{\calL_1(L_2(\RR^d))}\\
\leq&\|T\|_{\calL_1(L_2(\RR^d))}^{\frac{1}{2}}
\Big(\|(1-P_R)v|T|^{\frac{1}{2}}\|_{\calL_2(L_2(\RR^d))}
+\||T|^{\frac{1}{2}}(1-P_R)\|_{\calL_2(L_2(\RR^d))}\Big)
\to0.
\end{align*}
\fi
As $R\to \infty$, we obtain
\begin{align*}
\Tr(G^zx^z)
=&(2\pi)^{-d}\tau_\theta(x^z)\int_{\RR^d}g_z(t) dt\\
=&(2\pi)^{-d}\tau_\theta(x^z)\int_{\BS^{d-1}}h(s)^z dS(s)\int_0^\infty r^{d-1}(1+r^2)^{-\frac{\beta z}{2}} dr.\\
=&(2\pi)^{-d}\frac{\Gamma(\frac{d}{2})\Gamma((\beta z-d)/2)}{2\Gamma(\beta z/2)}\tau_\theta(x^z)
\int_{\BS^{d-1}}h(s)^z dS.
\end{align*}

Next, we verify the analytic extension directly. 
\iffalse
Let
$$M=\max\{1,\|x\|_{L_\infty(\RR^d_\theta)},\|h\|_{L_\infty(\BS^{d-1})}\}.$$
Fix $z$ with $\Re z>1$ and choose $0<\delta<\Re z-1$. For $0<u\leq M$ and $0<|s|<\delta$, one has
$$\Big|\frac{u^{z+s}-u^z}{s}-u^z\log u\Big|
=\Big|s\int_0^1(1-t)u^{z+ts}(\log u)^2 dt\Big|
\lesssim_{z,\delta,M}|s|u.$$
All logarithmic expressions are assigned the value zero at $u=0$. Thus
$$\Big\|\frac{x^{z+s}-x^z}{s}-x^z\log x\Big\|_{L_1(\RR^d_\theta)}\lesssim_{z,\delta,M}|s|\tau_\theta(x)\to0$$
and
$$\int_{\BS^{d-1}}\Big|\frac{h(s)^{z+s}-h(s)^z}{s}-h(s)^z\log h(s)\Big| dS(s)\lesssim_{z,\delta,M}|w|\int_{\BS^{d-1}}h(s) dS(s)\to0,$$
as $w\to0$. 
\fi
It is clear that on $\Re z>1$,
$$\frac{d}{dz}\tau_\theta(x^z)=\tau_\theta(x^z\log x)$$
and
$$\frac{d}{dz}\int_{\BS^{d-1}}h(s)^z dS(s)=\int_{\BS^{d-1}}h(s)^z\log h(s) dS(s).$$
Define
$$\Omega_*=\Big\{z\in\CC:\Re z>\max\Big\{1,\frac{d-2}{\beta}\Big\}\Big\}$$
and
$$
F(z)=\frac{(2\pi)^{-d}}{\beta}\tau_\theta(x^z)\int_{\BS^{d-1}}h(s)^z dS(s)\frac{\Gamma(d/2)\Gamma(1+(\beta z-d)/2)}{\Gamma(\beta z/2)}.
$$
For $z\in\Omega_*$, $\re z>1$, $\re\Big(1+\frac{\beta z-d}{2}\Big)>0$ and $\re\frac{\beta z}{2}>0$. Thus $F$ is holomorphic on $\Omega_*$ and
$$F(z_0)=\frac{(2\pi)^{-d}}{\beta}\tau_\theta(x^{d/\beta})
\int_{\BS^{d-1}}h(s)^{d/\beta} dS(s)=c.$$
Since
$$\frac{\Gamma(d/2)\Gamma((\beta z-d)/2)}{2\Gamma(\beta z/2)}
=\frac{1}{\beta(z-z_0)} \frac{\Gamma(d/2)\Gamma(1+(\beta z-d)/2)}{\Gamma(\beta z/2)},$$
one has from \eqref{eq:szeta} that
$$\Tr(G^zx^z)=\frac{F(z)}{z-z_0},\quad z\in \Omega_*.$$
For $\re z>z_0$, 
\begin{align*}
\Tr(G^zx^z)-\frac{c}{z-z_0}
=&\frac{F(z)-F(z_0)}{z-z_0}\\
=&\int_0^1 F'\big(z_0+t(z-z_0)\big) dt,
\end{align*}
which is holomorphic on $\Omega_*$. This completes the proof.
\end{proof}

The following theorem is taken from \cite{T2025}, while the condition below first appeared in \cite{SZ2023}. More precisely, the following statement is \cite[Theorem 6.4]{T2025}; its hypotheses are \cite[Condition 6.3]{T2025}, which is attributed there to \cite[Condition 5.4.1]{SZ2023}.

\begin{theorem}\label{thm:tr1}
Let $p>2$ and $0\le \mathbf A,\mathbf B\in\calB(H)$ satisfy the following conditions:
\begin{description}
\item[(i)] $\mathbf A^p \mathbf B \in \mathcal{L}_{1,\infty}(H)$;
\item[(ii)] $\mathbf A^{q-2}[\mathbf A,\mathbf B] \in \mathcal{L}_1(H)$ for every $q>p$;
\item[(iii)] $\mathbf B^{\frac12} \mathbf A \mathbf B^{\frac12} \in \mathcal{L}_{p,\infty}(H)$;
\item[(iv)] $\big[\mathbf A,\mathbf B^{\frac12}\big] \in \mathcal{L}_{\frac{p}{2},\infty}(H)$.
\end{description} If there is a constant $c>0$ such that the function
$$z\mapsto \operatorname{Tr}(\mathbf A^z \mathbf B^z)-\frac{c}{z-p},\quad z\in\mathbb{C},\quad \Re(z)>p,$$
admits a continuous extension to the closed half plane $\{z\in\mathbb{C}:\Re(z)\ge p\}$, then
$$\lim_{t\to\infty} t^{\frac{1}{p}}\mu(t,\mathbf A\mathbf B)=\Big(\frac{c}{p}\Big)^{\frac{1}{p}}.$$
\end{theorem}

The following lemma will be needed in the proof of spectral asymptotic formulae.

\begin{lemma}
Let $\alpha>0$ and $g\in L_\infty(\BS^{d-1})$. Then
\begin{equation}\label{eq:smult}
\|\widetilde g \omega^{-\alpha}\|_{L_{\frac{d}{\alpha},\infty}(\RR^d)}
+\|\widetilde g \omega^{-\alpha}\|_{\ell_{\frac{d}{\alpha},\infty}(L_\infty(\RR^d);Q)}
\lesssim_{d,\alpha}\|g\|_{L_\infty(\BS^{d-1})}.
\end{equation}
\end{lemma}
\begin{proof}
For $\lambda>0$, by polar coordinates, one has
\begin{align*}
|\{t\in\RR^d:|\widetilde g(t)|\omega(t)^{-\alpha}>\lambda\}|
=&\int_{\BS^{d-1}}
\Big(\Big(\frac{|g(\zeta)|}{\lambda}\Big)^{\frac{2}{\alpha}}-1\Big)_+^{\frac{d}{2}} d\zeta\\
\leq&\frac{1}{\lambda^{\frac{d}{\alpha}}}\int_{\BS^{d-1}}|g(\zeta)|^{\frac{d}{\alpha}} d\zeta
\lesssim_d\Big(\frac{\|g\|_{L_\infty(\BS^{d-1})}}{\lambda}\Big)^{\frac{d}{\alpha}}.
\end{align*}
This proves the $L_{\frac{d}{\alpha},\infty}(\RR^d)$-bound. Since $t\in n+Q$ implies
$$1+|n|\leq1+|t|+\sqrt d\leq\sqrt2(1+\sqrt d)\omega(t),$$
we have
$$a_n:=\|\widetilde g\omega^{-\alpha}\mathds{1}_{n+Q}\|_{L_\infty(\RR^d)}
\lesssim_{d,\alpha}\|g\|_{L_\infty(\BS^{d-1})}(1+|n|)^{-\alpha}.$$
This implies
$$\#\{n\in\ZZ^d:a_n>\lambda\}
\lesssim_{d,\alpha}\Big(\frac{\|g\|_{L_\infty(\BS^{d-1})}}{\lambda}\Big)^{\frac{d}{\alpha}}.$$
This completes the proof.
\end{proof}

At first, we establish the spectral asymptotic formulae for positive Schwartz symbols.

\begin{lemma}\label{lem:core}
Let $0<\beta<1$, $0\leq x,x^{\frac{1}{2}}\in\calS(\RR^d_\theta)$ and $0\leq h\in C^\infty(\BS^{d-1})$. Then
$$\lim_{t\to\infty}t^{\beta/d}\mu(t,x\pi_2(h)J^{-\beta})
=\big(d(2\pi)^d\big)^{-\beta/d}\|x\|_{L_{\frac{d}{\beta}}(\RR^d_\theta)}\|h\|_{L_{\frac{d}{\beta}}(\BS^{d-1})}.$$
\end{lemma}
\begin{proof}
If $x=0$ or $h=0$, Lemma \ref{lem:core} follows immediately. In the following, we assume that $x\ne0$ and $h\ne0$, so the constant $c$ in \eqref{eq:sres} is positive by the faithfulness of $\tau_\theta$ and the continuity of $h$. We will use Theorem \ref{thm:tr1}.

Now we verify the four conditions listed in Theorem \ref{thm:tr1} for $\mathbf A=\pi_2(h)J^{-\beta}$ and $\mathbf B=x$. By Theorem \ref{thm:cwl1} and \eqref{eq:smult},
$$x\big(\pi_2(h)J^{-\beta}\big)^{\frac{d}{\beta}}=x\pi_2(h^{\frac{d}{\beta}})J^{-d}\in\calL_{1,\infty}(L_2(\RR^d)),$$
which implies $\mathbf A^{\frac{d}{\beta}}\mathbf B\in \calL_{1,\infty}(L_2(\RR^d))$. For every $p>\frac{d}{\beta}$, $\pi_2(h^{p-2})$ is bounded. Then by Corollary \ref{cor:ang1},
\begin{align*}
\big(\pi_2(h)J^{-\beta}\big)^{p-2}\big[\pi_2(h)J^{-\beta},x\big]
=&\pi_2(h^{p-2})J^{-\beta(p-2)}[\pi_2(h)J^{-\beta},x]\\
\in&\calL_{\frac d{\beta(p-1)+1},\infty}(L_2(\RR^d))
\subseteq\calL_1(L_2(\RR^d)),
\end{align*}
since $\beta(p-1)+1>d+1-\beta>d$. This implies $\mathbf A^{p-2}[\mathbf A,\mathbf B] \in \mathcal{L}_1(L_2(\RR^d))$ for every $p>\frac{d}{\beta}$. Next, by Theorem \ref{thm:cw-hi} and \eqref{eq:smult},
\begin{align*}
\|\mathbf B^{\frac12} \mathbf A \mathbf B^{\frac12}\|_{\calL_{\frac{d}{\beta},\infty}(L_2(\RR^d))}=&\|x^{\frac{1}{2}}\pi_2(h)J^{-\beta}x^{\frac{1}{2}}\|_{\calL_{\frac{d}{\beta},\infty}(L_2(\RR^d))}\\
\lesssim&_{d,\beta}\|x^{\frac{1}{2}}\|_{L_{\frac{d}{\beta}}(\RR^d_\theta)}\|h\|_{L_\infty(\BS^{d-1})}\|x^{\frac{1}{2}}\|_{\calB(L_2(\RR^d))}<\infty.
\end{align*}
By Corollary \ref{cor:ang1},
$$\big[\mathbf A,\mathbf B^{\frac12}\big]=\big[\pi_2(h)J^{-\beta},x^{\frac{1}{2}}\big]
\in\calL_{\frac{d}{\beta+1},\infty}(L_2(\RR^d))
\subseteq\calL_{\frac{d}{2\beta},\infty}(L_2(\RR^d)).$$
Therefore, by Lemma \ref{lem:zeta} and Theorem \ref{thm:tr1},
$$\lim_{t\to\infty}t^{\beta/d}\mu(t,x\pi_2(h)J^{-\beta})=\big(d(2\pi)^d\big)^{-\beta/d}\|x\|_{L_{\frac{d}{\beta}}(\RR^d_\theta)}\|h\|_{L_{\frac{d}{\beta}}(\BS^{d-1})}.$$
This completes the proof.
\end{proof}

The following lemma will be used to extend Lemma \ref{lem:core} to $0<\beta\leq d$.

\begin{lemma}\label{lem:power}
Let $p\in\NN$, $0<\beta<1$ and $0\leq x,x^{\frac{1}{2}}\in\calS(\RR^d_\theta)$. If $0\leq h\in C^\infty(\BS^{d-1})$, then
$$
x^p\big(\pi_2(h)J^{-\beta}\big)^p
-\big(x^{\frac{1}{2}}\pi_2(h)J^{-\beta}x^{\frac{1}{2}}\big)^p
\in\calL_{\frac{d}{p\beta+1},\infty}(L_2(\RR^d)).
$$
\end{lemma}

\begin{proof}
By Theorem \ref{thm:cw-hi} and \eqref{eq:smult}, $x^{\frac12}\pi_2(h)J^{-\beta}x^{\frac12}\in\calL_{\frac{d}{\beta},\infty}(L_2(\RR^d))$. For $0\leq j\leq p-1$, we have $x^{p-j-\frac12}=x^{p-j-1}x^{\frac12}\in\calS(\RR^d_\theta)$. By Corollary \ref{cor:ang1},
\begin{align*}
&[x^{p-j-\frac{1}{2}},\pi_2(h)J^{-\beta}]
\big(\pi_2(h)J^{-\beta}\big)^{p-j-1}\\
=&\big([x^{p-j-\frac{1}{2}},\pi_2(h)J^{-\beta}]
J^{-\beta(p-j-1)}\big)\pi_2(h^{p-j-1})\in
\calL_{\frac{d}{(p-j)\beta+1},\infty}(L_2(\RR^d)).
\end{align*}
Thus by H\"older's inequality, we obtain for $1\leq j\leq p-1$ that
$$
\big(x^{\frac{1}{2}}\pi_2(h)J^{-\beta}x^{\frac{1}{2}}\big)^j
x^{\frac{1}{2}}[x^{p-j-\frac{1}{2}},\pi_2(h)J^{-\beta}]
\big(\pi_2(h)J^{-\beta}\big)^{p-j-1}\in\calL_{\frac{d}{p\beta+1},\infty}(L_2(\RR^d)),
$$
and it is clear that this estimate still holds for $j=0$ by the boundedness of $x^{\frac12}$. Finally, we expand the commutators to obtain
\begin{align*}
&\sum_{j=0}^{p-1}
\big(x^{\frac{1}{2}}\pi_2(h)J^{-\beta}x^{\frac{1}{2}}\big)^j
x^{\frac{1}{2}}[x^{p-j-\frac{1}{2}},\pi_2(h)J^{-\beta}]
\big(\pi_2(h)J^{-\beta}\big)^{p-j-1}\\
=&\sum_{j=0}^{p-1}\Big(
\big(x^{\frac{1}{2}}\pi_2(h)J^{-\beta}x^{\frac{1}{2}}\big)^j
x^{p-j}\big(\pi_2(h)J^{-\beta}\big)^{p-j}\\
&-\big(x^{\frac{1}{2}}\pi_2(h)J^{-\beta}x^{\frac{1}{2}}\big)^{j+1}
x^{p-j-1}\big(\pi_2(h)J^{-\beta}\big)^{p-j-1}
\Big)\\
=&x^p\big(\pi_2(h)J^{-\beta}\big)^p
-\big(x^{\frac{1}{2}}\pi_2(h)J^{-\beta}x^{\frac{1}{2}}\big)^p.
\end{align*}
Since each summand belongs to $\calL_{\frac{d}{p\beta+1},\infty}(L_2(\RR^d))$, this completes the proof.
\end{proof}

In what follows, we will show Theorem \ref{thm:allp}. We proceed with the proof for a special case of Theorem \ref{thm:allp} (involving only one term).

\begin{lemma}\label{lem:pos}
Let $\alpha\in(0,d]$ with $\frac{d}{\alpha}\ne2$, $x\in\calS(\RR^d_\theta)$ and $g\in C(\BS^{d-1})$. Then $x\pi_2(g)J^{-\alpha}\in\calL_{\frac{d}{\alpha},\infty}(L_2(\RR^d))$ and
\begin{align}\label{eq:sp1}
\lim_{t\to\infty}t^{\frac{\alpha}{d}}\mu(t,x\pi_2(g)J^{-\alpha})
=\big(d(2\pi)^d\big)^{-\frac{\alpha}{d}}\|x\|_{L_{\frac{d}{\alpha}}(\RR^d_\theta)}\|g\|_{L_{\frac{d}{\alpha}}(\BS^{d-1})}.
\end{align}
\end{lemma}

\begin{proof}
By Theorem \ref{thm:cw-hi}, Theorem \ref{thm:cwl1} and \eqref{eq:smult}, one has $x\pi_2(g)J^{-\alpha}\in\calL_{\frac{d}{\alpha},\infty}(L_2(\RR^d))$. To show \eqref{eq:sp1}, we first handle $|x|\pi_2(g)J^{-\alpha}$ with $g\geq0$, keeping $x\in\calS(\RR^d_\theta)$ and approximating $|x|$ by Lemma \ref{lem:root}. Choose an integer $m>\alpha$, and so $0<\frac{\alpha}{m}<1$ and $\frac{md}{\alpha}>d\geq2$. Then we take $g=h^m$ with $0\leq h\in C^\infty(\BS^{d-1})$. By Lemma \ref{lem:root},
$$|x|^2(|x|^2+\delta)^{\frac{1}{4m}-1},\quad
\big(|x|^2(|x|^2+\delta)^{\frac{1}{4m}-1}\big)^2
\in\calS(\RR^d_\theta).$$
By Corollary \ref{cor:ang1},
\begin{align*}
&|x|^2(|x|^2+\delta)^{\frac{1}{4m}-1}\pi_2(h)J^{-\frac{\alpha}{m}}|x|^2(|x|^2+\delta)^{\frac{1}{4m}-1}\\
&-\big(|x|^2(|x|^2+\delta)^{\frac{1}{4m}-1}\big)^2\pi_2(h)J^{-\frac{\alpha}{m}}\\
=&|x|^2(|x|^2+\delta)^{\frac{1}{4m}-1}\Big[\pi_2(h)J^{-\frac{\alpha}{m}},|x|^2(|x|^2+\delta)^{\frac{1}{4m}-1}\Big]\\
\in&\calL_{\frac{md}{\alpha+m},\infty}(L_2(\RR^d))
\subseteq\big(\calL_{\frac{md}{\alpha},\infty}(L_2(\RR^d))\big)_0.
\end{align*}
By Lemma \ref{lem:core} and Lemma \ref{lem:asst1}, we obtain
\begin{align*}
&\lim_{t\to\infty}t^{\frac{\alpha}{md}}\mu\Big(t,|x|^2(|x|^2+\delta)^{\frac{1}{4m}-1}\pi_2(h)J^{-\frac{\alpha}{m}}|x|^2(|x|^2+\delta)^{\frac{1}{4m}-1}\Big)\\
=&\big(d(2\pi)^d\big)^{-\frac{\alpha}{md}}\Big\|\big(|x|^2(|x|^2+\delta)^{\frac{1}{4m}-1}\big)^2\Big\|_{L_{\frac{md}{\alpha}}(\RR^d_\theta)}\|h\|_{L_{\frac{md}{\alpha}}(\BS^{d-1})}.
\end{align*}
This implies
\begin{align*}
&\lim_{t\to\infty}t^{\frac{\alpha}{d}}\mu\Big(t,\Big(|x|^2(|x|^2+\delta)^{\frac{1}{4m}-1}\pi_2(h)J^{-\frac{\alpha}{m}}|x|^2(|x|^2+\delta)^{\frac{1}{4m}-1}\Big)^m\Big)\\
=&\big(d(2\pi)^d\big)^{-\frac{\alpha}{d}}\Big\|\big(|x|^2(|x|^2+\delta)^{\frac{1}{4m}-1}\big)^{2m}\Big\|_{L_{\frac{d}{\alpha}}(\RR^d_\theta)}\|g\|_{L_{\frac{d}{\alpha}}(\BS^{d-1})}.
\end{align*}
By Lemma \ref{lem:power},
\begin{align*}
&\big(|x|^2(|x|^2+\delta)^{\frac{1}{4m}-1}\big)^{2m}
\pi_2(g)J^{-\alpha}\\
&-\Big(|x|^2(|x|^2+\delta)^{\frac{1}{4m}-1}\pi_2(h)J^{-\frac{\alpha}{m}}|x|^2(|x|^2+\delta)^{\frac{1}{4m}-1}\Big)^m\\
\in&\calL_{\frac{d}{\alpha+1},\infty}(L_2(\RR^d))
\subseteq\big(\calL_{\frac{d}{\alpha},\infty}(L_2(\RR^d))\big)_0.
\end{align*}
Applying Lemma \ref{lem:asst1} again, we obtain
\begin{equation}\label{eq:sapprox}
\begin{split}
&\lim_{t\to\infty}t^{\frac{\alpha}{d}}\mu\Big(t,\big(|x|^2(|x|^2+\delta)^{\frac{1}{4m}-1}\big)^{2m}\pi_2(g)J^{-\alpha}\Big)\\
=&\big(d(2\pi)^d\big)^{-\frac{\alpha}{d}}\Big\|\big(|x|^2(|x|^2+\delta)^{\frac{1}{4m}-1}\big)^{2m}\Big\|_{L_{\frac{d}{\alpha}}(\RR^d_\theta)}\|g\|_{L_{\frac{d}{\alpha}}(\BS^{d-1})}.
\end{split}
\end{equation}
By Lemma \ref{lem:root}, as $\delta\to0$,
$$\Big\|\big(|x|^2(|x|^2+\delta)^{\frac{1}{4m}-1}\big)^{2m}
-|x|\Big\|_{L_{\frac{d}{\alpha}}(\RR^d_\theta)}\to0.$$
If $\frac{d}{\alpha}>2$, Theorem \ref{thm:cw-hi} and \eqref{eq:smult} imply 
\begin{align}\label{eq:xg1}
\|x\pi_2(g)J^{-\alpha}\|_{\calL_{\frac{d}{\alpha},\infty}(L_2(\RR^d))}\lesssim_{d,\alpha}\|x\|_{L_{\frac{d}{\alpha}}(\RR^d_\theta)}
\|g\|_{L_\infty(\BS^{d-1})}
\end{align}
and as $\delta\to0$,
\begin{align*}
&\Big\|\Big(\big(|x|^2(|x|^2+\delta)^{\frac{1}{4m}-1}\big)^{2m}
-|x|\Big)\pi_2(g)J^{-\alpha}\Big\|_{\calL_{\frac{d}{\alpha},\infty}(L_2(\RR^d))}\\
\quad\lesssim&_{d,\alpha}\Big\|\big(|x|^2(|x|^2+\delta)^{\frac{1}{4m}-1}\big)^{2m}-|x|\Big\|_{L_{\frac{d}{\alpha}}(\RR^d_\theta)}
\|g\|_{L_\infty(\BS^{d-1})}\to0.
\end{align*}
When $1\leq\frac{d}{\alpha}<2$, by Proposition \ref{prop:col2}, $$\Big(\big(|x|^2(|x|^2+\delta)^{\frac{1}{4m}-1}\big)^{2m}-|x|\Big)P_0\in\calL_{\frac{d}{\alpha}}(L_2(\RR^d)).$$
It follows from Theorem \ref{thm:cwl1} and \eqref{eq:smult} that \begin{align}\label{eq:xg2}
\|x\pi_2(g)J^{-\alpha}\|_{\calL_{\frac{d}{\alpha},\infty}(L_2(\RR^d))}\lesssim_{d,\alpha}\|xP_0\|_{L_{\frac{d}{\alpha}}(\RR^d_\theta)}
\|g\|_{L_\infty(\BS^{d-1})}
\end{align}
and as $\delta\to0$,
\begin{align*}&\Big\|\Big(\big(|x|^2(|x|^2+\delta)^{\frac{1}{4m}-1}\big)^{2m}
-|x|\Big)\pi_2(g)J^{-\alpha}\Big\|_{\calL_{\frac{d}{\alpha},\infty}(L_2(\RR^d))}\\
\quad\lesssim&_{d,\alpha}\Big\|\Big(\big(|x|^2(|x|^2+\delta)^{\frac{1}{4m}-1}\big)^{2m}-|x|\Big)P_0\Big\|_{\calL_{\frac{d}{\alpha}}(L_2(\RR^d))}\|g\|_{L_\infty(\BS^{d-1})}
\to0,
\end{align*}
where we use \eqref{eq:slconv}. Applying Lemma \ref{lem:asapp} for $\delta=\frac{1}{n}$ to \eqref{eq:sapprox}, we obtain \eqref{eq:sp1} for $|x|\pi_2(g)J^{-\alpha}$ when $g=h^m$ with $0\leq h\in C^\infty(\BS^{d-1})$.

Next, we take $0\leq g\in C(\BS^{d-1})$ and choose $0\leq h_n\in C^\infty(\BS^{d-1})$ such that
$$\|h_n-g^{1/m}\|_{L_\infty(\BS^{d-1})}\to0.$$
This also implies $\|h_n^m-g\|_{L_\infty(\BS^{d-1})}\to0$. Note that \eqref{eq:xg1} and \eqref{eq:xg2} hold for any $g\in L_\infty(\BS^{d-1})$. Then one has
\begin{align*}
&\||x|\pi_2(h_n^m-g)J^{-\alpha}\|_
{\calL_{\frac{d}{\alpha},\infty}(L_2(\RR^d))}\\
\quad\lesssim&_{d,\alpha}
\begin{cases}
\|x\|_{L_{\frac{d}{\alpha}}(\RR^d_\theta)}
\|h_n^m-g\|_{L_\infty(\BS^{d-1})},
&\displaystyle \frac{d}{\alpha}>2,\\
\|xP_0\|_{\calL_{\frac{d}{\alpha}}(L_2(\RR^d))}
\|h_n^m-g\|_{L_\infty(\BS^{d-1})},
&\displaystyle 1\leq\frac{d}{\alpha}<2,
\end{cases}\\
\quad\to&0.
\end{align*}
Note that $\|h_n^m-g\|_{L_{\frac{d}{\alpha}}(\BS^{d-1})}\to0$ as well. Hence, by Lemma \ref{lem:asapp} again, we obtain \eqref{eq:sp1} for $|x|\pi_2(g)J^{-\alpha}$ with $g\geq0$.

Finally, for general $g\in C(\BS^{d-1})$, define $g/|g|=1$ on the zero set of $g$. Then $\pi_2(g/|g|)$ is unitary and
$$x\pi_2(g)J^{-\alpha}=x\pi_2(|g|)J^{-\alpha}\pi_2(g/|g|).$$
Moreover,
\begin{align*}
\big(x\pi_2(|g|)J^{-\alpha}\big)^*\big(x\pi_2(|g|)J^{-\alpha}\big)
=&J^{-\alpha}\pi_2(|g|)x^*x\pi_2(|g|)J^{-\alpha}\\
=&\big(|x|\pi_2(|g|)J^{-\alpha}\big)^*\big(|x|\pi_2(|g|)J^{-\alpha}\big).
\end{align*}
Therefore, we obtain
$$\mu(t,x\pi_2(g)J^{-\alpha})=\mu(t,|x|\pi_2(|g|)J^{-\alpha}).$$
This completes the proof.
\end{proof}

We next extend the result for one term to finite sums. We first assume that the functions $g_j$ have pairwise disjoint supports.

\begin{lemma}\label{lem:disj}
Let $\alpha\in(0,d]$ with $\frac{d}{\alpha}\ne2$ and $x_1,\ldots,x_M\in\calS(\RR^d_\theta)$. Let $g_1,\ldots,g_M\in C^\infty(\BS^{d-1})$ be nonnegative with pairwise disjoint supports. Then
\begin{equation}\label{eq:sdisj}
\lim_{t\to\infty}t^{\frac{\alpha}{d}}
\mu\Big(t,\sum_{k=1}^M x_k\pi_2(g_k)J^{-\alpha}\Big)
=\big(d(2\pi)^d\big)^{-\frac{\alpha}{d}}\Big(\sum_{k=1}^M
\|x_k\|_{L_{\frac{d}{\alpha}}(\RR^d_\theta)}^{\frac{d}{\alpha}}
\|g_k\|_{L_{\frac{d}{\alpha}}(\BS^{d-1})}^{\frac{d}{\alpha}}\Big)^{\frac{\alpha}{d}}.
\end{equation}
\end{lemma}

\begin{proof}
For each $k$ and $t>0$, recall that 
$$\mu\big(t,x_k\pi_2(g_k)J^{-\alpha}\big)
=\mu\big(t,|x_k|\pi_2(g_k)J^{-\alpha}\big).$$
Let
$$c_k=(d(2\pi)^d)^{-\frac{\alpha}{d}}
\|x_k\|_{L_{\frac{d}{\alpha}}(\RR^d_\theta)}
\|g_k\|_{L_{\frac{d}{\alpha}}(\BS^{d-1})}.$$
Define $d_T(\lambda)=\operatorname{Tr}(\mathds{1}_{(\lambda,\infty)}(|T|))$ for $T\in \mathcal{B}(H)$. By Lemma \ref{lem:pos}, we obtain
$$
\lim_{t\to\infty}t^{\frac{\alpha}{d}}
\mu\big(t,x_k\pi_2(g_k)J^{-\alpha}\big)=c_k, \quad \text{and} \quad
\lim_{\lambda\to0}\lambda^{\frac{d}{\alpha}}
d_{x_k\pi_2(g_k)J^{-\alpha}}(\lambda)=c_k^{\frac{d}{\alpha}}.
$$
For $k\ne l$, the disjointness of the supports yields
\begin{align*}
&\big|x_k\pi_2(g_k)J^{-\alpha}\big|^2
\big|x_l\pi_2(g_l)J^{-\alpha}\big|^2\\
=&J^{-\alpha}\pi_2(g_k)|x_k|^2\pi_2(g_k)J^{-\alpha}
J^{-\alpha}\pi_2(g_l)|x_l|^2\pi_2(g_l)J^{-\alpha}=0.
\end{align*}
Thus we obtain mutually orthogonal positive operators
$\big|x_k\pi_2(g_k)J^{-\alpha}\big|^2$. For
$$Z=\Big(\sum_k\big|x_k\pi_2(g_k)J^{-\alpha}\big|^2\Big)^{\frac12},$$
we have
$$d_Z(\lambda)=\sum_{k=1}^M d_{x_k\pi_2(g_k)J^{-\alpha}}(\lambda).$$
This implies
\begin{align*}
\lim_{t\to\infty}t^{\frac{\alpha}{d}}\mu(t,Z)
=\Big(\lim_{\lambda\to0}\lambda^{\frac{d}{\alpha}}d_Z(\lambda)\Big)^{\frac{\alpha}{d}}
=\Big(\sum_{k=1}^M c_k^{\frac{d}{\alpha}}\Big)^{\frac{\alpha}{d}}=\mathcal A,
\end{align*}
where $\mathcal A$ is the right-hand-side of \eqref{eq:sdisj}. By Lemma \ref{lem:pos}, both $x_k\pi_2(g_k)J^{-\alpha}$ and
$\pi_2(g_k)J^{-\alpha}x_k$ belong to
$\calL_{\frac{d}{\alpha},\infty}(L_2(\RR^d))$. By Corollary \ref{cor:ang1}, we have
\begin{align}\label{eq:xgJ1}
[x_k,\pi_2(g_k)J^{-\alpha}]\in\calL_{d/(\alpha+1),\infty}(L_2(\RR^d)).
\end{align}
By H\"older's inequality, we obtain
\begin{align*}
&\big|x_k\pi_2(g_k)J^{-\alpha}\big|^2
-\big|\pi_2(g_k)J^{-\alpha}x_k\big|^2\\
=&\big(x_k\pi_2(g_k)J^{-\alpha}\big)^*[x_k,\pi_2(g_k)J^{-\alpha}]
+[x_k,\pi_2(g_k)J^{-\alpha}]^*\big(\pi_2(g_k)J^{-\alpha}x_k\big)\\
\in&\calL_{d/(2\alpha+1),\infty}(L_2(\RR^d))
\subseteq\big(\calL_{d/(2\alpha),\infty}(L_2(\RR^d))\big)_0.
\end{align*}
Note that for $k\ne l$,
$$
\big(\pi_2(g_k)J^{-\alpha}x_k\big)^*
\big(\pi_2(g_l)J^{-\alpha}x_l\big)
=x_k^*J^{-\alpha}\pi_2(g_k)\pi_2(g_l)J^{-\alpha}x_l=0,
$$
and so
$$\Big|\sum_k\pi_2(g_k)J^{-\alpha}x_k\Big|^2
=\sum_k\big|\pi_2(g_k)J^{-\alpha}x_k\big|^2.$$
Since $Z^2=\sum_k\big|x_k\pi_2(g_k)J^{-\alpha}\big|^2$ and
$\mu(t,Z^2)=\mu(t,Z)^2$, we have
$$
\lim_{t\to\infty}t^{2\frac{\alpha}{d}}\mu(t,Z^2)=\mathcal A^2,
$$
and
\begin{align*}
Z^2-\Big|\sum_k\pi_2(g_k)J^{-\alpha}x_k\Big|^2
=&\sum_k\Big(\big|x_k\pi_2(g_k)J^{-\alpha}\big|^2
-\big|\pi_2(g_k)J^{-\alpha}x_k\big|^2\Big)\\
\in&\big(\calL_{d/(2\alpha),\infty}(L_2(\RR^d))\big)_0.
\end{align*}
By Lemma \ref{lem:asst1},
$$\lim_{t\to\infty}t^{\frac{\alpha}{d}}
\mu\Big(t,\sum_k\pi_2(g_k)J^{-\alpha}x_k\Big)=\mathcal A.$$
It then follows from Lemma \ref{lem:asst1} and \eqref{eq:xgJ1} that
$$
\lim_{t\to\infty}t^{\frac{\alpha}{d}}
\mu\Big(t,\sum_kx_k\pi_2(g_k)J^{-\alpha}\Big)=\mathcal A,
$$
which is \eqref{eq:sdisj}.
\end{proof}

We employ a standard approximation argument to remove the disjoint-support assumption in Lemma \ref{lem:disj}.

\begin{theorem}\label{thm:finite}
Let $d\geq2$ and $\alpha\in(0,d]$ with $\frac{d}{\alpha}\ne2$. For $x_1,\ldots,x_N\in\calS(\RR^d_\theta)$ and $g_1,\ldots,g_N\in C(\BS^{d-1})$, we have
\begin{equation}\label{eq:finite}
\lim_{t\to\infty}t^{\frac{\alpha}{d}}
\mu\Big(t,\sum_{j=1}^N\pi_1(x_j)\pi_2(g_j)J^{-\alpha}\Big)
=(d(2\pi)^d)^{-\frac{\alpha}{d}}
\|\sum_jx_j\otimes g_j\|_{L_{\frac{d}{\alpha}}(L_\infty(\RR^d_\theta)\overline\otimes L_\infty(\BS^{d-1}))}.
\end{equation}
\end{theorem}

\begin{proof}
We write
$$F=\sum_jx_j\otimes g_j\quad\text{and}\quad T=\sum_jx_j\pi_2(g_j)J^{-\alpha}.$$
For each $m\in \mathbb N$, take a finite Borel partition $\{E_{m,k}\}_{k=1}^{M_m}$ of $\BS^{d-1}$ such that there exist disjoint compact sets $K_{m,k}\subseteq E_{m,k}$ satisfying
$$\operatorname{diam}(E_{m,k})\le m^{-1},\quad
\Big|\BS^{d-1}\setminus\bigcup_{k=1}^{M_m}K_{m,k}\Big|<m^{-\frac{d}{\alpha}}.$$
After ignoring empty sets $K_{m,k}$, we can choose smooth functions $\varphi_{m,k}$ and points $s_{m,k}\in K_{m,k}$ such that
\begin{align*}
0\le\varphi_{m,k}\le1,\quad
\varphi_{m,k}|_{K_{m,k}}=1,\quad
\operatorname{diam}\big(\operatorname{supp}\varphi_{m,k}\big)\le3m^{-1},
\end{align*}
\begin{align*}
\operatorname{supp}\varphi_{m,k}\cap\operatorname{supp}\varphi_{m,l}=\varnothing\quad(k\ne l),
\end{align*}
and set
$$
g_{j,m}=\sum_k g_j(s_{m,k})\varphi_{m,k},\quad
\mathcal A_{m,k}=\sum_j g_j(s_{m,k})x_j.
$$
Let
$$
s_j(\delta)=\sup_{\substack{\xi,\eta\in\BS^{d-1}\\|\xi-\eta|\le\delta}}\big|g_j(\xi)-g_j(\eta)\big|, \quad \text{then} \ 
\quad s_j(\delta)\to0\quad(\delta\to0).
$$
Hence, as $m\to\infty$,
\begin{align*}
\|g_{j,m}-g_j\|_{L_{\frac{d}{\alpha}}(\BS^{d-1})}
\le |\BS^{d-1}|^{\frac{\alpha}{d}} s_j(3m^{-1})
+2m^{-1}\|g_j\|_{L_\infty(\BS^{d-1})}
\to0.
\end{align*}
Now define
$$F_m=\sum_k \mathcal A_{m,k}\otimes\varphi_{m,k}
=\sum_jx_j\otimes g_{j,m},\quad
T_m=\sum_k \mathcal A_{m,k}\pi_2(\varphi_{m,k})J^{-\alpha}.$$
By the triangle inequality, we obtain that as $m\to \infty$,
$$\|F_m-F\|_{L_{\frac{d}{\alpha}}(L_\infty(\RR^d_\theta)\overline\otimes L_\infty(\BS^{d-1}))}
\leq\sum_j\|x_j\|_{L_{\frac{d}{\alpha}}(\RR^d_\theta)}
\|g_{j,m}-g_j\|_{L_{\frac{d}{\alpha}}(\BS^{d-1})}\to0.$$
By Lemma \ref{lem:pos},
$$T,\ T_m\in\calL_{\frac{d}{\alpha},\infty}(L_2(\RR^d)).$$
Denote by
$$S_{j,m}=x_j\pi_2(g_{j,m}-g_j)J^{-\alpha},\quad
T_m-T=\sum_{j=1}^N S_{j,m}.$$
By Lemma \ref{lem:pos}, we have that for each $j$,
$$\limsup_{t\to\infty}t^{\frac{\alpha}{d}}\mu(t,S_{j,m})
=\big(d(2\pi)^d\big)^{-\frac{\alpha}{d}}
\|x_j\|_{L_{\frac{d}{\alpha}}(\RR^d_\theta)}
\|g_{j,m}-g_j\|_{L_{\frac{d}{\alpha}}(\BS^{d-1})}.$$
Fix $t>0$ and $\eta>0$. Then we can choose finite-rank operators $Y_{j,m}$ such that
$$\rank Y_{j,m}\leq\frac{t}{N},\quad
\|S_{j,m}-Y_{j,m}\|_{\calB(L_2(\RR^d))}\leq\mu(t/N,S_{j,m})+\eta.$$
This implies
$$\rank\Big(\sum_{j=1}^N Y_{j,m}\Big)\leq\sum_{j=1}^N\rank Y_{j,m}\leq t,$$
and hence
\begin{align*}
\mu(t,T_m-T)\leq&\Big\|\sum_{j=1}^N(S_{j,m}-Y_{j,m})\Big\|_{\calB(L_2(\RR^d))}\\
\leq&\sum_{j=1}^N\|S_{j,m}-Y_{j,m}\|_{\calB(L_2(\RR^d))}\\
\leq&\sum_{j=1}^N\mu(t/N,S_{j,m})+N\eta.
\end{align*}
By letting $\eta\to0$, we obtain
$$\mu(t,T_m-T)\leq\sum_{j=1}^N\mu(t/N,S_{j,m}).$$
Therefore, one has
\begin{align*}
\limsup_{t\to\infty}t^{\frac{\alpha}{d}}\mu(t,T_m-T)
\leq&\limsup_{t\to\infty}\sum_{j=1}^N t^{\frac{\alpha}{d}}\mu(t/N,S_{j,m})\\
\leq&\sum_{j=1}^N\limsup_{t\to\infty}t^{\frac{\alpha}{d}}\mu(t/N,S_{j,m})\\
=&N^{\frac{\alpha}{d}}\sum_{j=1}^N\limsup_{s\to\infty}s^{\frac{\alpha}{d}}\mu(s,S_{j,m}).
\end{align*}
Since $N$ is fixed and $g_{j,m}\to g_j$ in $L_{\frac{d}{\alpha}}(\BS^{d-1})$ for every $j$, we obtain that as $m\to\infty$,
$$\limsup_{t\to\infty}t^{\frac{\alpha}{d}}\mu(t,T_m-T)
\leq N^{\frac{\alpha}{d}}\big(d(2\pi)^d\big)^{-\frac{\alpha}{d}}
\sum_{j=1}^N\|x_j\|_{L_{\frac{d}{\alpha}}(\RR^d_\theta)}
\|g_{j,m}-g_j\|_{L_{\frac{d}{\alpha}}(\BS^{d-1})}\to0.$$
By disjointness of the supports, we have
$$\|F_m\|_{L_{\frac{d}{\alpha}}(L_\infty(\RR^d_\theta)\overline\otimes L_\infty(\BS^{d-1}))}^{\frac{d}{\alpha}}
=\sum_k\|\mathcal A_{m,k}\|_{L_{\frac{d}{\alpha}}(\RR^d_\theta)}^{\frac{d}{\alpha}}
\|\varphi_{m,k}\|_{L_{\frac{d}{\alpha}}(\BS^{d-1})}^{\frac{d}{\alpha}}.$$
By Lemma \ref{lem:disj},
$$\lim_{t\to\infty}t^{\frac{\alpha}{d}}\mu(t,T_m)
=(d(2\pi)^d)^{-\frac{\alpha}{d}}\|F_m\|_{L_{\frac{d}{\alpha}}(L_\infty(\RR^d_\theta)\overline\otimes L_\infty(\BS^{d-1}))}.$$
By Lemma \ref{lem:asapp}, in view of
$$\|F_m-F\|_{L_{\frac{d}{\alpha}}(L_\infty(\RR^d_\theta)\overline\otimes L_\infty(\BS^{d-1}))}\to0,\quad
\limsup_{t\to\infty}t^{\frac{\alpha}{d}}\mu(t,T_m-T)\to0,$$
we obtain
\begin{align*}
\lim_{t\to\infty}t^{\frac{\alpha}{d}}\mu(t,T)
=&\lim_{m\to\infty}\lim_{t\to\infty}t^{\frac{\alpha}{d}}\mu(t,T_m)\\
=&(d(2\pi)^d)^{-\frac{\alpha}{d}}
\lim_{m\to\infty}\|F_m\|_{L_{\frac{d}{\alpha}}(L_\infty(\RR^d_\theta)\overline\otimes L_\infty(\BS^{d-1}))}\\
=&(d(2\pi)^d)^{-\frac{\alpha}{d}}
\|F\|_{L_{\frac{d}{\alpha}}(L_\infty(\RR^d_\theta)\overline\otimes L_\infty(\BS^{d-1}))},
\end{align*}
which is exactly \eqref{eq:finite}. This completes the proof.
\end{proof}

Next, we turn to the remaining case $\frac{d}{\alpha}=2$.

\begin{theorem}\label{thm:p2}
For $x_1,\ldots,x_N\in\calS(\RR^d_\theta)$ and $g_1,\ldots,g_N\in C(\BS^{d-1})$, we have
\begin{equation}\label{eq:finite2}
\lim_{t\to\infty}t^{\frac{1}{2}}
\mu\Big(t,\sum_j\pi_1(x_j)\pi_2(g_j)J^{-\frac{d}{2}}\Big)
=(d(2\pi)^d)^{-\frac{1}{2}}
\Big\|\sum_jx_j\otimes g_j\Big\|_{L_2(L_\infty(\RR^d_\theta)\overline\otimes L_\infty(\BS^{d-1}))}.
\end{equation}
\end{theorem}

\begin{proof}
Fix $y\in\calS(\RR^d_\theta)$ and $v\in C(\BS^{d-1})$. Take $u\in C^\infty(\BS^{d-1})$. By Corollary \ref{cor:ang1},
$$
[y,\pi_2(\overline u)]J^{-d}\in\calL_{d/(d+1),\infty}(L_2(\RR^d))\subseteq\calL_1(L_2(\RR^d)).
$$
By the boundedness of $\pi_2(v)$,
$$
[\pi_2(\overline u),y]\pi_2(v)J^{-d}\in\calL_1(L_2(\RR^d)).
$$

Next, we take $u\in C(\BS^{d-1})$ and choose $u_n\in C^\infty(\BS^{d-1})$ with $\|u_n-u\|_{L_\infty(\BS^{d-1})}\to0$. By Proposition \ref{prop:col2}, we have $yP_0\in\calL_1(L_2(\RR^d))$. By Theorem \ref{thm:cwl1} and \eqref{eq:smult},
$$
yJ^{-d}\in\calL_{1,\infty}(L_2(\RR^d)).
$$
Note that
\begin{align*}
\big([\pi_2(\overline u),y]-[\pi_2(\overline{u_n}),y]\big)J^{-d}=\pi_2(\overline u-\overline{u_n})yJ^{-d}
-yJ^{-d}\pi_2(\overline u-\overline{u_n}).
\end{align*}
Consequently, as $n\to \infty$,
\begin{align*}
&\|\big([\pi_2(\overline u),y]-[\pi_2(\overline{u_n}),y]\big)J^{-d}\|_{\calL_{1,\infty}(L_2(\RR^d))}\\
\lesssim&\|u-u_n\|_{L_\infty(\BS^{d-1})}
\|yJ^{-d}\|_{\calL_{1,\infty}(L_2(\RR^d))}
\to0.
\end{align*}
For every $n$, we have $[\pi_2(\overline{u_n}),y]J^{-d}\in\calL_1(L_2(\RR^d))$. By the boundedness of $\pi_2(v)$, we obtain
$$
[\pi_2(\overline u),y]\pi_2(v)J^{-d}
\in\big(\calL_{1,\infty}(L_2(\RR^d))\big)_0,
\quad u,v\in C(\BS^{d-1}).
$$
By Proposition \ref{prop:cube} with symbol $h(t)=(1+|t|^2)^{-d/4}$, 
$$
[J^{-\frac{d}{2}},y]J^{-\frac{d}{2}}\in\calL_{d/(d+1),\infty}(L_2(\RR^d))\subseteq\calL_1(L_2(\RR^d)).
$$
By the boundedness of $\pi_2(\overline u),\pi_2(v)$,
$$
\pi_2(\overline u)[J^{-\frac{d}{2}},y]J^{-\frac{d}{2}}\pi_2(v)\in\calL_1(L_2(\RR^d))\subseteq\big(\calL_{1,\infty}(L_2(\RR^d))\big)_0.
$$
Therefore, we obtain
\begin{align*}
&J^{-\frac{d}{2}}\pi_2(\overline u)y\pi_2(v)J^{-\frac{d}{2}}
-y\pi_2(\overline u v)J^{-d}\\
=&[\pi_2(\overline u),y]\pi_2(v)J^{-d}
+\pi_2(\overline u)[J^{-\frac{d}{2}},y]J^{-\frac{d}{2}}\pi_2(v)\\
\quad\in&\big(\calL_{1,\infty}(L_2(\RR^d))\big)_0.
\end{align*}
This implies that
\begin{align*}
\sum_{j,k}J^{-\frac d2}\pi_2(\overline{g_j})x_j^*x_k\pi_2(g_k)J^{-\frac d2}
-\sum_{j,k}x_j^*x_k\pi_2(\overline{g_j}g_k)J^{-d}
\in\big(\calL_{1,\infty}(L_2(\RR^d))\big)_0.
\end{align*}
Since
\begin{align*}
&\sum_{j,k}J^{-\frac d2}\pi_2(\overline{g_j})x_j^*x_k\pi_2(g_k)J^{-\frac d2}\\
=&\Big(\sum_jx_j\pi_2(g_j)J^{-\frac d2}\Big)^*
\Big(\sum_kx_k\pi_2(g_k)J^{-\frac d2}\Big),
\end{align*}
one has
\begin{align}\label{eq:ssquare}
\begin{split}
&\Big(\sum_jx_j\pi_2(g_j)J^{-\frac d2}\Big)^*
\Big(\sum_kx_k\pi_2(g_k)J^{-\frac d2}\Big)-\sum_{j,k}x_j^*x_k\pi_2(\overline{g_j}g_k)J^{-d}\\
\in&\big(\calL_{1,\infty}(L_2(\RR^d))\big)_0.
\end{split}
\end{align}
By Theorem \ref{thm:finite} with $\alpha=d$, applied to the families $\{x_j^*x_k\}_{j,k}$ and $\{\overline{g_j}g_k\}_{j,k}$, and using
\begin{align*}
&\Big(\sum_jx_j\otimes g_j\Big)^*\Big(\sum_kx_k\otimes g_k\Big)
=\sum_{j,k}x_j^*x_k\otimes\overline{g_j}g_k,
\end{align*}
we obtain
\begin{align*}
&\lim_{t\to\infty}t \mu\Big(t,
\sum_{j,k}\pi_1(x_j^*x_k)\pi_2(\overline{g_j}g_k)J^{-d}\Big)\\
=&\frac1{d(2\pi)^d}
\Big\|\Big(\sum_jx_j\otimes g_j\Big)^*
\Big(\sum_kx_k\otimes g_k\Big)\Big\|_{L_1\big(L_\infty(\RR^d_\theta)\bar\otimes L_\infty(\BS^{d-1})\big)}\\
=&\frac1{d(2\pi)^d}
\Big\|\sum_jx_j\otimes g_j\Big\|_{L_2\big(L_\infty(\RR^d_\theta)\bar\otimes L_\infty(\BS^{d-1})\big)}^2.
\end{align*}
By Lemma \ref{lem:asst1} and \eqref{eq:ssquare},
\begin{align*}
&\lim_{t\to\infty}t^{\frac12}\mu\Big(t,
\sum_j\pi_1(x_j)\pi_2(g_j)J^{-\frac d2}\Big)\\
=&\Big(\lim_{t\to\infty}t \mu\Big(t,
\sum_{j,k}\pi_1(x_j^*x_k)\pi_2(\overline{g_j}g_k)J^{-d}\Big)\Big)^{\frac12}\\
=&\frac1{\sqrt{d(2\pi)^d}}
\Big\|\sum_jx_j\otimes g_j\Big\|_{L_2\big(L_\infty(\RR^d_\theta)\bar\otimes L_\infty(\BS^{d-1})\big)}.
\end{align*}
This completes the proof.
\end{proof}

Finally, we come to the proof of Theorem \ref{thm:allp}.

\begin{proof}[Proof of Theorem \ref{thm:allp}]
The case $\frac{d}{\alpha}\ne2$ follows from Theorem \ref{thm:finite}, and the case $\frac{d}{\alpha}=2$ follows from Theorem \ref{thm:p2}.
\end{proof}

\bigskip
\section{Proof of Theorem \ref{thm:m1}}\label{s:comm}
In this section, we complete the proof of Theorem \ref{thm:m1}. At first, we give the following spectral asymptotic of the commutator $[(-\Delta)^{\frac{\epsilon}{2}},x]$ for $x\in\calS(\RR^d_\theta)$.

\begin{lemma}\label{lem:smalaw}
Let $-\frac d2<\epsilon\neq 0<1$ and $x\in\calS(\RR^d_\theta)$. Then
$$[(-\Delta)^{\frac{\epsilon}{2}},x]\in\calL_{\frac{d}{1-\epsilon},\infty}(L_2(\RR^d)).$$
Moreover,
\begin{equation}\label{eq:slaw}
\lim_{t\to\infty}t^{\frac{1-\epsilon}{d}}\mu\big(t,[(-\Delta)^{\frac{\epsilon}{2}},x]\big)
=|\epsilon|(d(2\pi)^d)^{-\frac{1-\epsilon}{d}}\Big(\int_{\BS^{d-1}}\tau_\theta\Big(\Big|\sum_{j=1}^ds_j\partial_jx\Big|^{\frac{d}{1-\epsilon}}\Big)dS(s)\Big)^{\frac{1-\epsilon}{d}}.
\end{equation}
\end{lemma}

\begin{proof}
Proposition \ref{prop:DW1} yields
$$[(-\Delta)^{\frac{\epsilon}{2}},x]\in\calL_{\frac{d}{1-\epsilon},\infty}(L_2(\RR^d)).$$
Recall that $\varLambda_j(t)=t_j$ is the $j$-th coordinate function on $\BS^{d-1}$. Since $0<1-\epsilon<1+\frac{d}{2}\leq d$, $\partial_jx\in\calS(\RR^d_\theta)$ and $\varLambda_j\in C(\BS^{d-1})$, by Theorem \ref{thm:allp}, we have
\begin{equation}
\lim_{t\to\infty}t^{\frac{1-\epsilon}{d}}\mu\big(t,\epsilon\sum_{j=1}^d(\partial_jx)\pi_2(\varLambda_j)J^{-(1-\epsilon)}\big)
=|\epsilon|(d(2\pi)^d)^{-\frac{1-\epsilon}{d}}\Big(\int_{\BS^{d-1}}\tau_\theta\Big(\Big|\sum_{j=1}^ds_j\partial_jx\Big|^{\frac{d}{1-\epsilon}}\Big)dS(s)\Big)^{\frac{1-\epsilon}{d}}.
\end{equation}
By Proposition \ref{prop:hom},
$$
\lim_{t\to\infty}t^{\frac{1-\epsilon}{d}}
\mu\Big(t,[(-\Delta)^{\frac{\epsilon}{2}},x]
-\epsilon\sum_{j=1}^d(\partial_jx)\pi_2(\varLambda_j)J^{-(1-\epsilon)}\Big)=0.
$$
Therefore, Lemma \ref{lem:asst1} yields \eqref{eq:slaw}.
\end{proof}

Then one can use Lemma \ref{lem:smalaw} to obtain spectral asymptotics of $[(-\Delta)^{\frac{\epsilon}{2}},x]$ for $x\in\dW^1_{\frac{d}{1-\epsilon}}(\RR^d_\theta)$.

\begin{theorem}\label{thm:claw}
Let $-\frac d2<\epsilon\neq0<1$. For every $x\in\dW^1_{\frac{d}{1-\epsilon}}(\RR^d_\theta)$, one has
\begin{equation}\label{eq:asym0}
\lim_{t\to\infty}t^{\frac{1-\epsilon}{d}}\mu\big(t,[(-\Delta)^{\frac{\epsilon}{2}},x]\big)
=|\epsilon|(d(2\pi)^d)^{-\frac{1-\epsilon}{d}}\Big(\int_{\BS^{d-1}}\tau_\theta\Big(\Big|\sum_{j=1}^ds_j\partial_jx\Big|^{\frac{d}{1-\epsilon}}\Big)dS(s)\Big)^{\frac{1-\epsilon}{d}},
\end{equation}
and
$$
\|[(-\Delta)^{\frac{\epsilon}{2}},x]\|_{\calL_{\frac{d}{1-\epsilon},\infty}(L_2(\RR^d))}
\asymp_{d,\epsilon}\|x\|_{\dW^1_{\frac{d}{1-\epsilon}}(\RR^d_\theta)}.
$$
\end{theorem}

\begin{proof}
Choose a sequence $\{x_n\}\subseteq\calS(\RR^d_\theta)$ that converges to $x$ in $\dW^1_{\frac{d}{1-\epsilon}}(\RR^d_\theta)$ as in the proof of Proposition \ref{prop:clapp}. By Lemma \ref{lem:ge1},
\begin{align*}
&\Big|\Big(\int_{\BS^{d-1}}\tau_\theta\Big(\Big|\sum_js_j\partial_jx_n\Big|^{\frac{d}{1-\epsilon}}\Big)dS\Big)^{\frac{1-\epsilon}{d}}
-\Big(\int_{\BS^{d-1}}\tau_\theta\Big(\Big|\sum_js_j\partial_jx\Big|^{\frac{d}{1-\epsilon}}\Big)dS\Big)^{\frac{1-\epsilon}{d}}\Big|\\
\quad\leq&\Big(\int_{\BS^{d-1}}\tau_\theta\Big(\Big|\sum_js_j\partial_j(x_n-x)\Big|^{\frac{d}{1-\epsilon}}\Big)dS\Big)^{\frac{1-\epsilon}{d}}\\
\quad\lesssim&_{d,\epsilon}
\|x_n-x\|_{\dW^1_{\frac{d}{1-\epsilon}}(\RR^d_\theta)}\to0.
\end{align*}
Moreover, \eqref{eq:clapp} yields
\begin{align*}
&\limsup_{t\to\infty}t^{\frac{1-\epsilon}{d}}
\mu\big(t,[(-\Delta)^{\frac{\epsilon}{2}},x_n]-[(-\Delta)^{\frac{\epsilon}{2}},x]\big)\\
\leq&\|[(-\Delta)^{\frac{\epsilon}{2}},x_n]-[(-\Delta)^{\frac{\epsilon}{2}},x]\|_{\calL_{\frac{d}{1-\epsilon},\infty}(L_2(\RR^d))}\\
\lesssim&_{d,\epsilon}\|x_n-x\|_{\dW^1_{\frac{d}{1-\epsilon}}(\RR^d_\theta)}\to0.
\end{align*}
By Lemma \ref{lem:smalaw}, we obtain
\begin{align*}
&\lim_{n\to\infty}\lim_{t\to\infty}t^{\frac{1-\epsilon}{d}}
\mu\big(t,[(-\Delta)^{\frac{\epsilon}{2}},x_n]\big)\\
=&\lim_{n\to\infty}|\epsilon|(d(2\pi)^d)^{-\frac{1-\epsilon}{d}}
\Big(\int_{\BS^{d-1}}\tau_\theta\Big(\Big|\sum_js_j\partial_jx_n\Big|^{\frac{d}{1-\epsilon}}\Big)dS(s)\Big)^{\frac{1-\epsilon}{d}}\\
=&|\epsilon|(d(2\pi)^d)^{-\frac{1-\epsilon}{d}}\Big(\int_{\BS^{d-1}}\tau_\theta\Big(\Big|\sum_js_j\partial_jx\Big|^{\frac{d}{1-\epsilon}}\Big)dS(s)\Big)^{\frac{1-\epsilon}{d}}.
\end{align*}
Then \eqref{eq:asym0} follows from Lemma \ref{lem:asapp}. By Lemma \ref{lem:ge1} and \eqref{eq:asym0},
\begin{align*}
\|x\|_{\dW^1_{\frac{d}{1-\epsilon}}(\RR^d_\theta)}
\lesssim&_{d,\epsilon}\Big(\int_{\BS^{d-1}}\tau_\theta\Big(\Big|\sum_js_j\partial_jx\Big|^{\frac{d}{1-\epsilon}}\Big)dS\Big)^{\frac{1-\epsilon}{d}}\\
\asymp&_{d,\epsilon}\lim_{t\to\infty}t^{\frac{1-\epsilon}{d}}\mu\big(t,[(-\Delta)^{\frac{\epsilon}{2}},x]\big)\\
\leq&\|[(-\Delta)^{\frac{\epsilon}{2}},x]\|_{\calL_{\frac{d}{1-\epsilon},\infty}(L_2(\RR^d))}.
\end{align*}
The upper bound of $\|[(-\Delta)^{\frac{\epsilon}{2}},x]\|_{\calL_{\frac{d}{1-\epsilon},\infty}(L_2(\RR^d))}$ follows from Proposition \ref{prop:clapp}. This completes the proof.
\end{proof}

It remains to show that Theorem \ref{thm:m1} (i) implies Theorem \ref{thm:m1} (ii). To this end, we need the following proposition. Recall that $h_t$ is the heat kernel defined in Lemma \ref{lem:heat}. As in \cite[Lemma 7.2 and Lemma 7.3]{FSZ2024}, we use asymptotics of singular values to obtain a lower bound after multiplication on both sides. Here multiplication by $z\in\calS(\RR^d_\theta)$ and the approximation properties established in Section \ref{s:appr} allow us to derive this bound for $y_t=h_t*x$ with $x\in\calS'(\RR^d_\theta)$.

\begin{proposition}\label{prop:lasym}
Let $-\frac{d}{2}<\epsilon\neq 0<1$. Suppose that $x\in\calS'(\RR^d_\theta)$ satisfies $[(-\Delta)^{\frac{\epsilon}{2}},x]\in\calL_{\frac{d}{1-\epsilon},\infty}(L_2(\RR^d))$. Let $y_t=h_t*x$ for $t>0$. If $z=z^*\in\calS(\RR^d_\theta)$, then
\begin{equation}\label{eq:lasym}
\begin{split}
&\lim_{s\to\infty}s^{\frac{1-\epsilon}{d}}
\mu\big(s,z[(-\Delta)^{\frac{\epsilon}{2}},y_t]z\big)\\
=&|\epsilon|(d(2\pi)^d)^{-\frac{1-\epsilon}{d}}
\Big(\int_{\BS^{d-1}}\tau_\theta\Big(
\Big|\sum_{j=1}^ds_jz(\partial_jy_t)z\Big|^{\frac{d}{1-\epsilon}}
\Big)dS(s)\Big)^{\frac{1-\epsilon}{d}}.
\end{split}
\end{equation}
Moreover,
\begin{equation}\label{eq:lgrad}
\sum_{j=1}^d\|z(\partial_jy_t)z\|_{L_{\frac{d}{1-\epsilon}}(\RR^d_\theta)}
\lesssim_{d,\epsilon}\|z\|_{L_\infty(\RR^d_\theta)}^2
\|[(-\Delta)^{\frac{\epsilon}{2}},y_t]\|_{\calL_{\frac{d}{1-\epsilon},\infty}(L_2(\RR^d))}.
\end{equation}
\end{proposition}

The proof of Proposition \ref{prop:lasym} requires the following lemma. We compare $z[(-\Delta)^{\epsilon/2},y_t]z$ with $\epsilon\sum_j z(\partial_jy_t)z\pi_2(\varLambda_j)J^{-(1-\epsilon)}$ and prove that their difference belongs to $\big(\calL_{\frac{d}{1-\epsilon},\infty}(L_2(\RR^d))\big)_0$.

\begin{lemma}\label{lem:lprin}
Under the same assumption as in Proposition \ref{prop:lasym}, one has
\begin{equation}\label{eq:lprin}
z[(-\Delta)^{\frac{\epsilon}{2}},y_t]z
-\epsilon\sum_{j=1}^d z(\partial_jy_t)z\pi_2(\varLambda_j)J^{-(1-\epsilon)}
\in\big(\calL_{\frac{d}{1-\epsilon},\infty}(L_2(\RR^d))\big)_0.
\end{equation}
\end{lemma}

\begin{proof}
Fix $t>0$. By Lemma \ref{lem:sapp}, we can choose a sequence $(x_n)_{n\geq 1}\subseteq\calS(\RR^d_\theta)$ such that $x_n\to x$ strongly in $\calS'(\RR^d_\theta)$, and let $y_{t,n}=h_t*x_n$. By Lemma \ref{lem:heat} (i),
$$
\widehat{y_{t,n}}=e^{-t|\cdot|^2}\widehat{x_n}\in\calS(\RR^d),
\quad y_{t,n}\in\calS(\RR^d_\theta).
$$
By Proposition \ref{prop:hprod} and the continuous embedding
$\calS(\RR^d_\theta)\hookrightarrow L_\infty(\RR^d_\theta)$, one has
\begin{align}\label{eq:zyz1}
\begin{split}
z(\partial_jy_t)z
=&\lim_{n\to\infty}z(\partial_jy_{t,n})z\\
=&\lim_{n\to\infty}\Big(\partial_j(zy_{t,n}z)-(\partial_jz)y_{t,n}z-zy_{t,n}(\partial_jz)\Big)\\
=&\partial_j(zy_tz)-(\partial_jz)y_tz-(zy_t)(\partial_jz)
\end{split}
\end{align}
in $\calB(L_2(\RR^d))$.
By Lemma \ref{lem:heat} (iii), the commutator $[(-\Delta)^{\frac{\epsilon}{2}},y_t]$ is bounded. For $f,g\in\calS(\RR^d)$, we have $zf,zg\in\calS(\RR^d)$, and by Proposition \ref{prop:dw1} (i)
$$
\Phi_{\epsilon;V(zf),V(zg)}\in\calS(\RR^d_\theta).
$$
Using $z=z^*$ and \eqref{eq:cform0}, we obtain that as $n\to \infty$,
\begin{align*}
\big\langle z[(-\Delta)^{\frac{\epsilon}{2}},y_{t,n}]zf,\overline g\big\rangle
=&\big\langle[(-\Delta)^{\frac{\epsilon}{2}},y_{t,n}]zf,\overline{zg}\big\rangle\\
=&\langle y_{t,n},\Phi_{\epsilon;V(zf),V(zg)}\rangle\\
\to&\langle y_t,\Phi_{\epsilon;V(zf),V(zg)}\rangle\\
=&\big\langle z[(-\Delta)^{\frac{\epsilon}{2}},y_t]zf,\overline g\big\rangle.
\end{align*}
By Lemma \ref{lem:hprod}, Proposition \ref{prop:hprod}, Proposition \ref{prop:DW1} and the continuous embedding
$\calS(\RR^d_\theta)\hookrightarrow L_\infty(\RR^d_\theta)$, we have
\begin{equation}\label{eq:lprod}
\begin{split}
z[(-\Delta)^{\frac{\epsilon}{2}},y_t]z
=&\lim_{n\to\infty}z[(-\Delta)^{\frac{\epsilon}{2}},y_{t,n}]z\\
=&\lim_{n\to\infty}\Big([(-\Delta)^{\frac{\epsilon}{2}},zy_{t,n}z]
-[(-\Delta)^{\frac{\epsilon}{2}},z](y_{t,n}z)
-(zy_{t,n})[(-\Delta)^{\frac{\epsilon}{2}},z]\Big)\\
=&[(-\Delta)^{\frac{\epsilon}{2}},zy_tz]
-[(-\Delta)^{\frac{\epsilon}{2}},z](y_tz)
-(zy_t)[(-\Delta)^{\frac{\epsilon}{2}},z].
\end{split}
\end{equation}
in $\calB(L_2(\RR^d))$. By Proposition \ref{prop:hom},
\begin{align}\label{eq:z1}
[(-\Delta)^{\frac{\epsilon}{2}},zy_tz]
-\epsilon\sum_{j=1}^d\partial_j(zy_tz)\pi_2(\varLambda_j)J^{-(1-\epsilon)}
\in\big(\calL_{\frac{d}{1-\epsilon},\infty}(L_2(\RR^d))\big)_0
\end{align}
and
\begin{align}\label{eq:z2}
[(-\Delta)^{\frac{\epsilon}{2}},z]
-\epsilon\sum_{j=1}^d(\partial_jz)\pi_2(\varLambda_j)J^{-(1-\epsilon)}
\in\big(\calL_{\frac{d}{1-\epsilon},\infty}(L_2(\RR^d))\big)_0.
\end{align}
Since $y_tz$ and $zy_t$ are bounded, one has
\begin{align*}
[(-\Delta)^{\frac{\epsilon}{2}},z](y_tz)
-\epsilon\sum_{j=1}^d(\partial_jz)\pi_2(\varLambda_j)J^{-(1-\epsilon)}(y_tz)
\in\big(\calL_{\frac{d}{1-\epsilon},\infty}(L_2(\RR^d))\big)_0,
\end{align*}
and
\begin{align*}
(zy_t)[(-\Delta)^{\frac{\epsilon}{2}},z]
-\epsilon\sum_{j=1}^d (zy_t)(\partial_jz)\pi_2(\varLambda_j)J^{-(1-\epsilon)}
\in\big(\calL_{\frac{d}{1-\epsilon},\infty}(L_2(\RR^d))\big)_0.
\end{align*}

Now apply Corollary \ref{cor:ang1} with $a=b=0$, $\gamma=1-\epsilon$ and $g=\varLambda_j$. Since $y_tz\in\calS(\RR^d_\theta)$, we obtain
$$[\pi_2(\varLambda_j)J^{-(1-\epsilon)},y_tz]\in\calL_{\frac{d}{2-\epsilon},\infty}(L_2(\RR^d)).$$
Note that as $k\to\infty$,
\begin{align*}
&(k+1)^{\frac{1-\epsilon}{d}}\mu\big(k,(\partial_jz)[\pi_2(\varLambda_j)J^{-(1-\epsilon)},y_tz]\big)\\
\leq&\|\partial_jz\|_{L_\infty(\RR^d_\theta)}(k+1)^{\frac{1-\epsilon}{d}}\mu\big(k,[\pi_2(\varLambda_j)J^{-(1-\epsilon)},y_tz]\big)\\
\leq&\|\partial_jz\|_{L_\infty(\RR^d_\theta)}\|[\pi_2(\varLambda_j)J^{-(1-\epsilon)},y_tz]\|_{\calL_{\frac{d}{2-\epsilon},\infty}(L_2(\RR^d))}(k+1)^{-\frac1d}
\to0.
\end{align*}
This implies
\begin{align}\label{eq:z3}
(\partial_jz)[\pi_2(\varLambda_j)J^{-(1-\epsilon)},y_tz]
\in\big(\calL_{\frac{d}{1-\epsilon},\infty}(L_2(\RR^d))\big)_0.
\end{align}
By using \eqref{eq:lprod}, \eqref{eq:zyz1}, \eqref{eq:z1}, \eqref{eq:z2} and \eqref{eq:z3}, we obtain
\begin{align*}
&z[(-\Delta)^{\frac{\epsilon}{2}},y_t]z
-\epsilon\sum_{j=1}^d z(\partial_jy_t)z\pi_2(\varLambda_j)J^{-(1-\epsilon)}\\
=&\Big([(-\Delta)^{\frac{\epsilon}{2}},zy_tz]
-\epsilon\sum_{j=1}^d\partial_j(zy_tz)\pi_2(\varLambda_j)J^{-(1-\epsilon)}\Big)\\
&-\Big([(-\Delta)^{\frac{\epsilon}{2}},z]
-\epsilon\sum_{j=1}^d(\partial_jz)\pi_2(\varLambda_j)J^{-(1-\epsilon)}\Big)(y_tz)\\
&-(zy_t)\Big([(-\Delta)^{\frac{\epsilon}{2}},z]
-\epsilon\sum_{j=1}^d(\partial_jz)\pi_2(\varLambda_j)J^{-(1-\epsilon)}\Big)\\
&-\epsilon\sum_{j=1}^d(\partial_jz)[\pi_2(\varLambda_j)J^{-(1-\epsilon)},y_tz]\\
\in&\big(\calL_{\frac{d}{1-\epsilon},\infty}(L_2(\RR^d))\big)_0.
\end{align*}
This proves \eqref{eq:lprin}.
\end{proof}

\begin{proof}[Proof of Proposition \ref{prop:lasym}]
Note that $0<1-\epsilon<1+\frac{d}{2}\leq d$. By Lemma \ref{lem:hprod}, $z(\partial_jy_t)z\in\calS(\RR^d_\theta)$. By Lemma \ref{lem:lprin}, Theorem \ref{thm:allp} and Lemma \ref{lem:asst1}, we get \eqref{eq:lasym}.
\iffalse
Recall from Lemma \ref{lem:ge1} that
$$
z(\partial_ky_t)z
=\frac{d}{|\BS^{d-1}|}
\int_{\BS^{d-1}}s_k
\Big(\sum_{j=1}^ds_jz(\partial_jy_t)z\Big)dS(s).
$$
Minkowski's and H\"older's inequalities give, for each $k$,
\begin{align*}
&\|z(\partial_ky_t)z\|_{L_{\frac{d}{1-\epsilon}}(\RR^d_\theta)}\\
\leq&\frac d{|\BS^{d-1}|}\int_{\BS^{d-1}}|s_k|
\Big\|\sum_{j=1}^ds_jz(\partial_jy_t)z\Big\|_{L_{\frac{d}{1-\epsilon}}(\RR^d_\theta)}dS(s)\\
\leq&\frac d{|\BS^{d-1}|}
\Big(\int_{\BS^{d-1}}|s_k|^{\frac{d}{d-1+\epsilon}}dS(s)\Big)^{\frac{d-1+\epsilon}{d}}\\
&\times\Big(\int_{\BS^{d-1}}\tau_\theta\Big(
\Big|\sum_{j=1}^ds_jz(\partial_jy_t)z\Big|^{\frac{d}{1-\epsilon}}
\Big)dS(s)\Big)^{\frac{1-\epsilon}{d}}.
\end{align*}
\fi
In the same way as in Lemma \ref{lem:ge1}, and by using \eqref{eq:lasym}, we obtain
\begin{align*}
\sum_{j=1}^d\|z(\partial_jy_t)z\|_{L_{\frac{d}{1-\epsilon}}(\RR^d_\theta)}\lesssim&_{d,\epsilon}
\Big(\int_{\BS^{d-1}}\tau_\theta\Big(
\Big|\sum_{j=1}^ds_jz(\partial_jy_t)z\Big|^{\frac{d}{1-\epsilon}}
\Big)dS(s)\Big)^{\frac{1-\epsilon}{d}}\\
\lesssim&_{d,\epsilon}
\|z[(-\Delta)^{\frac{\epsilon}{2}},y_t]z\|_{\calL_{\frac{d}{1-\epsilon},\infty}(L_2(\RR^d))}\\
\leq&\|z\|_{L_\infty(\RR^d_\theta)}^2
\|[(-\Delta)^{\frac{\epsilon}{2}},y_t]\|_{\calL_{\frac{d}{1-\epsilon},\infty}(L_2(\RR^d))},
\end{align*}
where the second inequality follows from \eqref{eq:lasym}. This completes the proof.
\end{proof}

Now we can show that Theorem \ref{thm:m1} (i) implies Theorem \ref{thm:m1} (ii). The estimates for $z_\delta(\partial_jy_t)z_\delta$ obtained from Proposition \ref{prop:lasym} and Lemma \ref{lem:heat} are uniform in $\delta$ and $t$. We first let $\delta\to0$ and then $t\to0^+$ to obtain $\partial_jx\in L_{\frac{d}{1-\epsilon}}(\RR^d_\theta)$ from the weak Schatten class assumption.

\begin{theorem}\label{thm:clconv}
Let $-\frac{d}{2}<\epsilon\neq 0<1$ and $x\in\calS'(\RR^d_\theta)$. If $[(-\Delta)^{\frac{\epsilon}{2}},x]\in\calL_{\frac{d}{1-\epsilon},\infty}(L_2(\RR^d))$, then $x\in\dW^1_{\frac{d}{1-\epsilon}}(\RR^d_\theta)$ and
$$\|x\|_{\dW^1_{\frac{d}{1-\epsilon}}(\RR^d_\theta)}
\lesssim_{d,\epsilon}
\|[(-\Delta)^{\frac{\epsilon}{2}},x]\|_{\calL_{\frac{d}{1-\epsilon},\infty}(L_2(\RR^d))}.
$$\end{theorem}

\begin{proof}
Let $x\in\calS'(\RR^d_\theta)$ and $y_t=h_t*x$. Let $(z_\delta)_{\delta>0}$ denote the sequence from Lemma \ref{lem:ztest}. For each $j$, by Proposition \ref{prop:hprod},
$$z_\delta(\partial_jy_t)z_\delta\in\calS(\RR^d_\theta).$$
By Lemma \ref{lem:ztest}, as $\delta\to 0$, one has
\begin{equation}\label{eq:zdist}
\langle z_\delta(\partial_jy_t)z_\delta,a\rangle
=\langle\partial_jy_t,z_\delta az_\delta\rangle
\to\langle\partial_jy_t,a\rangle,
\quad a\in\calS(\RR^d_\theta).
\end{equation}
From Proposition \ref{prop:lasym} and Lemma \ref{lem:heat} (iii), we obtain
\begin{align*}
\sum_{j=1}^d\|z_\delta(\partial_jy_t)z_\delta\|_{L_{\frac{d}{1-\epsilon}}(\RR^d_\theta)}
\lesssim&_{d,\epsilon}\|z_\delta\|_{L_\infty(\RR^d_\theta)}^2
\|[(-\Delta)^{\frac{\epsilon}{2}},y_t]\|_{\calL_{\frac{d}{1-\epsilon},\infty}(L_2(\RR^d))}\\
\lesssim&_{d,\epsilon}\|[(-\Delta)^{\frac{\epsilon}{2}},x]\|_{\calL_{\frac{d}{1-\epsilon},\infty}(L_2(\RR^d))}.
\end{align*}
Fix $t>0$ and choose a sequence $(\delta_n)_{n\geq 1}\subseteq (0, \infty)$ such that $\delta_n\to0$ as $n\to \infty$. By the reflexivity of $L_{\frac{d}{1-\epsilon}}(\RR^d_\theta)$, one can assume that there exists a common subsequence and $v_1,\ldots,v_d\in L_{\frac{d}{1-\epsilon}}(\RR^d_\theta)$ such that
$$
z_{\delta_n}(\partial_jy_t)z_{\delta_n}\to v_j
\quad\text{in }L_{\frac{d}{1-\epsilon}}(\RR^d_\theta),
\quad 1\leq j\leq d.
$$
For every $a\in\calS(\RR^d_\theta)$, \eqref{eq:zdist} yields
\begin{align*}
\tau_\theta(v_ja)
=&\lim_{n\to\infty}\tau_\theta(z_{\delta_n}(\partial_jy_t)z_{\delta_n}a)\\
=&\lim_{n\to\infty}\langle\partial_jy_t,z_{\delta_n}az_{\delta_n}\rangle\\
=&\langle\partial_jy_t,a\rangle.
\end{align*}
This implies $\partial_jy_t=v_j\in L_{\frac{d}{1-\epsilon}}(\RR^d_\theta)$. For every $a\in L_{\frac{d}{d-1+\epsilon}}(\RR^d_\theta)$ with $\|a\|_{L_{\frac{d}{d-1+\epsilon}}(\RR^d_\theta)}\leq1$, H\"older's inequality yields
\begin{align*}
|\tau_\theta(v_ja)|
=\lim_{n\to\infty}|\tau_\theta(z_{\delta_n}(\partial_jy_t)z_{\delta_n}a)|
\leq\liminf_{n\to\infty}\|z_{\delta_n}(\partial_jy_t)z_{\delta_n}\|_{L_{\frac{d}{1-\epsilon}}(\RR^d_\theta)}.
\end{align*}
Taking the supremum over all such $a$, we obtain
\begin{align*}
\|v_j\|_{L_{\frac{d}{1-\epsilon}}(\RR^d_\theta)}
\leq\liminf_{n\to\infty}\|z_{\delta_n}(\partial_jy_t)z_{\delta_n}\|_{L_{\frac{d}{1-\epsilon}}(\RR^d_\theta)}.
\end{align*}
We conclude that
\begin{align*}
\|y_t\|_{\dW^1_{\frac{d}{1-\epsilon}}(\RR^d_\theta)}
=&\sum_{j=1}^d\|v_j\|_{L_{\frac{d}{1-\epsilon}}(\RR^d_\theta)}\lesssim_{d,\epsilon}
\|[(-\Delta)^{\frac{\epsilon}{2}},x]\|_{\calL_{\frac{d}{1-\epsilon},\infty}(L_2(\RR^d))}.
\end{align*}
By Lemma \ref{lem:heat} (i) and Lemma \ref{lem:wgrad}, one has 
$x\in\dW^1_{\frac{d}{1-\epsilon}}(\RR^d_\theta)$ and
$$
\|x\|_{\dW^1_{\frac{d}{1-\epsilon}}(\RR^d_\theta)}
\lesssim_{d,\epsilon}
\|[(-\Delta)^{\frac{\epsilon}{2}},x]\|_{\calL_{\frac{d}{1-\epsilon},\infty}(L_2(\RR^d))}.
$$
This completes the proof.
\end{proof}

Finally, we come to the proof of Theorem \ref{thm:m1}.

\begin{proof}[Proof of Theorem \ref{thm:m1}]
First, suppose that $x\in\dW^1_{\frac{d}{1-\epsilon}}(\RR^d_\theta)$. By \eqref{eq:transfer7},
$$
\mu\big(t,[(-\Delta_\theta)^{\frac{\epsilon}{2}},M_x]\big)
=\mu\big(t,[(-\Delta)^{\frac{\epsilon}{2}},x]\big),
\quad t\geq0.
$$
By Theorem \ref{thm:claw}, we have that Theorem \ref{thm:m1} (ii) implies Theorem \ref{thm:m1} (i). Conversely, suppose that $x\in\calS'(\RR^d_\theta)$ and $[(-\Delta_\theta)^{\frac{\epsilon}{2}},M_x]\in\calL_{\frac{d}{1-\epsilon},\infty}(L_2(\RR^d_\theta))$. By Theorem \ref{thm:clconv}, we obtain that $x\in\dW^1_{\frac{d}{1-\epsilon}}(\RR^d_\theta)$, and then Theorem \ref{thm:m1} (i) implies Theorem \ref{thm:m1} (ii). This completes the proof.
\end{proof}

\bigskip

\appendix
\section{Proof of Lemma \ref{lem:exp}}\label{sec:app1}
%\subsection{Proof of Lemma \ref{lem:exp}}
%\begin{proof}[Proof of Lemma \ref{lem:exp}]
We proceed directly with the proof of Lemma \ref{lem:exp}. For $g\in\calS(\RR^d)$,  by \eqref{eq:con1} and \eqref{eq:vmv1},
$$xg(s)=f *_\theta g(s) = \int_{\RR^d} e^{-\frac{\ii}{2} \langle s, \theta r\rangle} f(s- r) g(r) dr.$$
Note that
\begin{align*}
&s^\alpha\partial_s^\beta(f*_\theta g)(s)\\
=&\sum_{\gamma\leq\alpha, \delta\leq\beta}\binom{\alpha}{\gamma}\binom{\beta}{\delta}
\Big(\frac{\ii}{2}\Big)^{|\delta|}\int_{\RR^d}e^{-\frac{\ii}{2}\langle r,\theta s\rangle}r^\gamma(\theta r)^\delta f(r)(s-r)^{\alpha-\gamma}(\partial^{\beta-\delta}g)(s-r) dr\\
=&\sum_{\gamma\leq\alpha,\delta\leq\beta}\binom{\alpha}{\gamma}\binom{\beta}{\delta}\Big(\frac{\ii}{2}\Big)^{|\delta|}\Big[\big((\cdot)^\gamma(\theta \cdot)^\delta f\big)*_\theta\big((\cdot)^{\alpha-\gamma}\partial^{\beta-\delta}g\big)\Big](s).
\end{align*}
For $N\in\NN_0$, define
$$A_N=1+\sum_{|\alpha|+|\beta|\leq N}\sum_{\gamma\leq\alpha, \delta\leq\beta}\binom{\alpha}{\gamma}\binom{\beta}{\delta}2^{-|\delta|}\|s^\gamma(\theta s)^\delta f(s)\|_{L_1(\RR^d)}<\infty.$$
From Minkowski's inequality, one has
\begin{align*}
&q_N(f*_\theta g)\\
\leq&\sum_{|\alpha|+|\beta|\leq N}\sum_{\gamma\leq\alpha, \delta\leq\beta}\binom{\alpha}{\gamma}\binom{\beta}{\delta}2^{-\delta|}\|s^\gamma(\theta s)^\delta f(s)\|_{L_1(\RR^d)}\|s^{\alpha-\gamma}\partial_s^{\beta-\delta}g(s)\|_{L_2(\RR^d)}\\
\leq& A_Nq_N(g).
\end{align*}
Let $n\in\NN$. Let $f^{*_\theta1}:=f$ and $f^{*_\theta(n+1)}:=f*_\theta f^{*_\theta n}$. By induction, one has $q_N(f^{*_\theta n})\leq A_N^{n-1}q_N(f)$. Thus, for every $R>0$,
\begin{align*}
\sum_{n=1}^\infty\sup_{|t|\leq R}q_N\Big(\frac{(\ii t)^n}{n!}f^{*_\theta n}\Big)\leq\frac{q_N(f)}{A_N}\big(e^{A_NR}-1\big)<\infty
\end{align*}
and
\begin{align*}
\sum_{n=1}^\infty\sup_{|t|\leq R}q_N\Big(\frac{d}{dt}\Big(\frac{(\ii t)^n}{n!}f^{*_\theta n}\Big)\Big)\leq q_N(f)e^{A_NR}<\infty.
\end{align*}
Therefore, the map
$$t\mapsto u_t:=\sum_{n=1}^\infty \frac{(\ii t)^n}{n!}f^{*_\theta n}$$
is continuously differentiable from $\RR$ into $\calS(\RR^d)$.

By Proposition \ref{prop:iso2} and the Cauchy--Schwarz inequality,
\begin{align*}
\|U_\theta(g)\|_{L_\infty(\RR^d_\theta)}\lesssim&_d\|g\|_{L_1(\RR^d)}\\
=&\int_{\RR^d}\omega(s)^{-(d+1)} \omega(s)^{d+1}|g(s)| ds\\
\leq& \|\omega^{-(d+1)}\|_{L_2} \|\omega^{d+1}g\|_{L_2}\\
\lesssim&_d \|\omega^{d+1}g\|_{L_2}\\
\leq& q_{d+1}(g).
\end{align*}
By \eqref{eq:mul1}, $U_\theta(f^{*_\theta n})=x^n$. Therefore we obtain that with respect to the operator norm,
$$U_\theta(u_t)=\sum_{n=1}^\infty\frac{(\ii t)^n}{n!}x^n=e^{\ii tx}-1.$$
Note that $\frac{d}{dt}u_t=\ii f+\ii f*_\theta u_t=\ii f+\ii xu_t$ and $u_0=0$. This implies
\begin{align*}
\frac{d}{dt}\big(e^{-\ii tx}u_t\big)
=\ii e^{-\ii tx}f.
\end{align*}
Since $u_0=0$, we have
$$e^{-\ii tx}u_t=\ii\int_0^t e^{-\ii sx}f ds.$$
This yields
$$u_t=\ii\int_0^t e^{\ii(t-r)x}f dr.$$
By Proposition \ref{prop:iso1},
\begin{equation}\label{eq:q0}
\begin{split}
q_0(u_t)=&\|u_t\|_{L_2(\RR^d)}\\
\leq& \int_0^t \|e^{\ii(t-r)x}f\|_{L_2} dr\\
=&|t|q_0(f).
\end{split}
\end{equation}

For any $\alpha,\beta\in\NN_0^d$, since 
$$\ii e^{-\ii tx}\Big(f*_\theta\big((\cdot)^\alpha\partial^\beta u_t\big)\Big)
=\ii e^{-\ii tx}x\big((\cdot)^\alpha\partial^\beta u_t\big),$$
one has
\begin{align*}
\frac{d}{dt}\Big(e^{-\ii tx}(\cdot)^\alpha\partial^\beta u_t\Big)
=&-\ii e^{-\ii tx}x\big((\cdot)^\alpha\partial^\beta u_t\big)+\ii e^{-\ii tx}(\cdot)^\alpha\partial^\beta f\\
&+\ii e^{-\ii tx}\sum_{\gamma\leq\alpha, \delta\leq\beta}\binom{\alpha}{\gamma}\binom{\beta}{\delta}\Big(\frac{\ii}{2}\Big)^{|\delta|}\Big[\big((\cdot)^\gamma(\theta\cdot)^\delta f\big)*_\theta
\big((\cdot)^{\alpha-\gamma}\partial^{\beta-\delta}u_t\big)\Big]\\
=&\ii e^{-\ii tx}(\cdot)^\alpha\partial^\beta f\\
&+\ii\sum_{\substack{\gamma\leq\alpha, \delta\leq\beta\\
|\gamma|+|\delta|>0}}\binom{\alpha}{\gamma}\binom{\beta}{\delta}\Big(\frac{\ii}{2}\Big)^{|\delta|}e^{-\ii tx}U_\theta\big((\cdot)^\gamma(\theta \cdot)^\delta f\big)\big((\cdot)^{\alpha-\gamma}\partial^{\beta-\delta}u_t\big).
\end{align*}
Note that $u_0=0$, and then
\begin{align*}
\int_0^t\frac{d}{dr}\Big(e^{-\ii rx}(\cdot)^\alpha\partial^\beta u_r\Big) dr
=&e^{-\ii tx}(\cdot)^\alpha\partial^\beta u_t
-e^{0}(\cdot)^\alpha\partial^\beta u_0\\
=&e^{-\ii tx}(\cdot)^\alpha\partial^\beta u_t.
\end{align*}
We use Proposition \ref{prop:iso1} and Minkowski's inequality to get
\begin{align*}
\|(\cdot)^\alpha\partial^\beta u_t\|_{L_2(\RR^d)}
=&\Big\|\int_0^t\frac{d}{dr}\Big(e^{-\ii rx}(\cdot)^\alpha\partial^\beta u_r \Big)dr \Big\|_{L_2(\RR^d)}\\
\leq&|t|\|(\cdot)^\alpha\partial^\beta f\|_{L_2(\RR^d)}+\sum_{\substack{\gamma\leq\alpha, \delta\leq\beta\\
|\gamma|+|\delta|>0}}\binom{\alpha}{\gamma}\binom{\beta}{\delta}2^{-|\delta|}\|s^\gamma(\theta s)^\delta f(s)\|_{L_1(\RR^d)}\\
&\times\int_0^{|t|}\sup_{|v|\leq r}\|(\cdot)^{\alpha-\gamma}
\partial^{\beta-\delta}u_v\|_{L_2(\RR^d)} dr.
\end{align*}
Each summand with $|\alpha|+|\beta|\leq N$ satisfies $|\alpha-\gamma|+|\beta-\delta|\leq N-1$. Summing over all $|\alpha|+|\beta|\leq N$, we obtain that for $N\geq1$ and $R\geq0$
$$\sup_{|t|\leq R}q_N(u_t)\leq Rq_N(f)+A_N\int_0^R\sup_{|v|\leq r}q_{N-1}(u_v) dr.$$
Let $C_0=q_0(f)$ and $C_N=q_N(f)+\frac{A_N}{2}C_{N-1}$. When $N=0$, it follows from \eqref{eq:q0} that $q_0(u_t)\le C_0|t|$. Assume that the estimate holds at order $N-1$, namely
$$\sup_{|v|\leq r}q_{N-1}(u_v)\leq C_{N-1}r(1+r)^{N-1}.$$
Note that
\begin{align*}
\int_0^R r(1+r)^{N-1} dr\leq (1+R)^{N-1}\int_0^R r dr=\frac{R^2}{2}(1+R)^{N-1}
\leq \frac{R}{2}(1+R)^N.
\end{align*}
By the inductive hypothesis,
\begin{align*}
\sup_{|t|\leq R}q_N(u_t)
\leq& Rq_N(f)+A_NC_{N-1}\int_0^R r(1+r)^{N-1} dr\\
\leq& Rq_N(f)+\frac{A_NC_{N-1}}{2}R(1+R)^N\\
\leq& q_N(f)R(1+R)^N+\frac{A_NC_{N-1}}{2}R(1+R)^N\\
=&\Big(q_N(f)+\frac{A_NC_{N-1}}{2}\Big)R(1+R)^N\\
=&C_NR(1+R)^N.
\end{align*}
Hence, one can conclude
$$q_N(u_t)\leq C_N|t|(1+|t|)^N,$$
where $C_N$ depends only on $N,f,\theta,d$. This completes the proof.
%\end{proof}

\bigskip

\section{Cwikel estimates}\label{sec:app2}
This appendix contains the proofs of the Hilbert--Schmidt identity, the high-exponent Cwikel estimate, and the coefficient estimates used in Section \ref{s:cw}. All statements are valid for arbitrary real antisymmetric matrices $\theta$.
\subsection{Proof of Lemma \ref{lem:hs}}
%\begin{proof}[Proof of Lemma \ref{lem:hs}]
By Proposition \ref{prop:iso1}, we can write $x=U_\theta(f)$ with $f\in L_2(\RR^d)$. By \eqref{eq:defufh1} and change of variables,
$$
(U_\theta(f)M_g\xi)(t)
=\int_{\RR^d}f(t-s)g(s)e^{-\frac{\ii}{2}\langle t,\theta s\rangle}\xi(s) ds,\quad\xi\in L_2(\RR^d).
$$
This implies the kernel $K_{f,g}$ of $U_\theta(f)M_g$ is give by
$$ K_{f,g}(t,s)=f(t-s)g(s)e^{-\frac{\ii}{2}\langle t,\theta s\rangle}. $$
By Fubini's theorem,
\begin{align*}
\|K_{f,g}\|_{L_2(\RR^{2d})}^2
=&\int_{\RR^d}|g(s)|^2\int_{\RR^d}|f(t-s)|^2 dt ds\\
=&\|f\|_{L_2(\RR^d)}^2\|g\|_{L_2(\RR^d)}^2.
\end{align*}
This completes the proof.
%\end{proof}
\subsection{Proof of Theorem \ref{thm:cw-hi}}
%\begin{proof}[Proof of Theorem \ref{thm:cw-hi}]
Denote by
$$
(\mathcal A_1,\tau_1) = \big(L_\infty(\RR^d_\theta),\tau_\theta\big), \quad (\mathcal A_2,\tau_2) = \big(L_\infty(\RR^d),\int\big).
$$
Define normal $*$-representations on $L_2(\RR^d)$ by
$$
\pi_1(x)=x,\quad \pi_2(g)=M_g.
$$
Lemma \ref{lem:hs} yields \cite[Hypothesis 3.1]{LSZ2020}. Since $E$ interpolates the couple $(L_2(0,\infty),L_\infty(0,\infty))$,
$$
E(0,\infty)\subset L_2(0,\infty)+L_\infty(0,\infty).
$$
Hence, $x\otimes g\in E(L_\infty(\RR^d_\theta)\overline\otimes L_\infty(\RR^d))$ implies
$$
x\otimes g\in L_2(L_\infty(\RR^d_\theta)\overline\otimes L_\infty(\RR^d))+L_\infty(L_\infty(\RR^d_\theta)\overline\otimes L_\infty(\RR^d)).
$$
Therefore, one can apply \cite[Corollary 3.5]{LSZ2020} to obtain
$$
\|xM_g\|_{E(L_2(\RR^d))} \lesssim_E\|x\otimes g\|_{E(\mathcal N_\theta)}.
$$
In particular, if $E=L_p(0,\infty)$ with $2\leq p<\infty$, we obtain \eqref{eq:cw-q}. If $2<p<\infty$, we take $E=L_{p,\infty}(0,\infty)$ to obtain \eqref{eq:cw-wq}.
%\end{proof}

\subsection{Proof of Lemma \ref{lem:coef}}\label{appendixB:5.14}
%\begin{proof}[Proof of Lemma \ref{lem:coef}]
Let $\Delta_u=\sum_{j=1}^d\frac{\partial^2}{\partial u_j^2}$ be the classical Laplacian. Since $h_{r,v}$ is supported in the interior of $(-2,2)^d$ and
$(1-\Delta_u)^Le_m=\big(1+\frac{\pi^2}{4}|m|^2\big)^Le_m$, one has that by taking integration by parts and for every integer $L\geq0$,
$$c_{m,r}(v)=\Big(1+\frac{\pi^2}{4}|m|^2\Big)^{-L}
\int_{(-2,2)^d}(1-\Delta_u)^Lh_{r,v}(u)\overline{e_m(u)} du.$$
Thus
$$\|c_{m,r}\|_{L_\infty(Q)}
\lesssim_{L,N,f,\theta}(1+|m|)^{-2L}(1+|r|)^{-N}.$$
We can choose $2L\geq M$ to obtain the desired coefficient estimate. For $M>d$,
$$\sum_{m\in\ZZ^d}
\sup_{\substack{u\in(-2,2)^d\\v\in Q}}|c_{m,r}(v)e_m(u)|
\lesssim_{d,M,N,f,\theta}(1+|r|)^{-N}\sum_{m\in\ZZ^d}(1+|m|)^{-M}<\infty.$$
Therefore, the series above converges absolutely and uniformly. For fixed $v$, it equals $h_{r,v}$ in $L_2((-2,2)^d)$ because $(e_m)_m$ is an orthonormal basis. Since $e_m(u)$ is continuous with respect to the variable $u$, the series equals $h_{r,v}$ everywhere. In addition, since $\eta=1$ on $Q$, the stated identity follows immediately.
%\end{proof}

\subsection{Proof of Lemma \ref{lem:jtrace}}
%\begin{proof}[Proof of Lemma \ref{lem:jtrace}]
It is clear that the operator-valued function $J^{-(a+b)z}$ is differentiable with respect to the operator norm on the open strip $0<\Re z<1$, with its derivative $-(a+b)(\log J)J^{-(a+b)z}$. Let $z_n\to z_0$ in the closed strip $0\leq\Re z\leq1$. For every $\xi\in L_2(\RR^d)$, by the dominated convergence theorem,
\begin{align*}
&\|\big(J^{-(a+b)z_n}-J^{-(a+b)z_0}\big)\xi\|_{L_2(\RR^d)}^2\\
=&\int_{\RR^d}
\Big|(1+|t|^2)^{-\frac{(a+b)z_n}{2}}
-(1+|t|^2)^{-\frac{(a+b)z_0}{2}}\Big|^2
|\xi(t)|^2 dt\to0.
\end{align*}
The same argument also applies to the exponents $1-z_n$ and $1-z_0$. This implies that the following operator-valued function
$$z\mapsto J^{-(a+b)z}xJ^{-(a+b)(1-z)},
\quad 0\leq\Re z\leq1,$$
is holomorphic with respect to the operator norm on the open strip and strongly continuous on the closed strip. Moreover,
$$\|J^{-(a+b)z}xJ^{-(a+b)(1-z)}\|_{\calB(L_2(\RR^d))}
\leq\|x\|_{\calB(L_2(\RR^d))},
\quad 0\leq\Re z\leq1.$$
If $B\in\calB(L_2(\RR^d))$ is of finite rank, then the following function
$$z\mapsto\Tr\big(BJ^{-(a+b)z}xJ^{-(a+b)(1-z)}\big)$$
is holomorphic on $0<\Re z<1$. It remains to show the continuity of the above function on the closed strip. Let $m=\rank B$ and choose an orthonormal basis $(u_j)_{j\geq1}$ with $\operatorname{Ran}B=\operatorname{span}\{u_1,\ldots,u_m\}$. Then we can write
$$\Tr\big(BJ^{-(a+b)z}xJ^{-(a+b)(1-z)}\big)
=\sum_{j=1}^m\big\langle J^{-(a+b)z}xJ^{-(a+b)(1-z)}u_j,B^*u_j
\big\rangle_{L_2(\RR^d)}.$$
Note that each summand is continuous by the strong continuity. This completes the proof.
%\end{proof}

\subsection{Proof of Lemma \ref{lem:w11}}
%\begin{proof}[Proof of Lemma \ref{lem:w11}]
Fix $M,N\geq0$. We choose an integer $L\in\NN_0$ such that
$$L\geq\max\Big\{N+b+\gamma+2,\frac M2\Big\}.$$
Throughout this proof, implicit constants may depend on $b,c,\gamma,d,f,\sigma,\theta,\eta$ in addition to the parameters displayed in the subscripts, but are independent of $m,r,n,u,v$. For $r,n\in\ZZ^d$, $u\in(-2,2)^d$ and $v\in Q$
$$s=n+v,\quad\delta=r+u-v,\quad t=s+\delta=n+r+u.$$
For every $\alpha_2\in\NN_0^d$, it is clear that
$$|\partial_u^{\alpha_2}\omega(t)^{-b}|
\lesssim\omega(t)^{-b-|\alpha_2|}
\leq\omega(t)^{-b}\lesssim\omega(s)^{-b}\omega(\delta)^b.$$
For $0\leq\tau\leq1$, we also have
$$\omega(s)\leq\sqrt2 \omega(s+\tau\delta)\omega(\tau\delta)
\leq\sqrt2 \omega(s+\tau\delta)\omega(\delta).$$
Since
$$\sigma(s)-\sigma(t)
=-\int_0^1\sum_{j=1}^d\delta_j(\partial_j\sigma)(s+\tau\delta) d\tau$$
and
$$\partial_u^\nu\big(\sigma(s)-\sigma(t)\big)
=-(\partial^\nu\sigma)(t),\quad |\nu|\geq1, \nu\in\NN_0^d,$$
one has that for every $\nu\in\NN_0^d$,
$$|\partial_u^\nu\big(\sigma(s)-\sigma(t)\big)|
\lesssim\omega(s)^{-\gamma-1}\omega(\delta)^{\gamma+2}.$$
Now fix $\alpha\in\NN_0^d$ with $|\alpha|\leq2L$. By Leibniz's rule, we obtain
\begin{align*}
\big|\partial_u^\alpha h_{r,n,v}(u)\big|
\leq& \omega(s)^{-c}
\sum_{\alpha_1+\alpha_2+\alpha_3+\alpha_4=\alpha}
C_{\alpha_1\alpha_2\alpha_3\alpha_4}
\big|(\partial^{\alpha_1}\eta)(u)\big| 
\big|\partial_u^{\alpha_2}\omega(t)^{-b}\big| 
\big|(\partial^{\alpha_3}f)(\delta)\big|\\
&\times
\sum_{\beta+\nu=\alpha_4}C_{\beta\nu}
\Big|\Big(\frac{\ii}{2}\theta^{T}v\Big)^\beta\Big| 
\big|\partial_u^\nu\big(\sigma(s)-\sigma(t)\big)\big|\\
\lesssim&\omega(s)^{-c} \omega(s)^{-b}\omega(\delta)^b 
\omega(s)^{-\gamma-1}\omega(\delta)^{\gamma+2} 
(1+|r|)^{-L}\\
\lesssim&(1+|n|)^{-c-b-\gamma-1}(1+|r|)^{b+\gamma+2-L}.
\end{align*}
Therefore, the choice of $L$ yields
\begin{align*}
\sup_{u\in(-2,2)^d, v\in Q}|\partial_u^\alpha h_{r,n,v}(u)|
\lesssim&(1+|r|)^{-L+b+\gamma+2}(1+|n|)^{-b-c-\gamma-1}\\
\leq&(1+|r|)^{-N}(1+|n|)^{-\rho}.
\end{align*}
By using
$$(1-\Delta_u)^Le_m
=\Big(1+\frac{\pi^2}{4}|m|^2\Big)^Le_m,$$
and by taking integration by parts, we have
\begin{align*}
c_{m,r,n}(v)
=&\int_{(-2,2)^d} h_{r,n,v}(u) \overline{e_m(u)} du\\
=&\Big(1+\frac{\pi^2}{4}|m|^2\Big)^{-L}
\int_{(-2,2)^d} (1-\Delta_u)^L h_{r,n,v}(u) 
\overline{e_m(u)} du .
\end{align*}
This implies that for $v\in Q$,
$$|c_{m,r,n}(v)|\lesssim(1+|m|)^{-2L}(1+|r|)^{-N}(1+|n|)^{-\rho}.$$
Since $2L\geq M$, we obtain the desired coefficient estimates. If $M>d$, then
\begin{align*}
\sum_{m\in\ZZ^d}
\sup_{\substack{u\in(-2,2)^d\\v\in Q}}
|c_{m,r,n}(v)e_m(u)|\lesssim (1+|r|)^{-N}(1+|n|)^{-\rho}
\sum_{m\in\ZZ^d}(1+|m|)^{-M}<\infty.
\end{align*}
This yields the absolute and uniform convergence of the Fourier series on $(-2,2)^d\times Q$. By the same argument as in the last part of the proof of Lemma \ref{lem:coef} (see Appendix \ref{appendixB:5.14}), we get the desired Fourier identity everywhere.
%For each fixed $v\in Q$, the functions $(e_m)_{m\in\ZZ^d}$ form an orthonormal basis of $L_2((-2,2)^d)$, so the series converges to $h_{r,n,v}$ in that space. Its uniform sum therefore equals $h_{r,n,v}$ almost everywhere. Since both are continuous in $u$, equality holds everywhere, giving the asserted Fourier identity.
%\end{proof}

\bigskip

\noindent\textbf{AI Statement:} No artificial intelligence (AI) tools were used in the preparation of this paper. All ideas, proofs and writing were produced entirely by the authors.

\

\noindent\textbf{Acknowledgment:} The authors would like to thank Professor Fedor Sukochev and Professor Xiao Xiong for helpful comments.

\end{document}